\documentclass[onefignum,onetabnum]{siamart190516}

\usepackage{dsfont}

\usepackage{amsfonts,amssymb,stmaryrd}
\usepackage{amsmath}
\usepackage{mathrsfs,mathtools,enumitem}
\usepackage{algorithmic}
\usepackage{longtable,units,subdepth}
\Crefname{ALC@unique}{Line}{Lines} 

\usepackage{lipsum}
\usepackage{amsfonts}
\usepackage{graphicx}
\usepackage{epstopdf}
\usepackage{algorithmic}
\ifpdf
  \DeclareGraphicsExtensions{.eps,.pdf,.png,.jpg}
\else
  \DeclareGraphicsExtensions{.eps}
\fi

\newsiamremark{remark}{Remark}
\newsiamremark{hypothesis}{Hypothesis}
\crefname{hypothesis}{Hypothesis}{Hypotheses}
\newsiamthm{claim}{Claim}

\headers{DMP-preserving FEM}{G.R. Barrenechea, V. John, and P. Knobloch}

\title{Finite element methods respecting the discrete maximum principle for
convection-diffusion equations\thanks{Submitted to the editors DATE.
\funding{The research of all authors was supported by the programme {\it Research in pairs} of 
the Mathematisches Forschungsinstitut Oberwolfach (MFO),  grant No.~1937p.
The work of Gabriel R. Barrenechea has been partially funded by the Leverhume Trust via the Research Fellowship No. RF-2019-510. 
The work of Volker John has been supported by the Research Training 'Group
Differential Equation- and Data-driven Models in Life Sciences and Fluid
Dynamics' (DAEDALUS), RTG 2433, funded by the German Research Foundation (DFG).
The work of Petr Knobloch has been supported by the grant No.~22-01591S of the
Czech Science Foundation.}}}

\author{Gabriel R. Barrenechea\thanks{Department of Mathematics and Statistics, 
University of Strathclyde, 26 Richmond Street, Glasgow G1
1XH, Scotland. \email{gabriel.barrenechea@strath.ac.uk}}
\and Volker John\thanks{Weierstrass Institute for Applied Analysis and Stochastics, Leibniz Institute im
Forschungsverbund Berlin e.~V.~(WIAS), Mohrenstr.~29, 10117 Berlin, and
Freie Universit\"at Berlin, Department of Mathematics and Computer Science,
Arnimallee 6, 14195 Berlin, Germany.
 \email{john@wias-berlin.de}.}
\and Petr Knobloch\thanks{Department of Numerical Mathematics, Faculty of Mathematics and
Physics, Charles University, Sokolovsk\'a 83, 18675 Praha 8, Czech Republic.
\email{knobloch@karlin.mff.cuni.cz}}}

\usepackage{amsopn}

\ifpdf
\hypersetup{
  pdftitle={DMP preserving FEM},
  pdfauthor={G.R. Barrenechea, V. John, P. Knobloch}
}
\fi

\newcommand{\bb}{\boldsymbol b}
\newcommand{\bx}{\boldsymbol x}
\newcommand{\by}{\boldsymbol y}
\newcommand{\bv}{\boldsymbol v}
\newcommand{\bu}{\boldsymbol u}

\newcommand{\bn}{\boldsymbol n}
\newcommand{\bt}{\boldsymbol t}
\newcommand{\bJ}{\boldsymbol J}
\newcommand{\bJK}{{\boldsymbol J}_{\hspace*{-0.5mm}K}^{}}

\newcommand\calT{\mathscr{T}}
\newcommand\calF{\mathscr{F}}
\newcommand\calE{\mathscr{E}}
\newcommand\calI{\mathscr{I}}
\newcommand\calL{\mathscr{L}}

\newcommand{\bbR}{\mathbb{R}}

\newcommand\bbA{\mathbb{A}}
\newcommand\bbAd{\mathbb{A}_{\mathrm{d}}^{}}
\newcommand\bbAc{\mathbb{A}_{\mathrm{c}}^{}}
\newcommand\bbAN{\mathbb{A}_{\mathrm{N}}^{}}
\newcommand\bbAI{\mathbb{A}_{\mathrm{I}}^{}}
\newcommand\bbAB{\mathbb{A}_{\mathrm{B}}^{}}
\newcommand\bbAdI{\mathbb{A}_{\mathrm{d,I}}^{}}

\newcommand\bbD{\mathbb{D}}

\newcommand\bbP{\mathbb{P}}
\newcommand\bbQ{\mathbb{Q}}

\newcommand\bbB{\mathbb{B}}

\newcommand\bbMc{\mathbb M_{\mathrm{c}}}
\newcommand\bbMl{\mathbb M_{\mathrm{l}}}

\newcommand\fl{\ell}
\newcommand\cc{c}
\newcommand\mm{m}

\def\third{{\textstyle\frac13}}
\def\twothirds{{\textstyle\frac23}}
\def\sixth{{\textstyle\frac16}}
\def\MHAdK{{\mathbb{A}}_{\mathrm{d}}^K}
\def\MHAK{\hat{\mathbb{A}}_{\mathrm{c}}^K}
\def\MHTAK{\tilde{\mathbb{A}}_{\mathrm{c}}^K}
\def\MHTAKs#1{\tilde{\mathbb{A}}_{\mathrm{c}}^{K,#1}}
\def\MHAC{\hat{\mathbb{A}}_{\mathrm{c,MH}}}
\def\MHTAC{\tilde{\mathbb{A}}_{\mathrm{c,MH}}}
\def\MHaKij{\hat{c}^K_{ij}}

\newcommand{\bw}{\boldsymbol w}
\newcommand{\bff}{\boldsymbol f}

\newcommand{\mD}{\mathbb D}
\newcommand{\mL}{\mathbb L}
\newcommand{\mK}{\mathbb K}
\newcommand{\mBI}{\bbB_{\mathrm{I}}}
\newcommand{\mBB}{\bbB_{\mathrm{B}}}
\newcommand{\mKI}{\mK_{\mathrm{I}}}
\newcommand{\mKB}{\mK_{\mathrm{B}}}
\newcommand{\buI}{\bu_{\mathrm{I}}}
\newcommand{\buB}{\bu_{\mathrm{B}}}

\newsiamremark{assumption}{Assumption}

\begin{document}

\maketitle

\begin{abstract}
Convection-diffusion-reaction equations model the conservation of scalar quantities. From the analytic point 
of view, solution of these equations satisfy under certain conditions maximum principles, which represent 
physical bounds of the solution. That the same bounds are respected by numerical approximations 
of the solution is often of utmost importance in practice. The mathematical formulation of this
property, which contributes to  the physical consistency of a method,  
is called Discrete Maximum Principle (DMP). In many applications,
convection dominates diffusion by several orders of magnitude. 
It is well known
that  standard discretizations typically do not satisfy 
the DMP in this convection-dominated regime. In fact, in this case,
it turns out to be a challenging problem to construct discretizations that, on the one hand, 
respect the  DMP and, on the other hand, compute accurate solutions. 
This paper presents a survey on  finite element methods, with a main focus on the convection-dominated regime,  that satisfy a local or a global DMP. 
The concepts of the underlying numerical analysis are discussed.  
The survey reveals that for the steady-state problem there are only a few discretizations, all of them 
nonlinear, that at the same time satisfy the DMP and compute reasonably accurate solutions, e.g., algebraically 
stabilized schemes. Moreover, most of these discretizations have been developed in recent years, 
showing the enormous progress that has been achieved lately. 
Methods based on algebraic stabilization, nonlinear and linear ones, are currently as well 
the only finite element methods that combine the satisfaction of the global DMP and accurate numerical 
results for the evolutionary equations in the convection-dominated situation. 
\end{abstract}



\centerline{CONTENTS}
\medskip

\begin{center}
\begin{longtable}{llr}
\ref{Sec:Intro} & Introduction \dotfill & \pageref{Sec:Intro} \\
\ref{Sec:model} & The steady-state model problem, general notations \dotfill & \pageref{Sec:model}\\
\ref{Sec:steady_model} & The steady-state model problem \dotfill & \pageref{Sec:steady_model}\\
\ref{Sec:triangulations} & $\ldots$ Triangulations and finite element spaces \dotfill & \pageref{Sec:triangulations}\\
\ref{Sec:FEM-Matrices} & $\ldots$ Finite element matrices\dotfill  & \pageref{Sec:FEM-Matrices} \\
\ref{Sec:general_DMP} & General results on DMP satisfying discretizations \dotfill & \pageref{Sec:general_DMP}\\
\ref{Sec:DMP_linear_disc}& $\ldots$  Linear discretizations \dotfill & \pageref{Sec:DMP_linear_disc}\\
\ref{Sec:DMP_nonlinear_disc}& $\ldots$  Nonlinear discretizations\dotfill  & \pageref{Sec:DMP_nonlinear_disc}\\
\ref{Sec:cd_linear_methods} & Linear  discretizations of steady-state problems without convection\dotfill  & \pageref{Sec:cd_linear_methods}\\ 
\ref{Sec:cd_linear_methods_poisson} & $\ldots$ The Poisson problem \dotfill & \pageref{Sec:cd_linear_methods_poisson}\\
\ref{Sec:cd_linear_methods_reaction} & $\ldots$ The reaction-diffusion equation and mass lumping \dotfill & \pageref{Sec:cd_linear_methods_reaction} \\
\ref{Sec:cd_linear_methods_cd}  & Linear discretizations of the steady-state problem \dotfill & 
\pageref{Sec:cd_linear_methods_cd} \\
\ref{Sec:cd_linear_methods_gal} & $\ldots$  The Galerkin finite element method\dotfill & 
\pageref{Sec:cd_linear_methods_gal}\\
\ref{Sec:cd_linear_methods_iso_diff} & $\ldots$  Isotropic linear artificial diffusion \dotfill & 
\pageref{Sec:cd_linear_methods_iso_diff} \\
\ref{Sec:cd_linear_methods_upw} & $\ldots$ Upwind finite element methods \dotfill & 
\pageref{Sec:cd_linear_methods_upw}\\
\ref{sec:xu-zikatanov_fem} & $\ldots$ The edge-averaged  finite element method \dotfill& \pageref{sec:xu-zikatanov_fem} \\
\ref{Sec:cd_nonlinear_methods} & Nonlinear stabilized discretizations of the steady-state problem \dotfill& \pageref{Sec:cd_nonlinear_methods}\\
\ref{Sec:cd_nonlinear_methods_MH} & $\ldots$ The Mizukami--Hughes method\dotfill 
& \pageref{Sec:cd_nonlinear_methods_MH} \\
\ref{Sec:cd_nonlinear_methods_BE} &  $\ldots$ Burman--Ern Methods \dotfill & \pageref{Sec:cd_nonlinear_methods_BE}\\
\ref{Sec:cd_nonlinear_methods_afc} & $\ldots$ Algebraic Flux Correction methods \dotfill &  \pageref{Sec:cd_nonlinear_methods_afc} \\
\ref{Sec:cd_nonlinear_methods_edge} & $\ldots$ A monotone Local Projection Stabilized (LPS) method \dotfill & \pageref{Sec:cd_nonlinear_methods_edge}\\
\ref{sec:num_exam} & A numerical illustration \dotfill & \pageref{sec:num_exam} \\
\ref{sec:tcd} & Time-dependent problem \dotfill & \pageref{sec:tcd}\\
\ref{sec:tcd_cont_prob} & $\ldots$ The continuous problem \dotfill &  \pageref{sec:tcd_cont_prob}\\
\ref{sec:tcd_framework} & $\ldots$ Maximum principle, DMP, and positivity preservation \dotfill & \pageref{sec:tcd_framework} \\
\ref{sec:tcd_ext_methods_steady} &  $\ldots$ Linear methods \dotfill & \pageref{sec:tcd_ext_methods_steady} \\
\ref{sec:tcd_femfct} &  $\ldots$ FEM Flux-Corrected-Transport (FCT) schemes\dotfill & 
\pageref{sec:tcd_femfct}\\
\ref{sec:other_fems} & Other types of finite elements\dotfill & \pageref{sec:other_fems}\\
\ref{sec:other_fems_Q1} & $\ldots$  $\bbQ_1$ finite element \dotfill & \pageref{sec:other_fems_Q1}\\
\ref{sec:other_fems_H1} & $\ldots$  Higher order $H^1$-conforming finite elements\dotfill& \pageref{sec:other_fems_H1}\\
\ref{sec:other_fems_CR} & $\ldots$ Non-conforming finite elements of Crouzeix--Raviart type \dotfill & \pageref{sec:other_fems_CR}\\
\ref{sec:other_fems_DG} & $\ldots$ Discontinuous Galerkin finite element methods \dotfill&
\pageref{sec:other_fems_DG}\\
\ref{sec:hyperbolic_problems} & Brief comments on hyperbolic conservation laws \dotfill & 
\pageref{sec:hyperbolic_problems}\\
\ref{sec:summary} & Summary \dotfill & \pageref{sec:summary}\\
                  & References \dotfill & \pageref{sec:references}
\end{longtable}
\end{center}

\section{Introduction}\label{Sec:Intro}

Partial differential equations (PDEs) or systems of them are widely used for modeling processes from 
nature and industry. Usually, an analytic solution cannot
be obtained. In practice, numerical methods are utilized for computing
approximations of the solution. Such numerical methods consist of several
components, like discretizations with respect to different variables, 
approaches for solving nonlinear problems, and solvers for systems of
linear algebraic
equations. The actual choice of these components might be dictated by 
different goals, like efficiency, or accuracy with respect to quantities of
interest. A particular aspect of the second goal is the so-called physical 
consistency of a method, i.e., certain fundamental physical properties 
of the solution of the PDE should be inherited by
the numerical solution. For many practitioners, the physical consistency is 
an essential criterion for utilizing a numerical method. 

Classes of PDEs that can be found in many models from 
applications are elliptic linear second order equations
\begin{equation}\label{eq:elli_eq}
-\varepsilon\Delta u + \bb\cdot\nabla u + \sigma u =f \quad \mbox{in } \Omega,
\end{equation}
and their parabolic counterparts
\begin{equation}\label{eq:para_eq}
\partial_t u -\varepsilon\Delta u + \bb\cdot\nabla u + \sigma u =f \quad \mbox{in } (0,T] \times \Omega.
\end{equation}
In these equations $\Omega \subset \mathbb R^d$, $d \ge 1$, is a spatial domain, $(0,T]$
a time interval, and  $u$ is some scalar quantity like the temperature 
or a concentration. This scalar quantity is transported by molecular diffusion with 
the diffusion coefficient $\varepsilon\ [\unitfrac{m^2}{s}]$ and by convective transport
with the velocity field $\bb\ [\unitfrac{m}s]$. The zeroth order term in \eqref{eq:elli_eq}
and \eqref{eq:para_eq} is called reactive term with the reaction coefficient $\sigma\ 
[\unitfrac1s]$ and the term on the right-hand side describes sinks and sources of the 
scalar quantity. Both equations \eqref{eq:elli_eq} and \eqref{eq:para_eq}  have to be 
equipped with suitable boundary conditions at the boundary $\partial\Omega$ of $\Omega$
and \eqref{eq:para_eq} also  with an initial condition at $t=0$ in order to define  well-posed 
problems. Then, the analysis of \eqref{eq:elli_eq} and \eqref{eq:para_eq} is very well 
understood. In particular, it can be shown that under appropriate assumptions on
the data of the problems, so-called Maximum Principles (MP) are 
satisfied. That means, loosely speaking, that the solution at some point or in
some subdomain 
can be bounded a  priori, e.g., for a global MP by the values on $\partial\Omega$
and, for the evolutionary problem,  also on $\{0\}\times\Omega$. 
In case that the assumptions for the satisfaction of the MP are satisfied, it 
represents a fundamental physical property of solutions of \eqref{eq:elli_eq} and \eqref{eq:para_eq}.

A physically consistent discretization of \eqref{eq:elli_eq} and \eqref{eq:para_eq}
should satisfy discrete counterparts of the MP, the so-called Discrete Maximum Principle (DMP).
Discretizations that do not fulfill the DMP are prone to numerical solutions with unphysical 
values, so-called spurious oscillations. Usually, equations of type \eqref{eq:elli_eq} and \eqref{eq:para_eq}
are part of coupled problems and their numerical solution serves as input data for other 
equations. With spurious oscillations in this input, there is a high probability that also the 
numerical solutions of the remaining equations possess unphysical values and finally the 
numerical simulation of the coupled problem might blow up, as it is our own experience 
reported in \cite{JMRSTV09}. Consequently, the satisfaction of the DMP is essential for 
discretizations of \eqref{eq:elli_eq} and \eqref{eq:para_eq} to be useful for simulations
in applications. If this property is satisfied, then efficiency or the satisfaction of
other physical properties, like conservation properties, or
the accuracy with respect 
to quantities of interest, like norms in Sobolev spaces, are further criteria for 
selecting a method. 

The first proof of a maximum principle for a discretization of a PDE was 
presented by Gershgorin \cite{Ger30} already in 1930. A generalization of this 
result is given in the monograph by Collatz \cite{Col55} from 1955, whose English 
translation is \cite{Col60}. The consideration of discrete analogs of maximum principles 
can be found in papers by Bramble and Hubbard \cite{BH62,BH64} published in the early 1960s. 
In 1970, Ciarlet 
presented in \cite{Ciarlet70} necessary and sufficient conditions for a 
discretization to satisfy a DMP. In all these works, finite difference methods are considered. 
However, all arguments from linear algebra that were utilized in these papers can 
be applied analogously to linear systems of equations arising from other discretizations. 
The first work that studies the DMP explicitly for finite element methods was published in 
1973 by Ciarlet and Raviart \cite{CR73}. Since then, numerous papers appeared studying 
the DMP for different discretizations of elliptic and parabolic boundary value problems. 

Convection-diffusion-reaction equations \eqref{eq:elli_eq} and \eqref{eq:para_eq} possess 
a feature that makes the computation of a numerical solution challenging. In most 
applications, the convective transport by the velocity field strongly dominates
the diffusive
transport. Hence, the first order term in \eqref{eq:elli_eq} and \eqref{eq:para_eq} 
is dominant. Under appropriate conditions on the smoothness of the data, 
it can be shown that (weak) solutions of  \eqref{eq:elli_eq} and \eqref{eq:para_eq}
do not possess jumps, but they exhibit so-called layers. Layers are very thin regions 
where the norm of the gradient of the solution is very large. In the 
convection-dominated regime, the width of layer regions is 
much smaller than the affordable mesh width, apart from special cases when
anisotropic layer-adapted meshes can be constructed. Hence, in general,
layers cannot be resolved. Standard discretizations, like the Galerkin finite element 
method or central finite differences, cannot cope with this situation. In
general, numerical solutions computed with such discretizations are globally polluted with spurious 
oscillations. A well-known remedy consists in using so-called stabilized discretizations. 

Finite element methods are a popular approach for discretizing spatial 
derivatives. Major reasons include, but are not limited to, that unstructured 
meshes can be used easily, such that domains with complicated
boundaries can be coped with, and that for many problems they allow an error 
analysis. In a nutshell, finite element methods start with a weak formulation of the PDE, 
replace the infinite-dimensional function spaces with finite-dimensional ones, 
usually consisting of piecewise polynomial functions, and they might approximate, modify
or extend the forms (functionals, bilinear forms etc.) of the weak formulation. This procedure 
does not pay attention to physical consistency. The situation is different for other approaches, 
like finite volume methods, where a goal of the discretization process is to transfer 
conservation properties from the continuous to the discrete equation. However, in view of the 
attractive features of finite element methods, there has been a great interest in 
studying to which extent they lead to physically consistent discretizations and, in case 
of unsatisfactory findings, in developing modifications that possess the desired physical
consistency. 

The goal of the present paper consists in providing a survey on finite element 
methods that satisfy local or global DMPs for linear elliptic or parabolic 
problems. To keep the presentation focussed on the DMPs, other properties of 
the respective methods, like results from the finite element convergence theory,
will be discussed only in the form of brief comments.
On the one hand, many proofs concerning the DMPs use  
just basic tools from linear algebra and they will be presented such that main 
ideas of the numerical analysis become clear. But on the other hand, since this
survey is intended also for an audience without special knowledge in the 
mathematical analysis of the finite element method, it is referred to the 
literature for some other proofs, in particular for those which require many 
technical steps. Although the considered problems \eqref{eq:elli_eq} and
\eqref{eq:para_eq} are linear, both linear as 
well as nonlinear finite element methods for 
their discretization have been proposed. A nonlinear method contains stabilization terms 
whose parameters depend on the numerical solution. That such methods can be suitable 
becomes clear from the above described form of the solution: there are layers
and gently varying 
parts in the solution and an adequate discretization should treat both parts
differently. 

After formulating the steady-state problem and general notations
in Section~\ref{Sec:model}, the following
Section~\ref{Sec:general_DMP} will introduce general results concerning the DMP for both linear 
and nonlinear discretizations. Then, several sections follow that consider
discretizations of the steady-state
problem. First, problems without convection, in particular the Poisson problem, will be 
discussed in Section~\ref{Sec:cd_linear_methods}. Then, linear discretizations and 
finally nonlinear discretizations of convection-diffusion-reaction problems will be 
reviewed in Sections~\ref{Sec:cd_linear_methods_cd}
and~\ref{Sec:cd_nonlinear_methods}, respectively.
The theoretical considerations are illustrated by numerical results in
Section~\ref{sec:num_exam}.
In all these sections, only discretizations with conforming piecewise linear ($\bbP_1$) finite 
elements are considered, since most of the literature is for this case. 
Methods for parabolic problems, 
and $\bbP_1$ finite elements in space, will be reviewed in Section~\ref{sec:tcd}. 
The survey reveals that many finite element methods that satisfy the DMP for $\bbP_1$ finite 
elements transferred ideas from finite volume methods, like upwind techniques or the 
consideration of fluxes. Finite elements different than $\bbP_1$ are the topic
of Section~\ref{sec:other_fems}. The available results for the satisfaction of
the DMP
for other $H^1(\Omega)$-conforming finite elements, often even only
for the Poisson problem, pose usually very restrictive requirements on the
shape of the mesh cells, or they are even negative. Thus, it turns out that the restriction to the
$\bbP_1$ finite element 
in the literature (and the previous sections) has mathematical reasons. 
In addition, non-conforming finite elements are discussed. 
Then,  Section~\ref{sec:hyperbolic_problems} provides brief comments
on methods that satisfy the DMP for hyperbolic conservation laws.
Finally, a summary and an outlook are presented in Section~\ref{sec:summary}.

\section{The steady-state model problem, general notations}\label{Sec:model}

Let $\Omega \subset \bbR^d$, $d\in\{2,3\}$, be a bounded domain with polygonal
resp.~polyhedral and Lipschitz continuous boundary $\partial\Omega$. 
For a domain $D\subset\Omega$ we denote by $W^{m,p}(D)$ the space of functions 
in $L^p(D)$ with weak derivatives up to order $m$ belonging to $L^p(D)$, with 
the usual convention $W^{0,p}(D)=L^p(D)$. The notation $W^{m,p}_0(D)$ denotes 
the closure of $C^\infty_0(D)$ in $W^{m,p}(D)$. If $p=2$ and $m>0$, the usual 
notations $H^m(D)$ and $H_0^m(D)$ are used instead of $W^{m,p}(D)$ and
$W^{m,p}_0(D)$, respectively. The norm (seminorm) in $W^{m,p}(D)$ is denoted by 
$\|\cdot\|_{m,p,D}^{}$ ($|\cdot|_{m,p,D}^{}$), and whenever $p=2$, the index 
$p$ will be dropped from the notation, this is,
$\|\cdot\|_{m,D}^{}=\|\cdot\|_{m,2,D}^{}$. The inner product in $L^2(D)$ or
$L^2(D)^d$ is denoted by
$(\cdot,\cdot)_D^{}$, and the subindex will be dropped if $D=\Omega$.
The Euclidean norm of a vector is denoted by $|\cdot|$.
Finally, for
a number $a\in \bbR$, we define its positive and negative parts as follows:
\begin{equation*}
a^+:=\max\{ a,0\}\ge0 \qquad\textrm{and}\qquad a^-:= \min\{a, 0\}\le 0\,,
\end{equation*}
and the same notation is used to define the positive and negative parts of a real-valued function.

\subsection{The steady-state model problem}\label{Sec:steady_model}

Defining a characteristic length scale and a characteristic scale of the sought quantity, the 
steady-state equation \eqref{eq:elli_eq} can be transformed to a dimensionless problem, where 
we use for simplicity the same notations:
Find $u: \overline\Omega \to\mathbb{R}$ such that
\begin{equation}\label{steady-strong}
\begin{array}{rcll}
-\varepsilon\Delta u + \bb\cdot\nabla u + \sigma u &=& f & \textrm{in}\;\Omega\,, \\
u&=& g &\textrm{on}\;\partial\Omega\,. 
\end{array}
\end{equation}
For simplifying the following presentation, we will suppose that 
$\varepsilon>0$ and $\sigma\ge0$ are constants and  that $\bb$ is solenoidal.

Let $\bb \in W^{1,\infty}(\Omega)^d$, $f\in L^2(\Omega)$, 
and $g\in H^{1/2}(\partial\Omega)$, 
then the weak formulation of \eqref{steady-strong} reads as follows: Find $u\in H^1(\Omega)$ such that
$u|_{\partial\Omega}^{}=g$ and
\begin{equation}\label{steady-weak}
a(u,v)=(f,v)\qquad\forall\, v\in H^1_0(\Omega)\,,
\end{equation}
where $a(\cdot,\cdot)$ is the bilinear form given by
\begin{equation}\label{bilinear-a}
a(u,v)=\varepsilon\,(\nabla u,\nabla v)+(\bb\cdot\nabla u+\sigma u,v)\,.
\end{equation}
Under the stated assumptions on the smoothness of the data, the existence and uniqueness of a solution 
of \eqref{steady-weak} can be concluded from the Lax--Milgram theorem. The
weak maximum principle for 
a sufficiently regular solution reads as follows, e.g., see \cite[Chapter~3.1]{GT01}
or \cite[Chapter~6.4.1]{Eva10}.

\begin{theorem}[Weak maximum principle]
\label{thm:cd_weak_mp}Let  $u \in C^2(\Omega) \cap C(\overline{\Omega})$. Then
\[
\begin{array}{lcl}
-\varepsilon\Delta u + \bb\cdot\nabla u + \sigma u  \le  0 \quad\mbox{in } \Omega
\quad &\Longrightarrow &\quad
\displaystyle \max_{\bx \in \overline{\Omega}} u(\bx) \le\max_{\bx \in \partial\Omega} u^+(\bx),\\
-\varepsilon\Delta u + \bb\cdot\nabla u + \sigma u \ge 0 \quad\mbox{in } \Omega
\quad &\Longrightarrow &\quad
\displaystyle  \min_{\bx \in \overline{\Omega}} u(\bx) \ge \min_{\bx \in \partial\Omega} u^-(\bx).
\end{array}
\]
If $\sigma = 0$, then 
\[
\begin{array}{lcl}
-\varepsilon\Delta u + \bb\cdot\nabla u \le 0 \quad\mbox{in } \Omega
\quad &\Longrightarrow& \quad\displaystyle
\max_{\bx \in \overline{\Omega}} u(\bx) = \max_{\bx \in \partial\Omega} u(\bx),\\
-\varepsilon\Delta u + \bb\cdot\nabla u \ge 0 \quad\mbox{in } \Omega
\quad &\Longrightarrow &\quad\displaystyle
\min_{\bx \in \overline{\Omega}} u(\bx) = \min_{\bx \in \partial\Omega} u(\bx).
\end{array}
\]
\end{theorem}

\subsection{Triangulations and finite element spaces}\label{Sec:triangulations}

We denote by $\{\calT_h^{}\}_{h>0}^{}$ a family of  conforming and regular simplicial triangulations 
of $\Omega$ consisting of mesh cells 
$K$. Note that each mesh cell is the image of a fixed reference cell $\hat
K$ via an affine map. We use the notion of facet to denote an edge in 2d or a face
in 3d.
Let $h_G^{}= \mbox{diam}(G)$ be the diameter of a set $G$ and $h=\max \{h_K:K\in\calT_h^{}\}$.
For a mesh $\calT_h^{}$, the following notations are used:
\begin{enumerate}[leftmargin=*,label=\roman*)]
\item[$-$] internal vertices: $\{\bx_1^{},\ldots,\bx_M^{}\}$,
vertices on the boundary: $\{\bx_{M+1}^{},\ldots,\bx_N^{}\}$,
\item[$-$]  set of internal facets: $\calF_I^{}$, set of boundary facets:
$\calF_\partial^{}$, set of all facets:
$\calF_h^{}= \calF_I^{}\cup \calF_\partial^{}$,
\item[$-$] set of internal edges: $\calE_I^{}$, set of boundary edges:
$\calE_\partial^{}$, set of all edges:
$\calE_h^{}=\calE_I^{}\cup\calE_\partial^{}$, 
\item[$-$] for $K\in\calT_h^{}, F\in\calF_h^{}$, and a vertex $\bx_i^{}$, 
we define the sets
\[
\begin{array}{rclrcl}
\calF_K^{} &=&\{ F\in\calF_h^{}: F\subset K\}\,,\quad&
\calF_i^{} &=& \{ F\in \calF_h^{}: \bx_i^{}\in F\}\,,\\
\calE_K^{}&=& \{ E\in\calE_h^{}: E\subset K\}\,,\quad&
\calE_F^{} &=& \{E\in\calE_h^{}:E\subset F\}\,,
\end{array}
\]
\item[$-$] for $K\in\calT_h^{}, F\in\calF_h^{}$,  $E\in\calE_h^{}$, and a 
vertex $\bx_i^{}$, we define the following subsets of~$\overline\Omega$
\[
\begin{array}{rclrcl}
\omega_K^{} &=&\cup\{ K'\in\calT_h^{}:K\cap K'\not= \emptyset\}\,,\ &
\omega_F^{} &=&\cup\{ K\in\calT_h^{}:F\subset K\}\,,\\
\tilde\omega_F^{} &=&\cup\{ K\in\calT_h^{}:K\cap F\not= \emptyset\}\,,\ &
\omega_E^{} &=&\cup\{ K\in\calT_h^{}:E\subset K\}\,,\\
\omega_i^{} &=&\cup\{ K\in\calT_h^{}:\bx_i^{}\in K\}\,,&
\end{array}
\]
\item[$-$] for a vertex $\bx_i^{}$, we define the set of indices corresponding to
neighbor vertices by
\end{enumerate}
\begin{equation}
S_i^{} = \{j\in\{1,\ldots,N\}\setminus\{i\}\,:\,\bx_i^{}\;\textrm{and}\;\bx_j^{}\; \textrm{are endpoints of }\;E\in\calE_h^{}\}\,,\label{eq:def_S_i}
\end{equation}
\begin{enumerate}[leftmargin=*,label=\roman*)]
\item[$-$] for a facet $F\in\calF_I^{}$, we denote the jump of a function
across $F$ by $\llbracket\cdot\rrbracket_F^{}$. The orientation of the jump is
irrelevant, but fixed.
\end{enumerate}
Note that from the regularity of the triangulations a minimal angle condition follows, e.g., 
see \cite[Section~4.3]{BKK20}. In particular, the number of mesh cells in 
$\omega_K$, $\omega_E$, and $\omega_i$ is bounded
uniformly for all $K$, $E$, $i$, and $h$. In addition, the mesh regularity implies that there exists a positive constant 
$\rho$ such that
\begin{equation}\label{BE-4}
   h_K\le\rho\,h_F\qquad\forall\,\,K\subset\tilde\omega_F\,.
\end{equation}

\begin{figure}[t!]
\centerline{\includegraphics[width=0.5\textwidth]{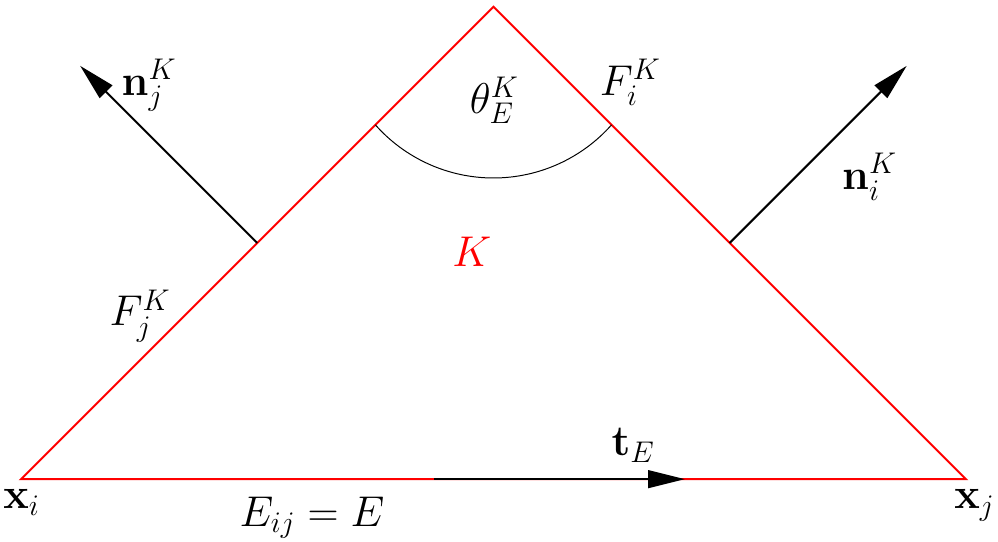}}
\caption{Notations for a triangle.}\label{fig:notations_tria}
\end{figure}

Let $\bx_i^{},\bx_j^{}$ be two vertices 
that are connected by an edge $E_{ij}^{}\in\calE_h^{}$ (or, simply $E$  when
there is no possible confusion) and $K\subset\omega_{E_{ij}^{}}^{}$, then, compare 
Figure~\ref{fig:notations_tria} for the two-dimensional situation, 
\begin{enumerate}[leftmargin=*,label=\roman*)]
\item[$-$] $F_i^K$ and $F_j^K$ are the facets of $K$ opposite $\bx_i^{}$ and
$\bx_j^{}$, respectively, with outer unit normals $\boldsymbol{n}_i^K$ and $\boldsymbol{n}_j^K$, respectively,
\item[$-$] $\theta_{E}^K$ is the angle formed by $F_i^K$ and $F_j^K$, or, more
precisely, $\theta_E^K$ is the dihedral angle given by (cf. \cite{BKKS09}))
\end{enumerate}
\begin{equation}\label{eq:dihedral_angles}
\cos\theta_E^K=-\boldsymbol{n}_i^K\cdot\boldsymbol{n}_j^{K}\,,
\end{equation}
\begin{enumerate}[leftmargin=*,label=\roman*)]
\item[$-$] $\kappa_E^K=F_i^K\cap F_j^K$; when $d=2$, we will adopt the convention $|\kappa_E^K|=1$, 
\item[$-$] $\bt_E^{}=(\bx_j^{}-\bx_i^{})/|\bx_j^{}-\bx_i^{}|$, where the
orientation of this tangent vector is irrelevant, but fixed,
\item[$-$] $\delta_E^{}v:=v(\bx_j^{})-v(\bx_i^{})$ for any function
$v\in C^0(\overline{\Omega})$ if the tangent vector $\bt_E^{}$ points from
$\bx_i$ to $\bx_j$, and $\delta_E^{}v:=v(\bx_i^{})-v(\bx_j^{})$ in the other
situation.
\end{enumerate}

Whether or not a discretization satisfies a DMP might depend on properties of the 
underlying mesh or family of meshes. Some relevant properties in two and three dimensions are defined 
next. 

\begin{definition}[Properties of meshes]
\label{Def:conditions-on-mesh}A mesh $\calT_h^{}$ will be said to be connected 
if, for any two vertices $\bx_i^{},\bx_j^{}$, there exists a path 
$j_0^{},\ldots,j_s^{}$ such that 
$E_{ij_0^{}}^{},E_{j_0^{}j_1^{}}^{},\ldots,E_{j_s^{}j}^{}$ are all edges in 
$\calE_h^{}$. In addition, the mesh $\calT_h^{}$ will be said to be:
\begin{enumerate}[leftmargin=*,label=\roman*)]
\item[$-$] weakly acute: if every internal dihedral angle  $\theta$ of the mesh satisfies $\theta\le \frac{\pi}{2}$,
\item [$-$] of  Xu--Zikatanov (XZ) type (cf. \cite{XZ99}): if, for every 
$E\in\calE_I^{}$, the following holds
\end{enumerate}
\begin{equation}\label{XZ-criterion}
\sum_{K\subset\omega_E} |\kappa_E^K|\cot\theta_E^K\ge 0\,,
\end{equation}
\begin{enumerate}[leftmargin=*,label=\roman*)]
\item[$-$] of Delaunay type: if the interior of the circumscribed sphere of 
any simplex from the mesh $\calT_h^{}$ does not contain any vertex of 
$\calT_h^{}$.
\end{enumerate}
\end{definition}

For $d=2$, the definition of a Delaunay mesh can be equivalently stated as 
follows: for every $E=K\cap K'\in\calE_I^{}$ there holds
\begin{equation*}
\theta_{E}^K+\theta_{E}^{K'}\le \pi\,.
\end{equation*}
In two dimensions, the XZ-criterion and the Delaunay property are equivalent.

\begin{definition}[Strictly acute and average acute families of meshes]  
A mesh family $\{\calT_h^{}\}_{h>0}^{}$ will be said to be strictly acute if 
there is a constant $\delta>0$ independent of $h$ such that every internal 
dihedral angle $\theta$ of any of the meshes satisfies 
\begin{equation}\label{eq:strictly_acute}
   \theta\le\frac{\pi}{2}-\delta\,.
\end{equation}
In two dimensions, a family $\{\calT_h^{}\}_{h>0}^{}$ will be said to be
average acute if, for every $h>0$ and every edge
$E=K\cap K'\in\calE_I^{}$, the following holds:
\begin{equation}\label{FEM-prelim-1}
\theta_{E}^K+\theta_{E}^{K'}\le \pi-\delta\,,
\end{equation}
where $\delta>0$ is independent of $h$.
\end{definition}

As already mentioned, most discretizations discussed in this survey are based on 
continuous piecewise linear finite elements. The corresponding finite element spaces and interpolation operators for this case will be defined next. 
Associated with the vertices $\{\bx_1^{},\ldots,\bx_N^{}\}$, the standard 
continuous piecewise linear basis functions 
$\phi_1^{},\ldots,\phi_N^{}$ are given by the property 
$\phi_i^{}(\bx_j^{})=\delta_{ij}^{}$ for $i,j\in\{1,\ldots,N\}$. Then, the 
corresponding conforming finite element spaces are
\begin{equation}\label{fem-space}
V_h^{}:=\mbox{span}\{\phi_1^{},\ldots,\phi_N^{}\}\quad\textrm{and}\quad
V_{h,0}^{}:= \mbox{span}\{\phi_1^{},\ldots,\phi_M^{}\}\,.
\end{equation}
Associated with $V_h^{}$,  the Lagrange interpolation operator is defined by 
\begin{equation*}
i_h^{}:C^0(\overline\Omega)\to  V_h^{}\,, \quad 
v\mapsto  i_h^{}v=\sum_{i=1}^Nv(\bx_i^{})\phi_i^{}\,.
\end{equation*}
We will also use the symbol $i_h^{}$ to interpolate functions with domain in the boundary of $\Omega$, this is,
$i_h^{}g=\sum_{i=M+1}^N g(\bx_i^{})\phi_i^{}$. 

\subsection{Finite element matrices}\label{Sec:FEM-Matrices}

In this section, the main finite element matrices are introduced.
The diffusion matrix $\bbAd$, the  convection matrix $\bbAc$, and the reaction matrix $\bbMc$, 
which is also called consistent mass matrix, are defined by 
\begin{alignat}{10}\label{def-dif-matrix}
\bbAd &= && (\fl_{ij}^{})_{i,j=1}^N \quad &&\textrm{where}\ &&  \fl_{ij}^{}=(\nabla\phi_j^{},\nabla\phi_i^{})\quad &&  \textrm{for}\; i,j=1,\ldots,N\,, \\
\bbAc &=&&  \,(\cc_{ij}^{})_{i,j=1}^N &&\textrm{where}&& \cc_{ij}^{}=(\bb\cdot\nabla\phi_j^{}, \phi_i^{})&& \textrm{for}\; i,j=1,\ldots,N\,, \label{def-conv-matrix}\\
\bbMc &=&&  \,(\mm_{ij}^{})_{i,j=1}^N &&\textrm{where} && \mm_{ij}^{}=(\phi_j^{}, \phi_i^{})&&\textrm{for}\; i,j=1,\ldots,N\,. \label{def-mass-matrix}
\end{alignat}
The entries of the matrices can be written as a sum of local entries, e.g.,
\[
\fl_{ij}^{} = \sum_{K\subset\omega_i \cap \omega_j} \fl_{ij}^K \quad \mbox{with } \fl_{ij}^K = (\nabla\phi_j,\nabla\phi_i)_K,
\]
and analogously for $\cc_{ij}^{}$ and $\mm_{ij}^{}$.

In the derivations made in the coming sections, having exact formulae for the diffusion and consistent mass matrices will be of much use. 
A basic tool in the derivations below is a formula relating the gradient of the barycentric coordinates
and the normal outward vector to $K$. Since the basis function $\phi_i^{}|_K^{}$ vanishes on $F_i^K$, 
its derivative in any direction tangent to $F_i^K$ vanishes. So, $\nabla\phi_i^{}|_K$ is 
proportional to the unit normal $\bn_i^K$. Consider the height vector $\boldsymbol h_i^{}$ from 
$F_i^K$ to $\bx_i$. This vector is parallel to $\bn_i^K$, pointing in the opposite direction, 
and the derivative of $\phi_i^{}|_K^{}$ in the direction of $\boldsymbol h_i^{}$ is the constant $1/|\boldsymbol h_i|$.
Hence, using the formula for the volume of the simplex $K$ leads to 
\begin{equation}\label{grad-basis-functions}
\nabla\phi_i^{}|_K= - \frac{1}{|\boldsymbol h_i|} \boldsymbol{n}_i^K = 
-\frac{|F_i^K|}{d|K|}\boldsymbol{n}_i^K\,.
\end{equation}
So, in view of \eqref{eq:dihedral_angles}, the local diffusion matrix is given by
\begin{equation}\label{basic-expression}
\fl_{ij}^{K} = (\nabla\phi_j^{},\nabla\phi_i^{})_K^{} = |K|\frac{|F_j^K|\,|F_i^K|}{d^2|K|^2}\boldsymbol{n}_j^K\cdot \boldsymbol{n}_i^K
= -\frac{|F_j^K|\,|F_i^K|}{d^2|K|}\cos\theta_E^K.
\end{equation}

Concerning the mass matrix and using the formula for the integral of a product of barycentric coordinates, see, e.g., \cite{VS18} where this is proven
in any space dimension, one gets
\begin{equation}\label{local-mass-exact}
m_{ij}^{K}= 
\begin{cases} \displaystyle \frac{2|K|}{(d+1)(d+2)} & i=j,
\\[1em]
\displaystyle \frac{|K|}{(d+1)(d+2)} & \textrm{else}\,. \end{cases}
\end{equation}

Both in the steady-state and time-dependent situations, mass lumping is a 
widely used technique to discretize terms without spatial derivatives. The 
derivation of mass lumping starts with the construction of a dual mesh, which 
is a technique from finite volume methods. 
For each node $\bx_i$, all mesh cells $K\subset\omega_i$ are considered. 
In each mesh cell, 
a polyhedral subset with volume $|K|/(d+1)$ assigned to $\bx_i$ is constructed.
The vertices of this subset are $\bx_i$, the barycenter of $K$, midpoints of 
edges of $K$ containing $\bx_i$, and, if $d=3$, also the barycenters of faces
of $K$ containing $\bx_i$.
Now, the dual mesh cell $D_i$ is defined by the union of these subsets from 
all $K\subset\omega_i$. Consequently, one has 
\begin{equation*}
|D_i^{}|=\frac{|\omega_i^{}|}{d+1}\,.
\end{equation*}
Piecewise constant basis functions, given by 
\begin{equation}\label{FEM-upwind-2}
\psi_i^{}(\bx)=\left\{\begin{array}{cl} 1 & \textrm{if}\; \bx\in D_i^{}\,,\\
0 & \textrm{else}\, , \end{array}\right. \quad  i=1,\ldots,N,
\end{equation}
are associated with this dual mesh.
With the help of these functions, the following lumping operator is defined 
\begin{equation}\label{FEM-upwind-3}
\calL\ : \ C(\overline\Omega)\to L^2(\Omega)\,,\quad
v \mapsto \calL v=\sum_{i=1}^Nv(\bx_i^{})\psi_i^{}\,.
\end{equation}
In addition, the lumped $L^2(\Omega)$ inner product 
$(\cdot,\cdot)_h:C(\overline\Omega)\times C(\overline\Omega) \to \mathbb R$ 
is given by 
\begin{equation}\label{eq:lumped_inner_prod}
(f,g)_h^{}=(\calL f,\calL g)\,.
\end{equation}
Since $\{\psi_i^{}\}_{i=1}^N$ is an orthogonal set in $L^2(\Omega)$ 
and $(\psi_i^{},\psi_i^{})=|D_i^{}|$, one obtains
\begin{equation*}
(f,g)_h^{}=\sum_{i,j=1}^Nf(\bx_j^{})g(\bx_i^{})(\psi_j^{},\psi_i^{})=\sum_{i=1}^N |D_i|f(\bx_i^{})g(\bx_i^{})\,.
\end{equation*}
Using the lumped inner product, the following seminorm is induced in $C(\overline\Omega)$, which is a norm 
in  $V_h^{}$, 
\begin{equation*}
|f|_h^{} := (f,f)_h^{1/2} = \left( \sum_{i=1}^N  |D_i^{}|\,|f(\bx_i^{})|^2\right)^{1/2}\,.
\end{equation*}
Finally, the lumped mass matrix, which is a diagonal matrix, is defined as follows
\begin{equation}\label{def-lumpedmass-matrix}
\bbMl=(\tilde{\mm}_{ij}^{})_{i,j=1}^N\quad\textrm{where}\quad \tilde{\mm}_{ij}^{} = (\phi_j^{},\phi_i^{})_h^{} =(\calL\phi_j^{},\calL\phi_i^{})= |D_i^{}|\delta_{ij}^{}\,.
\end{equation}
Utilizing an exact quadrature rule for linears and the fact that the basis 
functions of $V_h$ form a partition of unity yields 
\begin{equation}\label{eq:mass_matrix_l_c}
\tilde{\mm}_{ii}^{} = |D_i^{}| = \sum_{K\subset\omega_i^{}} \frac{|K|}{d+1} = \sum_{K\subset\omega_i^{}} (1,\phi_i^{})_K^{} = (1,\phi_i^{}) = \sum_{j=1}^N (\phi_j^{}, \phi_i^{}) = \sum_{j=1}^N m_{ij}^{}\,.
\end{equation}
So, the lumped mass matrix can be computed directly from the consistent mass matrix, without the need to build the dual mesh. 

\section{General results on DMP satisfying discretizations}\label{Sec:general_DMP}

This section provides conditions for the satisfaction of local and global DMPs 
that are based on special properties of matrices for general linear discrete 
problems, and of nonlinear forms for general nonlinear discretizations. The 
presentation of the theory for linear discretizations is based on the concept 
of matrices of non-negative type, instead on the traditional approach with 
monotone matrices or, more special, M-matrices. This concept enables also the 
consideration of local DMPs. 

\subsection{Linear discretizations} 
\label{Sec:DMP_linear_disc}

Let a matrix  $(a_{ij}^{})^{i=1,\dots,M}_{j=1,\dots,N}\in
\mathbb{R}^{M\times N}$ and real numbers $f_1^{},\ldots, f_M^{}, g_1^{},\ldots,
g_{N-M}^{}$ with $M<N$ be given. A linear discretization leads to a system of 
linear algebraic equations of the following form: 
Find $\bu=(u_1^{},\ldots,u_N^{})^T\in\mathbb{R}^N$ such that
\begin{align}
\sum_{j=1}^N a_{ij}^{}u_j^{} &= f_i^{}\qquad\textrm{for}\; i=1,\ldots, M\,,\label{system-1}\\
u_i^{} &= g_{i-M}^{}\qquad\textrm{for}\; i=M+1,\ldots, N\,. \label{system-2}
\end{align}

\begin{remark}
The system matrix of the system \eqref{system-1}-\eqref{system-2} is of the form
\begin{equation}\label{system-matrix}
\bbA =\begin{pmatrix} \bbAI & \bbAB\\ \mathbb O & \mathbb I \end{pmatrix}\,,
\end{equation}
where $\bbAI \in \mathbb R^{M\times M}$ is the matrix associated with the internal (or non-Dirichlet) degrees of freedom, 
$\bbAB \in \mathbb R^{M\times (N-M)}$ is the matrix that couples the 
boundary values to the values in the interior of the domain, $\mathbb I \in
\mathbb R^{(N-M)\times (N-M)}$ is the identity matrix and $\mathbb O \in
\mathbb R^{(N-M)\times M}$ a matrix consisting of zeros. In what follows,
$\bbA$ will always denote the matrix given by \eqref{system-matrix}.
\hspace*{\fill}$\Box$\end{remark}

\begin{definition}[Matrix of non-negative type]A matrix
$(a_{ij}^{})^{i=1,\dots,m}_{j=1,\dots,n}\in\mathbb{R}^{m\times n}$ 
($m,n\in\mathbb{N}$) will be said to be of  non-negative type if
\begin{align}
   a_{ij}^{}&\le 0 \qquad \forall\,\,i\neq j,\, 1\le i\le m,\,1\le j\le n\,,
   \label{NN-1}\\
   \sum_{j=1}^n\,a_{ij}^{} &\ge 0\qquad \forall\,\,1\le i\le m\,.\label{NN-2}
\end{align}
\end{definition}

One should notice that the notion of a matrix of non-negative type must not be 
confused with the notion of a non-negative matrix as it is studied, e.g., 
in \cite[Chapter~2]{Var00}.

\begin{remark} In some cases, e.g., when $\sigma=0$ in \eqref{steady-strong}, the matrix 
$\mathbb{A}$ will satisfy a stronger property than \eqref{NN-2}, namely
\begin{equation}
\sum_{j=1}^N a_{ij}^{} = 0\qquad \forall\; 1\le i\le M\,.\label{zero-row-sum}
\end{equation}
With this property, it will be possible to derive stronger statements for the DMP than 
with \eqref{NN-2}.
\hspace*{\fill}$\Box$\end{remark}

The next result is a local version of the results given in 
\cite{Ciarlet70,CR73}.
 
\begin{theorem}[Local DMP in the case of matrices of non-negative type]
\label{thm:algebraic-local-DMP}Let $a_{ii} >0$ for $i=1,\ldots,M$.
Then, any possible solution of \eqref{system-1}-\eqref{system-2} 
satisfies
\begin{equation}\label{eq:algebraic-local-DMP1}
f_i^{}\le 0\quad \Longrightarrow\quad u_i^{}\le \max\limits_{j\neq i,
a_{ij}^{}\not= 0} u_j^{+}\,, \qquad\qquad
f_i^{}\ge 0\quad\Longrightarrow\quad u_i^{}\ge \min\limits_{j\neq i, a_{ij}^{}\not= 0} u_j^{-}
\end{equation}
for all $i=1,\dots,M$ if and only if  $\mathbb{A}$ is of non-negative type.
The implications
\begin{equation}\label{eq:algebraic-local-DMP2}
f_i^{}\le 0\quad\Longrightarrow\quad u_i^{}\le \max\limits_{j\neq i,
a_{ij}^{}\not= 0} u_j^{}\,, \qquad\qquad
f_i^{}\ge 0\quad\Longrightarrow\quad u_i^{}\ge \max\limits_{j\neq i, a_{ij}^{}\not= 0} u_j^{}\,
\end{equation}
hold true for all $i=1,\dots,M$ if and only if  $\mathbb{A}$ is of non-negative 
type and satisfies in addition \eqref{zero-row-sum}.
\end{theorem}

\begin{proof}Consider any $i\in\{1,\dots,M\}$ and let $f_i\le0$. If
$\mathbb{A}$ is of non-negative type, then it follows from \eqref{system-1}, 
\eqref{NN-1}, and \eqref{NN-2} that
\begin{equation*}
   a_{ii}\,u_i=f_i-\sum_{j\neq i}\,a_{ij}\,u_j\le\sum_{j\neq i}\,(-a_{ij})\,
   \max\limits_{j\neq i, a_{ij}^{}\not= 0} u_j^{+}\le a_{ii}\,
   \max\limits_{j\neq i, a_{ij}^{}\not= 0} u_j^{+}\,,
\end{equation*}
which implies \eqref{eq:algebraic-local-DMP1}. If, in addition,
\eqref{zero-row-sum} holds, then \eqref{eq:algebraic-local-DMP2} follows from
\begin{equation*}
   a_{ii}\,u_i=f_i-\sum_{j\neq i}\,a_{ij}\,u_j\le\sum_{j\neq i}\,(-a_{ij})\,
   \max\limits_{j\neq i, a_{ij}^{}\not= 0} u_j=a_{ii}\,
   \max\limits_{j\neq i, a_{ij}^{}\not= 0} u_j\,.
\end{equation*}
The statements for $f_i\ge0$ follow analogously. The necessity of the conditions
on $\mathbb{A}$ can be proved by constructing appropriate counterexamples, see
\cite[Appendix]{BJK16}. 
\end{proof}

In the context of numerical approximation of PDEs,  Theorem~\ref{thm:algebraic-local-DMP} 
implies a local DMP.  It should be emphasized that 
for the local DMP the invertibility of $\mathbb{A}$ is not a necessary condition. In 
particular, it holds also for convection-diffusion equations \eqref{steady-strong}, without reactive
term, and with pure Neumann  boundary conditions as long as their discretization leads
to a system matrix of non-negative type and there is a solution.

Next, the global version of the DMP is shown. Its proof is based on a technique
developed in \cite{Knobloch10} and can be considered as a generalization of
\cite[Theorem~3]{Ciarlet70}. 

\begin{theorem}[Global DMP in the case of matrices of non-negative type]\label{thm:algebraic-global-DMP}
Let us suppose that $\mathbb{A}$ is of non-negative type and that the matrix $\bbAI = (a_{ij}^{})_{i,j=1}^M$ is
invertible. Then, system \eqref{system-1}-\eqref{system-2} possesses a unique solution. 
This solution satisfies
\begin{equation}\label{global-DMP-1}
\begin{array}{l}
f_i^{}\le 0\quad\forall\;\,i=1,\ldots,M\quad\Longrightarrow \quad\max\limits_{i=1,\ldots,N} u_i^{}\le \max\limits_{j=M+1,\ldots,N} u_j^{+}\,, \\
f_i^{}\ge 0\quad\forall\;\,i=1,\ldots,M\quad\Longrightarrow\quad \min\limits_{i=1,\ldots,N} u_i^{}\ge \min\limits_{j=M+1,\ldots,N} u_j^{-}\,.
\end{array}
\end{equation}
In addition, if $\mathbb{A}$ satisfies \eqref{zero-row-sum}, the following holds
\begin{equation}\label{global-DMP-2}
\begin{array}{l}
f_i^{}\le 0\quad\forall\;\,i=1,\ldots,M\quad\Longrightarrow \quad\max\limits_{i=1,\ldots,N} u_i^{}=\max\limits_{j=M+1,\ldots,N} u_j^{}\,,\\
f_i^{}\ge 0\quad\forall\;\,i=1,\ldots,M\quad\Longrightarrow \quad\min\limits_{i=1,\ldots,N} u_i^{}=\min\limits_{j=M+1,\ldots,N} u_j^{}\,.
\end{array}
\end{equation}
\end{theorem}

\begin{proof} Inserting the values from \eqref{system-2} in \eqref{system-1} leads 
to a linear system of equations for $u_1,\ldots,u_M$ with the matrix $\bbAI$.
From the assumed invertibility of this matrix, the existence of a unique solution of 
\eqref{system-1}-\eqref{system-2} follows.

Next, the first statement of \eqref{global-DMP-1} will be shown. The second 
statement of \eqref{global-DMP-1} follows by changing the signs of $\bu$ and of
the right-hand side of \eqref{system-1}-\eqref{system-2}. Let 
\begin{equation*}
s=\max_{i=1,\ldots,N}u_i^{}\quad\textrm{and}\quad J=\{ i\in \{1,\ldots, N\}:u_i^{}=s\}\,.
\end{equation*}
If $s\le 0$, then \eqref{global-DMP-1} holds trivially. So, consider $s>0$ and assume that $J\subset\{1,\ldots,M\}$. It will be shown that
\begin{equation}\label{Glob-Alg-DMP-1}
\exists k\in J\;\textrm{such that}\;\; \mu_k^{}:= \sum_{j\in J}a_{kj}^{} > 0\,.
\end{equation}
Let us suppose that \eqref{Glob-Alg-DMP-1} does not hold. Then, one concludes 
by combining \eqref{NN-1} and \eqref{NN-2} that 
\begin{equation*}
\sum_{j\in J}a_{ij}^{} = 0\qquad\forall\; i\in J\,.
\end{equation*}
Hence, the matrix $(a_{ij}^{})_{i,j\in J}^{}$ is singular because the sum of its columns is zero.
With $(a_{ij}^{})_{i,j\in J}^{}$, also its transposed $(a_{ji}^{})_{i,j\in J}^{}$
is singular. Hence, there exist numbers
$v_i^{}, i\in J$, not all zero, such that
\begin{equation}\label{eq:lem_global_dmp_00}
\sum_{i\in J} a_{ij}^{}v_i^{}=0\qquad\forall\; j\in J\,.
\end{equation}
In addition, applying  that  $\mathbb{A}$ is of non-negative type one finds that $a_{ij}^{}=0$ for all $i\in J$ and all $j\not\in J$. Using this property, \eqref{eq:lem_global_dmp_00}, and 
defining the vector $\tilde{\bv}=(\tilde{v}_i^{})_{i=1}^M$, where $\tilde{v}_i^{}=v_i^{}$ if $i\in J$,
and $\tilde{v}_i^{}=0$ otherwise, yields
\begin{equation*}
\sum_{i=1}^M a_{ij}^{}\tilde{v}_i^{}=\sum_{i\in J}a_{ij}^{}v_i^{}=0\,,
\end{equation*}
for all $j\in \{1,\ldots,M\}$. This implies that the matrix $\bbAI$ is singular,
which contradicts the hypothesis. So, \eqref{Glob-Alg-DMP-1} holds. 

Denoting now
\begin{equation*}
r = \max_{i\not\in J} u_i^{+}\,,
\end{equation*}
one obtains with $f_i^{}\le 0$ for all $i$, \eqref{NN-1}, and \eqref{NN-2}
\begin{eqnarray*}
s\mu_k^{} & = & \sum_{j\in J}a_{kj}^{}u_j^{} = f_k - \sum_{j\not\in J}a_{kj}^{}u_j^{}
\le - \sum_{j\not\in J}a_{kj}^{}u_j^{}
= \sum_{j\not\in J} (-a_{kj}^{})u_j^{} \le r\sum_{j\not\in J} (-a_{kj}^{})\\
&= &  r\left(\sum_{j=1}^N (-a_{kj}^{})  +  \sum_{j\in J} a_{kj}^{}\right)  \le r\mu_k\,.
\end{eqnarray*}
This implies that $s\le r$, which is a contradiction to the definition of $s$. Hence, $J\cap\{M+1,\ldots,N\}\not=\emptyset$ and \eqref{global-DMP-1} follows. 

The validity of \eqref{global-DMP-2} easily follows from \eqref{global-DMP-1}.
Since \eqref{zero-row-sum} holds, one can add a sufficiently large positive
constant $q>0$  to every $u_i^{}$ in such a way that all components of this
new vector $\tilde\bu$ are positive. Then, the first statement of 
\eqref{global-DMP-1} holds for $\tilde\bu$ without the positive parts, which
implies the first statement of \eqref{global-DMP-2}.
\end{proof}

\begin{remark}\label{rem:necessary_nonsingularity}
If the global DMP \eqref{global-DMP-1} holds and $\bu\in\mathbb{R}^N$ 
is such that $u_{M+1}^{}=\ldots=u_N^{}=0$ and $\buI:=(u_1^{},\ldots,u_M^{})^T$
satisfies $\bbAI\buI=0$, then $\max_{i=1,\ldots,N} u_i^{}\le0$ and 
$\min_{i=1,\ldots,N} u_i^{}\ge0$ so that $\bu=0$. Consequently,
the validity of the global DMP \eqref{global-DMP-1}
implies that the matrix $\bbAI$ is invertible. Thus, this additional assumption
(in comparison to the assumptions of Theorem~\ref{thm:algebraic-local-DMP} for
the local DMP) is necessary.
\hspace*{\fill}$\Box$\end{remark}

\begin{remark}
It is easy to construct a matrix $\bbA$ of non-negative type and a vector 
$\bu=(u_1^{},\ldots,u_N^{})^T$ such that the right-hand side of some of the
implications in Theorem~\ref{thm:algebraic-local-DMP} holds for all
$i=1,\dots,M$ but the corresponding right-hand side in
Theorem~\ref{thm:algebraic-global-DMP} is not satisfied. Thus, a global DMP
cannot be obtained as a consequence of the validity of the corresponding local
DMPs. On the other hand, it can also happen that the global DMP holds but the
local one not since the assumption that $\bbA$ is of non-negative type is not
necessary for the validity of the global DMP.
\hspace*{\fill}$\Box$\end{remark}

\begin{remark}
A situation considered sometimes in the literature is the case of homogeneous 
Dirichlet boundary values. In this case, the proof of 
Theorem~\ref{thm:algebraic-global-DMP} does not require any assumptions on 
the submatrix $\bbAB=(a_{ij}^{})^{i=1,\dots,M}_{j=M+1,\ldots,N}$. However, 
such assumptions are needed in the general case, and consequently considering 
homogeneous Dirichlet boundary conditions is only a particular situation.
\hspace*{\fill}$\Box$\end{remark}

\begin{remark}
From the previous theorems, it follows that both the local and global DMPs are 
satisfied if $\bbA$ is of non-negative type and $\bbAI$ is invertible. Since
$\det\bbA=\det\bbAI$, one observes that $\bbAI$ is invertible if and only if 
$\bbA$ is invertible. Moreover, a direct calculation shows that 
\begin{equation}\label{eq:matrix_inverse}
\bbA =\begin{pmatrix} \bbAI & \bbAB\\ \mathbb O & \mathbb I \end{pmatrix}
\quad \Longleftrightarrow \quad \bbA^{-1} = \begin{pmatrix}  \bbA_{\rm I}^{-1} 
& -\bbA_{\rm I}^{-1} \bbAB\\ \mathbb O & \mathbb I \end{pmatrix}\,.
\end{equation}
In addition, an interesting observation is that the proof of \eqref{global-DMP-1}
allows that $\bbAB=\mathbb O$. Hence, there is no connection between the 
degrees of freedom and the prescribed values on the boundary. In contrast, \eqref{zero-row-sum} in combination 
with the invertibility of $\bbAI$ requires that $\bbAB\neq\mathbb O$.
\hspace*{\fill}$\Box$\end{remark}

As discussed in the previous remark, the invertibility of $\bbAI$ is a 
necessary and sufficient condition for the well-posedness of the discrete 
problem and is also necessary for proving that a method satisfies a global DMP 
(cf.~Remark~\ref{rem:necessary_nonsingularity}). Then, under the assumptions of 
the previous theorems, the matrix $\bbAI$ is of non-negative type (since $\bbA$ 
is) and invertible.  It will be shown in Corollary~\ref{cor:m-matrix}
that these properties imply that the matrix $\bbAI$ belongs to the class of 
M-matrices defined next.

\begin{definition}[M-matrix, monotone matrix] A matrix $\mathbb Q = (q_{ij}^{})_{i,j=1}^n$ is an M-matrix 
if:
\begin{enumerate}[leftmargin=*,label=\roman*)]
\item The off-diagonal entries are non-positive, i.e., 
$q_{ij}^{} \le 0$, $i,j=1,\ldots,n, \ i\neq j$;
\item $\mathbb Q$ is non-singular; and \label{ass:m_matrix_ii}
\item $\mathbb Q^{-1} \ge 0$. \label{ass:m_matrix_iii}
\end{enumerate}
A matrix that satisfies conditions~\ref{ass:m_matrix_ii} and~\ref{ass:m_matrix_iii} is called monotone matrix.
\end{definition}

In the above definition, the condition $\mathbb Q^{-1} \ge 0$ means that all
entries of the matrix $\mathbb Q^{-1}$ are non-negative. In the following,
an analogous notation will be used also for vectors, e.g., $\bv\ge0$ means that 
all entries of the vector $\bv$ are non-negative. 

\begin{remark}\label{rem:monotonicity}
A monotone matrix $\mathbb Q$ can be equivalently characterized by the property 
that, for any $\bv\in\mathbb R^n$, the validity of ${\mathbb Q}\bv\ge0$ implies 
$\bv\ge0$.  Indeed, if this implication holds, then ${\mathbb Q}$ is
non-singular (since ${\mathbb Q}\bv=0$ implies both $\bv\ge0$ and $-\bv\ge0$)
and if $\bv$ is any column of $\mathbb Q^{-1}$, one has ${\mathbb Q}\bv\ge0$
and hence $\bv\ge0$ so that $\mathbb Q^{-1} \ge 0$. On the other hand, if
$\mathbb Q$ is monotone, then ${\mathbb Q}\bv\ge0$ implies that 
$\bv={\mathbb Q}^{-1}{\mathbb Q}\bv\ge0$.
\hspace*{\fill}$\Box$\end{remark}

\begin{theorem}[Equivalence of the monotonicity and the global DMP]
\label{thm:monotonicity-and-global-DMP}Let the row sums of the matrix $\bbA$
be non-negative. Then the global DMP \eqref{global-DMP-1} is satisfied
if and only if $\bbA$ is monotone.
\end{theorem}

\begin{proof}
If the global DMP holds, then, for any $\bv\in\mathbb{R}^N$ satisfying 
$\bbA\bv\ge0$, one has $v_i\ge\min_{j=M+1,\ldots,N} v_j^{-}=0$ for all
$i=1,\dots,N$ so that $\bbA$ is monotone due to 
Remark~\ref{rem:monotonicity}.  Reciprocally, let $\bbA$ be
monotone and let $\bu\in\mathbb{R}^N$ be the solution of
\eqref{system-1}-\eqref{system-2} with $f_i^{}\ge 0$, $i=1,\ldots,M$. Set
$c:=\min_{j=M+1,\ldots,N} u_j^{-}$ and define $\bv\in\mathbb{R}^N$ by 
$v_i=u_i-c$. Since $c\le 0$ and the row sums of $\bbA$ are non-negative, one 
has $\bbA\bv\ge0$. Then the monotonicity of $\bbA$ implies that $\bv\ge0$ and
hence $u_i\ge c$ for $i=1,\dots,N$. Thus the global DMP holds.
\end{proof}

\begin{corollary}[M-matrix property of $\bbA$]
\label{cor:m-matrix}
If the matrix $\bbA$ is invertible and of non-negative type, then both $\bbA$
and $\bbAI$ are M-matrices.
\end{corollary}

\begin{proof}
If $\bbA$ is invertible and of non-negative type, then, according to 
Theorem~\ref{thm:algebraic-global-DMP}, the global DMP \eqref{global-DMP-1} is 
satisfied and $\bbA$ is monotone in view of
Theorem~\ref{thm:monotonicity-and-global-DMP}. Consequently, $\bbA$ is an
M-matrix. In view of \eqref{eq:matrix_inverse}, $\bbAI$ is an M-matrix as well.
\end{proof}

\begin{remark}\label{rem:A-AI}
Using \eqref{eq:matrix_inverse}, it follows immediately that if $\bbA$ is an 
M-matrix \linebreak (monotone matrix) also
$\bbAI$ is an M-matrix (monotone matrix). Conversely, if $\bbAI$ is an M-matrix
(monotone matrix) and $\bbAB \le 0$ (in particular, if $\bbA$ is of
non-negative type), then $\bbA$ is an M-matrix (monotone matrix).
\hspace*{\fill}$\Box$\end{remark}

\begin{remark}
The analysis for linear discretizations was performed purely on the algebraic level. We 
like to emphasize that the results concerning the vector $\bu$ with respect to
the DMP can be 
transferred to the corresponding finite element function only in special cases,
like for the  
$\bbP_1$ finite element. Finite element spaces where such a transfer is not possible are
discussed in Section~\ref{sec:other_fems}.
\hspace*{\fill}$\Box$\end{remark}

\subsection{Nonlinear discretizations} 
\label{Sec:DMP_nonlinear_disc}

In this section we will deal with two types of nonlinear discretizations of 
\eqref{steady-strong} which will be considered in variational forms with the 
$\bbP_1$ finite element spaces \eqref{fem-space}:

\noindent \underline{Type~I:} Find $u_h^{}\in V_h^{}$ such that $u_h^{}|_{\partial\Omega}^{}=i_h^{}g$, and
\begin{equation}\label{nonlinear-generic-1}
a(u_h^{},v_h^{})+ j_h(u_h^{};v_h^{})
= (f,v_h^{})\qquad\forall\, v_h^{}\in V_{h,0}^{}\,,
\end{equation}
where $a(\cdot,\cdot)$ is the bilinear form given by \eqref{bilinear-a}, and
$j_h(\cdot;\cdot)$ is a nonlinear stabilizing term, linear in the second argument.

\noindent \underline{Type~II:} Find $u_h^{}\in V_h^{}$ such that $u_h^{}|_{\partial\Omega}^{}=i_h^{}g$, and
\begin{equation}\label{nonlinear-generic-2}
a(u_h^{},v_h^{})+ d_h^{}(u_h^{};u_h^{},v_h^{}) = (f,v_h^{})\qquad\forall\, v_h^{}\in V_{h,0}^{}\,,
\end{equation}
where $a(\cdot,\cdot)$ is the bilinear form given by \eqref{bilinear-a}, and 
$d_h^{}(\cdot;\cdot,\cdot)$ is nonlinear in the first argument and linear in
the remaining two arguments. We assume that $d_h^{}(\cdot;\cdot,\cdot)$ vanishes if the second 
argument is constant, i.e.,
\begin{equation}\label{eq:d_h_row_sums}
   d_h^{}(w_h^{};1,v_h^{})=0\qquad\forall\,\,w_h^{},v_h^{}\in V_h^{}
\end{equation}
and that,
for all $w_h^{}\in V_h^{}$, the bilinear form
$d_h^{}(w_h^{};\cdot,\cdot)$ is positive semidefinite, i.e.,
\begin{equation}\label{eq:pos_semidef_dh}
   d_h^{}(w_h^{};v_h^{},v_h^{})\ge 0\qquad\forall\,\,w_h^{},v_h^{}\in V_h^{}\,.
\end{equation}

Due to the nonlinear character of \eqref{nonlinear-generic-1} and \eqref{nonlinear-generic-2} the results presented in the last section cannot be applied. We present below two criteria for the satisfaction
of the DMP. In both cases the criteria are related to the following remark: in order to prove the DMP, the only argument used
concerns the entries of the row that corresponds to a node
where an extremum of a discrete solution is encountered. So, to prove the DMP, it is not necessary to
modify every equation, but only those associated with local extrema of a solution $u_h^{}$. Based on this idea, in \cite{BE05} 
a criterion was proposed in order to prove the DMP for a nonlinear 
discretization of Type~I.  Here, we present the following two variants of
this criterion.

\begin{definition}[Strong and weak DMP properties]
\label{def:DMP-criterion}The nonlinear form $j_h(\cdot;\cdot)$ is said to 
satisfy the strong DMP property if the following condition holds: If $u_h^{}$ 
attains a strict local minimum (maximum) at an interior node $\bx_i^{}$, then there 
exist constants $\alpha_F^{}> 0$, $F\in \calF_i^{}$, such that
\begin{equation*}
a(u_h^{},\phi_i^{})+ j_h(u_h^{};\phi_i^{})\le -\sum_{F\in\calF_i^{}}\alpha_F^{}\left|\llbracket \nabla u_h^{}\rrbracket_F^{}\right|\,,
\end{equation*}
(resp. $\ge \sum_{F\in\calF_i^{}}\alpha_F^{}|\llbracket \nabla u_h^{}\rrbracket_F^{}|$).
The form $j_h(\cdot;\cdot)$ is said to satisfy the weak DMP property if the same conclusion
holds under the extra assumption that the local minimum (maximum) satisfies $u_h^{}(\bx_i^{}) <0$ (resp. $u_h^{}(\bx_i^{}) >0$).
\end{definition}

\begin{definition}[Strong and weak DMP properties for non-strict extrema]
\label{def:DMP-criterion-non-strict}The nonlinear form $j_h(\cdot;\cdot)$ is
said to satisfy the strong or weak DMP property for non-strict extrema
if the conditions from Definition~\ref{def:DMP-criterion} hold not only in case
of a strict local minimum (maximum) but also in case of 
a non-strict local minimum (maximum) of $u_h^{}$ at the node $\bx_i^{}$.
\end{definition}

\begin{theorem}[Local and global DMPs for nonlinear discretizations of Type~I]
\label{Thm:Gen-Non-DMP-1}Let us suppose that $j_h(\cdot;\cdot)$ satisfies the weak DMP property. 
Then,  method \eqref{nonlinear-generic-1} satisfies the local DMP in the following sense:
\begin{align}
(f,\phi_i^{})\le 0\;\Longrightarrow\; \max_{\omega_i^{}}u_h^{} \le
\max_{\partial\omega_i^{}}u_h^+\,, \qquad
(f,\phi_i^{})\ge 0\;\Longrightarrow\; \min_{\omega_i^{}}u_h^{} \ge \min_{\partial\omega_i^{}}u_h^-\,, \label{Gen-Non-2}
\end{align}
for all $i=1,\dots,M$.
If $j_h(\cdot;\cdot)$ satisfies the strong DMP property,  \eqref{nonlinear-generic-1} satisfies the local DMP in the following sense:
\begin{align}
(f,\phi_i^{})\le 0\;\Longrightarrow\;\max_{\omega_i^{}}u_h^{} =
\max_{\partial\omega_i^{}}u_h^{}\,,\qquad
(f,\phi_i^{})\ge 0\;\Longrightarrow\;\min_{\omega_i^{}}u_h^{} = \min_{\partial\omega_i^{}}u_h^{}\,, \label{Gen-Non-4}
\end{align}
for all $i=1,\dots,M$.
In addition, the global DMP is also satisfied in the following form
\begin{align}
f\le 0\;\textrm{\rm in}\;\Omega\;\Longrightarrow\; \max_{\overline\Omega} u_h^{}\le
\max_{\partial\Omega}u_h^+\,,\qquad
f\ge 0\;\textrm{\rm in}\;\Omega\;\Longrightarrow\; \min_{\overline\Omega} u_h^{}\ge \min_{\partial\Omega}u_h^-\,, \label{Gen-Non-6}
\end{align}
if $j_h(\cdot;\cdot)$ satisfies the weak DMP property for non-strict extrema and in the form
\begin{align}
f\le 0\;\textrm{\rm in}\;\Omega\;\Longrightarrow\; \max_{\overline\Omega} u_h^{}=
\max_{\partial\Omega}u_h^{}\,,\qquad
f\ge 0\;\textrm{\rm in}\;\Omega\;\Longrightarrow\; \min_{\overline\Omega} u_h^{}= \min_{\partial\Omega}u_h^{}\,, \label{Gen-Non-8}
\end{align}
if $j_h(\cdot;\cdot)$ satisfies the strong DMP property for non-strict extrema.
\end{theorem}

\begin{proof}
The idea of the proof originates from \cite{BE05}.
Consider any $i\in\{1,\dots,M\}$ and let $(f,\phi_i^{})\le0$.  Since
${\max_{\omega_i}}u_h^{}$ is attained at a node, one has
${\max_{\omega_i}}u_h^{}=\max\{u_h^{}(\bx_i^{}),{\max_{\partial\omega_i}}u_h^{}\}\le\max\{u_h^{}(\bx_i^{}),{\max_{\partial\omega_i}}u_h^+\}$.
Thus, \eqref{Gen-Non-2} trivially holds if $u_h^{}(\bx_i^{})\le0$ and hence it
suffices to assume that $u_h^{}(\bx_i^{})>0$ or that the strong DMP property
holds.   Let us assume that
$u_h^{}(\bx_i^{})>{\max_{\partial\omega_i}}u_h^{}$.
Then $u_h$ attains a strict local maximum at $\bx_i^{}$ and hence the 
strong (weak) DMP property implies that
\begin{equation*}
0\ge (f,\phi_i^{})= a(u_h^{},\phi_i^{}) + j_h(u_h^{};\phi_i^{})\ge \sum_{F\in\calF_i^{}}\alpha_F^{}|\llbracket\nabla u_h^{}\rrbracket_F^{}|\,.
\end{equation*}
Thus, $\nabla u_h^{}$ is a constant in $\omega_i^{}$ and hence $u_h^{}$ is a 
$\bbP_1$ function in $\omega_i^{}$,  which is a contradiction since $u_h$ was
assumed to attain a strict local extremum in $\bx_i^{}$. Consequently, 
$u_h^{}(\bx_i^{})\le{\max_{\partial\omega_i}}u_h^{}$, which proves 
\eqref{Gen-Non-4} and also \eqref{Gen-Non-2}.
If $(f,\phi_i^{})\ge0$, one can proceed analogously.

For the global results \eqref{Gen-Non-6}, \eqref{Gen-Non-8}, let us suppose 
that $f\le 0$ in $\Omega$ and that the solution attains a global maximum at 
$\bx_i^{}$ with some $i\in \{1,\ldots,M\}$. If only the weak DMP property
holds, it is again sufficient to assume that $u_h(\bx_i^{})>0$. Then, 
analogously as for the local result, one deduces that $u_h^{}$ is a
$\bbP_1$ function in $\omega_i^{}$. Since $u_h$ attains an extremum at
$\bx_i^{}$, it has to be constant in $\omega_i^{}$, and thus the
global maximum is attained at a node $\bx_j^{}\in\partial\omega_i^{}$. If 
$\bx_j^{}\in\partial\Omega$, there is nothing more to prove. Otherwise,
we proceed as above and conclude that $u_h^{}$ is constant in 
$\omega_j^{}$ as well. Continuing in the same fashion, and using that
the mesh is connected, one can conclude that the global maximum is reached at a 
point on the boundary $\partial\Omega$. 
\end{proof}

To treat problems of Type~II, we introduce the following condition, 
reminiscent of \cite{Kno17} (see also \cite{BJKR18}).

\begin{definition}[Algebraic DMP property]
\label{def:algebraic_DMP_property}
We will say that $d_h^{}(\cdot;\cdot,\cdot)$ satisfies the algebraic DMP
property if the following condition holds: Consider any $u_h^{}\in V_h^{}$ and
any $i\in\{1,\dots,M\}$. If $u_h^{}(\bx_i^{})$ is a strict local extremum of 
$u_h^{}$ on $\omega_i^{}$, i.e.,
\begin{displaymath}
   u_h^{}(\bx_i^{})>u_h^{}(\bx)\quad\forall\,\,\bx\in\omega_i^{}\setminus\{\bx_i^{}\}
   \qquad\mbox{or}\qquad
   u_h^{}(\bx_i^{})<u_h^{}(\bx)\quad\forall\,\,\bx\in\omega_i^{}\setminus\{\bx_i^{}\}\,,
\end{displaymath}
then
\begin{equation}\label{Gen-Non-11}
   a(\phi_j^{},\phi_i^{})+  d_h^{}(u_h^{};\phi_j^{},\phi_i^{})\le 0\qquad
   \forall\,\,j\in S_i^{}
\end{equation}
and
\begin{equation}\label{eq:d_h_sparse}
   d_h^{}(u_h^{};\phi_j^{},\phi_i^{})=0\qquad\forall\,\,j\not\in S_i\cup\{i\}\,.
\end{equation}
\end{definition}

One can notice that, in essence, what \eqref{Gen-Non-11} states is that only
the $i^{\rm th}$ row in the nonlinear
system \eqref{nonlinear-generic-2} behaves like a matrix of non-negative type,
and not all the rows, in contrast to the case of linear discretizations.
The algebraic DMP property is sufficient for proving the local DMP. The proof
of the global DMP requires a sign condition also in case of non-strict extrema.

\begin{definition}[Algebraic DMP property for non-strict extrema]
\label{def:nonstrict_algebraic_DMP_property}
We will say that $d_h^{}(\cdot;\cdot,\cdot)$ satisfies the algebraic DMP
property for non-strict extrema if the following condition holds: Consider any 
$u_h^{}\in V_h^{}$ and any $i\in\{1,\dots,M\}$. If $u_h^{}(\bx_i^{})$ is a 
local extremum of $u_h^{}$ on $\omega_i^{}$, i.e.,
\begin{displaymath}
   u_h^{}(\bx_i^{})\ge u_h^{}(\bx)\quad\forall\,\,\bx\in\omega_i^{}
   \qquad\mbox{or}\qquad
   u_h^{}(\bx_i^{})\le u_h^{}(\bx)\quad\forall\,\,\bx\in\omega_i^{}\,,
\end{displaymath}
then
\begin{equation}\label{Gen-Non-12}
   a(\phi_j^{},\phi_i^{})+  d_h^{}(u_h^{};\phi_j^{},\phi_i^{})\le 0\qquad
   \forall\,\,j\in S_i^{}\,\,\,\mbox{\rm with}\,\,\,
   u_h^{}(\bx_j^{})\neq u_h^{}(\bx_i^{})
\end{equation}
and \eqref{eq:d_h_sparse} holds.
\end{definition}

\begin{theorem}[Local and global DMPs for nonlinear discretizations of Type~II]
\label{local_DMP}Let $u_h\in V_h$ be a solution of \eqref{nonlinear-generic-2} and let us suppose that $d_h^{}(\cdot;\cdot,\cdot)$
satisfies the algebraic DMP property. Then the local DMP \eqref{Gen-Non-2} holds for all $i=1,\dots,M$. 
If, in addition, $\sigma=0$, then also the stronger form \eqref{Gen-Non-4} of
the local DMP holds for all $i=1,\dots,M$.

If $d_h^{}(\cdot;\cdot,\cdot)$ satisfies the algebraic DMP property for non-strict extrema, then the
global DMP \eqref{Gen-Non-6} is satisfied. If, in addition, $\sigma=0$, then 
also the stronger form \eqref{Gen-Non-8} of the global DMP holds.
\end{theorem}

\begin{proof}
Denote $u_i=u_h(\bx_i)$ and $\tilde a_{ij}=a(\phi_j^{},\phi_i^{})+
d_h^{}(u_h^{};\phi_j^{},\phi_i^{})$ for $i,j=1,\dots,N$, and let us prove the
local versions of the DMP. Consider any 
$i\in\{1,\dots,M\}$ and let $(f,\phi_i^{})\le0$. If $\sigma>0$, it suffices to 
consider $u_i>0$ since otherwise \eqref{Gen-Non-2} trivially holds (cf.~the
beginning of the proof of Theorem~\ref{Thm:Gen-Non-DMP-1}). Let us 
assume that $u_i>u_j$ for all $j\in S_i$. If $d_h^{}(\cdot;\cdot,\cdot)$ satisfies the algebraic DMP
property, then it follows from \eqref{nonlinear-generic-2} and 
\eqref{eq:d_h_sparse} that
\begin{equation}\label{eq:8a}
   A_i\,u_i+\sum_{j\in S_i}\,\tilde a_{ij}\,(u_j-u_i)=(f,\phi_i^{})\,,
\end{equation}
where $A_i:=\sum_{j=1}^N\tilde a_{ij}=(\sigma,\phi_i^{})$ due to
\eqref{eq:d_h_row_sums}. Moreover, \eqref{Gen-Non-11} implies
that the sum in \eqref{eq:8a} is non-negative. If $\sigma=0$, then $A_i=0$ and
hence there is $j\in S_i$ such that $\tilde a_{ij}<0$ since
$\tilde a_{ii}\ge\varepsilon\,|\phi_i^{}|_{1,\Omega}^2>0$ (see
\eqref{eq:pos_semidef_dh}). This implies that the sum in \eqref{eq:8a} is 
positive. If $\sigma>0$, then $A_i\,u_i>0$.
Thus, in both cases, the left-hand side of \eqref{eq:8a} is positive, which is
a contradiction. Therefore, there is $j\in S_i$ such that $u_i\le u_j$, which
proves \eqref{Gen-Non-2} and \eqref{Gen-Non-4}. If $(f,\phi_i^{})\ge0$, one can
proceed analogously.

The proof of the global DMP can be carried out analogously as for
Theorem~\ref{thm:algebraic-global-DMP}, see also the proof of Theorem~3 in
\cite{BJKR18}.
\end{proof}

\section{Linear  discretizations of steady-state problems without convection}\label{Sec:cd_linear_methods}

This first section on linear discretizations is devoted to the special case of 
\eqref{steady-strong} where $\bb = \boldsymbol{0}$.
For all linear discretizations, the proofs of the DMP will consist of checking the hypotheses of 
Theorem~\ref{thm:algebraic-local-DMP}. It turns out that the DMP is satisfied
only under appropriate requirements on the mesh.

A careful inspection of the statements of the results from Section~\ref{Sec:DMP_linear_disc} reveals that one only needs to
show properties for the first $M$ rows of the coefficient matrix of system \eqref{system-1}-\eqref{system-2}, that is,
one only needs to worry about the equations associated with nodes interior to $\Omega$. This observation
motivates to define, for $\bbA\in\bbR^{N\times N}$, the matrix $(\bbA)^M\in \bbR^{M\times N}$
as the matrix containing only the first $M$ rows of $\bbA$. In fact, showing that $(\bbA)^M$ is of non-negative type
is what is needed to use Theorems~\ref{thm:algebraic-local-DMP} and
\ref{thm:algebraic-global-DMP} due to the expression \eqref{system-matrix} for
the matrix associated with
the system \eqref{system-1}-\eqref{system-2}.

\subsection{The Poisson problem}
\label{Sec:cd_linear_methods_poisson}

In this section we will discuss necessary and sufficient conditions for the satisfaction of the DMP for
the Poisson problem. The argument relies on proving that the diffusion matrix $(\bbAd)^M$, defined in \eqref{def-dif-matrix},
is of non-negative type.  For the finite element method the first result in this direction is given in \cite{CR73}. Since in that paper the partial differential equation
is a reaction-diffusion equation,  the mesh is supposed to be acute and fine
enough (see Section~\ref{Sec:cd_linear_methods_reaction} below). Later, for the
Poisson problem in 2d, it was noted that it is only needed for the mesh to
satisfy the Delaunay criterion, see \cite[p.~78]{SF73}. Extensions to three
space dimensions can be found in \cite{BKK20}.

We start noticing that using \eqref{basic-expression} leads to the first proof 
of the satisfaction of the DMP for the Poisson problem. In fact, if the mesh 
$\calT_h^{}$ is weakly acute, then, using \eqref{basic-expression}, one has
$\fl_{ij}^{}=\sum_{K\subset\omega_i \cap \omega_j}\fl_{ij}^{K}\le 0$ for $i\neq j$. This observation has been widely used in the literature and provides a sufficient condition
for the satisfaction of the DMP for the Poisson equation.
The proof we present next was first given in
\cite[Lemma~2.1]{XZ99} and has the advantage that it presents a necessary and
sufficient condition on the mesh to guarantee the satisfaction of the local DMP.

\begin{theorem}[Sufficient and necessary condition for $(\bbAd)^M$ to be of non-negative type, \cite{XZ99}]
\label{theo:Delaunay-mesh-laplacian}A sufficient condition for the matrix 
$(\bbAd)^M$ to be of non-negative type is that the mesh $\calT_h^{}$
satisfies the XZ-criterion \eqref{XZ-criterion}. If any internal edge of
$\calT_h^{}$ has at least one endpoint in $\Omega$, then this condition is
necessary. In addition, $(\bbAd)^M$  satisfies \eqref{zero-row-sum}.
\end{theorem}

\begin{proof}
Let $\bx_i^{},\bx_j^{}$ be two different nodes contained in the same mesh cell 
$K \in \calT_h^{}$. We recall the following formulas for the volume of a
simplex
\begin{equation*}
|K| =  \frac{|F_i^K||F_j^K|}{2} \sin\theta_{E_{ij}^{}}^K \quad \mbox{if } d=2,\qquad
|K| =  \frac{2|F_i^K||F_j^K|}{3 |\kappa_{E_{ij}^{}}^K|} \sin\theta_{E_{ij}^{}}^K \quad \mbox{if } d=3\,.
\end{equation*}
Inserting them in \eqref{basic-expression}, and using the convention that 
$|\kappa_{E_{ij}^{}}^K|=1$ if $d=2$ gives
\begin{equation}\label{eq:lK_ij}
\fl_{ij}^{K} = -\frac{1}{d(d-1)}|\kappa_{E_{ij}^{}}^K|\cot\theta_{E_{ij}^{}}^K \,.
\end{equation}
Thus, for $i\in\{1,\dots,M\}$ and $j\in S_i$, 
\begin{equation}\label{eq:l_ij}
\fl_{ij}^{}= \sum_{K\subset\omega_{E_{ij}}}  \fl_{ij}^{K} 
= - \sum_{K\subset\omega_{E_{ij}}} \frac{|\kappa_{E_{ij}^{}}^{K}|\cot\theta_{E_{ij}^{}}^K}{d(d-1)}\,,
\end{equation}
and then \eqref{NN-1} is satisfied if \eqref{XZ-criterion} holds. If the set 
$\calE_I^{}$ consists only of edges $E_{ij}^{}$ with $i\in\{1,\dots,M\}$ and 
$j\in S_i$, then \eqref{XZ-criterion} is necessary for the validity of
\eqref{NN-1}. Finally, since the basis functions form a partition of unity, 
one has
\begin{equation} \label{eq:A_d_row_sum}
\sum_{j=1}^N\,\fl_{ij}^{}=\sum_{j=1}^N\,(\nabla\phi_j^{},\nabla\phi_i^{})=
(\nabla 1,\nabla\phi_i^{})=0\,.
\end{equation}
So, \eqref{zero-row-sum} is satisfied, and in particular \eqref{NN-2}.
\end{proof}

\begin{remark}The statement of Theorem~\ref{theo:Delaunay-mesh-laplacian} 
implies, in connection with Theorem~\ref{thm:algebraic-local-DMP}, that the 
local DMP is satisfied if and only if the mesh is of XZ-type, with the slight
exception concerning edges whose endpoints are both on $\partial\Omega$. In 
addition, Theorems~\ref{theo:Delaunay-mesh-laplacian} and
\ref{thm:algebraic-global-DMP} show that the validity of the XZ-criterion
implies the global DMP. However, in this case, the XZ-criterion is not
necessary. Indeed, in \cite{DDS05} a two-dimensional example is constructed
where the global DMP is satisfied although the mesh is not of XZ-type. 
Nevertheless, in general, if the mesh is not of XZ-type,
then the global DMP might be violated as an example in \cite{BKKS09}
demonstrates. 
\hspace*{\fill}$\Box$\end{remark}

\begin{remark}\label{rem:diff_mat_non_sing} Let $\bbAdI \in \mathbb R^{M\times M}$ denote the $M\times M$ submatrix of the diffusion
matrix only considering the non-Dirichlet nodes, i.e.,  the analog of $\bbAI$ in \eqref{system-matrix}. Then,  
$\bbAdI$ is non-singular, since the corresponding bilinear form is elliptic on $H_0^1(\Omega)$.
\hspace*{\fill}$\Box$\end{remark}

\begin{remark} \label{rem:aniso_diff_poisson}
A Poisson problem with heterogeneous anisotropic 
diffusion is given by 
\begin{equation}\label{eq:het_ani_diff_prob}
\begin{array}{rcll}
-\nabla \cdot \left( \mathbb E(\bx) \nabla u \right)  &=& f & \textrm{in}\;\Omega\,, \\
u&=& g &\textrm{on}\;\partial\Omega\,,
\end{array}
\end{equation}
with the symmetric diffusion tensor $\mathbb E(\bx)$. The tensor $\mathbb E$ 
depends on the spatial variable $\bx$, which makes it heterogeneous, and in 
addition it is allowed to have different eigenvalues at a given $\bx$, making it
anisotropic. In any case, it will be assumed that $\mathbb E$ is symmetric and
strictly positive-definite in~$\Omega$. Numerous applications lead to
heterogeneous anisotropic diffusion, such as image processing \cite{TsiPet13}
and atmospheric modelling \cite{Stockie11}, just to name a few. 

Problem \eqref{eq:het_ani_diff_prob} was considered in 
\cite{LH10} for $\mathbb P_1^{}$ finite elements in two 
and three dimensions. The main condition on the mesh is the following: for
every element $K$ it is assumed that
\begin{equation}\label{eq:het_ani_diff_DMP}
\left(\bn_i^K\right)^T \mathbb E_K \bn_j^K \le 0 \quad \forall\ 
\bx_i^{},\bx_j^{}\in K,\,\, \bx_i^{}\neq\bx_j^{}, \quad
\forall\ K \in \calT_h\,,
\end{equation}
where $\mathbb E_K$ stands for an approximation of the integral of $\mathbb{E}$ in $K$ using quadrature.
By writing the global matrix as sum of local contributions it is proven that under this assumption the system matrix is 
of non-negative type, from which the validity of the 
DMP can be concluded using the 
results
presented in Section~\ref{Sec:DMP_linear_disc}. 
It can be readily seen that in the special case $\mathbb E_K = \mathbb I$, \eqref{eq:het_ani_diff_DMP}
reduces to the weakly acute angle condition from Definition~\ref{Def:conditions-on-mesh}.  A comprehensive 
interpretation of \eqref{eq:het_ani_diff_DMP} is provided in \cite{Hua11}. It turns out 
that \eqref{eq:het_ani_diff_DMP} is equivalent to the requirement that the angles are 
weakly acute with respect to an inner product induced by $\mathbb E_K^{-1}$. Condition 
\eqref{eq:het_ani_diff_DMP}  can be expressed in terms of the map from the reference cell 
to $K$. This formulation was utilized in \cite{LH10} for the construction of appropriate
meshes on which the numerical solution satisfies the global DMP. 

Later, in  \cite{Hua11}, the analysis from \cite{LH10} was refined for the
two-dimensional situation in order to
obtain a condition weaker than \eqref{eq:het_ani_diff_DMP}. The  numerical analysis studies the global stiffness matrix, 
in contrast to the analysis from  \cite{LH10}, and in the isotropic case
$\mathbb E_K = \mathbb I$ the resulting condition becomes that the mesh has to 
be  Delaunay.
\hspace*{\fill}$\Box$\end{remark}

\subsection{The reaction-diffusion equation and mass lumping}\label{Sec:cd_linear_methods_reaction}

So far the reaction was set to be zero to show the intrinsic link between the geometry of the mesh and
the properties of the matrix $\bbAd$. If reaction is added, the satisfaction of the DMP is in fact harder
than for the plain diffusion equation, as the next result shows.

\begin{lemma}[Sufficient condition for ${(\varepsilon\bbAd+\sigma\bbMc)^M}$ to be of non-negative type] 
\label{Diffusion-reaction}Let  
$\bbMc$ be the consistent mass matrix defined in \eqref{def-mass-matrix}.
Then, ${(\varepsilon\bbAd+\sigma\bbMc)^M}$ is of non-negative type if the
mesh family $\{\calT_h^{}\}_{h>0}^{}$ is strictly acute and $h$  satisfies 
\begin{equation}\label{RD-condition-h}
h^2 \le C\frac{\varepsilon}{\sigma}\,\cot\left(\frac\pi2 - \delta\right) = C\frac{\varepsilon}{\sigma}\,\tan\delta\,,
\end{equation}
 where $\delta$  is the angle from \eqref{eq:strictly_acute}, $C=12$ in 2d,
and $C$ depends only on the shape regularity of the mesh family
$\{\calT_h^{}\}_{h>0}^{}$ in 3d.
\end{lemma}

\begin{proof} The satisfaction of \eqref{NN-2} follows from
\eqref{eq:A_d_row_sum} and the fact that the row sum of the consistent mass 
matrix is positive, compare \eqref{eq:mass_matrix_l_c}.

Consider two nodes $\bx_i\neq\bx_j$ of a mesh cell $K\in\calT_h^{}$. The shape
regularity of the mesh implies that there is a constant $C_0$ such that
$|\kappa_{E_{ij}^{}}^K|\ge C_0h_K^{d-2}$ (note that one can set $C_0=1$ if
$d=2$).  Since $|K|\le h_K^d/(d(d-1))$, one
obtains using \eqref{eq:lK_ij}, the exact formula for the local mass matrix 
\eqref{local-mass-exact}, and the fact that the cotangent is monotonically
decreasing
\begin{align*}
\varepsilon\fl_{ij}^K+\sigma\mm_{ij}^K&=
-\varepsilon\frac{|\kappa_{E_{ij}^{}}^K|\cot\theta_{E_{ij}^{}}^K}{d(d-1)}+ \sigma\frac{|K|}{(d+1)(d+2)} \nonumber\\
&\le h_K^{d-2}\,\frac{(d-2)!}{(d+2)!}\left(
-\varepsilon\,C_0\,(d+1)(d+2)\cot(\frac{\pi}{2}-\delta)+
\sigma h_K^2\right).
\end{align*}
Hence, \eqref{RD-condition-h} with $C=C_0(d+1)(d+2)$ leads to $\varepsilon\fl_{ij}^{}+\sigma\mm_{ij}^{}\le 0$ for $i\neq j$, thus proving \eqref{NN-1}. 
\end{proof}

The last result shows that the presence of a positive reaction term makes the satisfaction of the DMP more difficult than for the
Poisson problem. In fact, the presence of the reaction imposes a restriction on the size of the mesh (cf. \eqref{RD-condition-h}) as well
as a stronger restriction on the geometry. While the need for a strictly acute
mesh family is clear from the proof, the restriction on the mesh
size has been slightly relaxed in, e.g., \cite{BKK08}, although some size restriction is always present as long as the consistent mass
matrix is used (see \cite{BKK08} for examples of non-satisfaction of the DMP if the mesh is not refined enough). So, we now move onto the 
presentation of a mass-lumping strategy that allows one to remove the size restriction without
affecting accuracy. The mass-lumped discretization of the reaction-diffusion equation
reads as follows: Find $u_h^{}\in V_h^{}$ such that $u_h^{}|_{\partial\Omega}^{}=i_h^{}g$, and
\begin{equation*}
\varepsilon(\nabla u_h^{},\nabla v_h^{})+\sigma(u_h^{},v_h^{})_h^{}= (f,v_h^{})\qquad\forall\, v_h^{}\in V_{h,0}^{}\,,
\end{equation*}
with $(\cdot,\cdot)_h^{}$ defined in \eqref{eq:lumped_inner_prod}. The
following result shows that the stiffness matrix of this modified Galerkin 
discretization is of non-negative type under the same conditions as the 
stiffness matrix of the pure diffusion problem. Thus, the modification
removes the restriction on the mesh size from Lemma~\ref{Diffusion-reaction}.

\begin{corollary}[Sufficient and necessary condition for 
${(\varepsilon\bbAd+\sigma\bbMl)^M}$ to be of non-negative 
type]\label{cor:diff_reac_non_neg_mat}Let $\bbMl$ be the lumped mass matrix 
defined in \eqref{def-lumpedmass-matrix}. Then, a sufficient condition for the 
matrix ${(\varepsilon\bbAd+\sigma\bbMl)^M}$ to be of non-negative type is that 
the mesh $\calT_h^{}$ is of XZ-type. If any internal edge of $\calT_h^{}$ 
has at least one endpoint in $\Omega$, then this condition is necessary.
\end{corollary}

\begin{proof} The proof follows by realizing that the lumping process removes
the positive off-diagonal entries of $\bbMc$, and then it becomes a
direct application of Theorem~\ref{theo:Delaunay-mesh-laplacian}.
\end{proof}

\begin{remark}
This section is finished with a brief discussion concerning the fact that an appropriate stabilized method for the reaction-diffusion equation
also helps relaxing the mesh conditions for the satisfaction of the DMP, even if it uses the consistent mass matrix.
This method, known as {Unusual Stabilized} finite element method (USFEM), was introduced in \cite{FF95} and
 reads as follows:  find $u_h^{}\in V_{h}^{}$ such that $u_h^{}|_{\partial\Omega}^{}=i_h^{}g$, and
\begin{eqnarray}\label{USFEM-1}
\lefteqn{
\varepsilon\, (\nabla u_h^{},\nabla v_h^{})+\sigma\,(u_h^{},v_h^{})-\sum_{K\in\calT_h^{}}\frac{h_K^2}{\sigma h_K^2+\varepsilon}(\sigma u_h^{},\sigma v_h^{})_K^{}}\\
&=& (f,v_h^{})-\sum_{K\in\calT_h^{}}\frac{h_K^2}{\sigma h_K^2+\varepsilon}(f,\sigma v_h^{})_K^{}
\quad \forall\ v_h^{}\in V_{h,0}^{}.\nonumber
\end{eqnarray}
The USFEM improves stability by subtracting a term of reaction type from both sides of the finite element equation. As a consequence, the corresponding  matrix 
$(\bbA)^M$ has the entries
\begin{equation*}
a_{ij}^{}=\varepsilon\, (\nabla \phi_j^{}, \nabla \phi_i^{})+ \sum_{K\in\calT_h^{}}\frac{\sigma\varepsilon}{\sigma h_K^2+\varepsilon} ( \phi_j^{},  \phi_i^{})_K^{}\,.
\end{equation*}
Following the same steps as in the proof of Lemma~\ref{Diffusion-reaction},
one can see that $a_{ij}^{}\le 0$ requires the mesh family to be strictly acute and 
\begin{equation}\label{USFEM-condition-h}
   \frac{\varepsilon}{\sigma h_K^2+\varepsilon}\,h_K^2
   \le C\frac{\varepsilon}{\sigma}\,\tan\delta\qquad\forall\ K\in \calT_h^{}\,,
\end{equation}
where $\delta$  is the angle from \eqref{eq:strictly_acute} and $C$ is the same
as in \eqref{RD-condition-h}. In the interesting case $\varepsilon\ll\sigma$, 
\eqref{USFEM-condition-h} is a much milder condition than 
\eqref{RD-condition-h}.  Moreover, \eqref{USFEM-condition-h} does not
restrict $h_K$ at all if $C\tan\delta\ge1$. 
Likewise important,   the sign of the right-hand side of \eqref{USFEM-1} is not affected, since
it can be written for every basis function $\phi_i^{}$  as
\begin{equation*}
\sum_{K\in\calT_h^{}}\frac{\varepsilon}{\sigma h_K^2+\varepsilon}(f, \phi_i^{})_K^{}\,.
\end{equation*}
Thus, for a uniform mesh with $h_K=h$ for any $K\in \calT_h^{}$, the USFEM
is equivalent to replacing $\varepsilon$ by $\varepsilon+\sigma h^2$ in the
standard Galerkin discretization so that it just adds isotropic linear
artificial diffusion of amount $\sigma h^2$, cf.~Section~\ref{Sec:cd_linear_methods_iso_diff}.

In summary, the USFEM \eqref{USFEM-1} preserves the DMP whenever \eqref{USFEM-condition-h} is satisfied.  
\hspace*{\fill}$\Box$\end{remark}

\section{Linear discretizations of the steady-state problem}  
\label{Sec:cd_linear_methods_cd} 

In this section the main ideas for a linear discretization of the
convection-diffusion equation \eqref{steady-strong} are given.
It should be kept in mind that the presentation of this
and the following sections focuses on the convection-dominated regime, even if 
this is not always explicitly stated, i.e., 
$\varepsilon$ has to be thought of being (very) small. 
First, to justify the need for stabilization we describe the standard Galerkin method and make it
explicit that, unless the mesh is acute, and prohibitively refined, the DMP cannot hold.
So, we then consider stabilized discretizations, where we review linear artificial diffusion, upwind methods,
and the edge-averaged finite element method.

\subsection{The Galerkin finite element method}
\label{Sec:cd_linear_methods_gal} 

The Galerkin finite element \linebreak meth\-od reads as follows: Find $u_h^{} \in V_{h}^{}$ such that
$u_h^{}=i_h^{}g$ on $\partial\Omega$ and
\begin{equation}\label{eq:galerkin_fem}
a(u_h^{},v_h^{}) = (f,v_h^{})
\quad \forall v_h^{} \in V_{h,0}^{} \,,
\end{equation} 
where $a(\cdot,\cdot)$ is defined in \eqref{bilinear-a}.
Following classical arguments (see, e.g. \cite{EG21-II}) one can derive {optimal order
error estimates, but with a constant that behaves like
$\|\bb\|_{0,\infty,\Omega}^{}h/\varepsilon$,
thus making these estimates not useful in practice, and somehow explaining why 
non-localized spurious oscillations
appear in the simulations. This feature is shared by all central
discretizations of the convective term (see, e.g., \cite{RST08} for extensive discussions
on this issue).}

To illustrate the restrictions of the Galerkin method with respect to the satisfaction of the DMP we focus on the special case where $d=2$ and $\sigma=0$. 
Here, the matrix associated with \eqref{eq:galerkin_fem} is
$(\bbA)^M=(\varepsilon\bbAd+\bbAc)^M$, compare \eqref{def-dif-matrix} and \eqref{def-conv-matrix}.
Since $\bb$ is solenoidal, $\bbAc$ satisfies 
\begin{equation}\label{semi-antisymmetry}
\cc_{ij}^{}=-\cc_{ji}^{}\qquad\textrm{for all}\; i,j=1,\ldots, M\;,
\end{equation}
i.e., there is a partial antisymmetry.

The next result states that the Galerkin method satisfies the DMP if the mesh
family $\{\calT_h^{}\}_{h>0}^{}$ is average acute and $h$ is sufficiently small.

\begin{theorem}[Conditions on the Galerkin method in 2d to satisfy the DMP]
Suppose \label{theo:fem-Gal-DMP}that $d=2$, $\sigma=0$, the mesh family
$\{\calT_h^{}\}_{h>0}^{}$ is average acute, and the data
and the mesh satisfy:
for all $E=K\cap K'\in \calE_I^{}$, it holds
\begin{equation}\label{FEM-Gal-3}
\frac{ (h_K+h_{K'})\|\bb\|_{0,\infty,\omega_E^{}}}
{3 \tan\frac{\delta}{2}} \le \varepsilon \,,
\end{equation}
where $\delta$ is the angle from \eqref{FEM-prelim-1}.
Then, the  matrix $(\varepsilon\bbAd+\bbAc)^M$ is of  non-negative type and satisfies 
\eqref{zero-row-sum}. 
\end{theorem}

\begin{proof}
Since the basis functions $\phi_1^{},\ldots,\phi_N^{}$ form a partition of unity, $(\varepsilon\bbAd+\bbAc)^M$  satisfies
\begin{equation}\label{eq:zero_row_diff_conv}
\sum_{j=1}^Na_{ij}^{}=\varepsilon\,(\nabla 1,\nabla \phi_i^{})+(\bb\cdot\nabla 1,\phi_i^{})
= 0\,, \quad i=1,\ldots,{M}\, ,
\end{equation}
which proves \eqref{zero-row-sum}.
It remains to show \eqref{NN-1}.
Let $E=K\cap K'\in\calE_I^{}$ with endpoints $\bx_i^{},\bx_j^{}$, $i\in \{1,\ldots,M\}$, 
$j\in \{1,\ldots,N\}$.
Using \eqref{eq:l_ij} and $|\kappa_{E_{ij}^{}}^{}|=1$ yields
\begin{align}
\fl_{ij}^{} &= (\nabla\phi_j^{},\nabla\phi_i^{})_K^{}+
(\nabla\phi_j^{},\nabla\phi_i^{})_{K'}^{}\label{Gal3}\\
&=-\frac{1}{2}\cot\theta_E^K-\frac{1}{2}\cot\theta_E^{K'} = -\frac{\sin(\theta_E^K+\theta_E^{K'})}{
2\sin\theta_E^K\sin\theta_E^{K'}}\,.\nonumber
\end{align}
In addition, since $\theta_E^K,\theta_E^{K'}\in (0,\pi)$, one has
\begin{align}
&\sin^2\left(\frac{\theta_E^K+\theta_E^{K'}}{2}\right)  =
\frac{1-\cos(\theta_E^K+\theta_E^{K'})}{2}\label{Gal4}\\
&= \frac{1-\cos\theta_E^K\cos\theta_E^{K'}}{2}+\frac{\sin\theta_E^K\sin\theta_E^{K'}}{2}
> \frac{\sin\theta_E^K\sin\theta_E^{K'}}{2} >0\,.\nonumber
\end{align}
Observing that the right-hand side of \eqref{Gal3} is negative, since the mesh
family is average acute and $\theta_E^K,\theta_E^{K'}\in (0,\pi)$, inserting
\eqref{Gal4} in \eqref{Gal3}, and using the monotonicity of the cotangent leads 
to 
\begin{align}
\fl_{ij}^{} &<
-\frac{\sin(\theta_E^K+\theta_E^{K'})}{4\sin^2\left(\frac{\theta_E^K+\theta_E^{K'}}{2}\right)} = 
-\frac{1}{2}\cot\frac{\theta_E^K+\theta_E^{K'}}{2}\label{Gal5}\\
&
\le -\frac{1}{2}\cot\left(\frac{\pi}{2}-\frac{\delta}{2}\right)
= - \frac{1}{2}\tan\frac{\delta}{2} < 0\,.\nonumber
\end{align}

Concerning the convective term, a direct calculation using \eqref{grad-basis-functions}, H\"older's inequality, and that the diameter of any facet of $K$ 
 is bounded by $h_K$, gives 
\begin{equation}\label{Gal6}
(\bb\cdot\nabla\phi_j,\phi_i)_K^{} = -\frac{|F_j^K|}{2|K|}(\bb\cdot\boldsymbol{n}_j^K,\phi_i)_K
\le \frac{h_K\|\bb\|_{0,\infty,K}}{2|K|}\frac{|K|}{3}\le \frac{h_K\|\bb\|_{0,\infty,K}}{6}\,.
\end{equation}

From \eqref{Gal5} and \eqref{Gal6}, one obtains the following upper bound for the off-diagonal matrix 
entries
\begin{equation}
a_{ij}^{}=\varepsilon\fl_{ij}^{} + \cc_{ij}^{} 
\le -\frac{\varepsilon}{2}\, \tan\frac{\delta}{2} 
+  \frac{(h_K^{}+h_{K'})\|\bb\|_{0,\infty,\omega_E^{}}}{6}\label{Gal7}
\end{equation}
and hence $a_{ij}^{}\le 0$ if \eqref{FEM-Gal-3} holds.
\end{proof}

\begin{remark}\label{Rem:acute-optimal} The geometrical  hypothesis on
$\calT_h^{}$ cannot be relaxed.  Indeed, suppose that 
$\{\calT_h^{}\}_{h>0}^{}$ is not average acute and choose an internal edge 
$E=K\cap K'\in\calE_I^{}$ with endpoints $\bx_i^{},\bx_j^{}$, 
$i ,j\in\{1,\ldots,M\}$, such that $\theta_E^K+\theta_E^{K'}=\pi$. Then, 
thanks to \eqref{Gal3}, it follows that $\fl_{ij}^{}= \fl_{ji}^{}=0$. So,
since $\bbAc$ satisfies \eqref{semi-antisymmetry}, then for any $\bb$ such that 
$\cc_{ij}^{}\not=0$, one has $\cc_{ij}^{}>0$ or $\cc_{ji}^{}>0$, which 
implies that $\bbA$ does not satisfy \eqref{NN-1}.
\hspace*{\fill}$\Box$\end{remark}

The discussion in this section  shows that the Galerkin method will not satisfy
the DMP in any practical situation. These observations were made as early
as \cite{Kik77}. On the other hand, supposing the mesh family
$\{\calT_h^{}\}_{h>0}^{}$ is average acute relaxes
the hypotheses made by \cite{Kik77,Codina93,BE02}, since in those works the 
results were proven for strictly acute mesh families.

\begin{remark} The analysis of  \cite{LH10} and \cite{Hua11} for heterogeneous 
anisotropic diffusion problems (cf.~Remark \ref{rem:aniso_diff_poisson})
was extended to convection-diffusion-reaction
problems in \cite{LHQ14}. Since a Galerkin discretization without
mass lumping was considered, a condition on the fineness of the mesh appears for the 
satisfaction of the DMP,  
cf.~Lemma~\ref{Diffusion-reaction} and Theorem~\ref{theo:fem-Gal-DMP}.
\hspace*{\fill}$\Box$\end{remark}

Concentrating for a brief discussion of an error estimate on the impact 
of diffusion and convection, i.e., considering $\sigma=0$ and homogeneous 
Dirichlet boundary conditions, one finds under the assumption that 
$u\in H^2(\Omega)$ that 
\begin{equation}\label{eq:gal_ee_est}
|u-u_h|_{1,\Omega}
\le Ch \left(1+\frac{\|\bb\|_{0,\infty,\Omega}\,h}{\varepsilon}\right) |u|_{2,\Omega},
\end{equation}
where $C$ comes from interpolation error estimates in the 
$L^2(\Omega)$ norm and in the $H^1(\Omega)$ seminorm. The term in the 
parentheses is very large in the convection-dominated case so that, although 
\eqref{eq:gal_ee_est} predicts first order error reduction, the error bound is
not useful as long as $h$ is not very small. In fact, large errors can be 
observed for the Galerkin method on coarse grids if the solution of
\eqref{steady-strong} possesses layers.

\subsection{Isotropic linear artificial diffusion}
\label{Sec:cd_linear_methods_iso_diff}

Restriction \eqref{FEM-Gal-3} can be circumvented by either refining the mesh or making the diffusion
of the discrete problem larger. This section will analyze a method that takes the latter approach and adds artificial
diffusion to the problem. It will turn out that the diffusion added needs to be  of a size proportional to the mesh size.
This method will also be supplemented with a mass lumping strategy  in order to avoid technical complications
due to the presence of reaction. 

The following finite element method with
added artificial diffusion will be studied: Find $u_h^{}\in V_h^{}$ 
such that $u_h^{}|_{\partial\Omega}^{}=i_h^{}g$, and 
\begin{equation}\label{FEM-LD-1}
a_h^{}(u_h^{},v_h^{})+ s_h(u_h^{},v_h^{})
= (f,v_h^{})\qquad\forall\, v_h^{}\in V_{h,0}^{}\,,
\end{equation}
where the bilinear form $a_h^{}(\cdot,\cdot)$ is given by
\begin{equation}\label{bilinear-a-lumped}
   a_h^{}(u,v)=\varepsilon\,(\nabla u,\nabla v)+(\bb\cdot\nabla u,v)
   +\sigma\,(u,v)_h^{}\,,
\end{equation}
with $(\cdot,\cdot)_h^{}$ being the mass-lumped inner product defined in \eqref{eq:lumped_inner_prod},
and the added linear artificial diffusion term is given by
\begin{equation*}
s_h(u_h^{},v_h^{})=\sum_{K\in\calT_h^{}}\tilde{\varepsilon}_K^{}(\nabla u_h,\nabla v_h)_K^{}\,, \quad 
\tilde{\varepsilon}_K\ge 0\,.
\end{equation*}
In this section we consider the following expression for the added diffusion  \cite{Kik77}:
\begin{equation}\label{FEM-LD-2}
\tilde{\varepsilon}_K^{}:=\max\left\{ c_0^{}\frac{ h_K\|\bb\|_{0,\infty,\Omega}}{
\tan\frac{\delta}{2}}-\varepsilon, 0 \right\}\,,
\end{equation}
where $\delta$ is the constant from \eqref{FEM-prelim-1} and $c_0^{}>0$ is a constant that is only linked to the shape regularity of the triangulation, 
see \eqref{LD-5} below. 
One notices the close relation between \eqref{FEM-LD-2} and \eqref{FEM-Gal-3}. In fact, the added diffusion is built in such a way that
once the mesh is sufficiently fine, \eqref{FEM-LD-1} reduces to  the standard
Galerkin method (up to the lumping of the reaction term).
Later  works proposed slightly different versions of $\tilde{\varepsilon}_K^{}$, e.g., see \cite{Codina93,BE02}.

The analysis of \eqref{FEM-LD-1} was carried out originally in \cite{Kik77}
under the assumption that the mesh families
are strictly acute. The analysis presented below is detailed for $d=2$, and relaxes this hypothesis and requires only
average acute mesh families (the case $d=3$ is discussed in Remark~\ref{Rem-LD-3D}).

\begin{theorem}[DMP for isotropic linear artificial diffusion in 2d]
\label{theo:fem-LD-DMP}Let
us suppose $d=2$, that the mesh family is average acute,
$\tilde{\varepsilon}_K^{}$ are defined
by \eqref{FEM-LD-2}, and $c_0^{}$ is large enough (see \eqref{LD-5}). Then, \eqref{FEM-LD-1} satisfies the DMP.
\end{theorem}

\begin{proof}
The proof consists in rewriting method~\eqref{FEM-LD-1} as follows: Find $u_h^{}\in V_h^{}$ such that $u_h^{}|_{\partial\Omega}^{}=i_h^{}g$, and
\begin{equation*}
\sum_{K\in\calT_h^{}}(\varepsilon+\tilde{\varepsilon}_K^{})(\nabla u_h^{},\nabla v_h^{})_K^{}
+(\bb\cdot\nabla u_h^{},v_h^{}) +\sigma\,(u_h^{},v_h^{})_h^{}=(f,v_h^{})\quad\forall\, v_h^{}\in V_{h,0}^{}\,.
\end{equation*}
Let $i\in\{1,\ldots,M\}$ and $i\not=  j\in\{1,\ldots,N\}$. Since
 the off-diagonal elements of the lumped mass matrix vanish, 
one gets
\begin{equation*}
a_{ij}^{}= \sum_{K\in\calT_h^{}}(\varepsilon+\tilde{\varepsilon}_K^{})(\nabla\phi_j^{},\nabla\phi_i^{})_K^{}+\cc_{ij}^{}\,.
\end{equation*}
Using the notation from the proof of Theorem~\ref{theo:fem-Gal-DMP} and
assuming that $(\nabla\phi_j^{},\nabla\phi_i^{})_K^{}\le0$ and
$(\nabla\phi_j^{},\nabla\phi_i^{})_{K'}\le0$, one can use the fact that
\begin{equation}\label{FEM-LD-4}
\varepsilon+\tilde{\varepsilon}_K^{} \ge
c_0^{}\frac{\|\bb\|_{0,\infty,\Omega}^{}h_K^{}}{\tan\frac{\delta}{2}}
\ge \frac{\|\bb\|_{0,\infty,\Omega}^{}(h_K^{}+h_{K'})}{3\tan\frac{\delta}{2}} \,,
\end{equation}
if
\begin{equation}\label{LD-5}
c_0^{}\ge \max_{K,K'\in\calT_h^{}:K\cap K'\in\calE_I^{}}\frac{h_K^{}+h_{K'}}{3\min\{h_K^{},h_{K'}\}}\,,
\end{equation}
which is a constant uniformly bounded thanks to the mesh regularity.  Then
an application of the techniques used to prove Theorem~\ref{theo:fem-Gal-DMP}
 shows that the system matrix of method~\eqref{FEM-LD-1} is of non-negative
type.  If, e.g., $(\nabla\phi_j^{},\nabla\phi_i^{})_{K'}^{}>0$, then 
$(\nabla\phi_j^{},\nabla\phi_i^{})_K^{}\le0$ since the mesh family is average
acute. Moreover, since $\theta_E^{K'}\ge\frac\pi2$, one has $h_{K'}=h_E\le
h_K$. Therefore, $\varepsilon+\tilde{\varepsilon}_{K'}
\le\varepsilon+\tilde{\varepsilon}_K^{}$ and hence
\begin{equation*}
(\varepsilon+\tilde{\varepsilon}_K^{})(\nabla\phi_j^{},\nabla\phi_i^{})_K^{}+
(\varepsilon+\tilde{\varepsilon}_{K'})(\nabla\phi_j^{},\nabla\phi_i^{})_{K'}\le
(\varepsilon+\tilde{\varepsilon}_K^{})\,\fl_{ij}^{}\,.
\end{equation*}
Now one can apply \eqref{FEM-LD-4} and conclude that $a_{ij}^{}\le0$
analogously as before.
For $\sigma=0$, the method satisfies \eqref{zero-row-sum}.
Finally, the theorem follows from the results of
Section~\ref{Sec:DMP_linear_disc}.
\end{proof}

\begin{remark} Once again, the hypothesis on the mesh family being average 
acute is sharp. In fact,  analogous considerations as made in 
Remark~\ref{Rem:acute-optimal} hold in this case. 
\hspace*{\fill}$\Box$\end{remark}

\begin{remark}\label{Rem-LD-3D}
We now briefly discuss the case $d=3$. For this case one needs to assume that 
the mesh family $\{\calT_h^{}\}_{h>0}^{}$
is strictly acute. Let $\delta>0$ be the angle from
\eqref{eq:strictly_acute}, and let the added diffusion be given by 
\begin{equation*}
\tilde{\varepsilon}_K^{} = \max\left\{c_0^{}\frac{ h_K\|\bb\|_{0,\infty,K}}{\tan\delta}-\varepsilon,0\right\}\,.
\end{equation*}
 Then, following
the same steps as to reach \eqref{FEM-LD-4} and using that $|\kappa_{E_{ij}}^K|\ge Ch_K^{}$ (thanks to the mesh regularity) one gets
\begin{align*}
a_{ij}^{} &= \sum_{K\in\calT_h^{}}(\varepsilon+\tilde{\varepsilon}_K^{})(\nabla\phi_j^{},\nabla\phi_i^{})_K^{}+\cc_{ij}^{}\nonumber\\
&\le
\sum_{K\subset\omega_i \cap \omega_j}\left\{-\frac{\varepsilon+\tilde{\varepsilon}_K^{}}{6}|\kappa_{E_{ij}}^K|\cot\theta_{E_{ij}^{}}^K+
\frac{h_K^2\|\bb\|_{0,\infty,K}^{}}{24}\right\}\nonumber\\
&\le \sum_{K\subset\omega_i \cap \omega_j}\left\{-C\, c_0^{}\frac{ h_K\|\bb\|_{0,\infty,K}}{6\tan\delta} h_K^{}\tan\delta+
\frac{h_K^2\|\bb\|_{0,\infty,K}^{}}{24}\right\}\nonumber\\
&= \sum_{K\subset\omega_i \cap \omega_j} h_K^2\|\bb\|_{0,\infty,K}^{} \left\{-\frac{C\,
c_0^{}}{ 6}+
\frac{1}{24}\right\}\,.
\end{align*}
By supposing $c_0^{}$ is large enough one concludes that $a_{ij}^{}\le 0$.
Thus, in three space dimensions the same result holds as in 2d under the
assumption of a strictly acute mesh family. 
\hspace*{\fill}$\Box$\end{remark}

The last theorem shows that method~\eqref{FEM-LD-1} satisfies the DMP under 
much milder assumptions than the Galerkin method. 

We finish this section with a short comment on an error estimate for method \eqref{FEM-LD-1}. We place ourselves in the same situation as in Section~\ref{Sec:cd_linear_methods_gal}, i.e., 
$\sigma =0$, $g=0$, and $u\in H^2(\Omega)$, and assuming $\tilde\varepsilon_K=\tilde\varepsilon$ for any 
$K\in\calT_h^{}$, gives the estimate 
\begin{align*}
|u-u_h|_{1,\Omega}&\le Ch \left(1+\frac{\|\bb\|_{0,\infty,\Omega}\,h}{\varepsilon+\tilde\varepsilon}\right) |u|_{2,\Omega}
+ \frac{\tilde\varepsilon}{\varepsilon + \tilde\varepsilon} |u|_{1,\Omega}\\
                  &\le Ch \left(1+\frac{\tan\delta}{c_0}\right) |u|_{2,\Omega}
+ \frac{\tilde\varepsilon}{\varepsilon + \tilde\varepsilon} |u|_{1,\Omega},
\end{align*}
where $C$ is again only linked to interpolation error estimates.
In contrast to the error estimate \eqref{eq:gal_ee_est} for the Galerkin 
method, the factor in front of $|u|_{2,\Omega}$ is of order ${\cal O}(1)$.
However, due to the consistency error estimated by the term including
$|u|_{1,\Omega}$, there is no reduction of the bound proportional to some power
of the mesh size as long as the P\'eclet number 
$\|\bb\|_{0,\infty,\Omega}^{}h/\varepsilon$ is large.
Note that this second term is strictly 
monotonically decreasing as $\tilde \varepsilon$ tends to zero and eventually it vanishes.

An extension of the linear isotropic diffusion method has recently been
proposed in \cite{BBK17a}. The interest in this extension
by itself is limited, but it opens the door for a LPS-based nonlinear discretization, to be presented
in Section~\ref{Sec:cd_nonlinear_methods_edge}.

\subsection{Upwind finite element methods} 
\label{Sec:cd_linear_methods_upw}

In this section, one of the earliest proposals for satisfying the DMP in the framework of 
finite element methods for convection-diffusion equations is reviewed. 
The basic  idea of this method consists in  discretizing the convective term in a finite volume manner
and utilizing an upwind technique. The first method of this type was 
developed in \cite{Tabata77}. An improved method is presented in \cite{BT81} and an extension
to non-conforming finite elements in \cite{OU84}, see Section~\ref{sec:other_fems_CR} for more details.
Although the methods from \cite{Tabata77,BT81} were originally proposed for
transient problems, compare Section~\ref{sec:tcd_ext_methods_steady}, we present
here their  steady-state versions as they contain the main ideas.
From the numerical experience reported in the literature, it is known that 
linear upwind methods lead to solutions with smeared layers, see also Section~\ref{sec:num_exam}.
This situation might explain that, to the best of our knowledge, the methods from \cite{Tabata77,BT81} are 
rarely used nowadays. So, their presentation will be kept brief, with an emphasis on 
the earlier method from \cite{Tabata77}.

In \cite{Tabata77}, a two-dimensional problem without reactive term is considered. These
assumptions will be relaxed below. In the first step of this method, one defines 
for an internal node $\bx_i^{}$ a so-called upwind simplex $K_i^{\rm up}$: 
$\bx_i^{}$ is a vertex of $K_i^{\rm up}$
and the straight half-line starting at $\bx_i^{}$ with direction $-\bb(\bx_i^{})$ 
intersects $K_i^{\rm up}$. If this line is parallel to a face (edge) $F$, then one chooses one element of $\omega_F^{}$ at random.
For nodes at the boundary, the construction is performed analogously. 
If $-\bb(\bx_i^{})$ points outside the domain, then $\bx_i^{}$ belongs to the inlet boundary, which means that
a Dirichlet condition is imposed at it, and, in turn, the test functions vanish at $\bx_i^{}$. This means that the 
upwind simplex can be chosen at random, as this choice will not affect the 
result. To simplify the presentation, we define the upwind simplex as the empty
set in this case.  If $\bb(\bx_i^{})=\bold0$, one uses an arbitrary
element of $\omega_i^{}$ as $K_i^{\rm up}$. The choice of the upwind 
element is motivated by the following observation. Let $\bx_j^{}, j\not= i$, be 
the other vertices of the simplex $K_i^{\rm up}$. By construction, it holds
that $|\sphericalangle (-\bb(\bx_i^{}),\bn_i)| < \pi/2$ and 
$\pi/2 \le |\sphericalangle (-\bb(\bx_i^{}),\bn_j)| < 3\pi/2$ for $j\neq i$,
where $\bn_i$ and $\bn_j$ are the outer unit normals to the facets of 
$K_i^{\rm up}$ opposite $\bx_i^{}$ and $\bx_j^{}$, respectively.
From \eqref{grad-basis-functions}, it follows that 
\begin{equation}\label{FEM-upwind-4}
\bb(\bx_i^{})\cdot\nabla\phi_i^{} |_{K_i^{\rm up}} \ge 0
\quad\textrm{and}\quad 
\bb(\bx_i^{})\cdot\nabla\phi_j^{} |_{K_i^{\rm up}}\le 0\quad \mbox{for}\,\,
j\neq i\,,
\end{equation}
which will be of major importance later. 
With these definitions, the upwind method reads as follows: Find $u_h^{}\in V_h^{}$ such that
$u_h^{}|_{\partial\Omega}^{}=i_h^{}g$, and
\begin{equation}\label{FEM-upwind-5}
\varepsilon(\nabla u_h^{},\nabla v_h^{})^{}+\sum_{j=1}^N
\left(\boldsymbol{b}(\bx_j^{})\cdot\nabla u_h^{}|_{K_j^{\rm up}}
\psi_j^{},\calL v_h^{}\right) +\sigma\,(u_h^{},v_h^{})_h^{}= (f,v_h^{})_h^{}\,,
\end{equation}
for all $ v_h^{}\in V_{h,0}^{}$, where $\psi_j^{}$ is the dual basis function defined in \eqref{FEM-upwind-2}, $\calL$ the lumping operator from \eqref{FEM-upwind-3}, 
and $(\cdot,\cdot)_h$ the lumped inner product defined in \eqref{eq:lumped_inner_prod}.
The term $\nabla u_h^{}|_{K_j^{\rm up}}$ is defined to be the zero vector if the 
upwind simplex is the empty set, otherwise it is a constant vector on $K_j^{\rm up}$.

The analysis of the method simplifies greatly if one rewrites the convective term. Noticing that
the dual basis functions $\psi_1^{},\ldots,\psi_N^{}$ are orthogonal in $L^2(\Omega)$
and using \eqref{eq:mass_matrix_l_c}, one can 
see that for every $v_h^{}\in V_h^{}$ the following holds
\begin{eqnarray*}\lefteqn{
\sum_{j=1}^N \left(\boldsymbol{b}(\bx_j^{})\cdot\nabla u_h^{}|_{K_j^{\rm up}} \psi_j^{},\calL v_h^{}\right)}\\
&=& \sum_{i,j=1}^N \bb(\bx_j^{})\cdot\nabla u_h^{}|_{K_j^{\rm up}}  v_h^{}(\bx_i^{}) (\psi_j^{},\psi_i^{})
= \sum_{i=1}^N \bb(\bx_i^{})\cdot\nabla u_h^{}|_{K_i^{\rm up}}  v_h^{}(\bx_i^{}) |D_i^{}| \nonumber\\
&=& \sum_{i=1}^N \bb(\bx_i^{})\cdot\nabla u_h^{}|_{K_i^{\rm up}}  v_h^{}(\bx_i^{}) (1,\phi_i^{}) 
= \sum_{i=1}^N(\bb(\bx_i^{})\cdot\nabla u_h^{}|_{K_i^{\rm up}} , \phi_i^{})\,v_h^{}(\bx_i^{})\,.
\end{eqnarray*}
Thus,  method \eqref{FEM-upwind-5} can be rewritten as follows: Find $u_h^{}\in V_h^{}$ such that
$u_h^{}|_{\partial\Omega}^{}=i_h^{}g$, and
\begin{equation*}
\varepsilon(\nabla u_h^{},\nabla v_h^{}) + 
\sum_{i=1}^N (\bb(\bx_i^{})\cdot\nabla u_h^{}|_{K_i^{\rm up}} , \phi_i^{})\,
v_h^{}(\bx_i^{})  +\sigma\,(u_h^{},v_h^{})_h^{}
= (f,v_h^{})_h^{}\,,
\end{equation*}
for all $v_h^{}\in V_{h,0}^{}$.  

The result below establishes well-posedness and the satisfaction of the
DMP. In addition,
this result also relaxes the hypotheses made on the mesh family from strictly acute to the XZ-criterion.

\begin{theorem}[DMP for the upwind finite element method]
\label{theo:upwind_tabata}Let us suppose that the mesh satisfies the 
XZ-criterion. Then, the matrix corresponding to the discrete problem
\eqref{FEM-upwind-5} is of non-negative type and hence the solution satisfies
the local DMP. In addition, the discrete problem \eqref{FEM-upwind-5} is well 
posed and then also the global DMP follows. 
\end{theorem}

\begin{proof} {We will show that 
$(\varepsilon\,\bbAd+\hat{\bbA}_{\rm c}^{} + \sigma \bbMl)^M$,
where 
\begin{equation*}
 \hat{\bbA}_{\rm c}^{}=(\hat{\cc}_{ij}^{})\quad \mbox{with} \quad \hat{\cc}_{ij}^{}:= (\bb(\bx_i^{})\cdot\nabla \phi_j^{}|_{K_i^{\rm up}} , \phi_i^{})\,,
\end{equation*}
is of non-negative type.}
From Corollary~\ref{cor:diff_reac_non_neg_mat} it is known that 
${(\varepsilon\,\bbAd+\sigma \bbMl)^M}$ is of non-negative type
if the mesh satisfies the XZ-criterion. 
Moreover, thanks to \eqref{FEM-upwind-4} and to the fact that the basis functions form 
a partition of unity on $K_i^{\rm up}$, one has for 
$i,j=1,\ldots, N$
\begin{equation*}
\hat{\cc}_{ii}^{}  \ge 0\;,\qquad \hat{\cc}_{ij}^{} \le 0\quad\textrm{for}\;i\not= j\;,\quad\textrm{and}\quad
\sum_{j=1}^N \hat{\cc}_{ij}^{} = 0\,.
\end{equation*}
Hence, $ \hat{\bbA}_{\rm c}^{}$ is also of non-negative type. It follows that
${(\varepsilon\,\bbAd+\hat{\bbA}_{\rm c}^{}+\sigma \bbMl)^M}$
is of non-negative type and since the diagonal entries of this matrix are positive, the method
satisfies the local DMP thanks to Theorem~\ref{thm:algebraic-local-DMP}.

Since $\varepsilon(\fl_{ij}^{})_{i,j=1}^M$  is of non-negative type and it is invertible (thanks to Remark~\ref{rem:diff_mat_non_sing}), 
 and $(\hat{\cc}_{ij}^{})_{i,j=1}^M,  (\sigma \tilde m_{ij})_{i,j=1}^M$ are of non-negative type,  an application of
\cite[Theorem~5.1]{Knobloch10} shows that $(\varepsilon\fl_{ij}^{}+\hat{\cc}_{ij}^{}+ \sigma \tilde  m_{ij}^{})_{i,j=1}^M$
is invertible, which, in turn, implies that  \eqref{FEM-upwind-5} has a unique solution. 
Finally,  an application of Theorem~\ref{thm:algebraic-global-DMP} leads to the satisfaction of  the global DMP. 
\end{proof}

Alternative versions of the upwind method for $\bbP_1$ finite elements have been proposed over the years. For example,   in \cite{BT81},
also for time-dependent convection-diffusion equations, a method was proposed motivated by the fact that the  exact solution satisfies a discrete analog 
of a mass conservation property if a special boundary 
condition is applied, see Section~\ref{sec:tcd_ext_methods_steady} for some details.
This is an additional feature compared with the method from \cite{Tabata77}. Domains $\Omega \subset \mathbb R^d$ and triangulations of weakly acute type
are considered in \cite{BT81}. Again, the barycentric cell $D_i^{}$ around a vertex $\bx_i^{}$ is constructed. 
Then, appropriate discrete fluxes $\beta_{ij}^{}$ across the individual parts of $\partial D_i^{}$
are defined, which is a technique from finite volume methods. 
The discrete convective term has the form 
\begin{equation*}
\sum_{i=1}^N \sum_{j \in S_i} \left(\beta_{ij}^+ u_h(\bx_i^{}) +\beta_{ij}^-
u_h(\bx_j^{}) \right) v_h(\bx_i^{}), 
\end{equation*}
with $S_i^{}$ defined in \eqref{eq:def_S_i}. The coefficients $\beta_{ij}^{}$ 
should satisfy several conditions and concrete choices are given in 
\cite{BT81}. The off-diagonal entries of the convection matrix are
always non-positive and, for a particular choice of the coefficients
$\beta_{ij}^{}$ specified in \cite{BT81},
the row sums of this matrix vanish and thus the convection matrix is of 
non-negative type. Under these assumptions, the statements of Theorem~\ref{theo:upwind_tabata} can be transferred
literally to the method from \cite{BT81}. 

One further upwind method, based on a slightly different choice of the 
domains for the dual basis, was presented in \cite{Kan78}. 
A proposal for partial upwinding can be found in \cite{Ikeda83}. For a unified presentation 
of upwind finite element methods and 
some numerical results we refer  to \cite{Knobloch06b}. 

The numerical analysis of several linear finite element upwind schemes can be found in 
\cite{Ikeda83}, in particular in \cite[Section~4.7]{Ikeda83} for the steady-state convection-diffusion 
equation ($\sigma = 0$) in two dimensions. The error analysis for one of 
the methods is presented in detail. For weakly acute triangulations, sufficiently small mesh width, and $u$ being regular enough, the estimate 
\[
\|u-u_h\|_{0,\infty,\Omega} \le Ch
\]
is proved, with $C$ being independent of $\varepsilon$. It is remarked that the same result holds 
true for the methods from \cite{BT81,Tabata77}. The $d$-dimensional convection-diffusion-reaction
equation 
is studied in \cite{BT81}, where the reaction coefficient is assumed to be constant and  
mass lumping is used for the reactive 
term. It is proved that there exists a positive constant $C$, which does not depend on $\varepsilon$, such 
that 
\[
\|i_hu-u_h\|_{0,\infty,\Omega} \le Ch \|u\|_{2,p,\Omega}, \quad p>d,
\]
if the reaction constant is sufficiently large.

\subsection{The edge-averaged  finite element method}\label{sec:xu-zikatanov_fem}

This section describes the method proposed in \cite{XZ99} and its main properties.

A part of the analysis will be performed under the assumption that the matrix
$\bbAdI$ is irreducible. Let us mention that if the mesh is connected (see 
Definition~\ref{Def:conditions-on-mesh}), then the diffusion matrix $\bbAd$ 
(including all boundary nodes) is irreducible, compare \cite[Rem.~2.3]{DDS05}. 
As shown in the same paper, this property does not necessarily imply the 
irreducibility of $\bbAdI$. Despite this, it needs to be considered that the 
example provided in \cite{DDS05} is rather pathological. In fact, in the same 
paper it is already noted that refining the mesh once removes the reducibility 
of $\bbAdI$. Thus, from the available experience, one might state that the 
reducibility of $\bbAdI$ is an exceptional situation that can be cured by mesh 
refinements (with the resulting mesh being still very coarse). For this reason, 
assuming that the matrix $\bbAdI$ is irreducible does not seem to be a big loss 
of generality.  
 
The following rewriting of the discrete Laplacian matrix $\bbAd$, which was at 
the heart of the proof of Theorem~\ref{theo:Delaunay-mesh-laplacian}, will 
be fundamental for the derivation of the method. Consider any $u_h^{},v_h^{}\in
V_h^{}$ and any $K\in\calT_h^{}$, and denote by $\calI_K$ the index set of
nodes contained in $K$. Since the local diffusion matrices are
symmetric and have zero row sums, a direct calculation using $\delta_E^{}$ 
defined in Section~\ref{Sec:triangulations} yields
\begin{eqnarray*}
&&(\nabla u_h^{},\nabla v_h^{})_K  =  
\sum_{i,j\in\calI_K} \fl_{ij}^K u_i^{}v_j^{} =
\sum_{i,j\in\calI_K} \fl_{ij}^K u_i^{}(v_j^{}-v_i^{})\\
&& \hspace*{13mm}=
\sum_{i,j\in\calI_K,i<j} \fl_{ij}^K (u_i^{}-u_j^{})(v_j^{}-v_i^{}) 
= - \sum_{i,j\in\calI_K,i<j}
  \fl_{ij}^K\delta_{E_{ij}^{}}u_h^{}\,\delta_{E_{ij}^{}}v_h^{}\,,
\end{eqnarray*}
where we use the notation $u_i=u_h(\bx_i^{})$, $v_i=v_h(\bx_i^{})$,
$i=1,\dots,N$.
This formula is a sum over the edges of $K$, where every edge appears exactly 
once. Hence, denoting
\begin{equation*}
   \lambda_E^K=  \frac{|\kappa_E^K|\cot\theta_E^K}{d(d-1)}\,,
\end{equation*}
it follows from \eqref{eq:lK_ij} that
\begin{equation*}
   (\nabla u_h^{},\nabla v_h^{})_K = \sum_{E \in \calE_K^{}} \lambda_E^K\,
   \delta_E^{}u_h^{}\,\delta_E^{}v_h^{}\,.
\end{equation*}
Consider any $\boldsymbol{a}\in\mathbb R^d$ and set 
$u_h^{}(\bx)=\boldsymbol{a}\cdot\bx$. Then $u_h^{}\in V_h^{}$, 
$\nabla u_h^{}=\boldsymbol{a}$, and 
$\delta_E^{}u_h^{}=h_E^{}\,\boldsymbol{a}\cdot\boldsymbol{t}_E^{}$ for any
$E\in\calE_h^{}$. Thus, the previous identity implies that
\begin{equation}\label{FEM-XZ-1}
   (\boldsymbol{a},\nabla v_h^{})_K 
   = \sum_{E \in \calE_K^{}} h_E^{}\,\lambda_E^K\,
   \boldsymbol{a}\cdot\boldsymbol{t}_E^{}\,\delta_E^{}v_h^{}\qquad
   \forall\,\,\boldsymbol{a}\in\mathbb R^d,\,v_h^{}\in V_h^{},\,
   K\in\calT_h^{}\,.
\end{equation}

Another fundamental ingredient in the derivation of the method is the 
consideration of a conservative form of the convective term. We will present, 
just for simplicity, the case $\sigma=0$, although the case $\sigma>0$ is also 
treated in \cite{XZ99} using a mass-lumping strategy. Then, applying 
integration by parts, the bilinear form $a(\cdot,\cdot)$ defined in
\eqref{bilinear-a} satisfies
\begin{equation}\label{FEM-XZ-3}
   a(u,v)=(\varepsilon\nabla u-\bb\, u, \nabla v)\qquad
   \forall\,\,u\in H^1(\Omega),\,v\in H^1_0(\Omega)\,.
\end{equation}
The quantity $\bJ(u)=\varepsilon\nabla u-\bb\, u$ is called total flux.

A further ingredient is a function $\chi_E^{}$ defined, for each edge 
$E\in\calE_h^{}$, by
\begin{equation*}
   \frac{\partial\chi_E^{}}{\partial\boldsymbol{t}_E^{}}=
   -\frac{\bb\cdot\boldsymbol{t}_E^{}}{\varepsilon}\,,
\end{equation*}
which determines $\chi_E^{}$ uniquely up to an additive constant. This
definition implies that, for $u\in C^1(\overline\Omega)$, one has
\begin{equation*}
   \frac{\partial(e^{\chi_E^{}}u)}{\partial\boldsymbol{t}_E^{}}
   =\frac1\varepsilon\,e^{\chi_E^{}}\,\bJ(u)\cdot\boldsymbol{t}_E^{}\,,
\end{equation*}
which leads to
\begin{equation*}
\delta_E^{}\left(e^{\chi_E^{}} u\right) =
\frac1{\varepsilon}\int_E e^{\chi_E^{}}\bJ(u)\cdot\boldsymbol{t}_E^{}\ ds\,.
\end{equation*}
Thus, approximating $\bJ(u)$ on $K\subset\omega_E$ by a constant vector 
$\bJK(u)$ leads to the relation
\begin{equation}\label{FEM-XZ-5}
   \bJK(u)\cdot\boldsymbol{t}_E^{}\approx\varepsilon\,
   \frac{\delta_E^{}(e^{\chi_E^{}}u)}{\int_E e^{\chi_E^{}}ds}\,.
\end{equation}

Now, using the approximations $\bJK(u)$ in \eqref{FEM-XZ-3} with 
$v=v_h^{}\in V_{h,0}^{}$ and applying \eqref{FEM-XZ-1} and \eqref{FEM-XZ-5} 
leads to
\begin{align*}
   a(u,v_h^{})&\approx \sum_{K\in\calT_h^{}} 
   (\bJK(u),\nabla v_h^{})_K =\sum_{K\in\calT_h^{}} 
   \sum_{E\in\calE_K^{}} h_E^{}\,\lambda_E^K\,\bJK(u)\cdot\boldsymbol{t}_E^{}\,
   \delta_E^{}v_h^{}\\
   &\approx\sum_{K\in\calT_h^{}}
   \sum_{E\in\calE_K^{}} \,\lambda_E^K\,\tilde{\varepsilon}_E^{}(\bb)\,
   \delta_E^{}(e^{\chi_E^{}}u)\,\delta_E^{}v_h^{}\,,
\end{align*}
where
\begin{equation*}
   \tilde{\varepsilon}_E^{}(\bb)= 
   \frac{\varepsilon\,h_E^{}}{\int_E e^{\chi_E^{}}\ ds}
\end{equation*}
is the harmonic average of $\varepsilon\,e^{-\chi_E^{}}$ on the edge $E$. This
suggests to introduce the bilinear form
\begin{equation*}
   a_h^{}(u_h^{},v_h^{})=\sum_{E\in\calE_h^{}} \Bigg(\sum_{K\subset\omega_E}
   \lambda_E^K\Bigg)\tilde{\varepsilon}_E^{}(\bb)\,
   \delta_E^{}(e^{\chi_E^{}}u_h^{})\,\delta_E^{}v_h^{}\,,
\end{equation*}
which leads to the following  Xu--Zikatanov, or edge-averaged, finite element 
method: Find $u_h^{}\in V_h^{}$, such that 
$u_h^{}|_{\partial\Omega}^{}=i_h^{}g$, and
\begin{equation}\label{XZ-method}
   a_h^{}(u_h^{},v_h^{}) = (f,v_h^{})\qquad\forall\, v_h^{}\in V_{h,0}^{}\,.
\end{equation}

It is worth stressing that if one replaces $\chi_E^{}$ by $\chi_E^{}+c$, $c\in\mathbb{R}$, then,
in exact arithmetic, 
the bilinear form $a_h^{}(\cdot,\cdot)$ is not affected. Thus, the fact that $\chi_E^{}$ is defined up to an additive constant has no
effect in method~\eqref{XZ-method}.  It is observed in \cite{BCC98} that in two dimensions 
the edge-averaged finite element method is equivalent to the Scharfetter--Gummel finite volume scheme.

For analyzing  \eqref{XZ-method}, first two properties  of its system matrix 
will be proven. More precisely, we define the matrix 
$(\bbA)^M=(a_{ij}^{})^{i=1,\dots,M}_{j=1,\dots,N}$ given by
$a_{ij}^{}=a_h^{}(\phi_j^{},\phi_i^{})$. Then, the following results hold.

\begin{lemma}[Properties of the system matrix of \eqref{XZ-method}]
\label{lem:X-Z-irreducible}If the matrix $\bbAdI$ is irreducible, then
the matrix  $\bbAI=(a_{ij}^{})_{i,j=1}^M$ is irreducible, too. In addition, if
the XZ-condition \eqref{XZ-criterion} is satisfied, the 
diagonal entries of $\bbAI=(a_{ij}^{})_{i,j=1}^M$ are positive.
\end{lemma}

\begin{proof} 
Consider any $i,j\in\{1,\dots,M\}$, $i\neq j$. If $\bx_i^{}$, $\bx_j^{}$ are
not endpoints of the same edge, then $a_{ij}^{}=0=\fl_{ij}^{}$. Otherwise, in
view of \eqref{eq:l_ij},
\begin{equation}\label{eq:XZ-a_ij}
   a_{ij}^{}=-\Bigg(\sum_{K\subset\omega_{E_{ij}}} \lambda_{E_{ij}^{}}^K\Bigg) 
   \tilde{\varepsilon}_{E_{ij}^{}}^{}(\bb)\,e^{\chi_{E_{ij}^{}}(\bx_j^{})}
   =\fl_{ij}^{}\,
   \tilde{\varepsilon}_{E_{ij}^{}}^{}(\bb)\,e^{\chi_{E_{ij}^{}}(\bx_j^{})}\,.
\end{equation}
The positivity of the last two factors implies that $a_{ij}^{}=0$ if and only if
$\fl_{ij}^{}=0$, which proves the first part of the lemma. Furthermore, again in
view of \eqref{eq:l_ij},
\begin{equation*}
   a_{ii}^{}=\sum_{E\in\calE_h^{}:\,\bx_i^{}\in E}\Bigg(\sum_{K\subset\omega_E}
   \lambda_E^K\Bigg)\,\tilde{\varepsilon}_E^{}(\bb)\,e^{\chi_E^{}(\bx_i^{})} 
   =-\sum_{j\in S_i}\,\fl_{ij}^{}\,
   \tilde{\varepsilon}_{E_{ij}^{}}^{}(\bb)\,e^{\chi_{E_{ij}^{}}(\bx_i^{})}
\end{equation*}
for any $i\in\{1,\dots,M\}$. If \eqref{XZ-criterion} holds, then 
\eqref{eq:l_ij} implies that $\fl_{ij}^{}\le0$ for all $j\neq i$ and since
$\fl_{ii}^{}=|\phi_i^{}|_{1,\Omega}^2>0$, it follows from 
\eqref{eq:A_d_row_sum} that $\fl_{ij}^{}<0$ for at least one index $j\neq i$.
Therefore, $a_{ii}^{}>0$, which finishes the proof.
\end{proof}

\begin{theorem}[M-matrix property of the system matrix of the edge-averaged FEM]
\label{thm:XZ-DMP}Let the mesh be of XZ-type and let the matrix
$\bbAdI$ be irreducible. Then the system matrix of the discretization 
\eqref{XZ-method} is an M-matrix.
\end{theorem}

\begin{proof} First, note that the matrix $\bbAI$ is irreducible
by Lemma~\ref{lem:X-Z-irreducible}. We extend the matrix $(\bbA)^M$ to an
$N\times N$ matrix by setting $a_{ij}^{}=a_h^{}(\phi_j^{},\phi_i^{})$ for all
$i,j=1,\dots,N$. Then the representation \eqref{eq:XZ-a_ij} holds if 
$j\in S_i$, and $a_{ij}^{}=0$ if $j\not\in S_i\cup\{i\}$. Since $\calT_h^{}$ 
satisfies the XZ-condition \eqref{XZ-criterion}, one observes 
immediately that $a_{ij}^{}\le 0$ if $j\neq i$ and $i\le M$ or $j\le M$.
Moreover, from the definition of $\delta_E$, it follows directly that 
\begin{equation*}
\sum_{i=1}^N a_{ij}^{} = a_h^{}(\phi_j^{},1)=0\,,\qquad j=1,\dots,N\,.
\end{equation*}
Since the matrix $\bbAd$ is irreducible, there is $\tilde i\in\{M+1,\dots,N\}$ 
and $\tilde j\in\{1,\dots,M\}$ such that $a_{\tilde i \tilde j}<0$, which
implies that at least one column sum of $\bbAI$ is strictly positive (while the
remaining ones are at least non-negative). Hence, $\mathbb{A}_{\mathrm{I}}^T$ 
is irreducibly diagonally dominant and then, according to 
\cite[Theorem~3.27]{Var00}, $\mathbb{A}_{\mathrm{I}}^T$ is an M-matrix.
Consequently, also $\bbAI$ is an M-matrix and the theorem follows from
Remark~\ref{rem:A-AI}.
\end{proof}

The last result generalizes the result presented in \cite[Lemma~6.2]{XZ99} 
where it is shown that the bilinear form  $a_h^{}(\cdot,\cdot)$ from 
\eqref{XZ-method} satisfies an inf-sup condition for sufficiently small $h$, 
and thus showing well-posedness of \eqref{XZ-method} for sufficiently refined 
meshes (although we should mention that this generalization is already hinted 
in \cite[Remark~6.1]{XZ99}).

\begin{remark}
The M-matrix property proved in Theorem~\ref{thm:XZ-DMP} immediately implies
the positivity preservation of the discrete problem \eqref{XZ-method}, i.e., if
the right-hand side $f$ and the boundary condition $g$ are non-negative, then
also the discrete solution $u_h$ is non-negative. However, the M-matrix
property does not imply the local or global DMP. The validity of the DMPs
follows from Theorems \ref{thm:algebraic-local-DMP} and
\ref{thm:algebraic-global-DMP} if the convection field $\bb$ is constant since 
then the validity of \eqref{zero-row-sum} can be shown. However, in general,
the validity of the local and global DMPs is open.
\hspace*{\fill}$\Box$\end{remark}

The discrete problem \eqref{XZ-method} is well-posed under the assumptions of 
Theorem~\ref{thm:XZ-DMP} since the system matrix is an M-matrix. In more general situations 
the well-posedness for sufficiently small mesh sizes is shown in \cite{XZ99}. That paper 
presents also an error estimate of the form 
\[
\|i_hu-u_h\|_{1,\Omega} \le Ch \left(\sum_{K\in\calT_h} |\bJ(u)|_{1,p,K}^2 + \sum_{K\in\calT_h} |\sigma u|_{1,r,K}^2\right)^{1/2},
\]
assuming that the terms on the right-hand side are well defined for 
sufficiently large values of $p$ and $r$, where the concrete values depend on the dimension.

\section{Nonlinear stabilized discretizations of the steady-state problem}\label{Sec:cd_nonlinear_methods}

One common feature of all the discretizations presented in the previous
section is that they add global stabilizing terms, that is, the methods
modify the formulation in the whole domain (equivalently, they modify
every row in the system matrix). As a consequence, linear stabilized
methods that respect the DMP provide, in general, very diffused solutions.
Now, as it was mentioned earlier, in order to prove the DMP, one only needs
to analyze the rows of the matrix associated with nodes where an extremum is
attained. So, ideally, a method should modify only these rows of the matrix 
in order to have a good performance. The selection of these rows depends on 
the solution itself, thus such a method is necessarily nonlinear. This is why 
in this section we present several nonlinear finite element methods for the 
convection-diffusion equation that respect the DMP. In contrast to linear
methods, some of the nonlinear approaches even satisfy the DMP on general 
meshes, i.e., without any assumptions on the angles in the meshes.

\subsection{The Mizukami--Hughes method}
\label{Sec:cd_nonlinear_methods_MH}

The Mizukami--Hughes method is a nonlinear Petrov--Galerkin method proposed in 
\cite{MH85} and improved and further developed in
\cite{Kno06,Kno07,Knobloch10}. The idea of the method is to create an upwind
effect by means of solution-dependent weighting functions which guarantee that
the approximate solution satisfies a linear system with a matrix of
non-negative type. Up to the best of our knowledge, this is the first nonlinear 
DMP-satisfying method proposed for
the numerical solution of \eqref{steady-strong}. We shall confine ourselves to 
the two-dimensional case and to $\sigma=0$. Extensions to $\sigma>0$ and to 
three space dimensions can be found in \cite{Kno06}.

For any interior node $\bx_i$, $i\in\{1,\dots,M\}$, we introduce the weighting 
function
\begin{displaymath}
   {\widetilde\phi}_i=\phi_i
   +\sum_{K\subset\omega_i}C^K_{i}\,\chi_K\,.
\end{displaymath}
Here $\chi_K$ denotes the characteristic functions of mesh cells $K$ (i.e., $\chi_K=1$ 
in $K$ and $\chi_K=0$ elsewhere) and $C^K_{i}$ are constants which will be
determined later. The discretization of the convection-diffusion equation
reads as follows: Find $u_h^{}\in V_h^{}$ such that
$u_h^{}|_{\partial\Omega}^{}=i_h^{}g$, and
\begin{equation}\label{eq:mh_discr}
   \varepsilon\,(\nabla u_h,\nabla \phi_i)
   +(\bb_h\cdot\nabla u_h,{\widetilde\phi}_i)=(f,{\widetilde\phi}_i)\,,\qquad 
   i=1,\dots,M\,,
\end{equation}
where $\bb_h$ is a piecewise constant approximation of $\bb$. We shall also
use the notation $\bb_K:=\bb_h|_K^{}$ for $K\in\calT_h^{}$. The simplest
choice is to set $\bb_K$ equal to the value of $\bb$ at the barycenter of $K$.

The definition of the constants $C^K_{i}$ is based on the requirement that the
local convection matrix $\MHAK$ with entries
\begin{equation}\label{eq:mhak}
   \MHaKij=(\bb_K\cdot\nabla\phi_j,\widetilde\phi_i)_K^{}\,,\qquad 
   i=1,\dots,M\,,\,\,\, j=1,\dots,N\,,\,\,\, \bx_i,\bx_j\in K\,,
\end{equation}
is of non-negative type. In \cite{MH85}, it was further required that
\begin{equation}\label{eq:mh_cond}
   C^K_{i}\ge-\third\qquad\forall\,\,i\in\{1,\dots,N\}\,,\,\,
                                             \bx_i\in K\,,\qquad\qquad
   \sum_{\parbox{4ex}{\centerline{$\scriptstyle i=1$}\vspace{-1ex}
                      \centerline{$\scriptstyle
\bx_i\in K$}}}^NC^K_{i}=0\,.
\end{equation}

\begin{figure}[t!]
\centerline{\includegraphics[width=0.4\textwidth]{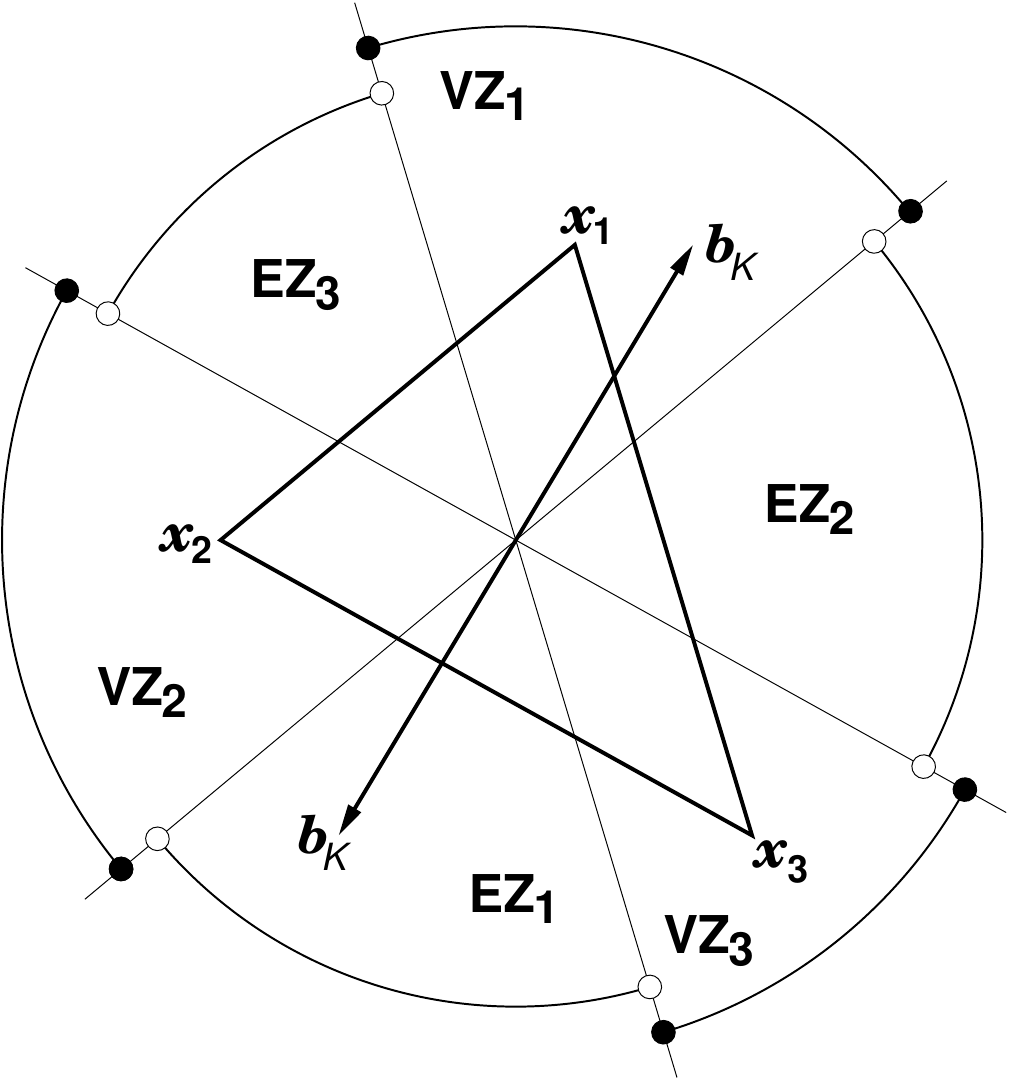}}
\caption{Definition of edge zones and vertex zones.}\label{fig:mh_fig}
\end{figure}

As we will see, the choice of the constants $C^K_{i}$ significantly depends on
the direction of the convection vector $\bb_K$ with respect to the edges of $K$.
To characterize the direction of $\bb_K$, we decompose any triangle $K$ into 
vertex zones and edge zones by
drawing lines parallel to the edges of $K$ which all intersect at the
barycenter of $K$, see Fig.~\ref{fig:mh_fig}. Denoting the vertices of $K$ by
$\bx_1$, $\bx_2$ and $\bx_3$, the set containing the vertex $\bx_i$, $i=1,2,3$, will be
called vertex zone ${\rm VZ}_i$. The remaining three sets are called edge zones
and the edge zone opposite the vertex $\bx_i$ will be denoted by ${\rm EZ}_i$.
The common part of the boundaries of two adjacent zones is included in the
respective vertex zone. The fact that the vector $\bb_K$ points from the
barycenter of $K$ into ${\rm VZ}_i$ or ${\rm EZ}_i$ will be shortly expressed
by $\bb_K\in{\rm VZ}_i$ or $\bb_K\in{\rm EZ}_i$, respectively. Without loss
of generality, one may assume that the vertices of $K$ are numbered in such a
way that $\bb_K\in{\rm VZ}_1$ or $\bb_K\in{\rm EZ}_1$ as depicted in
Fig.~\ref{fig:mh_fig}.

Using \eqref{grad-basis-functions}, it is easy to see that
\begin{align*}
     \bb_K\in {\rm VZ}_1\quad&\Longleftrightarrow\quad
     \bb_K\cdot\nabla\phi_1>0\,,\quad
     \bb_K\cdot\nabla\phi_2\le0\,,\quad
     \bb_K\cdot\nabla\phi_3\le0\,,\\
     \bb_K\in {\rm EZ}_1\quad&\Longleftrightarrow\quad
     \bb_K\cdot\nabla\phi_1<0\,,\quad
     \bb_K\cdot\nabla\phi_2>0\,,\quad
     \bb_K\cdot\nabla\phi_3>0\,,
\end{align*}
where we write $\nabla\phi_i$ instead of $\nabla\phi_i|_K^{}$ for simplicity.
Note that $\MHAK$ has always zero row sums, so that one has to assure only that
$\MHaKij\le0$ for $i\neq j$. Since
\begin{equation*}
   \MHaKij=\bb_K\cdot\nabla\phi_j|_K^{}\,|K|\,(\third+C^K_{i})\,,
\end{equation*}
one observes that, if $\bb_K\in {\rm VZ}_1$, this condition on $\MHAK$ can be 
easily satisfied by setting
\begin{equation}\label{eq:mh_vz1}
   C^K_{1}=\twothirds\,,\qquad C^K_{2}=C^K_{3}=-\third\,.
\end{equation}
However, if $\bb_K\in {\rm EZ}_1$, it is generally not possible to choose
the constants $C^K_1,C^K_2,C^K_3$ in such a way that \eqref{eq:mh_cond} holds and 
$\MHAK$ is of non-negative type.

Nevertheless, Mizukami and Hughes \cite{MH85} made the important observation 
that $u$ still solves the equation \eqref{steady-strong} if $\bb$ is replaced 
by any function $\tilde\bb$ such that $\tilde\bb-\bb$ is orthogonal to 
$\nabla u$. This suggests to define the constants $C^K_{i}$ in such a way that 
the matrix $\MHAK$ is of non-negative type for $\bb_K$ replaced by a function 
$\tilde\bb_K$ pointing into a vertex zone and preserving the product 
$\bb_K\cdot\nabla u_h|_K^{}$. Note that the local convection matrix $\MHAK$ will 
be still defined using $\bb_K$ and the vector $\tilde{\bb}_K$ is used only for 
defining the constants $C^K_{i}$. Since the constants $C^K_{i}$ depend through 
$\tilde{\bb}_K$ on the unknown discrete solution $u_h$, the resulting discrete 
problem is nonlinear.

Let us assume that $\bb_K\in{\rm EZ}_1$ and $\bb_K\cdot\nabla u_h|_K^{}\neq0$
and let $\bw\neq{\bf 0}$ be a vector orthogonal to $\nabla u_h|_K^{}$. We
introduce the sets
\begin{displaymath}
   V_k=\{\alpha\in{\mathbb{R}}\,;\,\,
                           \bb_K+\alpha\,\bw\in{\rm VZ}_k\}\,,\qquad k=2,3\,.
\end{displaymath}
The vectors $\bb_K+\alpha\,\bw$ play the role of $\tilde{\bb}_K$ mentioned
above. Is is easy to see that $V_2\cup V_3\neq\emptyset$.
Mizukami and Hughes showed that, depending on $V_2$ and $V_3$, the following
values of the constants $C^K_{i}$ should be used:
\begin{eqnarray}
   V_2\neq\emptyset\quad\&\quad V_3=\emptyset\quad\,\,
   &\Longrightarrow&\quad C^K_2=\twothirds\,,\quad C^K_1=C^K_3=-\third\,,
   \label{eq:mh_v2}\\[1ex]
   V_2  = \emptyset\quad\&\quad V_3\neq\emptyset\quad\,\,
   &\Longrightarrow&\quad C^K_3=\twothirds\,,\quad C^K_1=C^K_2=-\third\,,
   \label{eq:mh_v3}\\[1ex]
   V_2\neq\emptyset\quad\&\quad V_3\neq\emptyset\quad\,\,
   &\Longrightarrow&\quad C^K_1=-\third\,,\quad 
    C^K_2+C^K_3=\third\,,\label{eq:mh_v2v3}\\
   &&\quad C^K_2>-\third\,,\quad C^K_3>-\third\,.\nonumber
\end{eqnarray}
It was observed in \cite{Kno06} that the definition of $C^K_{i}$'s proposed in
\cite{MH85} for the case \eqref{eq:mh_v2v3} depends on the orientation of
$\bb_K$ and $\nabla u_h|_K^{}$ in a discontinuous way. This may deteriorate the
quality of the discrete solution and prevent the nonlinear iterative process
from converging. Therefore, another definition of these constants was
introduced in \cite{Kno06} for which the dependence on the orientation of 
$\bb_K$ and $\nabla u_h|_K^{}$ is continuous. To avoid technical digressions,
we refer to \cite{Kno06} for details.

It was also demonstrated in \cite{Kno06} that, in some cases, the solutions of 
the original Mizukami--Hughes method do not approximate boundary layers in a 
correct way. Therefore, if $\bb_K$ points into an edge zone, it was proposed to
set
\begin{equation}\label{eq:mh_bnd}
   C^K_{1}=C^K_{2}=C^K_{3}=-\third
\end{equation}
for any mesh cell $K\in\calT_h^{}$ having a node on $\partial\Omega$. 
Except for cases where these mesh cells form a strip along the boundary of an 
approximately constant width, the definition \eqref{eq:mh_bnd} is used also for
mesh cells whose all nodes are connected by edges to nodes on $\partial\Omega$.
The choice \eqref{eq:mh_bnd} suppresses the influence of the Dirichlet boundary
condition on the approximate solution inside $\Omega$, which may be important
if $K$ lies in the numerical boundary layer.

If $\bb_K\in{\rm EZ}_1$, $\bb_K\cdot\nabla u_h|_K^{}=0$ and
\eqref{eq:mh_bnd} is not used, then one sets
\begin{equation}\label{eq:mh_ez1}
   C^K_{1}=-\third\,,\qquad C^K_{2}=C^K_{3}=\sixth\,.
\end{equation}
Finally, one sets $C^K_{1}=C^K_{2}=C^K_{3}=0$ if $\bb_K=\boldsymbol0$.

Although the system matrix of \eqref{eq:mh_discr} is in general not of 
non-negative type, one can prove that, for meshes of XZ-type, the solution 
vector solves a linear
system of the form \eqref{system-1}--\eqref{system-2} with a non-singular 
matrix of non-negative type, which implies that the solution of the 
Mizukami--Hughes method satisfies local and global DMPs.

\begin{theorem}[Matrix of non-negative type for the Mizukami--Hughes method]
\label{thm:mh_nonnegative_type}Let the mesh $\calT_h^{}$ be of
XZ-type. Then the solution of the Mizukami--Hughes method
\eqref{eq:mh_discr} satisfies a linear system of the type
\eqref{system-1}--\eqref{system-2} with $f_i=(f,{\widetilde\phi}_i)$,
$i=1,\dots,M$, and $g_{i-M}=g(\bx_i)$, $i=M+1,\dots,N$, such that the
corresponding system matrix $\bbA$ given in \eqref{system-matrix} is 
of non-negative type and its block $\bbAI$ is invertible.
\end{theorem}

\begin{proof} 
Let $\bu$ be the coefficient vector corresponding to the solution of
\eqref{eq:mh_discr}. We shall show that, for any $K\in\calT_h^{}$, there is a 
matrix $\MHTAK$ of non-negative type such that
\begin{equation}
   \MHTAK\,\bu^K=\MHAK\,\bu^K\,,\label{eq:mh_97}
\end{equation}
where $\MHAK$ is defined by \eqref{eq:mhak} and $\bu^K$ consists of the 
components of $\bu$ corresponding to nodes of $K$. If $\bb_K=\boldsymbol0$ or 
$C^K_{i}$'s are defined in \eqref{eq:mh_vz1} or \eqref{eq:mh_bnd}, one can take
$\MHTAK=\MHAK$. In case of \eqref{eq:mh_ez1} which is used if
$\bb_K\cdot\nabla u_h|_K^{}=0$, one can set $\MHTAK=0$. It remains to 
define $\MHTAK$ in cases when the constants $C^K_{i}$ are defined by
\eqref{eq:mh_v2}--\eqref{eq:mh_v2v3}, which assumes that $\bb_K\in{\rm EZ}_1$
and $\bb_K\cdot\nabla u_h|_K^{}\neq0$. First, we introduce some auxiliary
notation. If, for some $k\in\{2,3\}$, the set $V_k$ is non-empty, we choose
$\alpha_k\in V_k$ and define the matrix $\MHTAKs{k}$ with entries
\begin{equation*}
   \tilde{c}^{K,k}_{ij}=(\bb_K+\alpha_k\,\bw)\cdot\nabla\phi_j|_K^{}\,|K|\,
   (\third+C^{K,k}_i)\,,\qquad i,j=1,2,3\,\,\, (\bx_i\in\Omega)\,,
\end{equation*}
where $C^{K,2}_i$ are defined as in \eqref{eq:mh_v2} and $C^{K,3}_i$ as in
\eqref{eq:mh_v3}. If $V_k=\emptyset$, we set $\MHTAKs{k}=0$. Then the 
matrices $\MHTAKs{2}$ and $\MHTAKs{3}$ are of non-negative type and
hence also
\begin{displaymath}
   \MHTAK:=(\third+{C}^K_2)\,\MHTAKs{2}+
                   (\third+{C}^K_3)\,\MHTAKs{3}
\end{displaymath}
is of non-negative type. Since $\bw\cdot\nabla u_h|_K=0$ and
\begin{displaymath}
   (\third+{C}^K_2)(\third+C^{K,2}_i)+
   (\third+{C}^K_3)(\third+C^{K,3}_i)=\third+C^K_{i}\,,\quad i=1,2,3\,,
\end{displaymath}
one obtains \eqref{eq:mh_97}.

The matrices $\MHAK$ and $\MHTAK$ are assembled to $M\times N$ matrices
$\MHAC$ and $\MHTAC$ for which $\MHAC\,\bu=\MHTAC\,\bu$ and $\MHTAC$ is of 
non-negative type. Since $\bu$ corresponds to the solution of 
\eqref{eq:mh_discr}, one has
$(\varepsilon\,(\bbAd)^M+\MHAC)\,\bu=\boldsymbol f$ with
${\boldsymbol f}=(f_1,\dots,f_M)$ introduced in the formulation of the 
theorem. As $\calT_h^{}$ is of XZ-type, the matrix $(\bbAd)^M$ is of
non-negative type according to Theorem~\ref{theo:Delaunay-mesh-laplacian}. 
Thus $\bu$ also satisfies the linear system 
$(\varepsilon\,(\bbAd)^M+\MHTAC)\,\bu=\boldsymbol f$
and the matrix $\mathbb{A}^M:= \varepsilon\,(\bbAd)^M+\MHTAC$ is of
non-negative type. Since the block $\bbAdI$ of $\bbAd$ is invertible
(cf.~Remark~\ref{rem:diff_mat_non_sing}), it
follows that also $\bbAI$ is invertible (see \cite[Theorem 5.1]{Knobloch10}). 
This finishes the proof.
\end{proof}

As discussed in \cite{Knobloch10}, the Mizukami--Hughes method corresponds to
the discretization of the convective term by standard upwind differencing. This 
is appropriate if the diffusion $\varepsilon$ is small in comparison to $\bb$.
However, if this is not the case, such a discretization leads to a low
accuracy since too much artificial diffusion is introduced. Therefore, in
\cite{Knobloch10}, the constants $C^K_{i}$ were defined in such a way that the 
matrix $\widetilde\varepsilon\,\MHAdK+\MHAK$ is of non-negative type, where
$\MHAdK$ is the local diffusion matrix and 
$\widetilde\varepsilon\in(0,\varepsilon)$ is close to $\varepsilon$. This does 
not change the method much in the convection-dominated case but improves the 
accuracy if $\varepsilon$ is not small.

To the best of our knowledge, there are no error estimates available for
the Mizukami--Hughes method. Also, the solvability of the nonlinear problem seems 
to be still an open problem.

\subsection{Burman--Ern Methods}\label{Sec:cd_nonlinear_methods_BE}

In this section we will present the finite element method, based on a continuous interior penalty idea,
presented in \cite{BE05}. The analysis of this method requires the mesh to be of XZ-type, so
we will assume that throughout this section. In the work \cite{BE05} the method is presented
with two stabilizations, namely, a linear one (e.g., SUPG or CIP), and the nonlinear
stabilizing term responsible for the DMP. To keep the discussion brief, we will start discussing
 the case of the reduced method, that is, the method only adds the nonlinear
stabilization to the Galerkin formulation. 
The proof of the local DMP  (cf. Theorem~\ref{Thm:Gen-Non-DMP-1}) is achieved by proving that the nonlinear
problem satisfies the weak DMP property (cf.
Definition~\ref{def:DMP-criterion}). 
So, as a motivation for the definition of the method we will now suppose that 
$u_h^{}(\bx_i^{})<0$, $i\in\{1,\dots,M\}$, is a local minimum in 
$\omega_i^{}$ and will bound $a(u_h^{},\phi_i^{})$. 
Thanks to the fact that the mesh is of XZ-type one has $\ell_{ij}^{}\le 0$ for
all $j\neq i$ (cf. Theorem~\ref{theo:Delaunay-mesh-laplacian}), and 
consequently
\begin{equation}\label{eq:BE-diff}
(\nabla u_h^{},\nabla \phi_i^{})  = \sum_{j\in S_i^{}} \ell_{ij}^{}(u_h^{}(\bx_j^{})- u_h^{}(\bx_i^{})) \le 0 \,.
\end{equation}
In addition, if the function $u_h^{}$ changes sign inside
$K\subset\omega_i^{}$, using a Taylor expansion at a zero of $u_h$, one gets
\begin{equation*}
(u_h^{},\phi_i^{})_K^{}\le \frac{|K|}{d+1}\,h_K^{}\,\big|\nabla u_h^{}|_K^{}\big|\,. 
\end{equation*}
If $u_h^{}\le 0$ in $K$ then one just bounds $(u_h^{},\phi_i^{})_K^{}\le 0$.
The convective term is bounded in a similar way leading to
\begin{equation*}
(\bb\cdot\nabla u_h^{}+ \sigma u_h^{},\phi_i^{}) \le \frac{1}{d+1}\sum_{K\subset\omega_i^{}}
\big( \|\bb\|_{0,\infty,K}^{}+\sigma\,h_K^{}\,\big)
|K|\,\big|\nabla u_h^{}|_K^{}\big|\,.
\end{equation*}

Next,  to bound the gradient of $u_h^{}$ in the last inequality one uses that 
$u_h^{}(\bx_i^{})$ is a local minimum and then the following bound holds 
(see \cite[Lemma~2.7]{BE05} for the proof):
\begin{equation*}
|\nabla u_h^{}|_K^{}|\le \sum_{F\in\calF_i^{}}|\llbracket\nabla
u_h^{}\rrbracket_F^{}|\qquad\forall\,\,K\subset\omega_i\,,
\end{equation*}
which leads to
\begin{align}\label{eq:a-bound}
a(u_h^{},\phi_i^{})&\le \frac{1}{d+1}\,\sum_{F\in\calF_i^{}}\,
\sum_{K\subset\omega_i^{}} 
\big( \|\bb\|_{0,\infty,K}^{}+\sigma\,h_K^{}\,\big)\,|K|
\,|\llbracket\nabla u_h^{}\rrbracket_F^{}|\\
&\le \frac1{d+1}\,
\sum_{F\in\calF_i^{}}\big( \|\bb\|_{0,\infty,\tilde\omega_F^{}}^{}+
\rho\,\sigma\,h_F^{}\big)
 |\omega_i|\,|\llbracket\nabla u_h^{}\rrbracket_F^{}|\,,\nonumber
\end{align}
where we used the fact that, in view of \eqref{BE-4}, one has
$h_K^{}\le\rho\,h_F^{}$ for any $K\subset\omega_i^{}$ and $F\in\calF_i^{}$.
Since $|\omega_i|\le\Omega_d\,(\max_{K\subset{\omega_i^{}}^{}}h_K^{})^d$, where
$\Omega_d$ is the measure of the unit ball in $\mathbb R^d$, one has
$|\omega_i|\le\Omega_d\,\rho^d\,h_F^d$ for any $F\in\calF_i^{}$. Using the mesh
regularity, one gets $|\omega_i|\le C\,\rho^d\,h_F\,|F|$, which gives
\begin{equation}\label{BE-a-bound}
   a(u_h^{},\phi_i^{})\le\frac{C\rho^d}{d+1}\,
   \sum_{F\in\calF_i^{}}\big( \|\bb\|_{0,\infty,\tilde\omega_F^{}}^{}+
   \rho\,\sigma\,h_F^{}\big)
   h_F^{}|F|\,|\llbracket\nabla u_h^{}\rrbracket_F^{}|\,.
\end{equation}

From the discussion above, one sees that in order to prove the DMP, one needs 
to control a term related to the jumps of the gradients of the discrete solution across the facets containing the local 
extrema. Motivated by this observation, in \cite{BE05} the following method is 
proposed:  Find $u_h^{}\in V_h^{}$
such that $u_h^{}|_{\partial\Omega}^{}=i_h^{}g$, and
\begin{equation}\label{BE-1}
a(u_h^{},v_h^{})+j_h^{}(u_h^{};v_h^{}) = (f,v_h^{})\qquad\forall\, v_h^{}\in
V_{h,0}^{}\,.
\end{equation}
Here,  $j_h^{}(\cdot;\cdot)$ is a stabilizing form given by
\begin{align}
j_h^{}(u_h^{};v_h^{})=& \,c_\rho^{}\sum_{F\in\calF_I^{}}\big(
\|\bb\|_{0,\infty,\tilde\omega_F^{}}^{}+
\rho\,\sigma\,h_F^{}\big)\,h_F^{}
\,\left(|\llbracket\nabla u_h^{}\rrbracket_F^{}|\, , 
b_F^{}(u_h^{};v_h^{})\,\right)_F,\label{BE-2}\\
b_F^{}(u_h^{};v_h^{}) =&
\,\sum_{E\in\calE_F^{}}h_E^{}\,{\rm sign}(\nabla u_h^{}\cdot\bt_E^{})\nabla v_h^{}\cdot\bt_E^{}\,.
\label{BE-3}
\end{align}
The parameter $c_\rho^{}>0$ depends on the mesh regularity through the quantity $\rho$. 
Using a regularized problem and Brouwer's fixed-point theorem in \cite{BE05} it is proven that \eqref{BE-1} admits at least one solution. 
Under the hypothesis that the mesh is of XZ-type, the following result
regarding the local DMP can be shown.

\begin{theorem}[Local DMP for the Burman--Ern method]\label{thm:BE05_DMP}Let 
us suppose that the mesh is of XZ-type. Then, if $c_{\rho}^{}$ is sufficiently 
large, the nonlinear form $j_h^{}(\cdot;\cdot)$ satisfies the weak DMP property 
if $\sigma>0$ and the strong DMP property if $\sigma=0$.  Consequently,
method \eqref{BE-1} satisfies the local DMP from
Theorem~\ref{Thm:Gen-Non-DMP-1}.
\end{theorem}

\begin{proof} Let us suppose that $u_h\in V_h$ has a strict local
minimum at the interior node  $\bx_i^{}$.  Then, for any $F\in\calF_i^{}$,
one has
\begin{equation}\label{eq:b_F-est}
b_F^{}(u_h^{};\phi_i^{})= \sum_{j\in S_i:\,E_{ij}\subset F}\,
	{\rm sign}\left(u_h(\bx_j)-u_h(\bx_i)\right)\left(\phi_i(\bx_j)-\phi_i(\bx_i)\right)=-(d-1)\,,
\end{equation}
since $\mbox{\rm card}\{j\in S_i:\,E_{ij}\subset F\}=d-1$. This
implies that
\begin{equation*}
j_h^{}(u_h^{};\phi_i^{}) \le -c_\rho^{}\,\sum_{F\in\calF_i^{}}\big(
\|\bb\|_{0,\infty,\tilde\omega_F^{}}^{}+
\rho\,\sigma\,h_F^{}\big)\,h_F^{}\,|F|\,|\llbracket\nabla u_h^{}\rrbracket_F^{}|\,.
\end{equation*}
Thus, combining this last bound with \eqref{BE-a-bound} (which was derived
for $u_h(\bx_i)<0$ but holds also for $u_h(\bx_i)\ge0$ if $\sigma=0$)  gives
\begin{equation*}
a(u_h^{},\phi_i^{})+j_h^{}(u_h^{};\phi_i^{})\le \left(\frac{C\rho^d}{d+1}-c_\rho^{}\right)
\sum_{F\in\calF_i^{}}\big( \|\bb\|_{0,\infty,\tilde\omega_F^{}}^{}+
\rho\,\sigma\,h_F^{}\big)\,h_F^{}\,|F|\,|\llbracket\nabla u_h^{}\rrbracket_F^{}|\,,
\end{equation*}
and the proof follows choosing $c_{\rho}^{}$ large enough provided that
$\|\bb\|_{0,\infty,\tilde\omega_F^{}}^{}+\rho\,\sigma\,h_F^{}>0$ for all
$F\in\calF_i^{}$. If this is not the case, one can employ the fact that the
previous inequality holds with the term $\varepsilon\,(\nabla u_h^{},\nabla
\phi_i^{})$ on the right-hand side. When deriving \eqref{eq:a-bound}, this term
was estimated by \eqref{eq:BE-diff}. However, since now $u_h^{}(\bx_i^{})$ is a
strict local minimum, it follows from \eqref{eq:BE-diff} that $(\nabla
u_h^{},\nabla \phi_i^{})$ is negative and hence can be estimated by
$-\sum_{F\in\calF_i^{}}\alpha_F^{}\left|\llbracket \nabla
u_h^{}\rrbracket_F^{}\right|$ with suitable positive constants $\alpha_F^{}$.
This finishes the proof.
\end{proof}

\begin{figure}[t!]
\centerline{\includegraphics[width=0.25\textwidth]{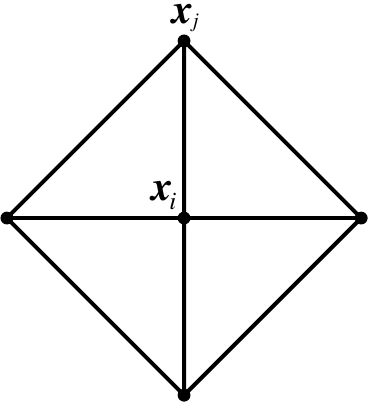}
\hspace*{10mm}
\includegraphics[width=0.25\textwidth]{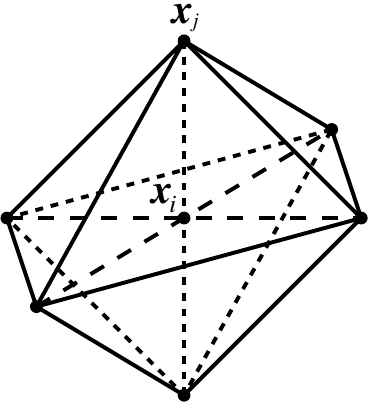}
}
\caption{Examples of patches $\omega_i^{}$ in 2d and 3d.}\label{fig:patches}
\end{figure}

\begin{remark}The validity of the global DMP seems to be open for method
\eqref{BE-1} since, in general, the stabilizing form $j_h(\cdot;\cdot)$ defined
in \eqref{BE-2} does not 
allow to prove the strong and weak DMP properties for non-strict extrema
formulated in Definition~\ref{def:DMP-criterion-non-strict}. To see this, let
us consider the patches $\omega_i^{}$ depicted in Fig.~\ref{fig:patches}. Let 
us decompose $\omega_i^{}$ into the sets
\begin{equation*}
   \omega_i^1=\cup\{ K\subset\omega_i^{}:\bx_j^{}\in K\}\,,\qquad
   \omega_i^2=\cup\{ K\subset\omega_i^{}:\bx_j^{}\not\in K\}\,.
\end{equation*}
Let $u_h^{}\in V_h$ be such that $u_h(\bx_j)\neq u_h(\bx_i)$ and 
$u_h(\bx_k)=u_h(\bx_i)$ for any vertex $\bx_k\in\omega_i^2$. Then
$u_h\in\bbP_1(\omega_i^1)$, $u_h$ is constant in $\omega_i^2$, it is not
constant in $\omega_i^{}$, and attains a local extremum at $\bx_i$. 
Consequently, $b_F^{}(u_h^{};v_h^{}) =0$ for any 
$F\subset\omega_i^2$ and any $v_h\in V_h$. On the other hand, for any $F\in
\calF_i^{}$ such that $\bx_j\in F$, one has 
$\llbracket\nabla u_h^{}\rrbracket_F^{}=\bold0$ since $\nabla u_h^{}$ is
constant in $\omega_i^1$. Thus, $j_h(u_h,\phi_i^{})=0$ which means that 
the term $j_h$ cannot be used to enforce the strong or weak DMP property for
non-strict extrema at the node $\bx_i$.
\hspace*{\fill}$\Box$\end{remark}

An alternative definition, hinted in \cite[Theorem~3.5]{BE05}, and developed further in \cite[Section~2.4]{BH04},
can be
obtained by replacing $|\llbracket\nabla u_h^{}\rrbracket_F^{}|$ in
\eqref{BE-2} by
\begin{equation*}
   m_F(u_h^{})=\max_{F^\prime\in\calF_I^{}:\,F^\prime\subset\omega_F^{}}\,
   |\llbracket\nabla u_h^{}\rrbracket_{F^\prime}|\,.
\end{equation*}
Then
\begin{equation}\label{BE-5}
j_h^{}(u_h^{};v_h^{})= \,c_\rho^{}\sum_{F\in\calF_I^{}}\big(
\varepsilon+\|\bb\|_{0,\infty,\tilde\omega_F^{}}^{}\,h_F^{}+
\rho\,\sigma\,h_F^2\big)
\left(m_F(u_h^{})\, , b_F^{}(u_h^{};v_h^{})\,\right)_F.
\end{equation}
For this stabilizing term, one can prove also the DMP properties for non-strict
extrema formulated in Definition~\ref{def:DMP-criterion-non-strict}.

\begin{theorem}[DMP for \eqref{BE-5}]\label{thm:BE06_DMP}Let 
us suppose that the mesh is of XZ-type. Then, if $c_{\rho}^{}$ is sufficiently 
large, the nonlinear form $j_h^{}(\cdot;\cdot)$ defined in \eqref{BE-5} 
satisfies the weak DMP property for non-strict extrema if $\sigma>0$ and the 
strong DMP property for non-strict extrema if $\sigma=0$. Consequently,
method \eqref{BE-1} with $j_h^{}(\cdot;\cdot)$ from \eqref{BE-5} satisfies both
the local and the global DMPs from Theorem~\ref{Thm:Gen-Non-DMP-1}.
\end{theorem}

\begin{proof} Let us suppose that $u_h^{}\in V_h$ has a local
minimum at the interior node  $\bx_i^{}$. Then, for any $F\in\calF_i^{}$,
one has
\begin{equation*}
   b_F^{}(u_h^{};\phi_i^{})=\sum_{j\in S_i:\,E_{ij}\subset F}\,
   {\rm sign}\left(u_h(\bx_j)-u_h(\bx_i)\right)
   \left(\phi_i(\bx_j)-\phi_i(\bx_i)\right)\le0\,.
\end{equation*}
Consider any $F\in\calF_i^{}$. If 
$\llbracket\nabla u_h^{}\rrbracket_F^{}\neq\bold0$, then there exists a vertex
$\bx_j^{}\in\omega_F^{}$ such that $u_h^{}(\bx_j^{})\neq u_h^{}(\bx_i^{})$. Let 
$F^{\prime\prime}\subset\omega_F^{}$ be a facet such that $\bx_i^{},\bx_j^{}\in
F^{\prime\prime}$. Then 
\begin{equation*}
   b_{F^{\prime\prime}}(u_h^{};\phi_i^{})\le
	-{\rm sign}\left(u_h(\bx_j)-u_h(\bx_i)\right)=-1\,.
\end{equation*}
Since $F\subset\omega_{F^{\prime\prime}}$, one gets
\begin{equation*}
   |\llbracket\nabla u_h^{}\rrbracket_F^{}|
   \le -m_{F^{\prime\prime}}(u_h^{})\,b_{F^{\prime\prime}}(u_h^{};\phi_i^{})\,.
\end{equation*}
If $\llbracket\nabla u_h^{}\rrbracket_F^{}=\bold0$, then this inequality
holds with any $F^{\prime\prime}\in\calF_i^{}$ satisfying
$F^{\prime\prime}\subset\omega_F^{}$ since the right-hand side is nonnegative.
Hence one finds that
\begin{align}\label{eq:bound_mF}
   \sum_{F\in\calF_i^{}}|\llbracket\nabla u_h^{}\rrbracket_F^{}|
   &\le-\sum_{F\in\calF_i^{}}m_{F^{\prime\prime}}(u_h^{})\,b_{F^{\prime\prime}}(u_h^{};\phi_i^{})\\
   &\le-(2\,d-1)\,\sum_{F\in\calF_i^{}} m_F(u_h^{})\,b_F(u_h^{};\phi_i^{})
   \nonumber
\end{align}
as the number of facets $F\in\calF_i^{}$ satisfying
$F\subset\omega_{F^{\prime\prime}}$ for a given $F^{\prime\prime}\in\calF_i^{}$
is $2\,d-1$. Using this estimate in the first inequality of \eqref{eq:a-bound}
(which was derived for $u_h(\bx_i)<0$ but holds also for $u_h(\bx_i)\ge0$ if
$\sigma=0$) and performing the same manipulations as used to derive
\eqref{BE-a-bound}, one obtains
\begin{equation*}
   a(u_h^{},\phi_i^{})\le-C\,\rho^d\,\frac{2\,d-1}{d+1}\,
   \sum_{F\in\calF_i^{}}\big( \|\bb\|_{0,\infty,\tilde\omega_F^{}}^{}+
   \rho\,\sigma\,h_F^{}\big)\,
   h_F^{}\left(m_F(u_h^{})\, , b_F^{}(u_h^{};\phi_i^{})\,\right)_F,
\end{equation*}
where $C$ is the same constant as in \eqref{BE-a-bound}. Thus,
$a(u_h^{},\phi_i^{})+j_h^{}(u_h^{};\phi_i^{})\le\frac12j_h^{}(u_h^{};\phi_i^{})$
if $c_\rho^{}\ge 2\,C\,\rho^d\,(2\,d-1)/(d+1)$. According to \eqref{eq:bound_mF}, one has
\begin{equation*}
j_h^{}(u_h^{};\phi_i^{})\le-\frac{c_\rho^{}}{2\,d-1}\
\min_{F\in\calF_i^{}}\left\{\big(\varepsilon+\|\bb\|_{0,\infty,\tilde\omega_F^{}}^{}\,h_F^{}+
\rho\,\sigma\,h_F^2\big)\,|F|\right\}\,
\sum_{F\in\calF_i^{}}|\llbracket\nabla u_h^{}\rrbracket_F^{}|\,,
\end{equation*}
which completes the proof. 
\end{proof}

\begin{remark} The methods just analyzed need the mesh to be of XZ-type. To
avoid this restriction, in \cite{BE04} the following method was proposed for 
the Poisson problem: Find $u_h^{}\in V_h^{}$ such that 
$u_h|_{\partial\Omega}^{}=i_hg$, and
\begin{equation}\label{BE04-1}
(\nabla u_h^{},\nabla v_h^{})+\delta\sum_{F\in\calF_I^{}} \left(|\llbracket
\nabla u_h^{}\rrbracket_F^{}|\,, b_F^{}(u_h^{};v_h^{})\,\right)_F =(f,v_h^{})\qquad\forall\, v_h^{}\in V_{h,0}^{}\,,
\end{equation}
where $b_F^{}$ is defined as in \eqref{BE-3} and $\delta>0$. Then, for $\delta >
\frac{1}{d(d-1)}$, method~\eqref{BE04-1} satisfies the strong DMP property
for any mesh. In fact, the main argument of the proof is the following 
observation from \cite{BE04}: regardless of the mesh,
\begin{equation}\label{BE04-2}
(\nabla u_h^{},\nabla \phi_i^{}) = \sum_{F\in \calF_i^{}} \big(
\llbracket\nabla u_h^{}\rrbracket_F^{}\cdot\bn_F,
\phi_i^{})_F^{}= \sum_{F\in \calF_i^{}}  \frac{|F|}{d}\llbracket\nabla u_h^{}\rrbracket_F^{} \cdot\bn_F\,,
\end{equation}
where $\bn_F$ is the unit normal vector to $F$ in the direction
corresponding to the orientation of the jump $\llbracket\cdot\rrbracket_F^{}$.
So, if $u_h^{}$ has a strict local minimum at an interior node 
$\bx_i^{}$,  it follows from \eqref{eq:b_F-est} that
\begin{equation*}
(\nabla u_h^{},\nabla \phi_i^{})
+\delta\sum_{F\in\calF_I^{}}\left(|\llbracket
\nabla u_h^{}\rrbracket_F^{}|\,, b_F^{}(u_h^{};\phi_i^{})\,\right)_F\le
\sum_{F\in \calF_i^{}}  \left(\frac{1}{d}
-\delta(d-1)\right)\,|F|\,|\llbracket\nabla u_h^{}\rrbracket_F^{}|\,.
\end{equation*}
Thus, for $\delta > \frac{1}{d(d-1)}$ 
\eqref{BE04-1} satisfies the strong DMP criterion. 

The main difference between \eqref{BE04-1} and \eqref{BE-1} resides on the size of the stabilization term. In fact, only considering
the powers of $h$ involved, the stabilization given in \eqref{BE04-1} is one size larger than the one from \eqref{BE-1}, as  \eqref{BE04-1} is designed
to match the behavior of the diffusion matrix given by \eqref{BE04-2}. So, even if this term is positive (as it would happen if a mesh that is
not of XZ-type is used), then the stabilization is large enough to compensate
for that. Even if in \cite{BE04} an extension to the 
convection-diffusion equation has been studied, this variant does not seem to have been applied to convection-dominated problems
in later years.
\hspace*{\fill}$\Box$\end{remark}

Method \eqref{BE-1} is the simplest form of a Burman--Ern method that respects the local DMP. 
In the presence of dominating convection, sometimes it is recommended to first add a linear stabilization term to stabilize
the convection, and only then to add a nonlinear stabilization to ensure the satisfaction of the DMP. With this objective in mind, 
this approach was pursued in \cite{BE05} by using a linear stabilization which can be given by the SUPG or CIP stabilization. We
now summarize briefly the results proven for the latter option. The CIP stabilizing term is defined as follows (see, e.g., \cite{BH04}) 
\begin{equation*}
s_h(u_h^{},v_h^{})=\sum_{F\in\calF_I^{}}\gamma_{\rm
cip}^{}\,\|\bb\|_{0,\infty,\Omega}^{}h_F^2\,(\llbracket\nabla
u_h^{}\rrbracket_F^{}, \llbracket \nabla v_h^{}\rrbracket_F^{})_F^{}\,,
\end{equation*}
where $\gamma_{\rm cip}^{}>0$. Using this stabilizing term, the following 
stabilized method is proposed in \cite{BE05}: Find $u_h^{}\in V_h^{}$
such that $u_h^{}|_{\partial\Omega}^{}=i_h^{}g$, and
\begin{equation}\label{BE-1-complete}
a(u_h^{},v_h^{})+s_h(u_h^{},v_h^{})+j_h^{}(u_h^{};v_h^{}) = (f,v_h^{})
\qquad{\forall\, v_h^{}\in V_{h,0}^{}}\,,
\end{equation}
with $j_h^{}(\cdot;\cdot)$ being a combination of \eqref{BE-2} and
\eqref{BE-5}.
The corresponding analogue of Theorem~\ref{thm:BE06_DMP} was proven for \eqref{BE-1-complete} in \cite[Theorem~3.5]{BE05}. For the diffusion-dominated regime, i.e., with the assumption
$ch\le \varepsilon$ for some appropriate constant $c$, the following error estimate
appears as a corollary of \cite[Theorem~3.10]{BE05}:
\begin{align}\label{BE-Error}
&\varepsilon^{\frac{1}{2}}|u-u_h^{}|_{1,\Omega}^{}+\sigma^{\frac{1}{2}}\|u-u_h^{}\|_{0,\Omega}^{}
+\|h^{\frac{1}{2}}\bb\cdot\nabla (u-u_h^{})\|_{0,\Omega}^{}\\
&\hspace*{10mm}+ s_h(u-u_h^{},u-u_h^{})^{\frac{1}{2}}
\le C \left(\varepsilon+\|\bb\|_{0,\infty,\Omega}^{}\,h
+\sigma\,h^2\right)^{\frac12}h\,\|u\|_{2,\Omega}^{}\,,\nonumber
\end{align}
where $C>0$ is independent of $h$ and all the physical parameters, provided that the exact solution $u$ belongs to $H^2(\Omega)$. 

The combination of linear and nonlinear stabilizations has two main effects in
this context. First, the addition of the linear stabilization term $s_h(\cdot,\cdot)$ 
allows for the extra control on the convective term appearing in \eqref{BE-1}, which is responsible for the estimate
\eqref{BE-Error}. This control is not possible to achieve if only the nonlinear stabilization $j_h^{}(\cdot,\cdot)$ is used. The second main effect
is computational. It  can be observed that, while the nonlinear stabilization
$j_h(\cdot,\cdot)$ is local (in the sense that it is active  mostly in the vicinity of extrema and layers),
the  linear stabilization term $s_h(\cdot,\cdot)$ is global, and thus it  helps dampening oscillations that appear away from the layers.

\begin{remark} Finally, it is worth mentioning that the works reviewed in this
section were not the first effort that was made in this direction by the
authors. In fact, in their previous paper \cite{BE02} the authors proposed a
nonlinear diffusion method that, under the assumption of acute meshes,
satisfies the global DMP. To improve the convergence of the nonlinear solver,
absolute values in the nonlinear terms were regularized, which however leads
to a violation of the DMP. Comprehensive numerical tests of three variants of 
the methods from \cite{BE02} can be also found in \cite{JK07,JK08}. In
particular, in \cite{JK08}, the authors did not succeed to solve the respective
nonlinear problems in a number of cases.
\hspace*{\fill}$\Box$\end{remark}

\subsection{Algebraic Flux Correction methods} 
\label{Sec:cd_nonlinear_methods_afc}

Algebraic flux correction (AFC) methods belong to the class of algebraically 
stabilized schemes which have been intensively developed in recent years, see,
e.g., 
\cite{BB17,BJK17,GNPY14,Kno21,Kuz06,Kuz07,Kuz09,Kuz12a,Kuz12,KS17,KT04,LKSM17}.
In contrast to the methods discussed in the previous sections, the 
stabilization is not introduced in a variational form but the starting point is
the system of linear algebraic equations corresponding to the Galerkin FEM 
discretization. Then, a nonlinear algebraic term is added to the linear system 
in order to enforce a DMP without an excessive smearing of the layers.

Let $\bbAN$ be the matrix corresponding to the standard Galerkin FEM 
\eqref{eq:galerkin_fem} with Neumann boundary conditions, i.e., 
\begin{equation}\label{eq:Galerkin_matrix}
   \bbAN=\varepsilon\bbAd+\bbAc+\sigma\bbMc\,.
\end{equation}
We will also consider a lumping of the reaction term in
\eqref{eq:galerkin_fem}, which leads to a matrix given by
\begin{equation}\label{eq:lumped_Galerkin_matrix}
   \bbAN=\varepsilon\bbAd+\bbAc+\sigma\bbMl\,.
\end{equation}
The discrete problem is then equivalent to the system \eqref{system-1},
\eqref{system-2}, where $f_i^{}=(f,\phi_i^{})$ for $i=1,\ldots,M$ and 
$g_{i-M}^{}=g(\bx_i^{})$ for $i=M+1,\ldots,N$. To derive an AFC scheme, first 
a symmetric artificial diffusion matrix $\bbD=(d_{ij})_{i,j=1}^N$ is introduced 
by
\begin{equation}
  d_{ij}=-\max\{0,a_{ij},a_{ji}\}\quad
  \mbox{for} \ i\ne j,\qquad\quad
  d_{ii}=-\sum_{j=1,j \neq i}^N d_{ij}\,.\label{eq:matrix_D}
\end{equation}
Hence $\mathbb D$ has zero row and column sums and the matrix $\bbAN+\bbD$ is of
non-negative type. Thus, replacing $\bbAN$ by $\bbAN+\bbD$ in \eqref{system-1},
one obtains the stabilized problem 
\begin{displaymath}
   (\bbAN+\bbD)^M\bu=\boldsymbol f
\end{displaymath}
satisfying the DMP (with ${\boldsymbol f}=(f_1,\dots,f_M)^T$). However,
like for the similar linear artificial diffusion method of
Section~\ref{Sec:cd_linear_methods_iso_diff}, the added artificial
diffusion is usually too large and leads to an excessive
smearing of layers. Therefore, it is necessary to restrict the artificial
diffusion to regions where the solution changes abruptly. Since these regions
are not known a priori, this will again lead to a nonlinear method.

The original derivation of the AFC method, e.g., in \cite{Kuz07}, is performed 
in such a way that first the term $(\bbD\bu)_i^{}$ is added to both sides of
\eqref{system-1} leading to
\begin{equation}\label{eq:equiv_system}
   (\bbAN+\bbD)^M\bu={\boldsymbol f}+\bbD^M\bu\,,
\end{equation}
and then the identity
\begin{displaymath}
   (\bbD\bu)_i^{}=\sum_{j=1}^N\,f_{ij}\qquad\mbox{with}\qquad
   f_{ij}=d_{ij}\,(u_j-u_i)
\end{displaymath}
is used.
The quantities $f_{ij}$ are called fluxes since they can be interpreted as
quantities which correspond to the intensity of the flow of $u$ between the 
nodes $\bx_i$ and $\bx_j$, see also the explanation of the concept of fluxes 
at the beginning of Section~\ref{sec:tcd_femfct}. It turns out that spurious 
oscillations in the approximate solution can be suppressed by damping the
above-introduced fluxes $f_{ij}$ appearing on the right-hand side of
\eqref{eq:equiv_system}. This damping is often called limiting and it
is achieved by multiplying the fluxes by solution-dependent correction factors
$\alpha_{ij}\in[0,1]$ called limiters. This leads to the nonlinear algebraic 
problem
\begin{align}
   \sum_{j=1}^N\,a_{ij}\,u_j
   +\sum_{j=1}^N\,(1-\alpha_{ij}(\bu))\,d_{ij}\,(u_j-u_i)
   &= f_i^{}\qquad\textrm{for}\; i=1,\ldots, M\,,\label{eq:afc-1}\\
   u_i^{} &= g_{i-M}^{}\qquad\textrm{for}\; i=M+1,\ldots, N\,. \label{eq:afc-2}
\end{align}
It is assumed that
\begin{equation}
   \alpha_{ij}=\alpha_{ji}\,,\qquad i,j=1,\dots,N\,,\label{eq:alpha_symm}
\end{equation}
and that, for any $i,j\in\{1,\dots,N\}$, the function
$\alpha_{ij}(\bu)(u_j-u_i)$ is a continuous function of
$\bu\in{\mathbb R}^N$. A theoretical analysis of the AFC scheme
\eqref{eq:afc-1}, \eqref{eq:afc-2} concerning the solvability, local DMP and
error estimation can be found in \cite{BJK16}; see also \cite{ABR17,Jha21} for
a posteriori error estimators.

The symmetry condition \eqref{eq:alpha_symm} is particularly important for 
several reasons. First, it guarantees that the resulting method is conservative.
Second, it implies that the matrix corresponding to the term arising from the
AFC is positive semidefinite. This shows that this term really enhances the
stability of the method and enables to estimate the error of the approximate
solution, see \cite{BJK16}. Finally, it was demonstrated in \cite{BJK15} that,
without the symmetry condition \eqref{eq:alpha_symm}, the nonlinear algebraic 
problem \eqref{eq:afc-1}, \eqref{eq:afc-2} is not solvable in general.

Recently, motivated by \cite{BB17}, a generalization of \eqref{eq:afc-1} was 
proposed in \cite{Kno21} by introducing the matrix 
$\bbB(\bu)=(b_{ij}(\bu))_{i,j=1}^N$ given by
\begin{align}
  &b_{ij}(\bu)=-\max\{0,(1-\alpha_{ij}(\bu))\,a_{ij},
                        (1-\alpha_{ji}(\bu))\,a_{ji}\}\quad
  \mbox{for} \ i\ne j,\label{eq:matrix_B1}\\
  &b_{ii}(\bu)=-\sum_{j=1,j \neq i}^N b_{ij}(\bu).\label{eq:matrix_B2}
\end{align} 
Then, instead of \eqref{eq:afc-1}, \eqref{eq:afc-2}, the following
algebraically stabilized problem is considered
\begin{align}
   \sum_{j=1}^N\,a_{ij}\,u_j+\sum_{j=1}^N\,b_{ij}(\bu)\,(u_j-u_i)
   &= f_i^{}\qquad\textrm{for}\; i=1,\ldots, M\,,\label{eq:asm-1}\\
   u_i^{} &= g_{i-M}^{}\qquad\textrm{for}\; i=M+1,\ldots, N\,.\label{eq:asm-2}
\end{align}
Under condition \eqref{eq:alpha_symm}, both algebraic problems,
\eqref{eq:afc-1}, \eqref{eq:afc-2} and \eqref{eq:asm-1}, \eqref{eq:asm-2}, are
equivalent. However, the advantage of \eqref{eq:asm-1}, \eqref{eq:asm-2} is
that the symmetry condition \eqref{eq:alpha_symm} is no longer necessary. Note
that the matrix $\bbB(\bu)$ is symmetric, has nonpositive off-diagonal entries
and has zero row and column sums. These properties imply that
\begin{equation*}
   \sum_{i,j=1}^N\,v_i\,b_{ij}(\bu)\,(v_j-v_i)
   =-\frac12\,\sum_{i,j=1}^N\,b_{ij}(\bu)\,(v_j-v_i)^2\ge0
   \quad\forall\,\,\bu,\bv\in{\mathbb R}^N\,.
\end{equation*}
Thus, the matrix ${\mathbb B}(\bu)$ is positive semidefinite for any 
$\bu\in{\mathbb R}^N$.

To write the above algebraic problem in a variational form, we denote
\begin{equation*}
   d_h(w;z,v)
   =\sum_{i,j=1}^N\,b_{ij}(w)\,z(\bx_j)\,v(\bx_i)
   \qquad\forall\,\,w,z,v\in C(\overline\Omega)\,,
\end{equation*}
with $b_{ij}(w):=b_{ij}(\{w(\bx_i)\}_{i=1}^N)$. Then 
\begin{equation}\label{eq:def_dh}
   d_h(w;\phi_j^{},\phi_i^{})=b_{ij}(w)
   \qquad\forall\,\,w\in C(\overline\Omega),\,i,j=1,\dots,N\,,
\end{equation}
and \eqref{eq:asm-1},
\eqref{eq:asm-2} is equivalent to problem \eqref{nonlinear-generic-2}, where
$a(\cdot,\cdot)$ is defined by \eqref{bilinear-a} in case of $\bbAN$ given
by \eqref{eq:Galerkin_matrix} and by \eqref{bilinear-a-lumped}
if $\bbAN$ given by \eqref{eq:lumped_Galerkin_matrix} is considered. 
The property \eqref{eq:matrix_B2} immediately implies the validity of
\eqref{eq:d_h_row_sums}. Since the matrix ${\mathbb B}(\bu)$ is positive 
semidefinite, the form $d_h$ also satisfies \eqref{eq:pos_semidef_dh}. Finally,
since $a_{ij}=a_{ji}=0$ if $j\not\in S_i\cup\{i\}$, one has
\begin{equation}\label{eq:d_h_sparse2}
   d_h^{}(w;\phi_j^{},\phi_i^{})=0\qquad\forall\,\,w\in C(\overline\Omega),\,
   j\not\in S_i\cup\{i\},\,i=1,\dots,N\,,
\end{equation}
so that \eqref{eq:d_h_sparse} always holds.

Of course, the properties of an algebraically stabilized scheme significantly
depend on the choice of the limiters $\alpha_{ij}$. Their design principles 
often originate from the time-dependent case where they should guarantee the 
positivity preservation, see Section~\ref{sec:tcd_femfct}. In the steady case,
a standard limiter is the Kuzmin limiter proposed in \cite{Kuz07} which was
thoroughly investigated in \cite{BJK16}. To define the limiter of
\cite{Kuz07}, one first computes, for $i=1,\dots,M$,
\begin{align}
   &P_i^+=\sum_{\mbox{\parbox{8mm}{\scriptsize\centerline{$j\in S_i$}
   \centerline{$a_{ji}\le a_{ij}$}}}}\,f_{ij}^+\,,\quad\,\,\,
   P_i^-=\sum_{\mbox{\parbox{8mm}{\scriptsize\centerline{$j\in S_i$}
   \centerline{$a_{ji}\le a_{ij}$}}}}\,f_{ij}^-\,,\label{eq:kuzmin_p}\\[1mm]
   &Q_i^+=-\sum_{j\in S_i}\,f_{ij}^-\,,\quad\,\,\,
   Q_i^-=-\sum_{j\in S_i}\,f_{ij}^+\,,\label{eq:kuzmin_q}
\end{align}
where $f_{ij}=d_{ij}\,(u_j-u_i)$, $f_{ij}^+=\max\{0,f_{ij}\}$, and 
$f_{ij}^-=\min\{0,f_{ij}\}$. We recall that $d_{ij}$ is defined in 
\eqref{eq:matrix_D} using the matrix $\bbAN$ from \eqref{eq:Galerkin_matrix} or
\eqref{eq:lumped_Galerkin_matrix}. Also the matrix entries appearing in
\eqref{eq:kuzmin_p} are taken from this matrix. Then, one defines
\begin{equation}\label{eq:R_i_plus_minus}
   R_i^+=\min\left\{1,\frac{Q_i^+}{P_i^+}\right\},\quad
   R_i^-=\min\left\{1,\frac{Q_i^-}{P_i^-}\right\},\qquad
   i=1,\dots,M\,.
\end{equation}
If $P_i^+$ or $P_i^-$ vanishes, one sets $R_i^+=1$ or $R_i^-=1$, respectively.
At Dirichlet nodes, these quantities are also set to be $1$, i.e.,
\begin{equation}\label{eq:R-Zalesak}
   R_i^+=1\,,\quad R_i^-=1\,,\qquad i=M+1,\dots,N\,.
\end{equation}
Furthermore, one sets
\begin{equation}\label{eq:alpha-definition}
        \widetilde\alpha_{ij}=\left\{
        \begin{array}{cl}
                R_i^+\quad&\mbox{if}\,\,\,f_{ij}>0\,,\\
                1\quad&\mbox{if}\,\,\,f_{ij}=0\,,\\
                R_i^-\quad&\mbox{if}\,\,\,f_{ij}<0\,,
        \end{array}\right.\qquad\qquad
        i,j=1,\dots,N\,.
\end{equation}
Finally, one defines
\begin{equation}\label{eq:kuzmin_symm_alpha}
   \alpha_{ij}=\alpha_{ji}=\widetilde\alpha_{ij}\qquad\mbox{if}\quad
   a_{ji}\le a_{ij}\,,\qquad i,j=1,\dots,N\,.
\end{equation}

\begin{theorem}[DMP for the AFC scheme with Kuzmin limiter]
\label{thm:kuzmin_dmp}Let
\begin{equation}\label{eq:assumption_min}
   \min\{a_{ij},a_{ji}\}\le0\qquad
   \forall\,\,i=1,\dots,M\,,\,\,j=1,\dots,N\,,\,\,i\neq j\,.
\end{equation}
Then the AFC scheme \eqref{eq:afc-1}, \eqref{eq:afc-2} with the Kuzmin limiter
defined by \eqref{eq:kuzmin_p}--\eqref{eq:kuzmin_symm_alpha} satisfies the 
algebraic DMP property formulated in
Definition~\ref{def:algebraic_DMP_property} and also the algebraic DMP property 
for non-strict extrema from
Definition~\ref{def:nonstrict_algebraic_DMP_property}.
\end{theorem}

\begin{proof}
Consider any $u_h^{}\in V_h^{}$, $i\in\{1,\dots,M\}$, and $j\in S_i$. Let
$\bu$ be the vector of nodal values of $u_h$ and assume that $u_i$ is a local 
extremum of $u_h$ on $\omega_i$ and that $u_i\neq u_j$. We want to prove that
\begin{equation}\label{aa}
   a_{ij}+(1-\alpha_{ij}(\bu))\,d_{ij}\le0\,.
\end{equation}
If $a_{ij}\le0$, then \eqref{aa} holds since
$(1-\alpha_{ij}(\bu))\,d_{ij}\le0$.
If $a_{ij}>0$, then $a_{ji}\le0$ due to \eqref{eq:assumption_min} and hence
$a_{ji}<a_{ij}$ and $d_{ij}=-a_{ij}<0$. Thus, if $u_i\ge u_k$ for all 
$k\in S_i$, then $f_{ij}>0$ and $f_{ik}\ge0$ for $k\in S_i$, so that 
$\alpha_{ij}=R_i^+=0$. Similarly, if $u_i\le u_k$ for all $k\in S_i$, then 
$f_{ij}<0$ and $f_{ik}\le0$ for $k\in S_i$, so that $\alpha_{ij}=R_i^-=0$. 
Since $a_{ij}+d_{ij}=0$, one concludes that \eqref{aa} holds.
\end{proof}

If the matrix \eqref{eq:lumped_Galerkin_matrix} with lumped reaction term is
considered, then the validity of \eqref{eq:assumption_min} is guaranteed if the 
triangulation $\calT_h^{}$ satisfies the XZ-criterion \eqref{XZ-criterion}.
The condition \eqref{eq:assumption_min} may be satisfied also if the
XZ-criterion is violated,
particularly, in the convection-dominated case, since the convection matrix is
skew-symmetric. However, in general, the validity of a DMP cannot be guaranteed
without the XZ-criterion. Moreover, if the matrix
\eqref{eq:lumped_Galerkin_matrix} is replaced by \eqref{eq:Galerkin_matrix},
then the validity of \eqref{eq:assumption_min} may be lost since some
off-diagonal entries of the matrix $\bbMc$ are positive.

It was shown in \cite{Kno17} that the DMP generally does not hold if
condition \eqref{eq:assumption_min} is not satisfied. This is due to the 
condition $a_{ji}\le a_{ij}$ used in \eqref{eq:kuzmin_symm_alpha} to symmetrize
the factors $\widetilde\alpha_{ij}$. Therefore, in \cite{Kno21}, it was
proposed to use the above limiter in the formulation \eqref{eq:asm-1},
\eqref{eq:asm-2} without the symmetry condition \eqref{eq:kuzmin_symm_alpha}.
To obtain a well defined problem satisfying a continuity assumption on 
$\alpha_{ij}(\bu)(u_j-u_i)$, the definition of $P_i^\pm$ was replaced by
\begin{equation}\label{eq:kuzmin_p_mod}
   P_i^+=\sum_{\mbox{\parbox{8mm}{\scriptsize\centerline{$j\in S_i$}
   \centerline{$a_{ij}>0$}}}}\,a_{ij}\,(u_i-u_j)^+\,,\quad\,\,\,
   P_i^-=\sum_{\mbox{\parbox{8mm}{\scriptsize\centerline{$j\in S_i$}
   \centerline{$a_{ij}>0$}}}}\,a_{ij}\,(u_i-u_j)^-\,.
\end{equation}
Then the DMP is satisfied without any additional condition on the matrix 
$\bbAN$, which means that it holds for any triangulation $\calT_h^{}$ and also 
without the lumping of the matrix $\bbMc$ in the Galerkin FEM. Note,
however, that if the reaction term is dominant, some lumping may be performed 
by the algebraic flux correction scheme.

\begin{theorem}[DMP for the algebraically stabilized scheme with modified 
Kuz\-min limiter]
\label{thm:asm_dmp}Let us consider the algebraically stabilized scheme 
\eqref{eq:asm-1}, \eqref{eq:asm-2} with $\alpha_{ij}=\widetilde\alpha_{ij}$ for 
$i,j=1,\dots,N$, where $\widetilde\alpha_{ij}$ is defined by 
\eqref{eq:kuzmin_p_mod} and \eqref{eq:kuzmin_q}--\eqref{eq:alpha-definition}. 
Then the algebraic DMP property and the algebraic DMP property for non-strict 
extrema are satisfied.
\end{theorem}

\begin{proof}
The proof is similar as for Theorem~\ref{thm:kuzmin_dmp}. Under the
assumptions made before \eqref{aa} we now want to prove that
\begin{equation}\label{aa2}
   a_{ij}-\max\{0,(1-\widetilde\alpha_{ij}(\bu))\,a_{ij},
                  (1-\widetilde\alpha_{ji}(\bu))\,a_{ji}\}\le0\,.
\end{equation}
Since this clearly holds if $a_{ij}\le0$, it suffices to investigate the case
$a_{ij}>0$. If $u_i\ge u_k$ for all $k\in S_i$, then $P_i^+\ge
a_{ij}\,(u_i-u_j)^+>0$, $f_{ij}>0$ and $f_{ik}\ge0$ for $k\in S_i$, so that
$\widetilde\alpha_{ij}=R_i^+=0$. If $u_i\le u_k$ for all $k\in S_i$, then
$P_i^-\le a_{ij}\,(u_i-u_j)^-<0$, $f_{ij}<0$ and $f_{ik}\le0$ for $k\in S_i$, 
so that $\widetilde\alpha_{ij}=R_i^-=0$. This implies \eqref{aa2}.
\end{proof}

If condition \eqref{eq:assumption_min} holds, then \eqref{eq:kuzmin_p} and
\eqref{eq:kuzmin_p_mod} are equivalent, and $b_{ij}(\bu)$ defined using the
modified Kuzmin limiter from Theorem~\ref{thm:asm_dmp} satisfies
$b_{ij}(\bu)=(1-\alpha_{ij}(\bu))d_{ij}$ with the Kuzmin limiter $\alpha_{ij}$
from \eqref{eq:kuzmin_symm_alpha}. Thus, under condition
\eqref{eq:assumption_min}, both approaches described above are equivalent.
The modified Kuzmin limiter was further improved and reformulated 
in \cite{JK22} leading to the Monotone Upwind-type Algebraically Stabilized 
(MUAS) method. The paper \cite{JK22} also contains a detailed analysis of
algebraically stabilized methods of the type \eqref{eq:asm-1},
\eqref{eq:asm-2}.  Further analytical and numerical studies of these 
approaches recently inspired the design of the Symmetrized Monotone Upwind-type 
Algebraically Stabilized (SMUAS) method in \cite{Kno22x}.

Another way how to construct a limiter leading to the DMP on arbitrary meshes
and without an explicit lumping of the matrix $\bbMc$ was proposed in \cite{BJK17}, using some 
ideas of \cite{Kuz12}. The definition of this limiter, which we call BJK
limiter, is inspired by the Zalesak algorithm that will be derived in
Section~\ref{sec:tcd_femfct} for the time-dependent case. It again relies on 
local quantities $P_i^+$, $P_i^-$, $Q_i^+$, $Q_i^-$ which are now computed for 
$i=1,\dots,M$ by
\begin{align}
   P_i^+&=\sum_{j\in S_i}\,f_{ij}^+\,,\qquad
   P_i^-=\sum_{j\in S_i}\,f_{ij}^-\,,\label{eq:BJK_p}\\
   Q_i^+&=q_i\,(u_i-u_i^{\rm max})\,,\qquad
   Q_i^-=q_i\,(u_i-u_i^{\rm min})\,,\label{eq:BJK_q}
\end{align}
where again $f_{ij}=d_{ij}\,(u_j-u_i)$ and
\begin{equation}
   u_i^{\max}= \max_{j\in S_i\cup\{i\}}\,u_j\,,\qquad
   u_i^{\min}= \min_{j\in S_i\cup\{i\}}\,u_j\,,\qquad
   q_i=\gamma_i\,\sum_{j\in S_i}\,d_{ij}\,,\label{eq:min_max_q}
\end{equation}
with fixed constants $\gamma_i>0$. Then one defines the factors
$\widetilde\alpha_{ij}$ by
\eqref{eq:R_i_plus_minus}--\eqref{eq:alpha-definition}.
Finally, the limiters are defined by
\begin{equation}\label{eq:BJK_symm_alpha}
   \alpha_{ij}=\min\{\widetilde\alpha_{ij},\widetilde\alpha_{ji}\}\,,
   \qquad i,j=1,\dots,N\,.
\end{equation}

\begin{theorem}[DMP for the AFC scheme with BJK limiter]
\label{thm:BJK_dmp}
The AFC scheme \eqref{eq:afc-1}, \eqref{eq:afc-2} with the BJK limiter
defined by \eqref{eq:BJK_p}--\eqref{eq:min_max_q},
\eqref{eq:R_i_plus_minus}--\eqref{eq:alpha-definition}, and
\eqref{eq:BJK_symm_alpha} satisfies the algebraic DMP property and also the 
algebraic DMP property for non-strict extrema.
\end{theorem}

\begin{proof}
The proof is similar as for Theorem~\ref{thm:kuzmin_dmp}. Under the
assumptions made before \eqref{aa} we now want to prove that
\begin{equation}\label{aaa}
   a_{ij}+(1-\min\{\widetilde\alpha_{ij}(\bu),\widetilde\alpha_{ji}(\bu)\})\,d_{ij}\le0\,.
\end{equation}
If $d_{ij}=0$, then $a_{ij}\le0$ and hence \eqref{aaa} holds. Thus, let us
assume that $d_{ij}<0$. If $u_i\ge u_k$ for all $k\in S_i$, then $f_{ij}>0$ and 
$u_i^{\max}=u_i$ so that $P_i^+>0$, $Q_i^+=0$ and
$\widetilde\alpha_{ij}=R_i^+=0$. Since $a_{ij}+d_{ij}\le0$, one obtains
\eqref{aaa}. If $u_i\le u_k$ for all $k\in S_i$, \eqref{aaa} follows 
analogously.
\end{proof}

It was proved in \cite{BJK17} that, for 
\begin{equation*}
   \gamma_i\ge\frac{\displaystyle\max_{\bx_j\in\partial\omega_i}\,|\bx_i-\bx_j|}
                 {\mbox{\rm dist}(\bx_i,\partial\omega_i^{\rm conv})}\,,
\end{equation*}
where $\omega_i^{\rm conv}$ is the convex hull of $\omega_i$,
the AFC scheme with the BJK limiter is linearity preserving, i.e., $\bbB(u)=0$
for $u\in\bbP_1(\mathbb{R}^d)$. This property may lead to improved convergence
results, see, e.g., \cite{BBK17b,BJKR18}. Note that large values of the 
constants $\gamma_i$ cause that more limiters
$\alpha_{ij}$ will be equal to 1 and hence less artificial diffusion is added,
which makes it possible to obtain sharp approximations of layers. On the other
hand, however, large values of $\gamma_i$'s also cause that the numerical
solution of the nonlinear algebraic problem becomes more involved.

\begin{remark}
The various limiters discussed above are inspired by techniques used in the
time-dependent case, where a classical approach is the above-mentioned Zalesak 
algorithm (cf.~Section~\ref{sec:tcd_femfct}). This algorithm cannot be
simply applied to the steady-state case since the quantities $Q_i^\pm$ are
defined using the mass matrix from the discretization of the time-derivative,
and a provisional solution of an explicit low-order scheme. The Kuzmin limiter
formulated in \eqref{eq:kuzmin_p}--\eqref{eq:kuzmin_symm_alpha} circumvents
this problem by defining $Q_i^\pm$ analogously as $P_i^\pm$ in the Zalesak
algorithm. The design of the BJK limiter is formally closer to the Zalesak
limiter and relies on a carefully selected multiplicative factor in the
definition of $Q_i^\pm$. The remaining approaches mentioned above use various
modifications of the Kuzmin limiter. As discussed above, the original Kuzmin
limiter satisfies the DMP only under the condition \eqref{eq:assumption_min}
whereas the other approaches satisfy the DMP without any condition on the
stiffness matrix. In addition the BJK limiter and the SMUAS limiter
\cite{Kno22x} are linearity preserving on arbitrary simplicial meshes.
Nevertheless, it is difficult to assess the quality of the resulting schemes
from these theoretical properties. Indeed, recent numerical results
\cite{JJ19,JJK22x,JKP22,Kno22x} reveal that depending on considered data and the 
used criterion (e.g., accuracy, efficiency or experimental convergence rate),
one can come to various conclusions concerning the quality of the
methods. For example, the BJK limiter often leads to sharp approximations of
layers but the nonlinear algebraic problems are difficult to solve and the
approximate solutions may be less accurate away from layers than for the Kuzmin
limiter.
\hspace*{\fill}$\Box$\end{remark}

Finally, let us present another way how to define the matrix $\bbB(\bu)$ in the
algebraically stabilized problem \eqref{eq:asm-1}, \eqref{eq:asm-2}, the
so-called BBK method proposed in \cite{BBK17b}. It is also referred to as 
smoothness-based viscosity and has its origin in the finite volume literature 
(see, e.g., \cite{JST81} and \cite{Jameson}).

Given $\bu\in\mathbb R^N$, one first defines the function $\xi_{\bu}\in V_h$ 
whose nodal values are given by
\begin{equation}\label{xi-definition}
\xi_{\bu}(\bx_i) = 
\begin{cases} \displaystyle
 \frac{\left|\sum_{j\in S_i} (u_i-u_j)\right|}{\sum_{j\in S_i}|u_i-u_j|}&\quad\mbox{if}\;\;
\displaystyle\sum_{j\in S_i} |u_i-u_j| \neq 0\, ,\\
    \qquad\qquad 0 &\quad\mbox{otherwise}\,,
 \end{cases}
\qquad i=1,\dots,N\,.
\end{equation}
Then, for any $i,j\in\{1,\dots,N\}$ such that there is an edge $E\in\calE_h^{}$ 
with endpoints $\bx_i,\bx_j$, one sets
\begin{equation}\label{B-BBK}
   b_{ij}(\bu)=-\gamma_0\,h_E^{d-1}\,
   \max_{\bx\in E} \big[ \xi_{\bu}(\bx)\big]^p\,,\qquad 
   p\in \; [1, +\infty)\,,
\end{equation}
where $\gamma_0$ is a fixed parameter, dependent on the data of 
\eqref{steady-strong}.
For other pairs of $i\neq j$, one sets $b_{ij}(\bu)=0$. Finally, the diagonal
entries of the matrix $\bbB(\bu)$ are again defined by \eqref{eq:matrix_B2}.
Then the corresponding form $d_h$ again satisfies \eqref{eq:d_h_row_sums},
\eqref{eq:pos_semidef_dh}, and \eqref{eq:d_h_sparse2}.

The value of $p$ determines the rate of decay of the numerical
diffusion with the distance to the critical points. A value closer to
$1$ adds more diffusion far away from layers and extrema, while a larger value
makes the diffusion vanish faster, but on the other hand, increasing $p$ may 
make the nonlinear system more difficult to solve. In our experience, values up 
to $p=20$ are considered safe to use
(see \cite{BBK17b} for a detailed discussion). Note also that, on symmetric
meshes, the method is linearity preserving. 

\begin{theorem}[DMP for the BBK method]
\label{thm:BBK_dmp}
Let the triangulation $\calT_h^{}$ satisfy the XZ-criterion 
\eqref{XZ-criterion}. Then there exist constants $C_0$ and $C_1$ depending
only on the shape regularity of $\calT_h^{}$
such that if $\gamma_0\ge C_0\|\bb\|_{0,\infty,\Omega}^{}+C_1\,\sigma\,h$, 
then the algebraically stabilized scheme \eqref{eq:asm-1}, \eqref{eq:asm-2} 
with $\bbB(\bu)$ defined by \eqref{xi-definition}, \eqref{B-BBK} satisfies the 
algebraic DMP property and also the algebraic DMP property for non-strict 
extrema.
\end{theorem}

\begin{proof}
We again start with the assumptions made in the proof of 
Theorem~\ref{thm:kuzmin_dmp} before \eqref{aa}. Then $\xi_{\bu}(\bx_i)=1$ and 
hence $b_{ij}(\bu)=-\gamma_0\,h_E^{d-1}$. In view of \eqref{eq:def_dh},
Theorem~\ref{theo:Delaunay-mesh-laplacian}, and the shape regularity of the 
mesh, one obtains
\begin{align*}
   a(\phi_j^{},\phi_i^{})+d_h^{}(u_h^{};\phi_j^{},\phi_i^{})
   &= \varepsilon\,(\nabla\phi_j,\nabla\phi_i)+(\bb\cdot\nabla\phi_j,\phi_i)
   + \sigma\,(\phi_j,\phi_i)-\gamma_0\,h_E^d\\
   &\le(C_0\,\|\bb\|_{0,\infty,\Omega}^{}+C_1\,\sigma\,h
   -\gamma_0)\,h_E^{d-1}
\end{align*}
and the result follows.
\end{proof}

Let us now briefly discuss different approaches to make the BBK method
linearity preserving on general meshes. The common point to all those alternatives is to introduce positive constants
$\beta_{ij}^{}$ for $j\in S_i^{}$ and modify slightly the definition \eqref{xi-definition} of $\xi_{\bu}^{}(\bx_i^{})$ as follows
\begin{equation*}
\xi_{\bu}(\bx_i) = 
\begin{cases} \displaystyle
 \frac{\left|\sum_{j\in S_i} \beta_{ij}^{}(u_i-u_j)\right|}{\sum_{j\in S_i} \beta_{ij}^{}|u_i-u_j|}&\quad\mbox{if}\;\;
\displaystyle\sum_{j\in S_i} |u_i-u_j| \neq 0\, ,\\
    \qquad\qquad 0 &\quad\mbox{otherwise}\,,
 \end{cases}
\qquad i=1,\dots,N\,.
\end{equation*}
In \cite[Remark~1]{BBK17b} a process to generate a linearity preserving method is described. It involves solving
local minimization problems in each node to determine the value of $\beta_{ij}$. An alternative approach is presented
in \cite[Section~4.3]{GP17}. If the support of the basis functions
$\phi_i^{}$ is convex,  then there exists a set of generalized barycentric
coordinates $(\omega_{ij})_{j\in S_i}$ 
such that its elements are non-negative functions, form a partition of unity,  
and $\bx=\sum_{j\in S_i}\omega_{ij}(\bx)\bx_j^{}$ for all $\bx\in\omega_i$.
A process to build these coordinates in higher dimensions is proposed in \cite{WSHD07} (see also \cite{floater_2015} for a comprehensive review
on the topic of generalized barycentric coordinates). Then, taking $\beta_{ij}^{}=\omega_{ij}(\bx_i)$, it can be proven that the
resulting method is linearity preserving.

We end this section again by discussing the solvability and error
estimates. It can be proven by means of Brouwer's fixed-point theorem that the
nonlinear algebraic problem \eqref{eq:asm-1}, \eqref{eq:asm-2} is solvable 
provided that the entries of the matrix $\bbB(\bu)$ are bounded functions of 
$\bu\in{\mathbb R}^N$ and, for any $i,j\in\{1,\dots,N\}$, the functions 
$b_{ij}(\bu)(u_j-u_i)$ are continuous, see, e.g., \cite{JK22}. This is the case 
for all the methods discussed in this section, cf.~\cite{BJKR18,JK22,Kno22x}.
A natural norm for estimating the errors of the solutions to the nonlinear
problems considered in this section is the solution-dependent norm proposed in
\cite{BJK16} given by
\begin{displaymath}
   \| v\|_h^{}:=\Big(\varepsilon\,|v|_{1,\Omega}^2+\sigma\,\|v\|_{0,\Omega}^2+
   d_h(u_h;v,v)\Big)^{1/2}\,.
\end{displaymath}
Then, if $u\in H^2(\Omega)$ and $\sigma>0$, one has (cf., e.g., \cite{BJK16})
\begin{displaymath}
   \| u-u_h\|_h^{}\le C\,(\varepsilon+\sigma^{-1}\,\|\bb\|_{0,\infty,\Omega}^2
   +\sigma)\,h\,\|u\|_{2,\Omega}^{}+ (d_h(u_h;i_h^{}u,i_h^{}u))^{1/2}\,,
\end{displaymath}
where $C$ is independent of $h$ and the data of the problem
\eqref{steady-strong}. The term \linebreak $(d_h(u_h;i_h^{}u,i_h^{}u))^{1/2}$ 
represents
an estimate of the consistency error induced by the algebraic stabilizations.
As its precise definition varies according to the choice of limiters, it is to 
be expected that different convergence orders may be proven for the different 
choices of limiters. A common feature of the analyses presented in 
\cite{BJK16,BBK17b} is the following: an ${\cal O}(h^{1/2})$ 
convergence can be proven for meshes of XZ-type. For non-XZ meshes, this
convergence order can be proven only in the convection-dominated case in
general since certain entries of the diffusion matrix may be positive.
Indeed, examples of non-convergence in the diffusion-dominated case are shown
for the Kuzmin limiter in \cite{BJK16}. Moreover, it was proven in
\cite{BBK17b,BJKR18} that the combination of Lipschitz continuity 
and linearity preservation leads to an ($\varepsilon$-dependent) improved  
error estimate of order ${\cal O}(h)$.

\subsection{A monotone Local Projection Stabilized (LPS) method}
\label{Sec:cd_nonlinear_methods_edge}

In this section we will review a LPS method that respects the DMP proposed in \cite{BBK17a}. Its motivation, already hinted in \cite{BBH11},
is to start with an optimal order stabilized method based on facets (e.g.~CIP), and to introduce a nonlinear switch that makes the method become
a first order linear artificial diffusion method in the vicinity of layers and 
extrema. 

The monotone LPS method is given by \eqref{nonlinear-generic-2} with
\begin{align}\label{LPS-nonlinear-stab}
   d_h^{}(w_h^{};u_h^{},v_h^{})  &=\sum_{F\in\calF_I^{}}\Big[
   \tau_F^{}\alpha_F^{}(w_h^{})(\nabla u_h^{},\nabla v_h^{})_{\omega_F^{}}^{}\\
   &\hspace*{5mm}+\,\gamma_F^{}\big(1-\alpha_F^{}\,(w_h^{})\big)
   (\nabla u_h^{}- {G}_F^{}\nabla u_h^{},
    \nabla v_h^{}- {G}_F^{}\nabla v_h^{})_{\omega_F^{}}^{}\Big].\nonumber
\end{align}
Here, for each $F\in\calF_I^{}$, the operator ${G}_F^{}$ provides a local
mean value defined by
\begin{equation*}
   {G}_F^{}q=\frac{(q,1)_{\omega_F}^{}}{|\omega_F^{}|}\,,\qquad 
   q\in L^1(\omega_F^{})\,,
\end{equation*}
which is computed component-wise in the case of vector-valued functions, and
$\tau_F^{}$, $\gamma_F^{}$ are stabilization parameters given by
\begin{equation}\label{LPS-param}
   \tau_F^{}= c_0^{}\|\bb\|_{0,\infty,\omega_F^{}}^{}h_F^{}
   \qquad\textrm{and}\qquad 
   \gamma_F^{} =\gamma_0^{}\min\left\{\|\bb\|_{0,\infty,\omega_F^{}}^{}h_F^{},
   \frac{h_F^2}{\varepsilon}\right\}\,,
\end{equation}
with positive constants $c_0^{}$ and $\gamma_0^{}$. 
The nonlinear switches $\alpha_F^{}$ need to be designed in such a way that 
they detect regions of extrema and large variations in the gradients, the 
latter indicating the possible presence of layers. For now, we will just 
assume that they satisfy the following two basic assumptions: 
\begin{enumerate}[leftmargin=2em,label=\roman*),ref=\roman*)]
\item \label{ass:LPS_i} $\alpha_F^{}: V_h^{}\to [0,1]$ are continuous functions; and 
\item $\alpha_F^{}(u_h^{})=1$ whenever $u_h^{}$ attains a local extremum at a 
node of a mesh cell containing $F$. \label{ass:LPS_ii}
\end{enumerate}
In \cite{BBK17a} it was proposed to define $\alpha_F^{}$ using regularized
versions of the Kuzmin limiter \eqref{eq:kuzmin_symm_alpha} or the 
smoothness-based indicator \eqref{xi-definition}.

The form $d_h^{}(\cdot;\cdot,\cdot)$ obviously
satisfies the assumptions \eqref{eq:d_h_row_sums} and
\eqref{eq:pos_semidef_dh}. In addition, since $(q-{G}_F^{}q,1)_{\omega_F}^{}=0$ 
for any $q\in L^1(\omega_F^{})$ and $F\in\calF_I^{}$, it can be also written as
\begin{align}\label{LPS-nonlinear-stab2}
   d_h^{}(w_h^{};u_h^{},v_h^{})  &=\sum_{F\in\calF_I^{}}\Big[
   \tau_F^{}\alpha_F^{}(w_h^{})(\nabla u_h^{},\nabla v_h^{})_{\omega_F^{}}^{}\\
   &\hspace*{5mm}+\,\gamma_F^{}\big(1-\alpha_F^{}\,(w_h^{})\big)
   (\nabla u_h^{}- {G}_F^{}\nabla u_h^{},\nabla v_h^{})_{\omega_F^{}}^{}\Big].
   \nonumber
\end{align}

\begin{remark}  A more natural way of writing \eqref{LPS-nonlinear-stab2}
would be to express the stabilizing term as follows
\begin{equation*}
   \sum_{F\in\calF_I^{}}\tilde{\tau}_F^{}(\nabla u_h^{}
   -\beta_F^{}(u_h^{}){G}_F^{}\nabla u_h^{}, 
   \nabla v_h^{})_{\omega_F^{}}^{}\,,
\end{equation*}
where $\tilde{\tau}_F^{}$ is a stabilization parameter, and 
$\beta_F^{}(u_h^{})=1-\alpha_F^{}(u_h^{})$. This writing does represent the 
idea of a method that includes transitions between low-order artificial 
diffusion and higher order local projection,  while at the same time
stressing the character of 
combining linear and nonlinear stabilization terms,
as it was
made in Section~\ref{Sec:cd_nonlinear_methods_BE} for the method given by \eqref{BE-1-complete}. Unfortunately,  numerical 
experimentation has shown that to obtain accurate results the stabilization 
parameters for the linear diffusion and local projection parts need to be of 
significantly different sizes. This has led to the (less natural)
writing \eqref{LPS-nonlinear-stab} for the stabilization term.

It is also worth mentioning that a similar strategy to the above monotone LPS
method, although using a local projection related to the Scott--Zhang 
interpolation operator, was used in \cite{BH14} to approximate the transport 
problem.
\hspace*{\fill}$\Box$\end{remark}

In \cite{BBK17a} it was proven that, under the assumptions i) and ii) on the 
limiters, the discrete problem has at least one solution.
Concerning the satisfaction of the DMP, we now report a proof slightly more 
specific than the one provided in \cite[\S~2.3]{BBK17a}. To avoid technical 
complications, we will present this result in two space dimensions and will 
suppose that $\sigma=0$.

\begin{theorem}[DMP for the monotone LPS method]
Let us suppose that $d=2$, the mesh family $\{\calT_h^{}\}_{h>0}^{}$ is weakly 
acute and average acute, $\sigma=0$ and the nonlinear switches $\alpha_F^{}$ 
satisfy \ref{ass:LPS_ii}. Then, there exists a constant $C>0$ 
depending
only on the shape regularity of the mesh family $\{\calT_h^{}\}_{h>0}^{}$
such that, if $c_0^{}$ from \eqref{LPS-param} satisfies
\begin{equation}\label{c0-large-LPS}
   c_0^{}\ge C\,\cot\frac{\delta}{2}\,,
\end{equation}
where $\delta$ is the angle appearing in \eqref{FEM-prelim-1}, then the form 
$d_h^{}(\cdot;\cdot,\cdot)$ defined in \eqref{LPS-nonlinear-stab} satisfies 
the algebraic DMP property and also the algebraic DMP property for non-strict 
extrema.
\end{theorem}

\begin{proof}
Consider any $u_h^{}\in V_h^{}$ and let us suppose that $u_h^{}$ attains a 
local extremum at an interior node $\bx_i^{}\in\Omega$. Consider any
$j\in\{1,\dots,N\}$. Since $\alpha_F^{}(u_h^{})=1$ for any $F\subset\omega_i$
and $\nabla\phi_i^{}|_{\omega_F^{}}^{}=0$ for any $F\not\subset\omega_i$, it
follows from \eqref{LPS-nonlinear-stab2} that
\begin{equation*}
   d_h^{}(u_h^{};\phi_j^{},\phi_i^{})  = 
   \sum_{F\in\calF_I^{},\,F\subset\omega_i}
   \tau_F^{}\,(\nabla \phi_j^{},\nabla \phi_i^{})_{\omega_F^{}}^{}\,,
\end{equation*}
which implies \eqref{eq:d_h_sparse}. Now consider any $j\in S_i$ and let us 
denote by $E=K\cap K'$ the edge connecting $\bx_i^{}$ and $\bx_j^{}$. Since the 
mesh is weakly acute, one has $(\nabla\phi_j^{},\nabla\phi_i^{})_K^{}\le 0$ for 
all $K\in\calT_h^{}$, which leads to $d_h^{}(u_h^{};\phi_j^{},\phi_i^{})\le
\tau_E^{}\,(\nabla\phi_j^{},\nabla\phi_i^{})_{\omega_E^{}}^{}
=\tau_E^{}\,\ell_{ij}$. Thus, applying \eqref{Gal7}, one arrives at
\begin{align*}
   a(\phi_j^{},\phi_i^{})+d_h^{}(u_h^{};\phi_j^{},\phi_i^{})&\le
   \tau_E^{}\,\ell_{ij}^{}+c_{ij}^{}
   =c_0^{}\,h_E^{}\,\|\bb\|_{0,\infty,\omega_E^{}}^{}\,\ell_{ij}^{}+c_{ij}^{}\\
   &\le-\frac{c_0^{}\,h_E^{}\,\|\bb\|_{0,\infty,\omega_E^{}}^{}}{2}\,
    \tan\frac{\delta}{2}+
   \frac{(h_K^{}+h_{K'}^{})\|\bb\|_{0,\infty,{\omega_E^{}}}^{}}{6}\,. 
\end{align*}
Thanks to the mesh regularity, one has 
$h_K^{}+h_{K'}^{}\le\tilde{C}\,h_E^{}$, where $\tilde{C}$ does not depend on 
the mesh size $h$. Hence, if \eqref{c0-large-LPS} holds with $C=\tilde{C}/3$,
one obtains \eqref{Gen-Non-11} and \eqref{Gen-Non-12}.
\end{proof}

We finish this section by summarizing the error estimates available for the method discussed in this section. Following standard estimates involving stability and asymptotic consistency,
an  $\mathcal O(h^{1/2})$ error estimate can be proven. In \cite[Section~2.4]{BBK17a}
a more refined analysis is carried out assuming that the functions $\alpha_F^{}$ decay with an appropriate rate away
from the layers, in other words, assuming that the nonlinear switch is active in only a small region of
the computational domain. More precisely, one starts  defining the region
\[ S_\alpha^{}:= \bigcup \left\{ K\in\calT_h^{}:\max_{F\in \calF_K^{}}\alpha_F^{}(u_h^{})> h^2\right\}\,, \]
and assumes  that $|S_\alpha^{}|=C\, h^s$ with $s>0$. In addition, for $r>0$ one defines the set
\[ S_{h,{\rm ext}}^{}:= \left\{\bx\in \Omega : |\nabla u(\bx)|\le C\,
h^r\,|u|_{2,\infty,\Omega}^{}\right\}\,,\]
and requires that 
\[ \sup_{\bx\in S_\alpha^{}}\inf_{\by\in S_{h,\textrm{ext}}}|\bx-\by|\le C\, h^r\,.
\]
Under these assumptions the following error estimate is proven in \cite[Lemma~2.6]{BBK17a} 
\begin{eqnarray*}
\lefteqn{\varepsilon^{\frac{1}{2}}|u-u_h^{}|_{1,\Omega}^{}+\sigma^{\frac{1}{2}}\|u-u_h^{}\|_{0,\Omega}^{}+d_h^{}(u_h^{};u-u_h^{},
u-u_h^{})^{\frac{1}{2}}}\nonumber\\
&\le& C\,\left(\varepsilon + \|\bb\|_{\infty,\Omega}^{}h + (\sigma+\sigma^{-1}\,|\bb|_{1,\infty,\Omega}^2)h^2\right)^{\frac{1}{2}}h\,|u|_{2,\Omega}^{} 
+ C\, h^{\frac{1+s}{2}}\big(h+h^r\big)\,|u|_{2,\infty,\Omega}^{}\,.
\end{eqnarray*}
Supposing in addition that $r+s/2\ge 1$ the improved estimate
\begin{eqnarray*}
\lefteqn{\varepsilon^{\frac{1}{2}}|u-u_h^{}|_{1,\Omega}^{}+\sigma^{\frac{1}{2}}\|u-u_h^{}\|_{0,\Omega}^{}+d_h^{}(u_h^{};u-u_h^{}, u-u_h^{})^{\frac{1}{2}}}\nonumber\\
&\le&\, C\,\left(\varepsilon + \|\bb\|_{\infty,\Omega}^{}h + (\sigma+\sigma^{-1}\,|\bb|_{1,\infty,\Omega}^2)h^2\right)^{\frac{1}{2}}h\,|u|_{2,\infty,\Omega}^{}\,
\end{eqnarray*}
is obtained.

\section{A numerical illustration}\label{sec:num_exam}

This section presents a brief numerical study that illustrates the behavior of several 
methods discussed in the previous chapters. 

\begin{figure}[t!]
\begin{center}
\centerline{\includegraphics[width=0.5\textwidth]{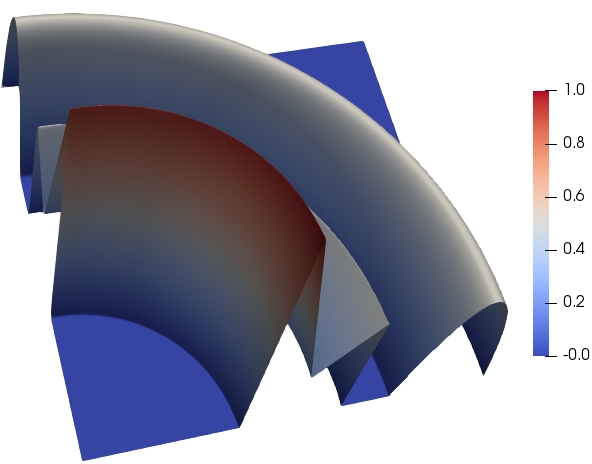}\hspace*{0.5em}
\includegraphics[width=0.45\textwidth]
{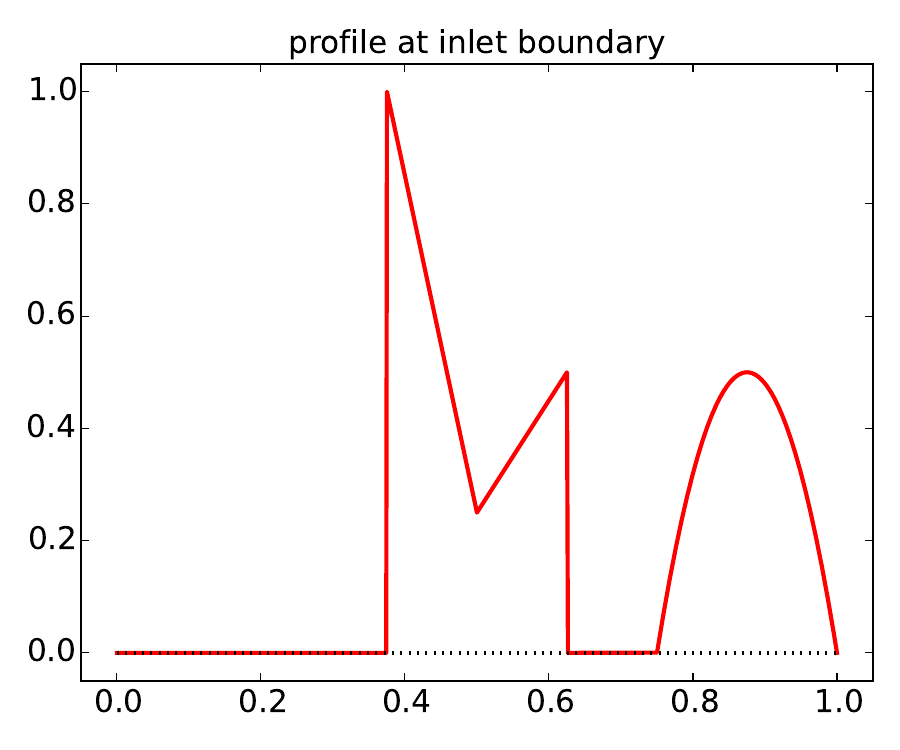}}
\caption{Numerical approximation of the solution (left) and profile at the inlet boundary (right).}\label{fig:three_int_layers_skewedM_bump_sol}
\end{center}
\end{figure}

In the considered example, a profile defined on the inlet boundary is transported through the 
domain $\Omega=(0,1)^2$.  The data of \eqref{steady-strong} are given by 
$\varepsilon = 10^{-5}$,  $\bb = (-y,x)^T$, and $\sigma=f=0$. Hence, the problem satisfies 
the conditions for the weak maximum principle from Theorem~\ref{thm:cd_weak_mp} for 
$\sigma=0$. The Dirichlet boundary condition at the inlet boundary $y=0$ 
is prescribed by 
\[
u(x,0) = \begin{cases}
\displaystyle \frac{x-0.375}{\xi} + 1 & \mbox{if } x \in [0.375-\xi, 0.375),\\[1em]
\displaystyle -0.75\frac{x-0.5}{0.125} + 0.25 & \mbox{if } x \in [0.375,0.5), \\[1em]
\displaystyle  0.25\frac{x-0.625}{0.125} + 0.5 & \mbox{if } x \in [0.5, 0.625),\\[1em]
\displaystyle -0.5\frac{x-0.625}{\xi} + 0.5& \mbox{if } x \in [0.625,0.625+\xi),\\[1em]
\displaystyle 32(x-0.75)(1-x)& \mbox{if } x \in [0.75,1],\\[1em]
0 & \mbox{else,}
\end{cases}
\]
with $\xi = 10^{-3}$.
A homogeneous Dirichlet boundary condition is prescribed at the boundary $x=1$ and homogeneous
Neumann conditions on the remaining part of the boundary. Figure~\ref{fig:three_int_layers_skewedM_bump_sol}
presents a numerical approximation of the solution and an illustration of the inlet
condition. 

\begin{figure}[t!]
\begin{center}
\begin{minipage}[c]{0.42\textwidth}
\includegraphics[width=\textwidth]
{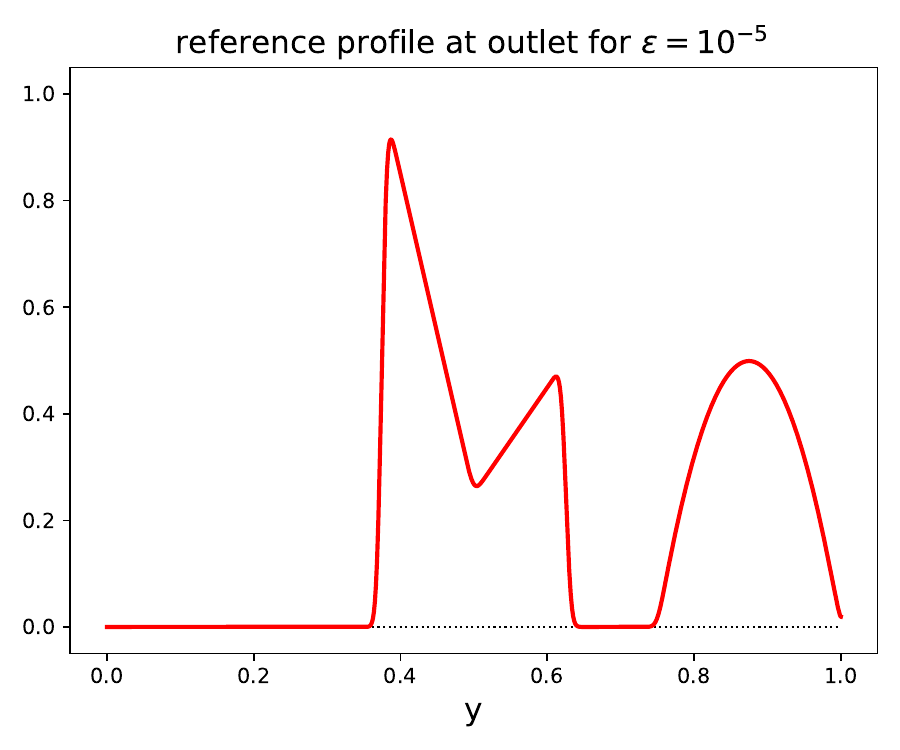}
\end{minipage}\hspace*{1ex}
\begin{minipage}[c]{0.53\textwidth}
\begin{tabular}{l|r}
quantity of interest & reference value\\
\hline
first maximum value&\hspace*{0.5em} 9.148468e-01\\
minimum value &2.642484e-01\\
second maximum value &4.699239e-01\\
width of the left profile & 2.628492e-01\\
maximum of the bump &4.989947e-01\\
width of the bump &2.367020e-01 \\
$u(0,1)$ & 1.914778e-02
\end{tabular}
\end{minipage}
\caption{Reference solution at the outlet boundary $x=0$ and corresponding reference values.}\label{fig:three_int_layers_skewedM_bump_out}
\end{center}
\end{figure}

For assessing the different methods, certain characteristic values of the solution at the 
outlet boundary $x=0$ are monitored. A reference solution was computed with the 
$\mathbb Q_2$ Galerkin FEM on a grid consisting of $4096 \times 4096$ squares
($67\,125\,249$ degrees of freedom, including Dirichlet nodes). 
Figure~\ref{fig:three_int_layers_skewedM_bump_out} depicts the reference solution at the outlet 
boundary.
For defining the reference 
values, the outlet boundary was decomposed into $100\,000$ intervals and the corresponding 
nodal values were used for computing the maximal and minimal values. The width of the 
left profile was defined by the condition $u(0,y)\ge0.1$ for $y \le 0.7$. For the width of the bump, 
also the condition $u(0,y)\ge0.1$ was used for computing the left point. Then, the width 
is defined by subtracting the $y$-coordinate of this point from $1$.
In all simulations, a linear interpolation was used for computing the widths. For 
the reference values, the above mentioned decomposition of the outlet boundary was used
and for the other simulations, an interpolation of the nodal values was applied. The 
reference values are provided in Figure~\ref{fig:three_int_layers_skewedM_bump_out}.

Simulations were performed for $\mathbb P_1$ finite elements. Initially, the domain was decomposed
into two triangles by using the diagonal from $(0,1)$ to $(1,0)$. Then, this decomposition was 
refined uniformly using red refinements. Linear systems of equations were solved with the sparse direct solver 
{\sc UMFPACK} \cite{Dav04} and nonlinear problems were solved with a simple fixed point iteration, 
e.g., see \cite{JK08} or 
the method \emph{fixed point rhs} from \cite{JJ19}, which has been proven to be the 
most efficient solver for AFC methods in the numerical studies of those papers. The iterations were stopped if the 
Euclidean norm of the residual vector was smaller than $10^{-10}$. Most of the computational results have been 
double checked with two codes, one of them {\sc ParMooN}, cf. \cite{GJM+16,WB_J17}.

From our numerical studies, only results will be presented where the numerical solution does not 
exhibit spurious oscillations, or 
more precisely, where the spurious oscillations are at most of the order of round-off
errors from floating point arithmetics or the stopping criterion for the iteration 
of a nonlinear discrete problem. There are many methods that compute solutions with 
small but still notable spurious oscillations, like some of the spurious
oscillations at layers diminishing (SOLD)
methods that can be found in the survey \cite{JK07}. However, such methods are not the topic of this review. 

The goal of computing oscillation-free numerical solutions could not be achieved for all methods 
presented in Section~\ref{Sec:cd_nonlinear_methods}. The proof of the DMP property for the edge stabilization
method of Burman and Ern from \cite{BE05} requires that the parameter $c_\rho$ from \eqref{BE-2} is 
sufficiently large, compare Theorem~\ref{thm:BE05_DMP}. In the numerical studies in \cite{BE05}, this 
parameter was set probably to $c_\rho = 5$ (this information is provided for an example with smooth solution
but not for an example with layers). But even with this parameter, notable spurious oscillations of the 
method are reported in \cite[Table~3]{BE05} for the case of a comparatively large diffusion coefficient. 
For the example studied here, we were able to solve the nonlinear problems (with two different codes)
for method \eqref{BE-1}--\eqref{BE-3} for parameters $c_\rho \precsim  0.005$. 
If a standard SUPG term is included, a numerical solution of the nonlinear problem 
was possible for $c_\rho \precsim 0.05$, which is the parameter choice for this 
method from \cite{JK07}. But in both cases and on all grids there are notable
undershoots of the computed solutions. This is the reason why we have not
reported the results from that method in this survey.

The precise definition of the constants $C^K_i$ used in the implementation of 
the Mizukami--Hughes method can be found in \cite[Fig.~8]{Kno06} or
\cite[Fig.~5]{Knobloch10}.
The algebraically stabilized method with BBK limiter was used with the
parameters $\gamma_0=0.75$ and $p=10$.

\begin{figure}[t!]
\begin{center}
\centerline{
\includegraphics[width=0.48\textwidth]{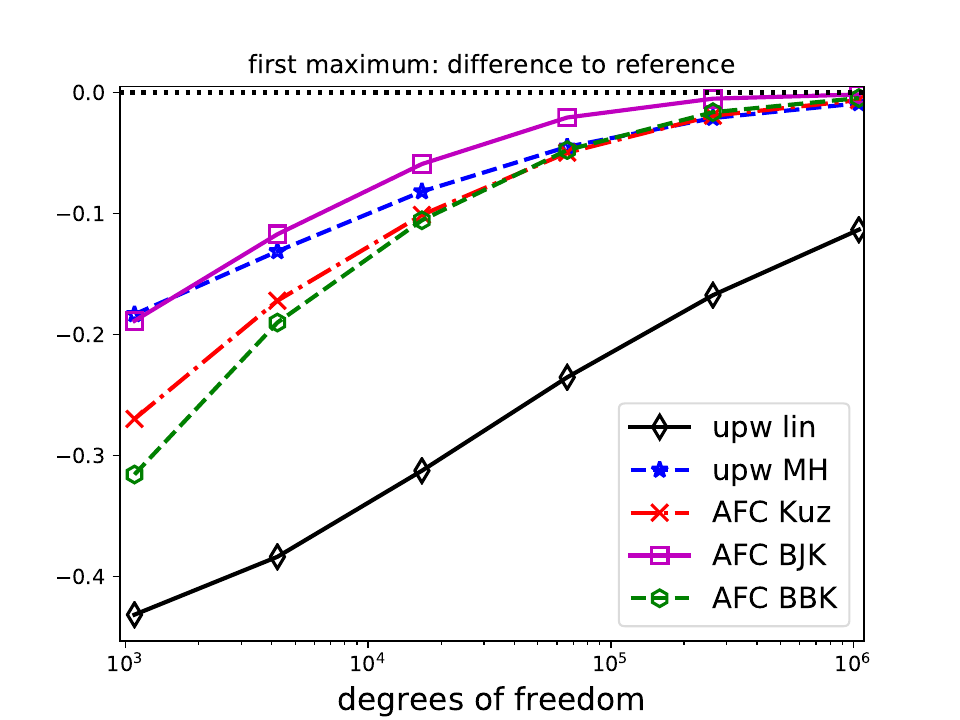}
\hspace*{1em}
\includegraphics[width=0.48\textwidth]{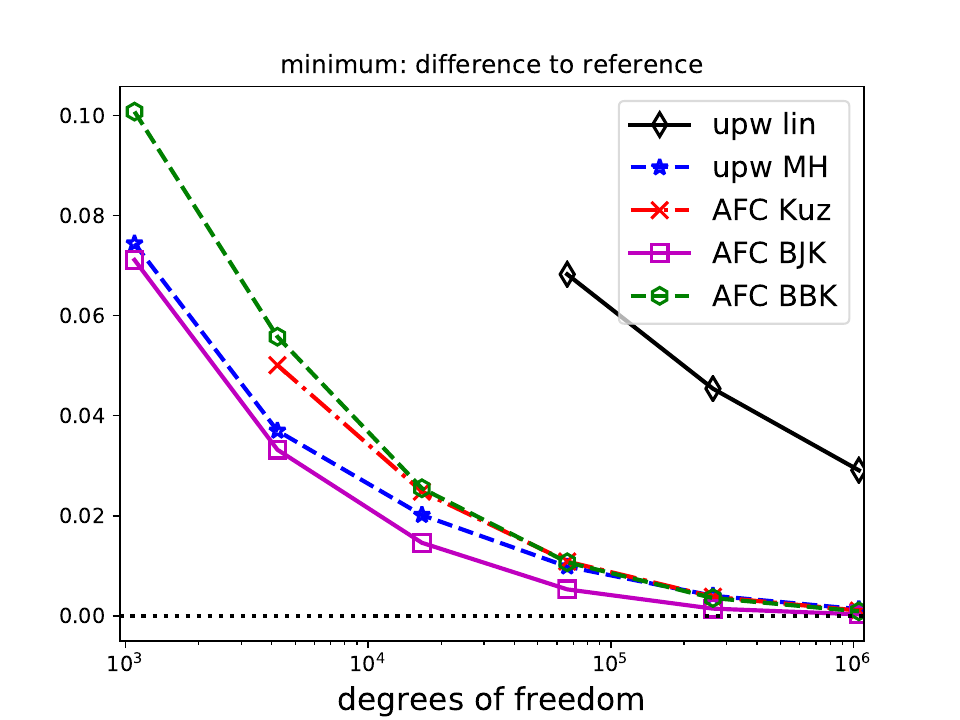}}
\centerline{
\includegraphics[width=0.48\textwidth]{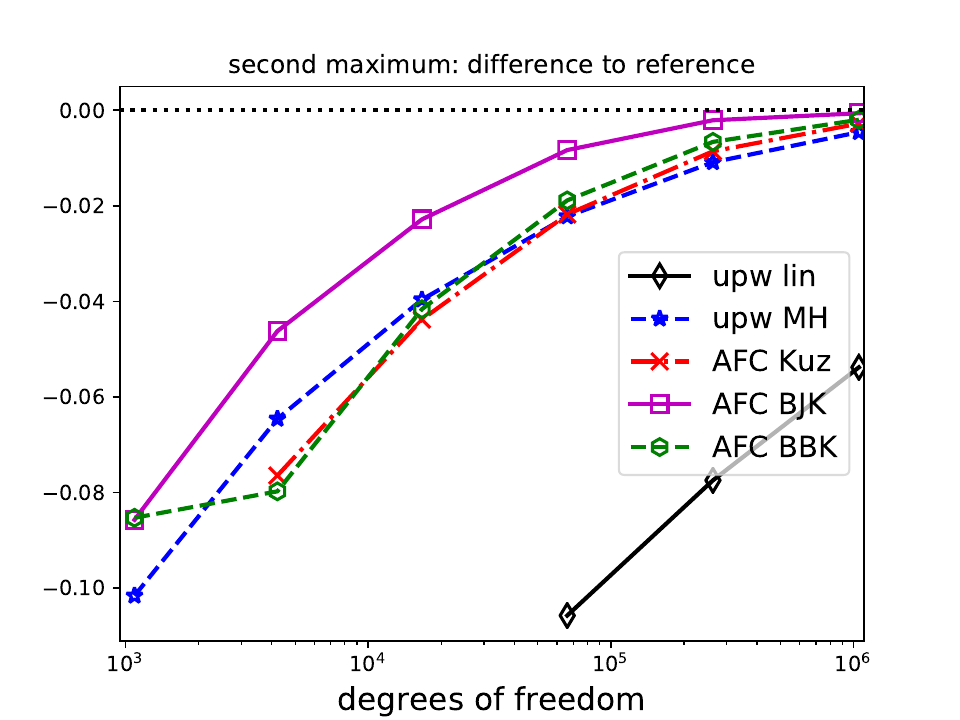}
\hspace*{1em}
\includegraphics[width=0.48\textwidth]{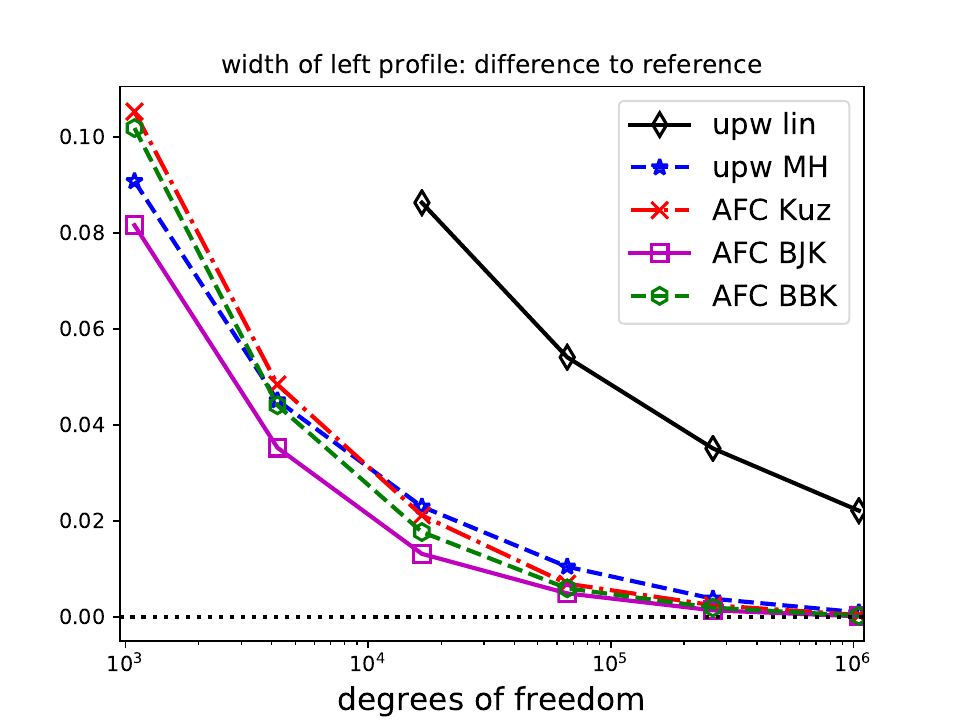}}
\centerline{
\includegraphics[width=0.48\textwidth]{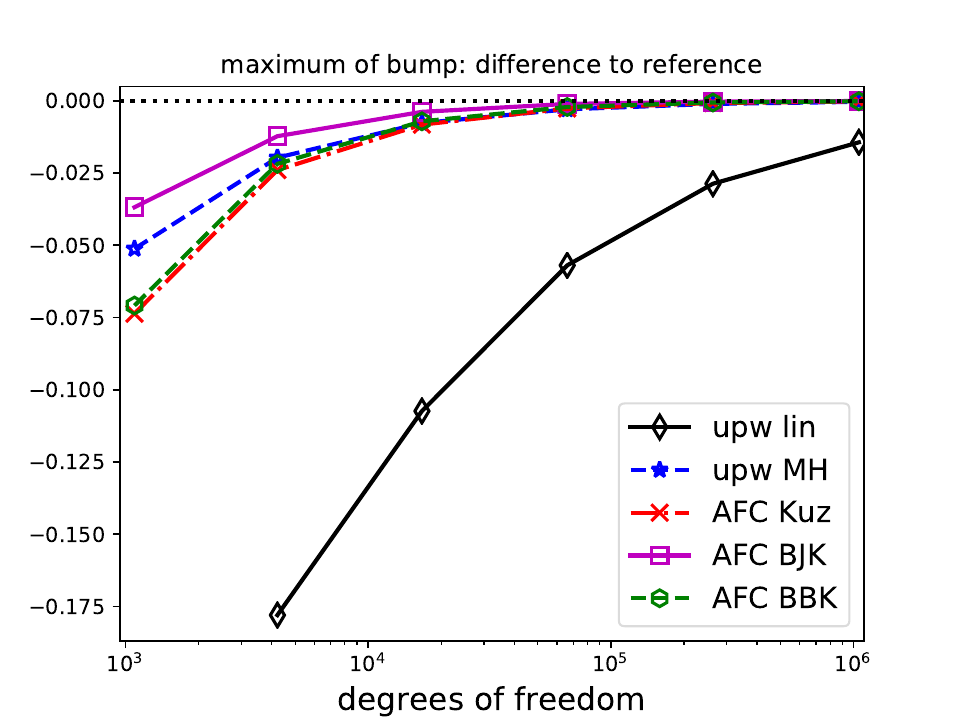}
\hspace*{1em}
\includegraphics[width=0.48\textwidth]{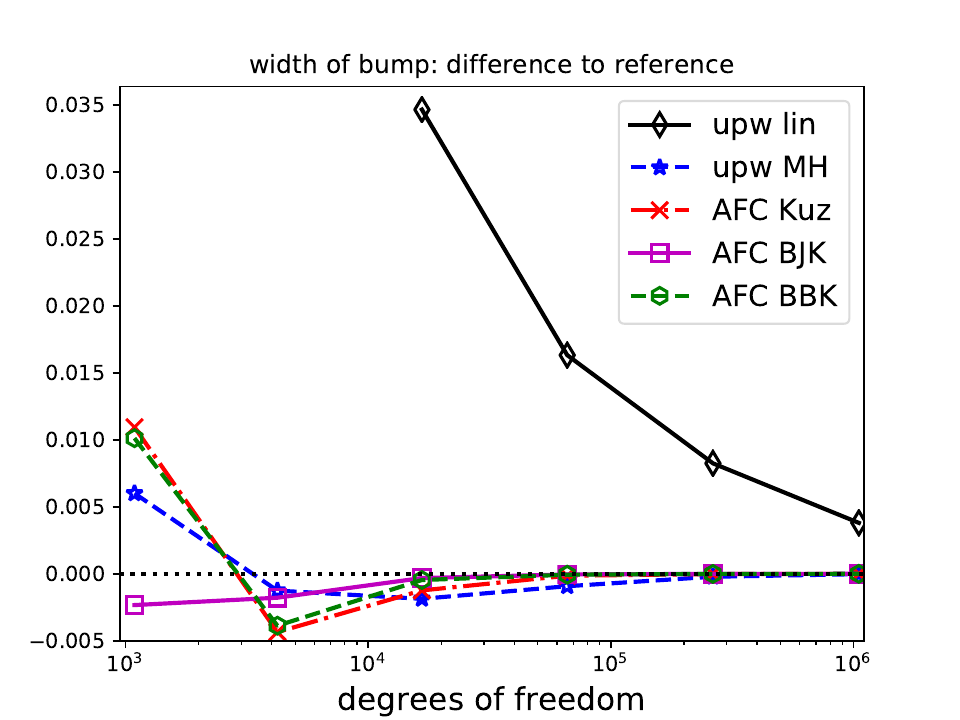}}
\centerline{
\includegraphics[width=0.48\textwidth]{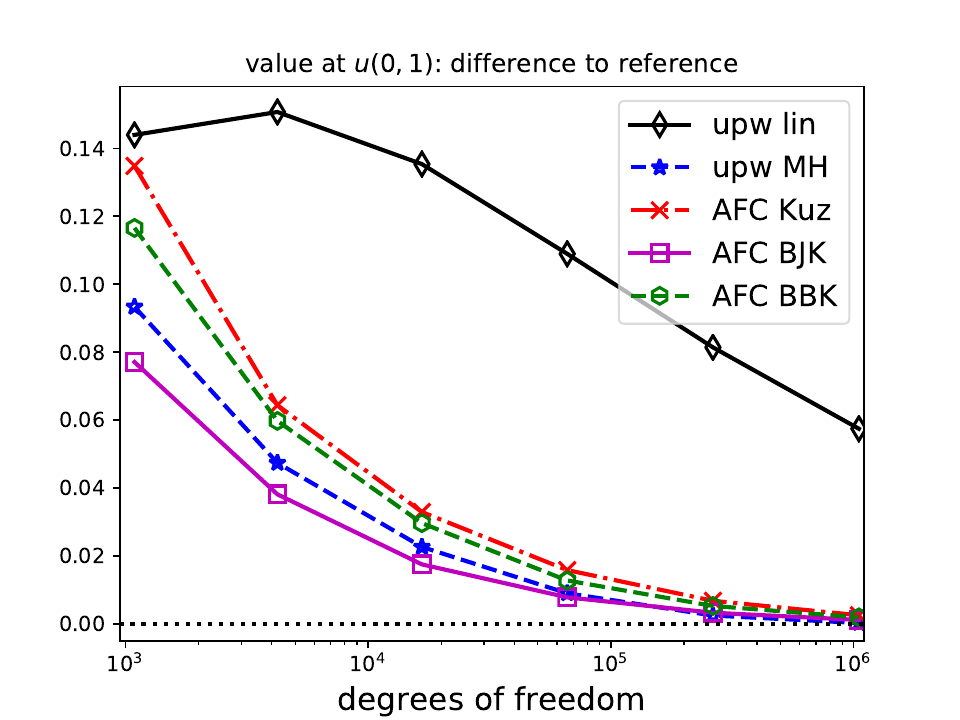}}
\caption{Differences of reference value and computed values for the quantities of interest.
}\label{fig:num_results_1em5}
\end{center}
\end{figure}

Figure~\ref{fig:num_results_1em5} presents the differences of the reference value and the values 
computed with the different methods for all quantities of interest. It can be seen that all nonlinear 
methods are much more accurate than the used linear method. The accuracy that is reached for the linear
upwind method with about $1\,000\,000$ degrees of freedom is usually achieved with the nonlinear methods 
already for about $4\,000$ or $16\,000$ degrees of freedom. One can also observe that there are some 
differences in the accuracy of the results computed with the different nonlinear discretizations, in 
particular on coarser grids. However, a comprehensive comparison of the different nonlinear methods, e.g., at other examples or with 
respect to the computational costs for solving the nonlinear problem, is outside the scope of this 
review. Some numerical comparisons of algebraically stabilized schemes can be found already in 
\cite{BJKR18,JJ19}.

In summary, the main messages that should be conveyed with this numerical study are that many nonlinear
discretizations which satisfy the DMP are much more accurate than linear discretizations with this 
property and that linear discretizations require prohibitively fine grids for computing accurate results
if the solution possesses layers.
This message is also supported by the recent paper \cite{JKP22} that 
contains results of comprehensive numerical studies not only for the methods
considered in this section but also for the edge-averaged method from
Section~\ref{sec:xu-zikatanov_fem}, the MUAS method~\cite{JK21} (see also
Section~\ref{Sec:cd_nonlinear_methods_afc}), and the monolithic convex limiting 
approach~\cite{Kuz20}.

\section{Time-dependent problem}\label{sec:tcd}

This section considers discretizations of time-dependent convection-diffusion-reaction 
equations, which use one-step $\theta$-schemes in time and finite element methods in space, and 
which satisfy a DMP. A few linear discretizations in space will be presented briefly and 
the class of FEM Flux-Corrected-Transport (FCT) schemes, which are usually nonlinear in space, 
will be discussed in detail.

\subsection{The continuous problem}\label{sec:tcd_cont_prob}

A time-dependent or evolutionary convection-diffusion-reaction initial-boundary value problem is given by 
\begin{equation}\label{eq:time_cdr} 
\begin{array}{rcll}
\partial_t u-\varepsilon\Delta u + \bb\cdot\nabla u + \sigma u &=& f &\textrm{in}\;
(0,T]\times\Omega\,,\\
u&=& g & \textrm{on}\;(0,T]\times \partial\Omega\,,\\
u(0,\cdot) &=& u_0 & \textrm{in}\;\Omega\,,
\end{array}
\end{equation}
where for the data of the problem, the same notations are used as in the
steady-state case. For simplicity, we will again suppose that
$\varepsilon>0$ and $\sigma\ge0$ are constants and that $\bb$ is solenoidal.
In \eqref{eq:time_cdr}, $T$ is the final time and $u_0=u_0(\bx)$ is a given 
initial condition. The velocity field $\bb$, the right-hand side $f$, and 
the boundary condition $g$ might depend on time and space. For brevity, the 
notation $\Omega_T = (0,T]\times \Omega$ is introduced and the parabolic 
boundary is denoted by $\Gamma_T = \overline{\Omega}_T\setminus{\Omega_T}$. 
Note that if $\sigma<0$, then a change of variable 
$\check u(t,\bx) = u(t,\bx) \exp(-\kappa t)$ leads to an evolutionary 
convection-diffusion-reaction equation for $\check u$ with the same terms for 
diffusion and convection, but the coefficient of the reactive term becomes 
$\sigma+\kappa$, such that $\sigma+\kappa\ge0$ holds for sufficiently 
large $\kappa$. In this way, many results obtained for $\sigma\ge0$ can be
extended to $\sigma<0$.

Consider for the moment a problem with $g=0$ on $(0,T]\times \partial\Omega$. 
Then, the definition and the analysis of a weak solution of \eqref{eq:time_cdr}
can be found, e.g., in \cite[Chapter~7.1]{Eva10}. For 
$\bb\in L^\infty(0,T;L^\infty(\Omega))$,
$f\in L^2(\Omega_T)$, and $u_0 \in L^2(\Omega)$, a function 
$u\in L^2(0,T;H_0^1(\Omega))$ with $\partial_t u \in L^2(0,T;H^{-1}(\Omega))$
is a weak solution of the convection-diffusion-reaction initial-boundary value problem if $u(0)=u_0$ and 
\[
\langle \partial_t u, v\rangle + \varepsilon(\nabla u, \nabla v)
+ (\bb\cdot\nabla u + \sigma u, v) = (f,v) \quad \forall\ v \in H_0^1(\Omega)
\]
almost everywhere in $[0,T]$, where $\langle\cdot,\cdot\rangle$ denotes the
duality pairing between $H^{-1}(\Omega)$ and $H_0^1(\Omega)$. The existence of 
a weak solution of \eqref{eq:time_cdr} can be proven with the Galerkin method, 
see also \cite{Eva10}. For proving uniqueness, it suffices to show that the 
fully homogeneous problem ($f=0$, $g=0$, $u_0=0$) possesses only the trivial 
solution, because the problem is linear. This statement can be proven using
the Gronwall lemma. Note that the condition $\sigma\ge0$ is not 
needed for these results. If $g$ does not vanish and it 
is sufficiently smooth, which will be assumed from now on, 
a problem with homogeneous boundary conditions can be constructed 
in the usual way by using a lifting of $g$ into $\Omega$ for each time and 
considering a problem for the difference of $u$ and the lifting. 

If $\sigma=0$, problem \eqref{eq:time_cdr} can be equivalently written in the
form
\begin{equation}\label{eq:time_cdr_cons} 
\begin{array}{rcll}
\partial_t u+ \nabla\cdot\left(-\varepsilon\nabla u + \bb u\right) &=& f &\textrm{in}\;(0,T]\times\Omega\,,\\
u&=& g & \textrm{on}\;(0,T]\times \partial\Omega\,,\\
u(0,\cdot) &=& u_0 & \textrm{in}\;\Omega\,,
\end{array}
\end{equation}
which is called conservative form and results from modeling the conservation of 
physical quantities. In \eqref{eq:time_cdr_cons}, $-\varepsilon\nabla u$ is 
called diffusive flux and $\bb u$ convective flux.

\subsection{Maximum principle, DMP, and positivity preservation}
\label{sec:tcd_framework}

It will be assumed in this section that $\bb\in C(\overline{\Omega}_T)$, 
such that this function is in particular bounded. From the practical point 
of view, the following weak maximum principle is of importance. Its proof can
be found in \cite[Chapter~7.1.4]{Eva10}, where also a strong maximum principle
is proven.

\begin{theorem}[Weak maximum principle]\label{thm:tcd_weak_mp}
Let $u \in C^2(\Omega_T) \cap C(\overline{\Omega}_T)$. Then
\begin{align}
&\partial_t u-\varepsilon\Delta u + \bb\cdot\nabla u + \sigma u \le 0 \quad\mbox{in } \Omega_T
\quad \Longrightarrow \quad
\max_{(t,\bx) \in \overline{\Omega}_T} u(t,\bx) \le \max_{(t,\bx) \in \Gamma_T}
u^+(t,\bx).\label{eq:tcd_mp_2}\\
&\partial_t u-\varepsilon\Delta u + \bb\cdot\nabla u + \sigma u \ge 0 \quad\mbox{in } \Omega_T
\quad \Longrightarrow \quad
\min_{(t,\bx) \in \overline{\Omega}_T} u(t,\bx) \ge \min_{(t,\bx) \in
\Gamma_T} u^-(t,\bx).\label{eq:tcd_mp_4}
\end{align}
If $\sigma = 0$, then
\begin{align}
&\partial_t u-\varepsilon\Delta u + \bb\cdot\nabla u \le 0 \quad\mbox{in } \Omega_T
\quad \Longrightarrow \quad
\max_{(t,\bx) \in \overline{\Omega}_T} u(t,\bx) = \max_{(t,\bx) \in \Gamma_T}
u(t,\bx).\label{eq:tcd_mp_1}\\
&\partial_t u-\varepsilon\Delta u + \bb\cdot\nabla u \ge 0 \quad\mbox{in } \Omega_T
\quad \Longrightarrow \quad
\min_{(t,\bx) \in \overline{\Omega}_T} u(t,\bx) = \min_{(t,\bx) \in \Gamma_T}
u(t,\bx).\label{eq:tcd_mp_3}
\end{align}
\end{theorem}

Consider problem \eqref{eq:time_cdr} with $\sigma=0$ and $f=0$. For a 
sufficiently smooth solution, it follows from \eqref{eq:tcd_mp_1} and
\eqref{eq:tcd_mp_3} that
\begin{equation}\label{eq:tcd_mp_both}
\min_{(t,\bx) \in \Gamma_T} u(t,\bx) \le u(t,\bx) \le \max_{(t,\bx) \in \Gamma_T} u(t,\bx)
\quad \forall\ (t,\bx)\in \Omega_T.
\end{equation}

Physical quantities whose behavior is modeled with 
convection-diffusion-reaction equations are often by definition non-negative, 
like concentrations or the temperature (in Kelvin). The mathematical 
formulation of this property is the so-called positivity preservation. Let the 
data of \eqref{eq:time_cdr} be non-negative, i.e., $f\ge 0$ in $\Omega_T$ (no 
sinks), $g\ge 0$ on $(0,T]\times \partial\Omega$, and $u_0 \ge 0$ in $\Omega$. 
Then it follows from \eqref{eq:tcd_mp_4} that
$u\ge 0$ in $\Omega_T$. If $\sigma<0$, then as already explained in 
Section~\ref{sec:tcd_cont_prob}, one can transform problem \eqref{eq:time_cdr} 
to an equivalent problem for $\check u(t,\bx) = u(t,\bx) \exp(-\kappa t)$ with 
non-negative reaction coefficient and non-negative data on the right-hand
sides. Then \eqref{eq:tcd_mp_4} implies that $\check u\ge 0$ in
$\Omega_T$ whence also $u\ge 0$ in $\Omega_T$. Thus, independently of the sign
of $\sigma$, the non-negativity of the data $f$, $g$, $u_0$ is sufficient for 
obtaining a non-negative solution.  Therefore,
besides the local and global DMP,  also the positivity preservation
of discretizations of the time-dependent 
problem is often studied in the literature.

Consider from now on the case that the right-hand side of \eqref{eq:time_cdr}
is identically zero. Moreover, for simplicity, we assume that the boundary 
condition $g$ is independent of time.
Let the time interval be decomposed by $0=t^0 < t^1< \ldots < t^J = T$. 
After having applied a one-step $\theta$ scheme in time and a linear discretization in space to \eqref{eq:time_cdr}, one arrives at time instant 
$t^{n+1}$ at an algebraic problem of the form 
\begin{equation}\label{eq:tcd_disc_prob}
\bbB \bu^{n+1} = \mK \bu^n,
\end{equation}
where $\bu^{n+1}$ is the sought solution vector at $t^{n+1}$ and $\bu^n$ is the
solution at time $t^n$. The matrices $\bbB$ and $\mK$ have the form
\eqref{system-matrix} so that the last $N-M$ equations of 
\eqref{eq:tcd_disc_prob} set the Dirichlet boundary conditions for $\bu^{n+1}$;
we recall that the last $N-M$ entries of $\bu^n$ and $\bu^{n+1}$ contain the 
boundary values. We assume that the matrices $\bbB$ and $\mK$ possess the
typical sparsity pattern corresponding to discretizations with $\bbP_1$ finite
elements, i.e.,
\begin{equation}\label{eq:tcd_S_i}
   b_{ij}=k_{ij}=0\qquad\forall\,\,j\not\in S_i^{}\cup\{i\}\,,\,\,
   1\le i\le M\,,
\end{equation}
where $S_i^{}$ is defined by \eqref{eq:def_S_i}.

Since the right-hand side of \eqref{eq:time_cdr} is identically zero, all cases of the maximum principle from Theorem~\ref{thm:tcd_weak_mp} apply. 
Now, conditions on the matrices $\bbB$ and $\mK$ will be derived such 
that a discrete version of \eqref{eq:tcd_mp_both} holds. 

\begin{lemma}[Local DMP] \label{lem:tcd_local_DMP}
Consider any $n\in\{0,\dots,J-1\}$ and denote
\begin{equation*}
u_i^{\mathrm{min}} =\min \left\{ \min_{j\in S_i \cup \{i\}} u_j^n,
\min_{j\in S_i} u_j^{n+1}\right\}, \quad u_i^{\mathrm{max}} = \max \left\{ \max_{j\in S_i \cup \{i\}} u_j^n,
\max_{j\in S_i} u_j^{n+1}\right\}
\end{equation*}
for $i=1,\ldots,M$. Assume that \eqref{eq:tcd_disc_prob} holds with
\eqref{eq:tcd_S_i} and
\begin{equation}\label{eq:tcd_mat_coeff}
b_{ii} > 0, \,\,\, k_{ii} \ge 0, \,\,\, b_{ij}\le 0, \,\,\, k_{ij} \ge 0 
\qquad \forall\ j\in S_i\,,\,\,1\le i\le M\,.
\end{equation}
If 
\[
\sum_{j\in S_i \cup \{i\}} b_{ij} = \sum_{j\in S_i \cup \{i\}} k_{ij}\,,
\qquad 1\le i\le M\,,
\]
then it follows that 
\begin{equation*}
u_i^{\mathrm{min}} \le u_i^{n+1} \le u_i^{\mathrm{max}}\,,\qquad 1\le i\le M\,.
\end{equation*}
\end{lemma}

\begin{proof} The proof will be given for the upper bound, the statement for 
the lower bound can be derived analogously. Consider any $i\in\{1,\dots,M\}$. 
Let $w_j = u_j^{n+1} - u_i^{\mathrm{max}}$ and 
$v_j = u_j^n - u_i^{\mathrm{max}}$ for $j=1,\dots,N$. Then $w_j \le 0$ for all $j\in S_i$ and 
$v_j \le 0$ for all $j \in S_i \cup \{i\}$. A direct calculation, utilizing the 
assumption on the row sums, reveals that 
\begin{equation*}
b_{ii} w_i = k_{ii} v_i + \sum_{j\in S_i} \left(k_{ij} v_j - b_{ij} w_j\right).
\end{equation*}
By construction and assumption \eqref{eq:tcd_mat_coeff}, the coefficient on the 
left-hand side is positive and the right-hand side is non-positive. Hence, one 
obtains $w_i \le 0$, which is equivalent to $u_i^{n+1} \le u_i^{\mathrm{max}}$.
\end{proof}

For studying global properties, it is convenient to write \eqref{eq:tcd_disc_prob} 
without the (trivial) equations for the values on the Dirichlet boundary:
\begin{equation}\label{eq:tcd_disc_prob_inner}
(\mBI | \mBB) \begin{pmatrix}\buI^{n+1}\\ \buB^{n+1}\end{pmatrix} = 
(\mKI | \mKB) \begin{pmatrix} \bu_{\mathrm{I}}^n\\ \bu_{\mathrm{B}}^n \end{pmatrix},
\end{equation}
with 
$\mBI, \mKI \in \mathbb R^{M\times M}$, $\mBB, \mKB\in \mathbb R^{M\times (N-M)}$,
$\buI^{n+1}, \bu_{\mathrm{I}}^n \in \mathbb R^M$, and $\buB^{n+1}, \bu_{\mathrm{B}}^n \in \mathbb R^{N-M}$.
It will be assumed that $\mBI$ is invertible. 
Note that from setting the Dirichlet boundary conditions, $\buB^{n+1} = \bu_{\mathrm{B}}^n$, but 
for the following considerations, these vectors might be even different. 

\begin{definition}[Positivity preservation] Method \eqref{eq:tcd_disc_prob_inner}
is said to be positivity preserving if the inequality $\buI^{n+1}\ge 0$ is valid for 
all non-negative vectors $\buB^{n+1}$, $\bu_{\mathrm{I}}^n$, $\bu_{\mathrm{B}}^n$.
\end{definition}

\begin{theorem}[Necessary and sufficient conditions for positivity preservation]\label{thm:tcd_pos_pres_disc_gen}
Method \eqref{eq:tcd_disc_prob_inner} is positivity preserving if and only if the two conditions
\begin{eqnarray}\label{eq:DMP_cond_rhs}
\mBI^{-1} (\mKI | \mKB) &\ge& 0,\\
\label{eq:DMP_cond_bdry}
-\mBI^{-1} \mBB &\ge& 0,
\end{eqnarray}
hold.
\end{theorem}

\begin{proof} 
The statement of the theorem follows immediately from 
the following representation
\begin{equation*}
\buI^{n+1} = \mBI^{-1} (\mKI | \mKB)\begin{pmatrix}\bu_{\mathrm{I}}^n\\ \bu_{\mathrm{B}}^n \end{pmatrix} - 
\mBI^{-1} \mBB \buB^{n+1},
\end{equation*}
which is obtained from \eqref{eq:tcd_disc_prob_inner}.
\end{proof}

\begin{definition}[Global DMP]\label{def:tcd_global_dmp} Method \eqref{eq:tcd_disc_prob_inner} is said 
to satisfy the (global) DMP if 
\begin{equation}\label{eq:tcd_dmp_def}
\min\left\{\buB^{n+1}, \bu_{\mathrm{I}}^n, \bu_{\mathrm{B}}^n \right\} \le u_i^{n+1} \le \max\left\{ \buB^{n+1}, \bu_{\mathrm{I}}^n, \bu_{\mathrm{B}}^n\right\}, \quad 1 \le i \le M,
\end{equation}
for each choice $\buB^{n+1},\bu_{\mathrm{I}}^n, \bu_{\mathrm{B}}^n$, where 
$(u_i^{n+1})_{i=1}^M=\buI^{n+1}$.
\end{definition}

In the following, a vector of length $k \in \mathbb N$ where all 
entries are $1$ is denoted by $\boldsymbol 1_k$. 

\begin{theorem}[Necessary and sufficient conditions for the global DMP]
\label{thm:tcd_dmp_thm_gen}Method \eqref{eq:tcd_disc_prob_inner} satisfies the 
global DMP if and only if \eqref{eq:DMP_cond_rhs}, \eqref{eq:DMP_cond_bdry}, 
and 
\begin{equation}\label{eq:DMP_cond_rows}
(\mBI | \mBB) \boldsymbol 1_N = 
(\mKI | \mKB) \boldsymbol 1_N
\end{equation}
hold, i.e., the $i$th row sums of $(\mBI | \mBB)$ and $(\mKI | \mKB)$ are identical, $i=1,\ldots,M$.
\end{theorem}

\begin{proof} The proof follows \cite{FH06}.

\emph{i) DMP $\Longrightarrow$ \eqref{eq:DMP_cond_rhs},
\eqref{eq:DMP_cond_bdry}, \eqref{eq:DMP_cond_rows}.} If $\buB^{n+1}$, 
$\bu_{\mathrm{I}}^n$, and $ \bu_{\mathrm{B}}^n$ are arbitrary non-negative 
vectors, then the left-hand inequality of \eqref{eq:tcd_dmp_def} states that 
$\buI^{n+1}$ is also non-negative. Hence, the method is positivity preserving 
and it follows from Theorem~\ref{thm:tcd_pos_pres_disc_gen} that 
\eqref{eq:DMP_cond_rhs} and \eqref{eq:DMP_cond_bdry} are satisfied. 

Choosing in \eqref{eq:tcd_dmp_def} $\buB^{n+1} = \boldsymbol 1_{N-M}$,
$\bu_{\mathrm{I}}^n = \boldsymbol 1_M$, and $\bu_{\mathrm{B}}^n = \boldsymbol 1_{N-M}$ yields $\buI^{n+1} = \boldsymbol 1_M$.
Inserting these vectors in \eqref{eq:tcd_disc_prob_inner} shows that \eqref{eq:DMP_cond_rows} is satisfied. 

\emph{ii) \eqref{eq:DMP_cond_rhs},
\eqref{eq:DMP_cond_bdry}, \eqref{eq:DMP_cond_rows} $\Longrightarrow$ DMP.} 
Denoting $u_{\mathrm{max}}^n =
\max\{\buB^{n+1},\bu_{\mathrm{I}}^n,\bu_{\mathrm{B}}^n\}$ and using 
\eqref{eq:DMP_cond_rhs}, \eqref{eq:DMP_cond_rows}, and
\eqref{eq:DMP_cond_bdry}, gives 
\begin{eqnarray*}
\buI^{n+1} &=& - \mBI^{-1} \mBB \buB^{n+1} +
\mBI^{-1} (\mKI | \mKB)\begin{pmatrix}\bu_{\mathrm{I}}^n\\ \bu_{\mathrm{B}}^n \end{pmatrix} \\
& \le & - \mBI^{-1} \mBB \buB^{n+1} + u_{\mathrm{max}}^n
\mBI^{-1} (\mKI | \mKB) \boldsymbol 1_N\\
& = & - \mBI^{-1} \mBB \buB^{n+1} + u_{\mathrm{max}}^n
\mBI^{-1} (\mBI | \mBB) \boldsymbol 1_N\\
& = & - \mBI^{-1} \mBB(\buB^{n+1}-u_{\mathrm{max}}^n\boldsymbol 1_{N-M}) +
u_{\mathrm{max}}^n \boldsymbol 1_M \le u_{\mathrm{max}}^n \boldsymbol 1_M\,,
\end{eqnarray*}
which is equivalent to the right-hand inequality in \eqref{eq:tcd_dmp_def}.
The left-hand inequality is proven similarly. 
\end{proof}

The concepts of positivity preservation and of the global DMP can be extended 
to non-vanishing right-hand sides, see \cite{FH06}. The necessary and 
sufficient requirements on the matrices for the satisfaction of these 
properties are the same as given in Theorems~\ref{thm:tcd_pos_pres_disc_gen} 
and~\ref{thm:tcd_dmp_thm_gen}.

\begin{corollary}[Positivity preservation and global DMP for monotone matrices]
\label{cor:tcd_pospres_DMP_mon_mat}
Let the matrix 
\[
\bbB = \begin{pmatrix} \mBI & \mBB \\ \mathbb O & \mathbb I \end{pmatrix}
\]
be monotone and let $\mK \ge 0$. Then method \eqref{eq:tcd_disc_prob_inner} 
is positivity preserving. If, in addition, the $i$th row sums of $\bbB$ and 
$\mK$ are identical, $i=1,\ldots,M$, then method \eqref{eq:tcd_disc_prob_inner}
satisfies the global DMP.
\end{corollary}

\begin{proof} 
From computing the inverse of $\bbB$, compare \eqref{eq:matrix_inverse}, it
follows that $\mBI^{-1}\ge 0$ and $-\mBI^{-1}\mBB\ge 0$. Since $\mK\ge 0$, the
conditions \eqref{eq:DMP_cond_rhs} and \eqref{eq:DMP_cond_bdry} are satisfied.
Thus, the corollary follows from Theorems~\ref{thm:tcd_pos_pres_disc_gen} 
and~\ref{thm:tcd_dmp_thm_gen}.
\end{proof}

\begin{remark}
Note that if $\bbB$ is a monotone matrix, $\mK \ge 0$, and $\bu^n\ge 0$,
then it immediately follows that the solution of \eqref{eq:tcd_disc_prob}
satisfies $\bu^{n+1} \ge 0$.
\hspace*{\fill}$\Box$\end{remark}

Another property that is often studied for discretizations of scalar evolutionary transport problems is the 
local extremum diminishing (LED) property. Considering a method that is only semi-discrete in space, the 
LED condition is as follows: if $u_i$ is a local maximum in space, then $du_i/dt \le 0$ and if $u_i$ is a local 
minimum in space, then $du_i/dt \ge 0$, i.e., a local maximum does not increase and a local minimum does not decrease. 
For a fully discrete method, discretized with a one-step $\theta$-scheme, the LED property states that if 
$u_i^{n+\theta}= \theta u_i^{n+1} + (1-\theta) u_i^n$ is a local maximum in space, then 
$u_i^{n+1}\le u_i^n$ and similarly for a local minimum, e.g., see \cite{BBH17}.

Section~\ref{sec:tcd_femfct} will discuss a class of nonlinear discretizations 
in some detail. A motivation for considering such discretizations for the 
convection-dominated regime is provided by a study of the limit case of 
\eqref{eq:time_cdr} with respect to small diffusion, i.e., the transport 
equation where $\varepsilon = 0$. Consider this case with constant convection 
$b\neq0$ and $\sigma=f=0$ in one dimension on the infinite domain 
$\Omega = (-\infty, \infty)$. The domain is decomposed using an equidistant 
grid with mesh width $h$ and the nodes $x_i$, $i\in \mathbb Z$. Then, the 
application of 
an explicit one-step $\theta$-scheme
leads to a problem of the form 
\begin{equation}\label{eq:one-step-scheme}
   u_j^{n+1} = \sum_{i=-S}^S \gamma_i u_{j+i}^{n}, \qquad j \in \mathbb Z\,,
\end{equation}
where $S$ is determined by the width of the stencil.
For this kind of problem there exists the notion of a monotonicity preserving scheme: for all monotone discrete initial conditions $u^0$, the solution 
$u^n$ possesses the same monotonicity for all $n\ge 1$. It can be shown that the scheme is 
monotonicity preserving if and only if $\gamma_i \ge 0$
for all $i\in\{-S,\dots,S\}$. Then, Godunov's order barrier theorem \cite{God59} states that if 
$C_{\mathrm{CFL}} = {|b|\tau}/h \not\in \mathbb N$, a linear monotonicity
preserving method of form 
\eqref{eq:one-step-scheme} cannot compute solutions exactly that are polynomials of degree~$2$. 
Hence, a linear monotonicity preserving method has to be of low order. For a more recent presentation 
of this topic see \cite{Wes01}. Using an implicit one-step scheme or a linear 
multi-step scheme instead of an explicit one-step scheme does not solve this 
issue, see \cite[Thm.~9.2.4]{Wes01}. 

The condition on the non-negativity of $\gamma_i$ resembles condition \eqref{eq:DMP_cond_rhs}, which is necessary 
for the positivity preservation and the satisfaction of the DMP. Thus, one can expect that for \eqref{eq:time_cdr}, 
in the convection-dominated regime, a linear discretization that possesses these properties
will be only of low-order. There is no mathematical proof of this expectation but computational evidence. This issue 
motivates the construction of nonlinear discretizations to obtain accurate schemes for \eqref{eq:time_cdr} that are positivity preserving 
and satisfy the DMP. 

\subsection{Linear methods}\label{sec:tcd_ext_methods_steady}

Utilizing a one-step $\theta$-scheme in combination with the Galerkin 
or some stabilized finite element method for the discretization of 
\eqref{eq:time_cdr} with $f=0$ leads to an algebraic system of the form 
\begin{equation}\label{eq:tcd_alg_system}
\left(\bbMc + \theta \tau \bbA_1 \right)^M \bu^{n+1} = \left(\bbMc - (1-\theta) \tau \bbA_2 \right)^M \bu^n,
\quad u_{i}^{n+1} = g_{i-M}^{n+1}, 
\end{equation}
$i=M+1,\ldots,N$,
where $\bbMc$ is the consistent mass matrix defined in \eqref{def-mass-matrix}, 
$\bbA_1$, $\bbA_2$ are stiffness matrices, and $\tau= t^{n+1}-t^n$ is the current time step. 
Consider a uniform spatial grid with mesh width $h$. Then, for standard Lagrangian finite element spaces, 
$\bbMc$ possesses positive off-diagonal entries of order $\mathcal O(h^d)$, compare \eqref{local-mass-exact} for $\bbP_1$ finite elements. 
Consequently, $\bbMc$ is not an 
M-matrix and as can be checked easily, e.g., for a one-dimensional problem, $\bbMc$ is not a monotone
matrix. The off-diagonal entries of $\tau\bbA_1$ are of order $\mathcal O(\tau h^{d-2})$ for the diffusive term 
and $\mathcal O(\tau h^{d-1})$ for the convective term. Hence, if $\tau$ is sufficiently small, the 
system matrix of \eqref{eq:tcd_alg_system} cannot be an M-matrix. In particular, any finite element 
analysis that considers the so-called continuous-in-time situation, i.e., only a semi-discretization 
in space, cannot apply the concept of M-matrices. It is shown in 
\cite{TW08} that a standard continuous-in-time finite element 
discretization of the heat equation cannot be positivity preserving and it cannot satisfy the global DMP. 
One can only hope 
for non-positive off-diagonal entries of the system matrix of \eqref{eq:tcd_alg_system} if 
$\tau$ is of order $\max\{h,h^2\}$. In fact, for the heat equation, discretized with a one-step 
$\theta$-scheme and the Galerkin FEM, sufficient conditions for the satisfaction
of the DMP were derived in \cite{FH06} that include a lower and an upper bound for the length of the 
time step, which are both of order $\mathcal O(h^2)$.

Note that this issue does not appear for finite volume and finite difference methods, where the 
temporal discretization leads to a diagonal matrix with positive diagonal entries. Studying 
positivity preservation and the DMP with the concept of M-matrices for finite element methods, the common way consists
in applying mass lumping, which is presented in Section~\ref{Sec:FEM-Matrices}. Utilizing a lumped mass matrix, the positivity
preservation can be proven for the heat equation in two dimensions, $\bbP_1$ finite elements, and 
under certain additional assumptions, see \cite{STW10}. An extension of this
result to three dimensions is also possible.

In \cite{FKK12} a class of problems was studied which includes the linear 
convection-diffusion-reaction equation as a special case. The considered
discretization was a one-step $\theta$-scheme combined with the Galerkin FEM.
The DMP is proven under a number of assumptions. 
Because of using the Galerkin FEM, the mesh width has to be sufficiently small, 
compare \cite[Thm.~5.2 (ii)]{FKK12},
in particular the bound for the mesh width tends to zero as $\varepsilon\to 0$. For a
sufficiently small mesh width, there is a lower bound for the time step of order $\mathcal O(h^2)$. 

As already mentioned in Section~\ref{Sec:cd_linear_methods_upw}, the upwind 
finite element method proposed in \cite{Tabata77} was formulated and studied 
for a two-dimensional time-dependent equation. The analysis is performed for 
the forward 
Euler scheme, where a lumped mass matrix is utilized, so that the
discretization of the time derivative corresponds to a finite difference or
finite volume one.
The key ingredient of this method 
is the discretization of the convective term, which is described in Section~\ref{Sec:cd_linear_methods_upw}.
From the proof presented in \cite{Tabata77}, it can be seen that the
assumptions of Corollary~\ref{cor:tcd_pospres_DMP_mon_mat} are satisfied
under an appropriate CFL condition, hence the method satisfies the DMP.
In the final part of \cite{Tabata77}, it is mentioned that the analysis can be 
extended to the (mass lumped) backward Euler scheme and to time-dependent convection fields. 

The upwind method proposed and analyzed in \cite{BT81} was also already presented in 
Section~\ref{Sec:cd_linear_methods_upw}. In \cite{BT81}, it was studied 
for the conservative form \eqref{eq:time_cdr_cons} of the convection-diffusion equation. 
In contrast with the method from \cite{Tabata77}, it satisfies a discrete analog 
of a mass conservation property if \eqref{eq:time_cdr_cons} is equipped with so-called free boundary condition
\[
\varepsilon\frac{\partial u}{\partial\bn} - \bb\cdot \bn\,u = 0\quad \textrm{on}\;(0,T]\times \partial\Omega.
\]
The upwind method is analyzed for this boundary condition, steady-state 
convection fields, and the mass lumped forward Euler scheme so that an 
appropriate CFL condition becomes necessary throughout the analysis. A brief 
description of the discretization of the convective term, leading to a 
convection matrix $\tilde\bbA_{\mathrm{c}}$, is already provided in 
Section~\ref{Sec:cd_linear_methods_upw}. Thus, the discretization of
\eqref{eq:time_cdr_cons} with $f=0$ and the free boundary condition is of the 
form 
\[
 (\bbMl)^M\bu^{n+1} =  (\bbMl)^M\bu^n
 -\tau(\varepsilon\bbA_{\mathrm{d}}+\tilde\bbA_{\mathrm{c}})^M\bu^n.
\]
The construction of $\tilde\bbA_{\mathrm{c}}$ assures that its row sums vanish.
The row sums of $\bbA_{\mathrm{d}}$ also vanish, see 
\eqref{eq:A_d_row_sum}, and hence the positivity preservation and the 
satisfaction of the global DMP for this upwind method can be inferred from 
Corollary~\ref{cor:tcd_pospres_DMP_mon_mat}.

\begin{remark} The techniques of 
\cite{LH10,Hua11} developed for problems with heterogeneous anisotropic
diffusion, see Remark~\ref{rem:aniso_diff_poisson}, were applied to study also the DMP for the 
heat equation in \cite{LH13}. The $\mathbb P_1^{}$ finite element in space is 
combined with a one-step $\theta$-method in time. Concerning the spatial mesh, 
the same conditions apply as for the steady-state diffusion problem. Using a 
lumped mass matrix, one obtains a restriction for the length of the time step,  
which is of the form 
\begin{equation*}
   \tau \le C \min_{K \in \calT_h}\min_{F\in\calF_K^{}}
   \frac{h_{K,F}^2}{\lambda_{\mathrm{max}} (\bar{\mathbb E}_K)}, 
\end{equation*}
where $h_{K,F}$ is the height from the facet $F\subset K$ to the vertex of
$K$ opposite $F$ and $\bar{\mathbb E}_K$ is defined to be the integral mean of
the diffusion tensor $\mathbb E$ on $K$.
\hspace*{\fill}$\Box$\end{remark}

\subsection{FEM Flux-Corrected-Transport (FCT) schemes}\label{sec:tcd_femfct}

A physical quantity is called extensive if it scales with the size of the physical 
problem. Examples are mass, momentum, or energy. Fluxes are quantities of an 
extensive variable that moves from one location in space to another one. That means, 
the amount of the variable that is removed from the first location is added at the 
second location. If numerical methods are formulated in terms of fluxes, 
they are called conservative if the same principle is applied as mentioned above: 
what is removed from one degree of freedom is added to another one. 
The conservation
of physical quantities in a numerical method contributes to the physical consistency of this method and 
thus, it helps that the method becomes accepted by practitioners.

The usual starting point for the construction of numerical methods based on fluxes is the 
conservative form \eqref{eq:time_cdr_cons} of the convection-diffusion equation. 
Natural discretizations for this form are finite difference and finite volume methods. 

For illustrating the concept of numerical fluxes, consider a finite difference method
for the one-dimensional analog of \eqref{eq:time_cdr_cons}
\begin{equation}\label{eq:time_cdr_cons_1d} 
\begin{array}{rcll}
\partial_t u+ \partial_x\left(-\varepsilon \partial_x u + b u\right) &=& 0 &\textrm{in}\;
(0,T]\times \Omega\,,\\
u&=& 0 & \textrm{on}\;(0,T]\times \partial\Omega\,,\\
u(\cdot,0) &=& u_0 & \textrm{in}\;\Omega\,,
\end{array}
\end{equation}
with $\Omega = (\xi_{\mathrm{l}}, \xi_{\mathrm{r}})$, 
$\xi_{\mathrm{l}} < \xi_{\mathrm{r}}$. Let $\overline\Omega$ be triangulated 
using an equidistant grid with mesh width $h$ and nodes
$\{x_i\}_{i=1}^N$, $x_1 = \xi_{\mathrm{l}}$, $x_N = \xi_{\mathrm{r}}$,
$x_i<x_{i+1}$.
Consider the step from time instant $t^n$ to $t^{n+1}$. 
A finite difference approximation of \eqref{eq:time_cdr_cons_1d} is said to be of 
conservative form, if it can be written for inner nodes in the form 
\[
u_i^{n+1} = u_i^n + \frac{\tau}{\frac12 (x_{i+1}-x_{i-1})}\left( f_{i-1/2} - f_{i+1/2}\right),
\]
where
$f_{i+1/2}$ and $f_{i-1/2}$ are numerical fluxes depending on diffusion and convection 
at one or several time levels. Utilizing the explicit Euler scheme for 
discretizing \eqref{eq:time_cdr_cons_1d} in time, the 
standard 3~point stencil for the discretization of the second derivative and 
a central finite difference defined on the points $x_{i+1/2} = (x_{i+1}+x_i)/2$
and $x_{i-1/2} = (x_{i}+x_{i-1})/2$ for the convective term yields
\begin{eqnarray*}
u_i^{n+1} & = & u_i^n + \tau \left[ \varepsilon \frac{u_{i+1}^n - 2u_i^n + u_{i-1}^n}{h^2} - \frac{b_{i+1/2}^nu_{i+1/2}^n - b_{i-1/2}^nu_{i-1/2}^n}{h}\right] \\
&=& u_i^n + \frac{\tau}{h} \left[-\varepsilon \frac{u_{i}^n - u_{i-1}^n}h + b_{i-1/2}^nu_{i-1/2}^n - \left(-\varepsilon \frac{u_{i+1}^n - u_i^n}h + b_{i+1/2}^nu_{i+1/2}^n\right) \right].
\end{eqnarray*}
Hence, the numerical analog of the fluxes of the continuous
problem, see the end of Section~\ref{sec:tcd_cont_prob}, is given by 
\[
 f_{i+1/2} = -\varepsilon \frac{u_{i+1}^n - u_i^n}h + b_{i+1/2}^nu_{i+1/2}^n\,,
\]
where the first term on the right-hand side is the numerical diffusive flux and 
the second term the numerical convective flux.  Usually, the values 
$u_{i\pm1/2}^n$ at $x_{i\pm1/2}$ are approximated using the values at the 
neighboring nodes with the aim to obtain a stable discretization. A classical 
example is the one-sided upwind approximation. 

The first development and implementation of a FCT scheme was performed for a finite difference method in one dimension 
in \cite{BB73}. Consider the step from one discrete time level to the next one, then 
the basic approach is as follows:
\begin{enumerate}[leftmargin=*]
\item\label{enu:low_order} A (linear) scheme is needed that guarantees that no nonphysical values are 
computed. Such a scheme has to utilize low-order fluxes, which possess a large 
amount of numerical diffusion. 
\item\label{enu:high_order} A second (linear) scheme with high-order fluxes is used, which is highly 
accurate for smooth regions of the solution. This scheme has only a small amount 
of numerical diffusion and its solution has spurious oscillations in a vicinity of
layers or shocks. 
\item\label{enu:anti_diff_flux} So-called antidiffusive fluxes are defined by the difference of the high and 
low-order fluxes from the two schemes. 
\item\label{enu:limiter} The solution at the new time level is obtained by adding appropriately 
weighted (limited) antidiffusive fluxes to the solution of the low-order scheme. 
The limiting process has to ensure that no unphysical values are created in this step. 
For smooth parts of the solution, the high-order scheme should be recovered. 
\end{enumerate}
FCT schemes were then transferred to one-dimensional finite volume methods. 
It turned out that the limiter for one-dimensional problems proposed in \cite{BB73} 
does not work properly in multiple dimensions. Thus, the next milestone in the 
development of FCT schemes was the proposal of a new limiter that works in 
multiple dimensions in \cite{Zal79}, the nowadays so-called Zalesak limiter. 
This limiter will be described within the presentation of the FEM-FCT methods. 
A good survey of the motivations for deriving FCT schemes and their main 
design principles can be found in the paper \cite{Zal12}, which concentrates on finite 
volume schemes on structured grids. 

The development of FCT schemes for finite element methods was driven by the goal to 
apply the FCT methodology on unstructured grids. To this end, a concept that resembles 
fluxes was introduced in finite element methods, the so-called algebraic fluxes. 
Algebraic fluxes are quantities $f_{ij}$ between adjacent degrees of freedom
$i$ and $j$ that 
are derived from algebraic quantities like matrices and vectors and for which $f_{ij}
= - f_{ji}$ (the flux property) holds. The vast majority of FEM-FCT methods 
have been developed for $\bbP_1$ and $\bbQ_1$ finite elements, where the 
degrees of freedom are function values at the vertices of the mesh cells. The 
first FEM-FCT schemes were proposed in 
\cite{LMPV87,PC86}.
Since then, FEM-FCT schemes have been improved and further developed, e.g., in
\cite{KT02,Kuz06,Kuz09,Kuz12,LKSM17}, see also the surveys in \cite{Kuz12a}
and \cite[Chapters~6.3, 7.5, 7.6]{KH15}.
Nevertheless, theoretical results on FEM-FCT schemes for time-dependent
convection-diffusion-reaction problems are far less developed than for the
related algebraically stabilized methods proposed for the steady-state problem
and discussed in Section~\ref{Sec:cd_nonlinear_methods_afc}. In particular, we
are not aware of any error estimates.

Whereas the FCT methodology is used in finite difference and finite volume schemes directly to define a discretization 
of the convection and diffusion operators with the goal to satisfy the DMP
locally, its application in the FEM is more indirect. There, the Galerkin FEM
discretization is reformulated equivalently such that the system matrix 
becomes an M-matrix and then the FCT methodology is utilized to modify the 
right-hand side such that the M-matrix property of the system matrix allows to 
satisfy the global DMP and the positivity preservation. 

In the following, a FEM-FCT scheme will be presented in detail, thereby explaining the derivation 
and application of the Zalesak limiter. The starting point is now problem
\eqref{eq:time_cdr} and it is again assumed that the right-hand side vanishes.
Moreover, for simplicity, we assume that the velocity field $\bb$ does not
depend on time.

The high-order method from Step~\ref{enu:high_order} of the basic FCT approach is 
the standard Galerkin FEM. Using a one-step $\theta$-scheme as 
temporal discretization,  $\theta \in (0,1]$, leads to the linear algebraic 
system
\begin{equation}\label{eq:fct_gal}
\left(\bbMc + \theta \tau \bbAN\right)^M \bu^{n+1} = \left(\bbMc - (1-\theta)
\tau \bbAN\right)^M \bu^n, 
\end{equation}
where the matrix $\bbAN$ is defined by \eqref{eq:Galerkin_matrix}. The system
\eqref{eq:fct_gal} has to be supplemented by Dirichlet boundary conditions for 
$\bu^{n+1}$. Like for the algebraic flux correction in the steady case, we 
define the matrix $\bbD=(d_{ij})_{i,j=1}^N$ by \eqref{eq:matrix_D} using the
entries of $\bbAN$. In addition, we introduce the matrix 
$\mL=(l_{ij})_{i,j=1}^N$ defined by
\begin{equation*}
\mL = \bbAN+\mD. 
\end{equation*}
As discussed in Section~\ref{Sec:cd_nonlinear_methods_afc}, the matrix $\mL$ is
of non-negative type and $\bbD$ is positive semidefinite.

Next, the low-order scheme from Step~\ref{enu:low_order} of the basic FCT 
algorithm is given by 
\begin{equation}\label{eq:fct_low_order}
\left(\bbMl + \theta \tau \mL \right)^M \tilde\bu = \left(\bbMl - (1-\theta) \tau \mL \right)^M \bu^n,
\quad \tilde u_{i} = g_{i-M}^{n+1}, \ i=M+1,\ldots,N,
\end{equation}
where the lumped mass matrix $\bbMl$ is defined in 
\eqref{def-lumpedmass-matrix}. Due to the assumptions on the data of
\eqref{eq:time_cdr}, the matrix $(\bbAN)_{\mathrm{I}}^{}$ is positive definite
and hence also the matrix $(\bbMl + \theta \tau \mL)_{\mathrm{I}}^{}$ is
positive definite. Consequently, the system matrix of \eqref{eq:fct_low_order},
defined by extending the matrix $\left(\bbMl + \theta \tau \mL \right)^M$ by
the lower blocks of \eqref{system-matrix}, is invertible. Since it is of 
non-negative type, Corollary~\ref{cor:m-matrix} implies that the system matrix
of \eqref{eq:fct_low_order} is an M-matrix. Thus, in view of 
Corollary~\ref{cor:tcd_pospres_DMP_mon_mat}, method \eqref{eq:fct_low_order} is
positivity preserving if 
\begin{equation}\label{eq:tcd_femfct_00}
\left(\bbMl - (1-\theta) \tau \mL\right)^M \ge 0.
\end{equation}
To simplify the presentation, we denote the diagonal entries of $\bbMl$ by
$m_i$ instead of $\tilde{\mm}_{ii}^{}$ considered in
\eqref{def-lumpedmass-matrix}. Since $\mL$ is of non-negative type and
$\mL_{\mathrm{I}}^{}$ is positive definite, one has $l_{ii}>0$ and 
$l_{ij}\le 0$ for $j\neq i$, $i=1,\dots,M$. Hence \eqref{eq:tcd_femfct_00}
holds if and only if $(1-\theta) \tau l_{ii}\le m_i$ for all $i=1,\ldots,M$,
which is satisfied if $\theta=1$ or if
\begin{equation}\label{eq:tcd_fct_cfl}
   \tau\le\frac{m_i}{(1-\theta)l_{ii}}, \qquad i=1,\ldots,M\,. 
\end{equation}
This is a CFL condition which can be checked easily in simulations. 

Although the solution of \eqref{eq:fct_low_order} does not possess unphysical 
values under the CFL condition \eqref{eq:tcd_fct_cfl}, it is usually very 
inaccurate. In the FEM-FCT methodology, a correction 
term $\tau \overline \bff$ is added, which leads to a method of the form 
\begin{equation}\label{eq:fct_low_order_correction}
\left(\bbMl + \theta \tau \mL \right)^M \bu^{n+1} = \left(\bbMl - (1-\theta) \tau \mL \right)^M \bu^n + \tau \overline \bff.
\end{equation}
If the solution is smooth in the whole domain, then \eqref{eq:fct_low_order_correction}
should recover the high-order method. A direct calculation, subtracting \eqref{eq:fct_gal} from 
\eqref{eq:fct_low_order_correction}, shows that in this case 
\begin{equation*}
\tau \overline \bff = \left(\bbMl-\bbMc\right)^M\big(\bu^{n+1}-\bu^n\big) + \tau \left(\mD\right)^M \big(\theta \bu^{n+1} + (1-\theta) \bu^n\big)
\end{equation*}
is the appropriate correction. The expression on the right-hand side can be written in 
terms of algebraic fluxes. Using the definition \eqref{def-lumpedmass-matrix} of 
the lumped mass matrix and that the row sums of $\mD$ are zero, one obtains by a 
straightforward calculation 
\begin{eqnarray*}
\tau \left(\overline\bff\right)_i & = & 
\sum_{j=1}^N \left[- m_{ij}\big(u_j^{n+1} - u_i^{n+1}\big) + m_{ij}
\big(u_j^n - u_i^n\big)\right]\\
& &+ \tau \sum_{j=1}^N \left[\theta d_{ij} \big(u_j^{n+1} - u_i^{n+1}\big) 
+ (1-\theta)d_{ij} \big(u_j^n - u_i^n\big)
\right]. 
\end{eqnarray*}
For computing the right-hand side, again the matrices without having 
imposed Dirichlet boundary conditions are used.  Thus, the antidiffusive 
fluxes from Step~\ref{enu:anti_diff_flux} of the basic FCT algorithm 
are given by 
\begin{eqnarray}\label{eq:fct_anti_diff_fluxes}
\hspace*{5mm}
f_{ij} &=& \frac1{\tau} \left[ -m_{ij}\big(u_j^{n+1} - u_i^{n+1}\big) + m_{ij}
\big(u_j^n - u_i^n\big)\right] \\
&& +\left[ \theta d_{ij} \big(u_j^{n+1} - u_i^{n+1}\big) 
+ (1-\theta)d_{ij} \big(u_j^n - u_i^n\big)
\right], \quad i,j = 1,\ldots,N.\nonumber
\end{eqnarray}
Because $\bbMc$ and $\mD$ are symmetric matrices, one has $f_{ij} = - f_{ji}$. Note that the 
fluxes depend on (unknown) values of the numerical solution at time level
$t^{n+1}$.

Now, following Step~\ref{enu:limiter} of the basic FCT algorithm, the solution for the inner nodes at the 
next time level is defined by 
\begin{equation}\label{eq:fct_low_order_limiter}
\left(\bbMl + \theta \tau \mL \right)^M \bu^{n+1} = \left(\bbMl - (1-\theta) \tau \mL \right)^M \bu^n + \tau \left(\sum_{j=1}^N \alpha_{ij}f_{ij}\right)_{i=1}^M,
\end{equation}
where the limiters $\alpha_{ij} = \alpha_{ji} \in [0,1]$ have to be chosen appropriately. 

In order to apply the framework presented in Section~\ref{sec:tcd_framework}, the nonlinear 
problem \eqref{eq:fct_low_order_limiter} is written in the following way:
\begin{eqnarray}
\left(\bbMl\right)^M \overline\bu & = & \left(\bbMl - (1-\theta) \tau \mL \right)^M \bu^n,
\label{eq:tcd_femfct_2}\\
\left(\bbMl\right)^M \tilde\bu & = & \left(\bbMl\right)^M \overline\bu + \tau \left(\sum_{j=1}^N \left( \alpha_{ij}f_{ij}\right)^{[n+1]}\right)_{i=1}^M,\label{eq:tcd_femfct_1}\\
\left(\bbMl + \theta \tau \mL \right)^M \bu^{n+1} &= & \left(\bbMl\right)^M \tilde\bu, \label{eq:tcd_femfct_0}
\end{eqnarray}
where the superscript $[n+1]$ indicates that the fluxes and limiters depend on the solution 
at time instant $t^{n+1}$.
The function $\overline\bu$, which is equipped with the boundary conditions at
$t^{n+1-\theta}$, 
has to be computed only in the first step. This function is needed 
because it enters the 
definition of a lower and an upper bound in the limiting process, see \eqref{eq:tcd_uimin_max_fct} below. 
Then, solving \eqref{eq:tcd_femfct_1}--\eqref{eq:tcd_femfct_0} has to be performed with an 
iterative process, where the boundary conditions at $t^{n+1}$ are utilized in \eqref{eq:tcd_femfct_0}.

First, positivity preservation will be discussed. Let $\bu^n \ge 0$. Assuming
the validity of the CFL condition \eqref{eq:tcd_fct_cfl}, one has
\eqref{eq:tcd_femfct_00} and hence $\overline\bu\ge0$ since $\bbMl$ is a
diagonal matrix with positive diagonal entries. In the next step, 
$u_i^{\mathrm{min}} \ge 0$, $i=1,\ldots,M$, are chosen and the limiters 
are determined such that $\tilde u_i \ge u_i^{\mathrm{min}}$,
$i=1,\ldots, M$, 
in \eqref{eq:tcd_femfct_1}. Finally, since $\bbMl\ge 0$ and the system matrix 
of \eqref{eq:tcd_femfct_0} equipped with Dirichlet boundary conditions (which 
are assumed to be non-negative) is an M-matrix, it follows from 
Corollary~\ref{cor:tcd_pospres_DMP_mon_mat} that $\bu^{n+1}\ge 0$.

For studying the satisfaction of the global DMP 
(cf.~Definition~\ref{def:tcd_global_dmp}), the computation of the limiters has 
to be explained in detail. Let $\bu^{(m)}$ be an approximation of $\bu^{n+1}$ 
after the $m$th iteration for solving 
\eqref{eq:tcd_femfct_1}--\eqref{eq:tcd_femfct_0}. Then, the algebraic fluxes 
defined in \eqref{eq:fct_anti_diff_fluxes} are approximated using 
$\bu^{(m)}$ instead of $\bu^{n+1}$, leading to fluxes $f_{ij}^{(m)}$. Consider 
any $i\in\{1,\dots,M\}$ and define
\begin{equation}\label{eq:tcd_uimin_max_fct}
   \overline u_i^{\mathrm{min}} = \min_{j\in S_i\cup\{i\}} \overline u_j\,,
   \qquad 
   \overline u_i^{\mathrm{max}} = \max_{j\in S_i\cup\{i\}} \overline u_j\,,
\end{equation}
with $S_i$ given by \eqref{eq:def_S_i}. Then the limiters $\alpha_{ij}^{(m)}$, 
where the superscript indicates that they depend on $f_{ij}^{(m)}$, are 
computed such that 
\begin{equation}\label{eq:tcd_femfct_01}
\overline u_i^{\mathrm{min}} \le \tilde u_i\le \overline u_i^{\mathrm{max}},
\end{equation}
where $\tilde\bu$ is the solution of \eqref{eq:tcd_femfct_1} with the fluxes 
$f_{ij}^{(m)}$ and the limiters $\alpha_{ij}^{(m)}$. Consider the upper bound 
and introduce non-negative numbers $R_i^+$ such that 
$\alpha_{ij}^{(m)}\le R_i^+$ if $f_{ij}^{(m)}>0$. Then
\begin{eqnarray*}
\tilde u_i &=& \overline u_i + \frac{\tau}{m_i} \sum_{j=1}^N \alpha_{ij}^{(m)}f_{ij}^{(m)}
\le \overline u_i + \frac{\tau}{m_i} \sum_{j=1}^N
\alpha_{ij}^{(m)}\left(f_{ij}^{(m)}\right)^+ \\
&\le& \overline u_i + \frac{\tau}{m_i} R_i^+
\sum_{j=1}^N\left(f_{ij}^{(m)}\right)^+.
\end{eqnarray*}
Thus, to satisfy the upper bound in \eqref{eq:tcd_femfct_01}, it suffices to
require that
\begin{equation}\label{eq:tcd_R_i_+_der}
   R_i^+ \le \frac{m_i}{\tau} \left(\overline u_i^{\mathrm{max}}-\overline u_i\right) \left(\sum_{j=1}^N\left(f_{ij}^{(m)}\right)^+\right)^{-1},
\end{equation}
where the right-hand side is non-negative thanks to the definition 
\eqref{eq:tcd_uimin_max_fct} of $\overline u_i^{\mathrm{max}}$. Note that if
$\big(f_{ij}^{(m)}\big)^+=0$ for all $j=1,\dots,N$, then the upper bound in 
\eqref{eq:tcd_femfct_01} always holds and $R_i^+$ can be defined arbitrarily.
Similarly, to satisfy the lower bound in \eqref{eq:tcd_femfct_01}, it suffices
to require that $\alpha_{ij}^{(m)}\le R_i^-$ if $f_{ij}^{(m)}<0$ with
\begin{equation}\label{eq:tcd_R_i_-_der}
   R_i^- \le \frac{m_i}{\tau} \left(\overline u_i^{\mathrm{min}}-\overline u_i\right) \left(\sum_{j=1}^N\left(f_{ij}^{(m)}\right)^-\right)^{-1}.
\end{equation}
Like in the previous case, if $\big(f_{ij}^{(m)}\big)^-=0$ for all 
$j=1,\dots,N$, then the lower bound in \eqref{eq:tcd_femfct_01} always holds 
and $R_i^-$ can be defined arbitrarily. Since the limiters need to belong to 
$[0,1]$ by definition, one has to require $R_i^+\le1$ and $R_i^-\le1$ besides
the conditions \eqref{eq:tcd_R_i_+_der} and \eqref{eq:tcd_R_i_-_der}. In 
addition, one has to take into account that the flux property is maintained 
after having applied the limiters, i.e., 
$\alpha_{ij}^{(m)}f_{ij}^{(m)} = -\alpha_{ji}^{(m)}f_{ji}^{(m)}$, which
requires $\alpha_{ij}^{(m)}=\alpha_{ji}^{(m)}$ since 
$f_{ij}^{(m)} = - f_{ji}^{(m)}$. Thus, one has to take the smaller value of the
above-derived bounds for $\alpha_{ij}^{(m)}$ and $\alpha_{ji}^{(m)}$.
Summarizing all these considerations leads to the algorithm for the Zalesak 
limiter from \cite{Zal79}, where for the sake of clarity the iteration index is 
neglected in its presentation:
\begin{enumerate}[leftmargin=*]
\item Compute
\[
P_i^+ = \sum_{j=1,j\neq i}^N f_{ij}^+, \qquad 
P_i^- = \sum_{j=1,j\neq i}^N f_{ij}^-.
\]
\item\label{alg:tcd_zal2} Compute
\[
Q_i^+ = \frac{m_i}{\tau} \left(\overline u_i^{\mathrm{max}}-\overline
u_i\right), \qquad
Q_i^- = \frac{m_i}{\tau} \left(\overline u_i^{\mathrm{min}}-\overline u_i\right).
\]
\item \label{alg:tcd_zal3} Compute
\[
R_i^+ = \min \left\{ 1, \frac{Q_i^+}{P_i^+} \right\}, \qquad
R_i^- = \min \left\{ 1, \frac{Q_i^-}{P_i^-} \right\}.
\]
If the denominator is zero, set the value equal to $1$. In addition, both values
are set to be $1$ at Dirichlet nodes.
\item Compute
\[
\alpha_{ij} = \begin{cases}
\min\{R_i^+,R_j^-\} & \mbox{if } f_{ij} > 0, \\
1 & \mbox{if } f_{ij} = 0,\\
\min\{R_i^-,R_j^+\} & \mbox{if } f_{ij} <0.
\end{cases}
\]
Note that the value for $f_{ij}=0$ does not possess any impact. 
\end{enumerate}
It should be emphasized that, like in the steady-state case, the fluxes and limiters are computed on the 
basis of the matrices for Neumann boundary conditions. 

The nonlinear discretization \eqref{eq:fct_low_order_limiter}, or equivalently \eqref{eq:tcd_femfct_2}--\eqref{eq:tcd_femfct_0}, together with a limiter of the form of Zalesak's limiter and fluxes depending on 
$\bu^{n+1}$ is 
called nonlinear FEM-FCT scheme. The standard approach for computing an approximation to the solution, which is already 
sketched above, is summarized in Algorithm~\ref{alg:nonlin_femfct}. The following theorem shows that, under 
appropriate conditions, all iterates satisfy the global DMP. 

\begin{algorithm}[h!]\caption{Iterative scheme for computing an approximation of the solution of the nonlinear FEM-FCT problem. 
Let $\bu^{(0)}=\bu^n$ and let $\mathrm{tol} > 0$ and a damping factor $\rho \in (0,1]$ be given.}\label{alg:nonlin_femfct}
\begin{algorithmic}[1]
\STATE{Solve \eqref{eq:tcd_femfct_2}.}
\FOR{$m=0,1,\ldots$}
\STATE{Compute the algebraic fluxes $f_{ij}^{(m)}$ as in \eqref{eq:fct_anti_diff_fluxes}
with $\bu^{n+1}$ replaced by $\bu^{(m)}$ and the corresponding limiters $\alpha_{ij}^{(m)}$
by Zalesak's algorithm, such that $(\bbMl)^M\tilde\bu $ can be computed from
\eqref{eq:tcd_femfct_1}.}
\IF{$\left|\left(\bbMl + \theta \tau \mL \right)^M \bu^{(m)}-(\bbMl)^M \tilde \bu\right| \le \mathrm{tol}$}\label{alg:nonlin_femfct_stop}
\STATE{Set $\bu^{n+1} := \bu^{(m)}$, break.}
\ENDIF
\STATE{Solve \eqref{eq:tcd_femfct_0} with the right-hand side 
$(\bbMl)^M \tilde\bu$ and Dirichlet boundary conditions at $t^{n+1}$. Denote 
the solution $\hat\bu$ and set 
$\bu^{(m+1)} = \bu^{(m)} + \rho \big(\hat \bu -\bu^{(m)}\big)$ for the inner 
nodes and $\bu^{(m+1)} = \hat \bu$ for the boundary nodes.}
\ENDFOR
\end{algorithmic}
\end{algorithm}

\begin{theorem}[Global DMP for the iterates of
Algorithm~\ref{alg:nonlin_femfct}] \label{thm:tcd_dmp_fct}Denote
\begin{eqnarray}
u^{\mathrm{min}} &=& \min\left\{ u_1^n,\ldots,u_N^n, g_1^{n+1-\theta},\ldots, g_{N-M}^{n+1-\theta},
g_1^{n+1},\ldots, g_{N-M}^{n+1} \right\}\,,\label{eq:tcd_min} \\
u^{\mathrm{max}} &=& \max \left\{ u_1^n,\ldots,u_N^n, g_1^{n+1-\theta},\ldots, g_{N-M}^{n+1-\theta},
g_1^{n+1},\ldots, g_{N-M}^{n+1} \right\}\,.\label{eq:tcd_max}
\end{eqnarray}
Let $\theta=1$ or the CFL condition \eqref{eq:tcd_fct_cfl} be satisfied and
let $\bu^{(0)} = \bu^n$ in Algorithm~\ref{alg:nonlin_femfct}. Let all row sums 
of $(\mL)^M$ vanish and let the Zalesak algorithm be applied to compute the flux 
limiters. Then all iterates $\bu^{(m)}$, $m=0,1,\ldots$, satisfy
$u^{\mathrm{min}} \le u_i^{(m)} \le u^{\mathrm{max}}$, $i=1,\ldots,N$. 
\end{theorem}

\begin{proof} Note that the boundary values of $\overline\bu$ are 
$g_1^{n+1-\theta},\ldots, g_{N-M}^{n+1-\theta}$. The CFL condition implies that 
\eqref{eq:tcd_femfct_00} holds. Thus, if all row sums of $(\mL)^M$ vanish, then 
the matrices of equation \eqref{eq:tcd_femfct_2} satisfy the assumptions of
Corollary~\ref{cor:tcd_pospres_DMP_mon_mat}. Hence it follows that
$u^{\mathrm{min}} \le \overline u_i \le u^{\mathrm{max}}$, $i=1,\ldots,N$.
Since the Zalesak limiter is constructed in such a way that the solution of 
\eqref{eq:tcd_femfct_1} satisfies \eqref{eq:tcd_femfct_01}, one also has
\begin{equation*}
   u^{\mathrm{min}}\le\tilde u_i\le u^{\mathrm{max}}, \qquad i=1,\ldots,M.
\end{equation*}
As already mentioned above, the matrix on the left-hand side of 
\eqref{eq:tcd_femfct_0}, extended by the rows for the Dirichlet conditions, is 
an M-matrix. Since the row sums of $(\mL)^M$ vanish, the matrices in 
\eqref{eq:tcd_femfct_0} satisfy the assumptions of 
Corollary~\ref{cor:tcd_pospres_DMP_mon_mat} and hence
\begin{alignat*}{2}
   u^{\mathrm{min}}&\le
   \min\left\{\tilde{u}_1^{},\ldots,\tilde{u}_M^{},g_1^{n+1},\ldots, g_{N-M}^{n+1}\right\}
   &&\le \hat u_i\,,\\
   \hat u_i&\le 
   \max\left\{\tilde{u}_1^{},\ldots,\tilde{u}_M^{},g_1^{n+1},\ldots, g_{N-M}^{n+1}\right\}
   &&\le u^{\mathrm{max}}\,,
\end{alignat*}
for $i=1,\ldots,N$.
Finally, from $\bu^{(m+1)} = (1-\rho) \bu^{(m)} + \rho \hat \bu$ for the inner nodes, it can be inferred that 
$u^{\mathrm{min}} \le u_i^{(m+1)} \le u^{\mathrm{max}}$, $i=1,\ldots,N$.
\end{proof}

Note that the statement of Theorem~\ref{thm:tcd_dmp_fct} does not depend on the 
form of the algebraic fluxes. 

Now, one has to study under which conditions the row sums of $(\mL)^M$ vanish.
Since the row sums of $\mD$ are zero by construction, the row sums of $(\mL)^M$ 
vanish if and only if the row sums of the matrix $(\bbAN)^M$ vanish. In view of 
\eqref{eq:zero_row_diff_conv}, this is the case if and only if $\sigma=0$.
The assumption that $\sigma=0$ has to be expected since it appears already for 
the continuous version \eqref{eq:tcd_mp_both} of the maximum principle. 

\begin{remark}
The group finite element 
method is an alternative assembling routine of the convective term for $\bbP_1$ and $\bbQ_1$ finite elements that is based on matrix-vector 
multiplications instead on numerical quadrature. It introduces a consistency error, see \cite{BK17}
for a numerical analysis of the method, but it is usually considerably more efficient than 
the standard discretization, see \cite{JN12}. The $i$th row sum of the matrix for the convective 
term reads as follows \cite{BK17,JN12} for $i=1,\dots,M$ 
\[
\sum_{j=1}^N \left(\sum_{k=1}^d \left(\partial_k \phi_j,\phi_i\right) b_k(\bx_j)\right), 
\]
where $b_k(\bx_j)$ is the value of the $k$th component of $\bb$ at the node 
$\bx_j$. With the same argument as for the standard discretization, one finds 
that this row sum vanishes if $\bb$ is constant with respect to space, i.e., 
$b_k(\bx_j) = b_k$. But for general convection fields, the row sums do not 
vanish and hence, for the group finite element method, the satisfaction of the 
global DMP can be inferred from Theorem~\ref{thm:tcd_dmp_fct} only for very 
special (academic) convection fields. 
\hspace*{\fill}$\Box$\end{remark}

\begin{lemma}[Local DMP for both substeps of the FEM-FCT scheme] 
\label{lem:tcd_fct_local_dmp}Let the assumptions of Theorem~\ref{thm:tcd_dmp_fct} be satisfied,
then the substeps of the FEM-FCT scheme satisfy the following local DMPs:

i)\,\, The solution $\overline \bu$ of \eqref{eq:tcd_femfct_2} satisfies
\begin{equation}\label{eq:sub_dmp1}
   \min_{j\in S_i\cup\{i\}} u_j^n\le \overline u_i \le \max_{j\in S_i\cup\{i\}} u_j^n\,,
   \qquad 1\le i\le M\,.
\end{equation}

ii)\,\, The solution $\bu^{n+1}$ of \eqref{eq:tcd_femfct_0} satisfies
\begin{equation}\label{eq:sub_dmp2}
\min \left\{\overline u_i^{\mathrm{min}}, \min_{j\in S_i} u_j^{n+1}\right\}
\le u_i^{n+1} \le 
\max \left\{ \overline u_i^{\mathrm{max}},\max_{j\in S_i} u_j^{n+1}\right\}\,,
\qquad 1\le i\le M\,.
\end{equation}
\end{lemma}

\begin{proof} 
Consider any $i\in\{1,\dots,M\}$. We will prove only the upper bounds in 
\eqref{eq:sub_dmp1} and \eqref{eq:sub_dmp2} since the proofs of the lower 
bounds proceed along the same lines.

Denote by $u_i^{\mathrm{max}}$ the right-hand side of \eqref{eq:sub_dmp1} and
set $\mK=\bbMl - (1-\theta) \tau \mL$. Then $(\mK)^M\ge0$ due to
\eqref{eq:tcd_femfct_00}. Using the notation $\mK=(k_{ij})_{i,j=1}^N$ and
the row sum property of $(\mL)^M$, the solution of \eqref{eq:tcd_femfct_2} 
satisfies
\begin{equation*}
m_i\overline u_i=\sum_{j\in
S_i\cup\{i\}}k_{ij}\,\big(u^n_j-u_i^{\mathrm{max}}\big)+m_i\,u_i^{\mathrm{max}}
\le m_i\,u_i^{\mathrm{max}}\,,
\end{equation*}
which implies the upper bound in \eqref{eq:sub_dmp1}.

Now denote by $u_i^{\mathrm{max}}$ the right-hand side of \eqref{eq:sub_dmp2}.
Then the $i$th row of \eqref{eq:tcd_femfct_0} can be written in the form
\begin{equation}\label{eq:sub_dmp3}
   (m_i + \theta \tau l_{ii})\big(u_i^{n+1}-u_i^{\mathrm{max}}\big)=
   m_i\,\big(\tilde u_i-u_i^{\mathrm{max}}\big) 
   - \theta\tau \sum_{j\in S_i} l_{ij}\,
   \big(u_j^{n+1} - u_i^{\mathrm{max}}\big)\,.
\end{equation}
Since $l_{ij}\le0$ for $j\in S_i$ and the Zalesak limiter is constructed in 
such a way that $\tilde\bu$ satisfies \eqref{eq:tcd_femfct_01}, the right-hand
side of \eqref{eq:sub_dmp3} is non-positive. As discussed above, the matrix
$\mL_{\mathrm{I}}^{}$ is positive definite and hence $l_{ii}>0$. Thus,
\eqref{eq:sub_dmp3} implies the upper bound in \eqref{eq:sub_dmp2}.
\end{proof}

Summarizing the statements of Lemma~\ref{lem:tcd_fct_local_dmp}, one finds that
the solution of the nonlinear problem \eqref{eq:fct_low_order_limiter}
satisfies
\begin{equation*}
u_i^{n+1} \le \max \left\{\overline u_i, \max_{j\in S_i} \overline u_j, \max_{j\in S_i} u_j^{n+1} \right\} 
\le \max \left\{\max_{j\in S_i \cup \{i\}} u_j^n, \max_{j\in S_i} \overline u_j, \max_{j\in S_i} u_j^{n+1} \right\}. 
\end{equation*}
Consequently, one cannot conclude that a local DMP of the form formulated in
Lemma~\ref{lem:tcd_local_DMP} is satisfied for \eqref{eq:fct_low_order_limiter}
since the values of the intermediate solution $\overline \bu$ might determine 
the maximum on the right-hand side of the above estimate. Likewise, one cannot 
prove the LED property for the fully discrete problem, but only for both 
substeps individually. For instance, if $u_i^n$ is a local maximum, it cannot 
be excluded that $\overline u_j > \overline u_i$ for some $j\in S_i$. In this 
case, it is $\overline u_i^{\mathrm{max}}\neq \overline u_i$
and the LED property of the second substep does not provide information on the value of $u_i^{n+1}$. 
That the local DMP and the LED property, which are usually stated in the literature for the semi-discrete problem, cannot be transferred to the fully discrete 
problem is already mentioned in \cite[Ex.~4.56]{Loh20}.

In \cite{JK21}, the existence of a solution of \eqref{eq:tcd_femfct_2}--\eqref{eq:tcd_femfct_0} is proven for 
arbitrary time steps. The existence and uniqueness of a solution for sufficiently small time steps is shown 
in \cite{JKK21}.

We like to mention that there are in practice a couple of algorithmic issues and variations of the 
FEM-FCT scheme, like prelimiting. Since this topic is outside the scope of this survey, we refer to 
\cite{Kuz12a} or \cite[Chapters~7.5, 7.6]{KH15} for detailed presentations. Note
that the global DMP is still satisfied as long as the fluxes are modified before the 
application of the Zalesak limiter. 

Method \eqref{eq:tcd_femfct_2}--\eqref{eq:tcd_femfct_0} with the fluxes \eqref{eq:fct_anti_diff_fluxes}
and the bounds for the limiter \eqref{eq:tcd_uimin_max_fct} is a nonlinear scheme. As shown in Theorem~\ref{thm:tcd_dmp_fct}, an accurate solution of the nonlinear problem is not 
necessary in order to satisfy the global DMP, since it is satisfied for each iterate, but the 
accuracy of the numerical solution depends on how accurately the nonlinear problems are solved. 
However, in practice, it might be of advantage to use 
a linear version of a FEM-FCT scheme
for the sake of high efficiency, thereby accepting some loss of accuracy. Note that already the first FEM-FCT scheme proposed in \cite{PC86} is a linear scheme.
Linear FEM-FCT schemes are systematically derived in \cite{Kuz09}. 

The source of nonlinearity of a nonlinear FEM-FCT scheme is the definition 
\eqref{eq:fct_anti_diff_fluxes} of the algebraic fluxes. A linear FEM-FCT 
scheme can be also considered in the form \eqref{eq:fct_low_order_limiter}, 
however, the fluxes $f_{ij}$ are independent of the solution $\bu^{n+1}$ at 
the new time level. To define these fluxes, the values of $\bu^{n+1}$ in the 
formula \eqref{eq:fct_anti_diff_fluxes} are approximated by the solution of an
appropriate problem, e.g., the high-order method \eqref{eq:fct_gal} or the
low-order method \eqref{eq:fct_low_order}, or by extrapolating the solution 
$\overline\bu$ of the explicit scheme \eqref{eq:tcd_femfct_2} to the time level 
$t^{n+1}$. For $\theta=1/2$, such extrapolation was considered in \cite{JN12},
leading to the approximation of $\bu^{n+1}$ by $2\,\overline\bu-\bu^n$. Then
the fluxes are given by
\begin{equation}\label{eq:tcd_fctlin_fluxes}
f_{ij} = -m_{ij} \left(\hat u_j-\hat u_i\right) + d_{ij}\left(\overline u_j - \overline u_i\right)
\end{equation}
with $\hat\bu=2(\overline\bu-\bu^n)/\tau$. Note that
\begin{equation}\label{eq:low_order_explicit}
   (\bbMl)^M \hat\bu = -(\mL)^M\bu^n\,,
\end{equation}
i.e., $\hat\bu$ is an approximation of the time derivative of $u$ corresponding
to the low-order scheme \eqref{eq:fct_low_order} with $\theta=0$.
Independently of how the algebraic fluxes are defined, the
limiting procedure remains the same as for the nonlinear FEM-FCT scheme. In
particular, the bounds \eqref{eq:tcd_uimin_max_fct} for the limiter are defined
using the solution of \eqref{eq:tcd_femfct_2}. Thus, one obtains the following
analog of Theorem~\ref{thm:tcd_dmp_fct}.

\begin{corollary}[Global DMP for the linear FEM-FCT scheme with Zalesak limiter]
Let the algebraic fluxes be defined by \eqref{eq:fct_anti_diff_fluxes} with 
$\bu^{n+1}$ approximated using the solution of a problem depending on $\bu^n$
such that the fluxes are independent of $\bu^{n+1}$. Let $\theta=1$ or the CFL 
condition \eqref{eq:tcd_fct_cfl} be satisfied, and let the bounds of the 
limiter be defined by \eqref{eq:tcd_uimin_max_fct} with $\overline\bu$ from
\eqref{eq:tcd_femfct_2}. Let all row sums of $(\mL)^M$ vanish and let the 
Zalesak algorithm be applied for computing the flux limiters. Then the solution
of the linear scheme \eqref{eq:fct_low_order_limiter} satisfies
$u^{\mathrm{min}} \le u_i^{n+1} \le u^{\mathrm{max}}$, $i=1,\ldots,N$, where
$u^{\mathrm{min}}$ and $u^{\mathrm{max}}$ are defined by
\eqref{eq:tcd_min} and \eqref{eq:tcd_max}, respectively.
\end{corollary}

\begin{proof} The proof proceeds along the lines of the corresponding proof 
for the nonlinear FEM-FCT scheme. It was already noted that the concrete form 
of the fluxes does not play any role. 
\end{proof}

Another linearization strategy proposed in \cite{Kuz09} is a
predictor-corrector approach directly based on the basic FCT algorithm. In the 
first step, an intermediate solution $\overline\bu$ at time level $t^{n+1}$ is 
computed, e.g., by solving a problem of form \eqref{eq:fct_low_order}. In this 
step, one has to ensure that $\overline\bu$ satisfies a global DMP, which will 
give rise to a CFL condition, like \eqref{eq:tcd_fct_cfl}. The solution
$\overline\bu$ is used for computing the algebraic fluxes and the bounds 
\eqref{eq:tcd_uimin_max_fct} for the limiter. Then the flux limiters are
computed in the same way as for the nonlinear FEM-FCT method and a corrected 
solution is defined by 
\begin{equation}\label{eq:fct_predictor_corrector}
   (\bbMl)^M \bu^{n+1} = (\bbMl)^M\overline\bu + \tau \left(\sum_{j=1}^N
   \alpha_{ij}f_{ij}\right)_{i=1}^M
\end{equation}
and Dirichlet boundary conditions at $t^{n+1}$. The algebraic fluxes can be
defined by the formula \eqref{eq:fct_anti_diff_fluxes} with $\bu^{n+1}$
replaced by $\overline\bu$, as considered in \cite{Loh20}. In \cite{Kuz09},
the formula \eqref{eq:fct_anti_diff_fluxes} is considered with $\theta=1$,
leading to \eqref{eq:tcd_fctlin_fluxes}, where $\hat\bu$ is again an 
approximation of the discrete time derivative $(\bu^{n+1}-\bu^n)/\tau$ which
can be defined by \eqref{eq:low_order_explicit}, see \cite{Kuz09,Kuz12a} for
alternative proposals.

\begin{theorem}[Global DMP for the predictor-corrector FEM-FCT scheme with 
Zalesak limiter]
Let $\overline\bu$ be the solution of \eqref{eq:fct_low_order} and let the
bounds of the limiter be defined by \eqref{eq:tcd_uimin_max_fct} using this
$\overline\bu$. Let the algebraic fluxes be defined by an approximation of 
\eqref{eq:fct_anti_diff_fluxes} such that they are independent of $\bu^{n+1}$
and let the Zalesak algorithm be applied for computing the flux limiters. 
Let $\theta=1$ or the CFL condition \eqref{eq:tcd_fct_cfl} be satisfied, and 
let all row sums of $(\mL)^M$ vanish.
Then the corrected solution defined by \eqref{eq:fct_predictor_corrector}
satisfies $u^{\mathrm{min}} \le u_i^{n+1} \le u^{\mathrm{max}}$, 
$i=1,\ldots,N$, where 
$u^{\mathrm{min}}=\min\{u_1^n,\ldots,u_N^n, g_1^{n+1},\ldots, g_{N-M}^{n+1}\}$ and 
$u^{\mathrm{max}}=\max\{u_1^n,\ldots,u_N^n, g_1^{n+1},\ldots, g_{N-M}^{n+1}\}$.
\end{theorem}

\begin{proof} Since the matrices in \eqref{eq:fct_low_order} satisfy all the 
assumptions of Corollary~\ref{cor:tcd_pospres_DMP_mon_mat}, the solution 
$\overline\bu$ of \eqref{eq:fct_low_order} satisfies $u^{\mathrm{min}} \le
\overline u_i \le u^{\mathrm{max}}$, $i=1,\ldots,N$. The Zalesak limiter
is constructed in such a way that the corrected solution satisfies 
$\overline u_i^{\mathrm{min}} \le u^{n+1}_i\le \overline u_i^{\mathrm{max}}$,
$i=1,\dots,M$, which implies the theorem.
\end{proof}

For a comprehensive evaluation of the gain of efficiency and loss of accuracy 
in using a linear scheme for several academic problems, we refer to the 
numerical studies in \cite{JN12}. In that paper, one can find also comparisons 
with a linear upwind finite element method and an example where some
shortcomings of the FEM-FCT method are presented.

\section{Other types of finite elements}
\label{sec:other_fems}

This section discusses results concerning the DMP and corresponding methods for 
finite elements other than continuous piecewise linears. It turns out that the results are
often negative, at least in dimensions higher than one, and that there are only few 
methods for which a DMP can be proven. This situation justifies the
concentration on
the $\bbP_1$ finite element in the previous sections. 

\subsection{$\bbQ_1$ finite element}
\label{sec:other_fems_Q1}

Triangulations made of quadrilaterals in two dimensions or hexahedra in three dimensions are 
widely used for problems from fluid dynamics. The lowest order continuous finite element space on 
such triangulations is the space $\bbQ_1$ consisting of piecewise $d$-linear functions. Strictly 
speaking, one has to distinguish between two types of such spaces, namely mapped 
and unmapped $\bbQ_1$ finite elements. For the mapped version the local space is defined on 
a reference cell $\hat K$, e.g., $\hat K  = [-1,1]^d$.
Then, the finite element space on a physical mesh cell $K$ is given by the reference map from $\hat K$ to $K$. 
For the unmapped version the local functions  are defined directly on the physical
mesh cells. Both definitions coincide if the reference map is affine, i.e., if $K$ is a parallelepiped. 
If this is not the case, the image of a $d$-linear function defined on $\hat K$ will not 
be a $d$-linear function on $K$. 

Concerning $\bbQ_1$ finite elements, investigations of the DMP have been
concentrated so far on meshes whose 
cells are Cartesian products of intervals, sometimes called blocks in the literature. For the Poisson
equation in two dimensions, it had been observed already in \cite{CH84} that the DMP is violated if the 
aspect ratio, i.e., the ratio of the lengths of the longest edge and the
shortest edge of the cell, 
becomes too large. Based on the tensor-product representation of the basis functions by one-dimensional
basis functions, one can derive with a straightforward calculation a formula 
for the local entries $\fl_{ij}^{K}$ of the diffusion matrix, compare
\cite[Sec.~4.6]{Vej11}. If the corresponding nodes $\bx_i$ and $\bx_j$ share a
common edge $E_1$, then one finds in particular that
\begin{equation*}
\fl_{ij}^K = -\frac{|K|}{3^{d-1}} \left(\frac1{h_{E_1}^2} -\sum_{k=2}^d \frac1{2 h_{E_k}^2}\right),
\end{equation*}
where $E_1, \ldots, E_d$ are mutually orthogonal edges of $K$.
Thus, for $d=2$, one obtains a non-positive entry, which is condition \eqref{NN-1} for a 
matrix of non-negative type, if the aspect ratio is lower than or equal to
$\sqrt2$. For $d=3$, the mentioned entries are non-negative
only for cubes, namely $\fl_{ij}^K=0$, see also \cite{KKK07}. Considering the 
relaxed requirement that the diffusion matrix should be monotone, then numerical studies 
in \cite{KV10} reveal that the aspect ratios might be larger, at least on sufficiently 
fine grids, about $2.16$ for $d=2$ and $1.05$ for $d=3$. An extension of the analysis
to reaction-diffusion equations can be found in \cite[Sec.~4.6]{Vej11}.

\subsection{Higher order $H^1$-conforming finite elements}
\label{sec:other_fems_H1}

\begin{figure}[t!]
\centerline{\includegraphics[width=0.48\textwidth]{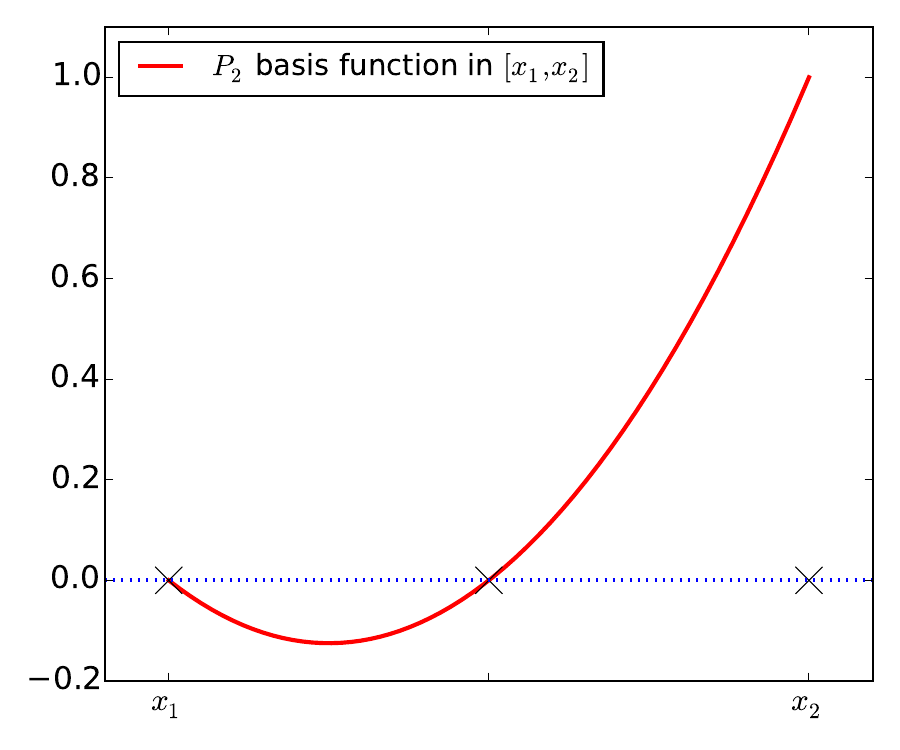}}
\caption{Basis function for $\bbP_2$ in the interval $[x_1,x_2]$. The degrees of freedom are indicated
with black crosses. The function is non-negative at the degrees of freedom, but takes negative
values in the interval.}\label{fig:p2_basis}
\end{figure}

Concerning the investigation of the DMP, a major difference between higher order $H^1$-conforming finite element
functions and $\bbP_1$ functions is as follows. Whereas local extrema are attained for $\bbP_1$ 
functions only in the degrees of freedom, i.e., geometrically at the vertices of the mesh 
cells, this is not the case for higher order finite element functions. As simple example, 
a one-dimensional standard $\bbP_2$ basis function is depicted in Figure~\ref{fig:p2_basis}, 
which takes its minimum between the locations of the degrees of freedom. 

A local DMP whose definition is restricted to the degrees of freedom has been studied for the Poisson 
equation in two dimensions already in \cite{Lor77,HM81}. It is shown in \cite{HM81} that such 
a DMP is satisfied for $\bbP_2$ finite elements only in special situations: on triangulations 
with equilateral triangles and on meshes consisting of squares in which the squares are divided by arbitrary
diagonals.  Note that these special triangulations impose severe restrictions 
on admissible forms of the domain. A more recent numerical study in
\cite{Vej10} shows that for $\bbP_2$ elements also triangulations with `nearly' 
equilateral triangles lead to a satisfaction of the DMP with respect to the 
degrees of freedom and that such a DMP is not satisfied for finite elements of 
degree three and higher. In addition, it is discussed in \cite{HM81} that even 
on special grids a DMP for the degrees of freedom is not valid for $\bbP_3$ 
finite elements. 

A proposal for extending an algebraically stabilized scheme to $\bbP_2$ finite 
elements such that the DMP for the nodal values is satisfied can be found 
in \cite{Kuz08}.

Already in \cite{HM81}, an example is given that the DMP for the degrees of freedom
does not imply a DMP for the finite element function. This issue might be crucial in coupled 
problems, when the $\bbP_2$ finite element solution is a coefficient in other equations and 
sufficiently accurate quadrature rules have to be utilized for assembling the finite element 
terms of the other equations. Usually, the nodes of such quadrature rules do not coincide 
with the geometric positions of the degrees of freedom of the $\bbP_2$ finite element function. 

In \cite{Lor77}, the special case of a triangulation consisting of squares that are divided
by diagonals which have all the same direction is studied. The proof of the DMP relies on a 
sufficient condition for the system matrix to be monotone. This condition is based, interestingly, 
on an additive decomposition of the system matrix, in its diagonal, a term that contains all 
positive off-diagonal entries, and a term that contains all negative off-diagonal entries. Then, it is 
assumed that the last term admits another additive decomposition that satisfies appropriate 
properties. A way that might be successful for deriving such a decomposition is provided. 
For details, it is referred
to \cite{Lor77,LZ20}.

At least for one-dimensional problems, some progress concerning the validation
of the DMP 
has been achieved, e.g., in \cite{VS07a,VS07b}. These results will not be discussed here 
since they do not generalize to higher dimensions. Another direction of research, 
inspired by \cite{Sch80}, consists in proving a so-called weak DMP, i.e.,
in showing that 
$\|u_h\|_{\infty,\Omega} \le C \|u_h\|_{\infty, \partial\Omega}$, where $C$ is independent of the 
mesh width, e.g., see \cite{LL21} for a recent contribution. Although mathematically 
certainly of interest, the weak DMP does not ensure the physical consistency of the 
numerical solution, even for $C=1$, e.g., if the solution is a concentration that should take 
values in $[0,1]$ in $\Omega$ and equals $1$ at some part of $\partial\Omega$, then negative 
values can still appear in a corresponding numerical solution.  A further direction of research consists in applying finite difference 
techniques for deriving a discrete problem for $\bbQ_2$ finite elements, e.g., see \cite{LZ20}
for a recent paper, which studies reaction-diffusion equations in two dimensions. Such methods 
possess the usual restriction of finite difference methods to simple domains. Results 
presented in \cite{LZ20} include the satisfaction of the global DMP on uniform meshes for the 
Poisson equation. If the uniform mesh is sufficiently fine, then the global DMP is also 
satisfied for the reaction-diffusion equation. 

\begin{remark}\emph{Bernstein finite element methods.} The presentation
of the FCT schemes in Section~\ref{sec:tcd_femfct} is 
completely algebraic, it did not exploit any special property of $\bbP_1$ finite elements. Only some 
general properties were used, like that the finite element basis forms a partition of unity and that the 
off-diagonal entries of $\bbMc$ are non-negative in order to obtain a well-defined lumped mass matrix. 
These two properties are also satisfied if the finite element basis consists of local Bernstein polynomials
of some degree. The finite element solution can be represented as 
a linear combination of these basis functions, which are non-negative, with so-called Bernstein coefficients. 
However, even in points that are degrees of freedom, the value of the solution usually does not coincide
with one of the Bernstein coefficients, in contrast to Lagrangian basis
functions. All statements proved in Section~\ref{sec:tcd_femfct} can be 
transferred to a FEM-FCT scheme with Bernstein polynomials, where everywhere 
the 
solution $u$ has to be replaced by the Bernstein coefficients, because they appear in the algebraic problems. Such a scheme for scalar transport equations 
is studied in \cite{LKSM17}.
\hspace*{\fill}$\Box$\end{remark}

\subsection{Non-conforming finite elements of Crouzeix--Raviart type}
\label{sec:other_fems_CR}

Consider a simplicial triangulation $\calT_h$ of $\Omega$. Then, the lowest order non-conforming
finite element space of Crouzeix--Raviart-type, proposed in \cite{CR73}, is defined by 
\begin{eqnarray*}
\bbP_1^{\mathrm{nc}} &=& \big\{v_h^{}\in L^2(\Omega)\ :\ v_h^{}|_K \in
\bbP_1(K)\,\,\,\forall\,\,K\in\calT_h,\ v_h^{}
\ \mbox{is continuous at the}\nonumber \\
&& \ \mbox{barycenters of all facets}\big\}.
\end{eqnarray*}
Functions from $\bbP_1^{\mathrm{nc}}$ are usually discontinuous across facets, so  
$\bbP_1^{\mathrm{nc}}\not\subset H^1(\Omega)$. The degrees of freedom are assigned to the 
facets. Consequently, the support of each nodal basis function consists of not more than two 
mesh cells. This property results in a small communication overhead in simulations on parallel
computers. Furthermore, the localized support leads to quite sparse matrices for many discretizations.

An upwind method for $\bbP_1^{\mathrm{nc}}$ was proposed in \cite{OU84}. To this end, a dual domain 
or lumping domain for each degree of freedom is considered. Since the degrees of freedom 
are assigned to the facets, the construction of the dual domain is much easier than for $\bbP_1$. 
For each degree of freedom, it is the polytope whose vertices are the vertices of the 
corresponding facet and the barycenter(s) of 
the mesh cell(s) where the facet belongs to. Integration by parts on the dual grid is applied 
to the convective term and then the fluxes across the facets of the dual mesh cells are 
approximated by an upwind technique.
The construction of the upwind fluxes leads on triangulations of acute type to a 
convection matrix that is of non-negative type. Also the diffusion matrix for $\bbP_1^{\mathrm{nc}}$
is of non-negative type on acute grids. Its restriction to the degrees of freedom that are 
not on the Dirichlet boundary is invertible, since the corresponding bilinear form is coercive 
with respect to a piecewise defined $H^1(\Omega)$ seminorm, which is a norm in
the subspace of $\bbP_1^{\mathrm{nc}}$ consisting of functions vanishing at
barycenters of facets contained in the Dirichlet boundary. Thus, from
\cite[Theorem~5.1]{Knobloch10} 
one can conclude the existence of a unique solution of the discrete problem and from 
Theorems~\ref{thm:algebraic-local-DMP} and~\ref{thm:algebraic-global-DMP} the satisfaction 
of the local and global DMP for the degrees of freedom, respectively, on acute triangulations. 

To the best of our knowledge, this upwind method is nowadays rarely used for the numerical 
solution of convection-diffusion-reaction equations. However, it gained some usefulness in the 
construction of multigrid methods for incompressible flow problems. For such problems, the pair 
$\bbP_1^{\mathrm{nc}}/\bbP_0$ satisfies a discrete inf-sup condition and applying the upwind technique
from \cite{OU84} leads to a convection-stabilized discretization of the incompressible Navier--Stokes 
equations. It was proposed in \cite{JM01} to utilize this discretization on lower levels of a 
multigrid method, leading to the so-called multiple discretization multilevel (MDML) method. 
A more recent comparison of solvers for the incompressible Navier--Stokes equations that
includes the MDML method can be found in \cite{ABJW18}.

The upwind technique from \cite{OU84} can be extended in a straightforward way to 
non-conforming rotated bilinear finite elements of lowest order for quadrilaterals and hexahedra
proposed in \cite{RT92}, see \cite{Tur94}. 

\subsection{Discontinuous Galerkin finite element methods}\label{sec:other_fems_DG}

Discontinuous Ga\-lerkin (DG)
methods were already proposed in \cite{RH73} for first order 
hyperbolic problems. They started to become also popular for discretizing second order 
elliptic equations in the 1990s. Meanwhile, a number of monographs are available, e.g., 
\cite{Riv08,DE12,DF15}. 

In DG methods, the finite element space consists of piecewise polynomials that
are completely discontinuous across facets of the mesh cells. Thus, a DG
finite element function is usually not contained in $H^1(\Omega)$. 

For DG methods, the notion of `satisfying a DMP' has to be revisited. 
In several papers on time-dependent 
transport and convection-diffusion equations, e.g., \cite{ZS10,ZZS13},  the
fact that DG allows to use the cell averages in natural way has been used to 
restrict the DMP to these quantities, and then the following criterion has been
proposed: let the cell-wise averages of the DG solution $u^n$ 
at time instant $t^n$ be in $[u^{\mathrm{min}},u^{\mathrm{max}}]$, then the 
DG method satisfies a DMP if the averages of $u^{n+1}$ at $t^{n+1}$ are also contained in this 
interval. For a detailed discussion of such methods, it is referred to the respective literature, 
e.g., \cite{Shu16}.  An alternative approach, more algebraic and based in 
the concept of invariant sets and domains for hyperbolic problems, is followed 
in \cite{GPT19,Hajduk21,Pazner21}.

 In here, we will detail an approach proposed for the 
convection-diffusion equation in \cite{BH15}. We start by defining the first order discontinuous space on a simplicial grid,
that is\footnote{Strictly speaking, the functions of $\bbP_1^{\mathrm{disc}}$
are well-defined only on the interiors of the mesh cells, since the limits to
the same point at the boundaries of mesh cells, approached from different mesh
cells, are usually different. To simplify the presentation, we will
nevertheless speak of values on facets or at vertices and mean always the limit
from the corresponding mesh cell.
}
\[
\bbP_1^{\mathrm{disc}} = \left\{ v_h \in L^2(\Omega) \ : \ \left. v_h \right|_K
\in \bbP_1(K) \,\,\,\forall\,\,K\in\calT_h\right\}.
\]
This space is equipped with the basis $\{\phi_i^K\}$, where for a mesh cell $K$ 
and a node $i$ such that $\bx_i$ is a vertex of $K$, the function $\phi_i^K$ is
defined as follows: $\phi_i^K$ is linear in $K$, $\phi_i^K|_K^{}(\bx_i) = 1$, 
$\phi_i^K|_K^{} = 0$ at all other vertices of $K$, and $\phi_i^K$ vanishes outside of $K$. 
The restriction of $v_h \in \bbP_1^{\mathrm{disc}}$ to a mesh cell $K$ is denoted by $v_h^K$. 

\begin{figure}[t!]
\centerline{\includegraphics[width=0.6\textwidth]{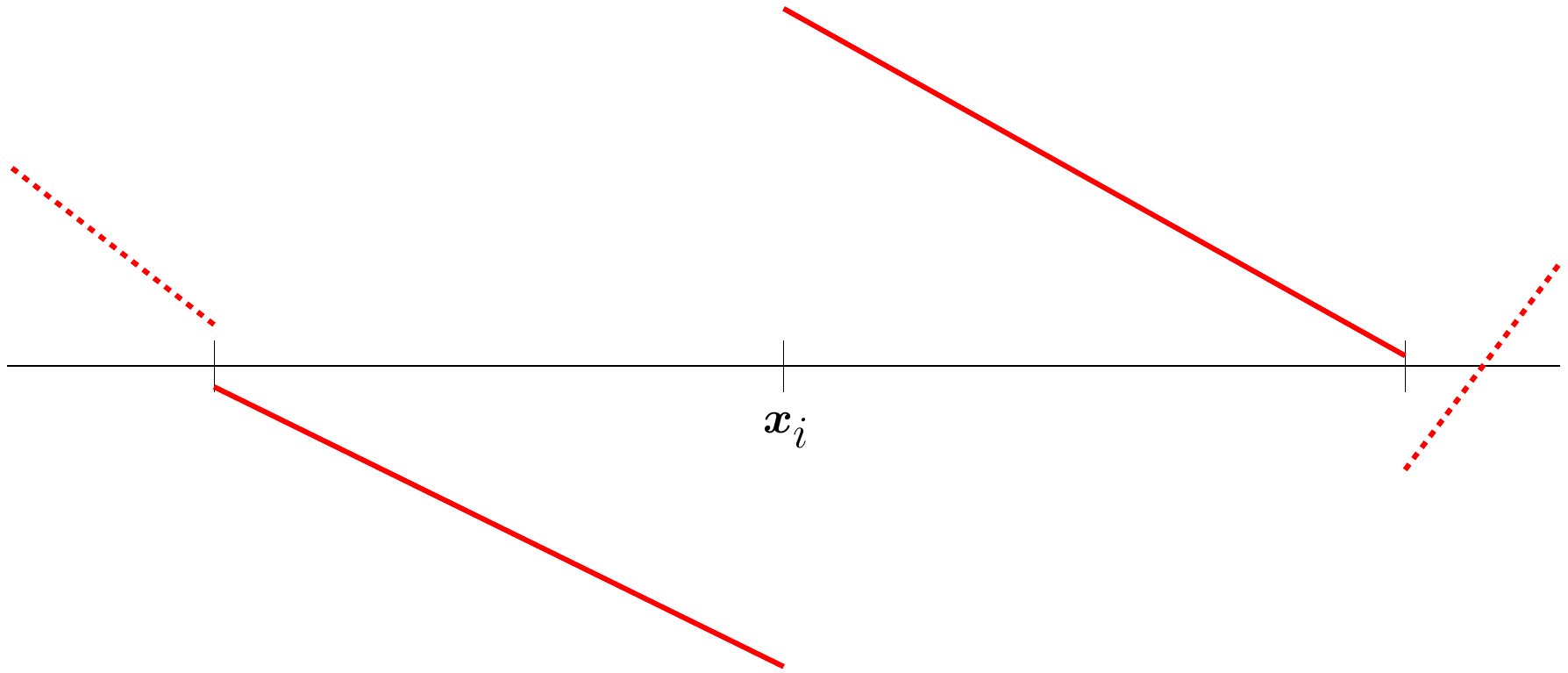}}
\caption{$\bbP_1^{\mathrm{disc}}$ function (in red) with local minimum and local maximum at $\bx_i$.}\label{fig:dg_loc_extr}
\end{figure}

The first observation is that even the notion of a local extremum is not clear
for functions from $\bbP_1^{\mathrm{disc}}$, compare
Fig.~\ref{fig:dg_loc_extr}, where the values at $\bx_i$ are both a strict
local minimum and a strict local maximum. To this end, the following definition
was introduced in \cite{BH15}.

\begin{definition}[Local discrete extremum for $\bbP_1^{\mathrm{disc}}$] 
\label{def:dg_loc_extr}The function $u_h \in \bbP_1^{\mathrm{disc}}$ has a 
local discrete minimum (resp. maximum) at the vertex $\bx_i$ in $K$ if $u_h^K(\bx_i) \le u_h(\bx)$
(resp. $u_h^K(\bx_i) \ge u_h(\bx)$) for all $\bx \in \omega_i$. 
\end{definition}

Then, a definition of a DMP for methods using $\bbP_1^{\mathrm{disc}}$ is given in \cite{BH15}, 
which is inspired by Definition~\ref{def:DMP-criterion-non-strict} for nonlinear forms with $\bbP_1$ functions. 

\begin{definition}[DMP for $\bbP_1^{\mathrm{disc}}$] 
\label{def:dmp_dg}Let $a_h \ : \ \bbP_1^{\mathrm{disc}} \times 
\bbP_1^{\mathrm{disc}} \to \mathbb R$ be a bilinear form.
This bilinear form is said to possess 
the DMP property if for all $u_h \in  \bbP_1^{\mathrm{disc}}$ and for all interior vertices 
$\bx_i$ where $u_h$ is locally minimal (resp. maximal) at $\bx_i$ in $K$, there
exist constants
$\alpha_F^{}> 0$ and $\zeta_K > 0$ such that 
\begin{equation}\label{eq:dg_dmp}
a_h \left(u_h,\phi_i^K\right) \le - \sum_{F\in\calF_i^{}\cap\calF_K^{}} \frac{\alpha_F}{h_F}
\int_F \left|\llbracket u_h^{}\rrbracket_F \right| \ d\boldsymbol s - \frac{\zeta_K}{h_K} 
\int_K \left| \nabla u_h^K\right| \ d\bx, 
\end{equation}
(resp. $ a_h \left(u_h,\phi_i^K\right) \ge
\sum_{F\in \calF_i^{}\cap\calF_K^{}} \frac{\alpha_F}{h_F}
\int_F \left|\llbracket u_h^{}\rrbracket_F \right| \ d\boldsymbol s + \frac{\zeta_K}{h_K} 
\int_K \left| \nabla u_h^K\right| \ d\bx$).
\end{definition}

Next, the consistency of the preceding definitions will be shown. 

\begin{lemma}[Consequences of the satisfaction of the DMP] Let $a_h \ : \ \bbP_1^{\mathrm{disc}} \times 
\bbP_1^{\mathrm{disc}} \to \mathbb R$ be a bilinear form that satisfies the DMP
property from 
Definition~\ref{def:dmp_dg} and consider the problem to find $u_h \in \bbP_1^{\mathrm{disc}}$ such 
that $a_h(u_h,v_h) = (f,v_h)$ for all $v_h \in \bbP_1^{\mathrm{disc}}$. 
\begin{enumerate}[leftmargin=*]
\item[i)] If $ f \ge 0$ (resp. $f\le 0$), then $u_h$ does not possess a strict local discrete minimum 
(resp. maximum), see Definition~\ref{def:dg_loc_extr},  at any interior point. 
\item[ii)] If $ f \ge 0$ (resp. $f\le 0$), then $u_h$ attains its global minimum (resp. maximum)
at the boundary $\partial \Omega$. 
\end{enumerate}
\end{lemma}

\begin{proof}
i) Assume that $u_h$ has a strict local discrete minimum at the interior node $\bx_i$
in the mesh cell $K$. Since $a_h(\cdot,\cdot)$ satisfies the DMP property, it follows from \eqref{eq:dg_dmp} that 
$a_h(u_h,\phi_i^K) \le 0$. On the other hand, one has $(f,\phi_i^K) \ge 0$ and then
$a_h(u_h,\phi_i^K) = 0$ holds. From \eqref{eq:dg_dmp}, one infers that then 
$\nabla u_h^K = \boldsymbol 0$ and hence $u_h^K$ is constant so that
the minimum is not strict.

ii) If $u_h^K(\bx_i)$ is a global minimum for some mesh cell $K$ and some
interior node $\bx_i\in K$, then it is also a local minimum and from the proof
of i), one gets that $u_h^K$ is constant. Moreover, it follows from the DMP
property that $\llbracket u_h^{}\rrbracket_F =0$ for all 
$F\in\calF_i^{}\cap\calF_K^{}$. Let $K'\subset\omega_i$ be a mesh cell that 
shares a common facet $F$ with $K$. As the jump $\llbracket u_h^{}\rrbracket_F$
vanishes, then $u_h^K(\bx) =u_h^{K'}(\bx)$  for all $\bx\in F$, and in 
particular $u_h^K(\bx_i) =u_h^{K'}(\bx_i)$. Thus, $u_h^{K'}(\bx_i)$ also is a
global minimum and it follows that $u_h^{K'}$ is constant. By induction, one 
finds that $u_h |_{\omega_i} = u_h^K(\bx_i)$ is a constant. Then, again 
by induction, it follows that $u_h$ is 
constant in $\Omega$ and in particular that $u_h|_{\partial\Omega} = u_h^K(\bx_i)$. Hence, the 
global minimum is attained at the boundary of $\Omega$. 
\end{proof}

One type of equations studied in \cite{BH15} is a steady-state convection-diffusion equation 
with conservative form of the convective term and solenoidal convection field. For the DG 
discretization of the diffusive term, the standard incomplete interior penalty (IIP) method is
used. This choice is motivated by the analysis of one-dimensional diffusion
problems that are discretized with DG methods, see \cite{HM13}. The convective term is 
integrated by parts and then an upwind discretization at interior facets is utilized. 
In addition, and this is the major algorithmic proposal of \cite{BH15}, a nonlinear, locally defined
artificial diffusion term built with the help of a shock detector
is added. 
For a one-dimensional problem, the DMP, according to Definition~\ref{def:dmp_dg}, is proven.
There are no analytic results for multiple dimensions. The main obstacle for such 
results is that a DMP is not available 
already for the usual interior penalty discretizations of the diffusion term. 
In the numerical studies presented in \cite{BH15}, small violations of the 
DMP can be observed for a simulation performed on an acute mesh in two dimensions.

A method that addresses the above mentioned issue of the DMP for interior penalty discretizations
of the diffusive term is proposed in \cite{BBH17}. This method augments the symmetric 
interior penalty method with a nonlinear discrete diffusion operator
related to the AFC/FCT schemes described in previous sections.
Then, it is shown in \cite{BBH17} that the proposed scheme 
for the steady-state convection-diffusion problem satisfies
a local DMP if the right-hand side of the  equation vanishes identically. 
This statement holds for arbitrary admissible grids and $\bbP_1^{\mathrm{disc}}$ finite elements
on simplices and discontinuous piecewise $d$-linear elements on quadrilaterals or hexahedra.
For the time-dependent case, a semi-discrete problem in space is considered and it is 
shown that the discrete 
scheme is LED, again in case that the right-hand side of the problem 
is identically zero.

High-order DG schemes based on algebraic flux correction were recently 
developed in \cite{Hajduk21,Pazner21} for hyperbolic conservation laws. 
While in \cite{Hajduk21} monolithic convex limiting with subcell flux
limiters is used, in \cite{Pazner21} an FCT-type predictor-corrector algorithm
is advocated.
The bound preserving DG scheme of \cite{Hajduk21} employs Bernstein polynomials to
facilitate the use of very high order spatial approximations. The limiting 
strategy of \cite{Pazner21} is tailor-made for Legendre-Gauss-Lobatto
DG bases, and makes
use of a novel sparse low-order invariant domain preserving method
whose stencil does not grow with the polynomial degree of the corresponding 
high-order method. The invariant domain preservation is proved under a CFL
condition.

\section{Brief comments on hyperbolic conservation laws}
\label{sec:hyperbolic_problems}

The aim of this section is to discuss briefly results on the satisfaction of 
the DMP for transport equations and nonlinear hyperbolic conservation laws. Presenting 
in detail the amount and variety of works devoted to hyperbolic  problems requires a review on its
own and it is clearly outside the scope of the present survey.   In particular, 
in this section we will only focus  on continuous finite element methods,
since for discontinuous Galerkin approaches there exist several well documented reviews (e.g., \cite{ZS13,Shu16}).
In addition, in recent years there has been an increasing interest in seeking suitable conforming approximations
for hyperbolic problems, since conforming approximations do not have a built-in stability, and hence
the challenge of finding structure-preserving stabilizing terms is different from the discontinuous counterparts.

The model problem considered in this section is  the extreme case $\varepsilon=0$, this is, the transport equation,
or, more generally, conservation equations of the form
\begin{equation}\label{transport-transient}
\partial_t^{} u + \textrm{div} \boldsymbol{f}( u) = 0\qquad \textrm{in}\; \Omega\,,
\end{equation} 
where $\boldsymbol{f}(u)$ is the flux function, provided with appropriate (inlet) boundary and initial conditions. 
If $\boldsymbol{f}(u)=\boldsymbol{b}u$, then \eqref{transport-transient}  reduces to the linear transport equation.

\begin{remark}
It is worth mentioning that the  case $\varepsilon=0$ allows to propose methods that respect the DMP on general
meshes in a more natural way. In fact, the added viscosity methods only need to deal with compensating for the wrong signed terms in
the convection matrix, and not with the possibly positive terms in the diffusion matrix, which are of a different order in terms of the
mesh size.  For example, in \cite{Bur15} an appropriate combination of upwinding and FCT-related techniques is
used to propose a nonlinear stabilized scheme that preserves the DMP for the linear transport equation. In addition, the time discretization is based 
on an explicit method, so the overhead of using a nonlinear discretization is minimal.
On the other hand, it is important to mention that the discrete maximum principle is not, in general, enough to
prove the convergence of a numerical scheme to the entropy solution of \eqref{transport-transient}, as it has been mentioned in, e.g.,
\cite{GP17}, where the authors show that, in order to converge to the entropy solution, the scheme needs also to control the maximum
wave speed. More precisely, in Lemma~4.6 in that reference, it is shown that the FCT algorithm, equipped with a limiter related to
the Zalesak one, might not converge for certain nonlinear fluxes, which is then
confirmed in the numerical experiments for Burgers' equation.
\hspace*{\fill}$\Box$\end{remark}

We start by mentioning that most of the references quoted in Section~\ref{sec:tcd_femfct} were, in fact, works developed for the transport, or Euler, equations. 
So,  this section will be devoted to  describing some of the more recent
developments of DMP-preserving schemes for this problem.  In 
\cite{Kuz20}, using the framework of algebraic flux correction and invariant 
domain preserving schemes, a monolithic approach to convex limiting is 
introduced for hyperbolic conservation laws. The convex limiting is thoroughly 
discussed for both scalar conservation laws (including the transport equation)
and hyperbolic systems.
In the context of the enriched finite element method
(proposed originally in \cite{BBHL04}), a FCT scheme for the transport equation is proposed in \cite{KHR20}  where the DMP
is proven (under appropriate CFL conditions)  for both the continuous and discontinuous parts of the solution.

In \cite{GN14} a first order added diffusion/viscosity method with an explicit time discretization is proposed for \eqref{transport-transient}.
The DMP for the resulting scheme is proven under a CFL condition.
On uniform meshes, the bilinear form of the first order diffusion used
in \cite{GN14} corresponds to the matrix $M_{\rm C}-M_{\rm L}$ used in
\cite{LMPV87}.    Later, in \cite{GPY17} the authors show
that it is impossible to propose an explicit continuous finite element method that is stabilized
with artificial viscosity and satisfies the DMP if the time derivative
is approximated using the consistent mass matrix.
In the  paper \cite{GNPY14} the authors propose a different technique: first,  a higher order added viscosity (defined as the 
minimum between the first order viscosity and the entropy residual) is added. The DMP cannot be proven for the resulting scheme, so they use a
technique related to the FCT method (linked to the graph-Laplacian writing of the added viscosity), supplied with  flux limiters related to those 
described in  Section~\ref{sec:tcd_femfct} 
(based on the Zalesak algorithm), and an approximation of the inverse of the consistent  mass matrix to correct the scheme.   
The combination of these techniques allows for the proof of the DMP.
Later, in \cite{GP17} a method, again related to the FCT family, is  proposed,  equipped with three different limiters, namely the Zalesak limiter,
the smoothness-based indicator, and a greedy viscosity algorithm. In addition, the satisfaction of the DMP and the convergence
to the entropy solution are shown. Some comparisons
in terms of robustness and reliability are also carried out in \cite{GP17}. 
Another work devoted to stabilization by the nonlinear diffusion operator (also 
referred to as graph Laplacian in some papers) is the work
 \cite{BB17}, where a regularization of the definition of the limiters is proposed in order to obtain twice differentiable limiters and to make
the discretization amenable to the use of Newton's method to solve the algebraic system.

In the context of the Burgers equation, in \cite{Bur07} numerical viscosity is added  to satisfy the DMP and prove convergence to the entropy solution of the hyperbolic
equation. In one space dimension the method consists of adding a numerical diffusion of the form $(\varepsilon(u_h^{})\partial_x^{}u_h^{},\partial_x^{}v_h^{})$ where $\varepsilon(u_h^{})$
is  designed to satisfy several hypotheses. These conditions imply the Lipschitz continuity of the stabilization and the fact that
the problem satisfies the strong DMP property (similar to those in Section~\ref{Sec:DMP_nonlinear_disc}). Under these assumptions, the finite element method is proven to 
converge to the entropy solution of Burgers' equation. Later, in \cite{BH14}, essentially the same assumptions are imposed on the coefficient
of the added diffusion, with the difference that in this case the diffusion is of the form of a local projection stabilization method. The method is
proven to be LED and  to converge to the entropy solution.

We next comment on the possibility of using both linear and nonlinear stabilizing terms in conservation laws. In fact, as it was mentioned in
previous sections, it has been observed in several works that the use of a nonlinear stabilization (e.g., FCT) alone does not suffice to build a convergent method. For example,
in \cite{GP17} it is shown that using nonlinear stabilization alone leads, in certain cases, to failure in convergence of the scheme. So, the authors take a different approach
by first adding an entropy viscosity to a method by using the consistent mass matrix, thus violating the DMP, and then applying a FCT technique as a post-processing to 
produce a DMP-preserving approximate solution.
In addition, in the work \cite{EG13} a combined use of linear (edge-based) stabilization and a nonlinear entropy viscosity is advocated. It is shown in that reference
that the addition of linear stabilization, if not weighted properly, can actually hinder the satisfaction of the DMP and increase the entropy violations, and even in some extreme cases, make a 
convergent method converge to the wrong weak solution. So, a nonlinear weight is introduced to balance the influence of 
the stabilizing terms and secure convergence to the entropy solution. We should, nevertheless, mention that even if the entropy viscosity method is claimed to satisfy a \textit{weakened} maximum
principle, there is no proof of DMP-satisfaction (or weakened DMP) available, although the authors show numerical evidence supporting the claim that the weighted method does satisfy
a weakened DMP. 

We finish this short section by mentioning two relatively recent works where DMP-preserving methods are introduced and that use LPS-related methods as linear stabilization.
In \cite{KBS17} a linear stabilizing term is first introduced. This term penalizes the fluctuations between the discrete solution and its local average (thus inspired by the LPS idea, but
departing from the classical LPS approaches). This method 
preserves the DMP but provides inaccurate results,
so the target function, that is, the function with respect to which the fluctuation is 
computed, is modified by adding to it an approximation of its gradient. This approximation is then limited using limiters that guarantee the LED property and linearity preservation (on general
meshes) of the resulting scheme. The authors claim that the linearity preserving limiter introduced in \cite[Section~7]{KBS17} can also be applied in different contexts, e.g.,
the AFC and FCT schemes.
The resulting method is tested in steady-state and time-dependent schemes showing that the combination of the gradient approximation as high order stabilization with the
LED limiter localizes the stabilization enough as to reduce the oscillations around the shocks without smearing the profiles in excess.
 Finally, in \cite{MSK18}  the authors present a nonlinear stabilization through discrete artificial diffusion supplemented by a monotone local projection operator based on limiting at the semi-discrete level.
The resulting method respects the DMP and is linearity preserving. 
The impact of the local projection operator is studied in the numerical experiments where it is shown that its addition (that acts as a high order background dissipation) helps to reduce the
terracing (and even eliminates it in some cases).

\section{Summary}
\label{sec:summary}

For convection-dominated convection-diffusion problems it is a challenging task to construct discretizations 
that  at the same time satisfy the DMP and compute accurate solutions. Enormous efforts have been spent 
since the 1980s in the development of schemes that enrich traditional stabilized finite element methods 
with extra terms to reduce the size of spurious oscillations, leading to the class of SOLD methods. 
However, this development turned out to be only little successful with respect to designing methods for 
which the DMP can be proven rigorously, since only the Mizukami--Hughes method satisfies this property. 
In the 2000s, a different class of methods was started to be developed, namely algebraically stabilized 
finite element methods. In that decade, FEM-FCT schemes for the time-dependent problem were proposed 
and at the end of that decade, the first AFC method for the steady-state problem. Then, in recent years, 
the analysis for AFC methods have been developed and further methods for the steady-state problem, 
like  modifications and extensions of algebraic stabilizations, have been developed. 
For all of these schemes, the DMP can be proven, sometimes under conditions on
the data or the mesh. 
In summary, there are meanwhile several, but still surprisingly few, finite element
methods available that satisfy the DMP and compute simultaneously quite accurate results. 

For the steady-state problem, all DMP-respecting finite element schemes with accurate solutions are nonlinear. 
It can be seen in the numerical example from Section~\ref{sec:num_exam} that, on the one hand, there are 
differences concerning the accuracy of the computed solutions, but on the other hand, the differences are 
not large. For the practical use of these methods, also aspects like the efficiency for solving 
the nonlinear problems and the efforts for implementing the methods in three dimensions are important.
Concerning the first issue, whose investigation is outside the scope of this survey, a comprehensive 
comparison of two algebraically stabilized schemes can be found in \cite{JJ19}. Simulations of three-dimensional 
problems with various algebraic stabilizations can be found in \cite{BJKR18,JJ19}.
Note the many algebraic stabilizations do work only with the matrices and vectors such that their implementation 
can be carried out independently of the dimension of the problem. 

There is a similar situation for the time-dependent problem: algebraically stabilized schemes are the 
currently best available finite element methods that satisfy the global DMP and compute accurate solutions. 
Here, also a linear variant is available which showed in several applications a good balance of accuracy
and efficiency. 

\bigskip

\textbf{Acknowledgment.} This work was initiated at a stay of the three authors
at the Mathematisches Forschungsinstitut  Oberwolfach (MFO) within the {\it
Research in pairs} programme, grant No.~1937p. 

\label{sec:references}

\bibliographystyle{siamplain}

\begin{thebibliography}{100}

\bibitem{ABJW18}
{\sc N.~Ahmed, C.~Bartsch, V.~John, and U.~Wilbrandt}, {\em An assessment of
  some solvers for saddle point problems emerging from the incompressible
  {N}avier-{S}tokes equations}, Comput. Methods Appl. Mech. Engrg., 331 (2018),
  pp.~492--513, \url{https://doi.org/10.1016/j.cma.2017.12.004}.

\bibitem{ABR17}
{\sc A.~Allendes, G.~R. Barrenechea, and R.~Rankin}, {\em Fully computable
  error estimation of a nonlinear, positivity-preserving discretization of the
  convection-diffusion-reaction equation}, SIAM J. Sci. Comput., 39 (2017),
  pp.~A1903--A1927, \url{https://doi.org/10.1137/16M1092763}.

\bibitem{BT81}
{\sc K.~Baba and M.~Tabata}, {\em On a conservative upwind finite element
  scheme for convective diffusion equations}, RAIRO Anal. Num\'er., 15 (1981),
  pp.~3--25, \url{https://doi.org/10.1051/m2an/1981150100031}.

\bibitem{BB17}
{\sc S.~Badia and J.~Bonilla}, {\em Monotonicity-preserving finite element
  schemes based on differentiable nonlinear stabilization}, Comput. Methods
  Appl. Mech. Engrg., 313 (2017), pp.~133--158,
  \url{https://doi.org/10.1016/j.cma.2016.09.035}.

\bibitem{BBH17}
{\sc S.~Badia, J.~Bonilla, and A.~Hierro}, {\em Differentiable
  monotonicity-preserving schemes for discontinuous {G}alerkin methods on
  arbitrary meshes}, Comput. Methods Appl. Mech. Engrg., 320 (2017),
  pp.~582--605, \url{https://doi.org/10.1016/j.cma.2017.03.032}.

\bibitem{BH14}
{\sc S.~Badia and A.~Hierro}, {\em On monotonicity-preserving stabilized finite
  element approximations of transport problems}, SIAM J. Sci. Comput., 36
  (2014), pp.~A2673--A2697, \url{https://doi.org/10.1137/130927206}.

\bibitem{BH15}
{\sc S.~Badia and A.~Hierro}, {\em On discrete maximum principles for
  discontinuous {G}alerkin methods}, Comput. Methods Appl. Mech. Engrg., 286
  (2015), pp.~107--122, \url{https://doi.org/10.1016/j.cma.2014.12.006}.

\bibitem{BCC98}
{\sc R.~E. Bank, W.~M. Coughran~jr., and L.~C. Cowsar}, {\em The finite volume
  {Scharfetter}-{Gummel} method for steady convection diffusion equations},
  Comput. Vis. Sci., 1 (1998), pp.~123--136,
  \url{https://doi.org/10.1007/s007910050012}.

\bibitem{BBK17a}
{\sc G.~R. Barrenechea, E.~Burman, and F.~Karakatsani}, {\em Blending low-order
  stabilised finite element methods: a positivity-preserving local projection
  method for the convection-diffusion equation}, Comput. Methods Appl. Mech.
  Engrg., 317 (2017), pp.~1169--1193,
  \url{https://doi.org/10.1016/j.cma.2017.01.016}.

\bibitem{BBK17b}
{\sc G.~R. Barrenechea, E.~Burman, and F.~Karakatsani}, {\em Edge-based
  nonlinear diffusion for finite element approximations of convection-diffusion
  equations and its relation to algebraic flux-correction schemes}, Numer.
  Math., 135 (2017), pp.~521--545,
  \url{https://doi.org/10.1007/s00211-016-0808-z}.

\bibitem{BJK15}
{\sc G.~R. Barrenechea, V.~John, and P.~Knobloch}, {\em Some analytical results
  for an algebraic flux correction scheme for a steady convection-diffusion
  equation in one dimension}, IMA J. Numer. Anal., 35 (2015), pp.~1729--1756,
  \url{https://doi.org/10.1093/imanum/dru041}.

\bibitem{BJK16}
{\sc G.~R. Barrenechea, V.~John, and P.~Knobloch}, {\em Analysis of algebraic
  flux correction schemes}, SIAM J. Numer. Anal., 54 (2016), pp.~2427--2451,
  \url{https://doi.org/10.1137/15M1018216}.

\bibitem{BJK17}
{\sc G.~R. Barrenechea, V.~John, and P.~Knobloch}, {\em An algebraic flux
  correction scheme satisfying the discrete maximum principle and linearity
  preservation on general meshes}, Math. Models Methods Appl. Sci., 27 (2017),
  pp.~525--548, \url{https://doi.org/10.1142/S0218202517500087}.

\bibitem{BJKR18}
{\sc G.~R. Barrenechea, V.~John, P.~Knobloch, and R.~Rankin}, {\em A unified
  analysis of algebraic flux correction schemes for convection-diffusion
  equations}, SeMA J., 75 (2018), pp.~655--685,
  \url{https://doi.org/10.1007/s40324-018-0160-6}.

\bibitem{BK17}
{\sc G.~R. Barrenechea and P.~Knobloch}, {\em Analysis of a group finite
  element formulation}, Appl. Numer. Math., 118 (2017), pp.~238--248,
  \url{https://doi.org/10.1016/j.apnum.2017.03.008}.

\bibitem{BBH11}
{\sc R.~Becker, E.~Burman, and P.~Hansbo}, {\em A finite element time
  relaxation method}, C. R. Math. Acad. Sci. Paris, 349 (2011), pp.~353--356,
  \url{https://doi.org/10.1016/j.crma.2010.12.010}.

\bibitem{BBHL04}
{\sc R.~Becker, E.~Burman, P.~Hansbo, and M.~G. Larson}, {\em A reduced
  $\mathbb{P}^1$-discontinuous {G}alerkin method.}, Chalmers Finite Element
  Center Preprint 2003-13, Chalmers University of Technology, G\"oteborg,
  Sweden, 2003.

\bibitem{BB73}
{\sc J.~P. {Boris} and D.~L. {Book}}, {\em {Flux-corrected transport. I:
  SHASTA, a fluid transport algorithm that works.}}, {J. Comput. Phys.}, 11
  (1973), pp.~38--69, \url{https://doi.org/10.1016/0021-9991(73)90147-2}.

\bibitem{BH62}
{\sc J.~H. Bramble and B.~E. Hubbard}, {\em On the formulation of finite
  difference analogues of the {D}irichlet problem for {P}oisson's equation},
  Numer. Math., 4 (1962), pp.~313--327,
  \url{https://doi.org/10.1007/BF01386325}.

\bibitem{BH64}
{\sc J.~H. Bramble and B.~E. Hubbard}, {\em New monotone type approximations
  for elliptic problems}, Math. Comp., 18 (1964), pp.~349--367,
  \url{https://doi.org/10.1090/S0025-5718-1964-0165702-X}.

\bibitem{BKK20}
{\sc J.~Brandts, S.~Korotov, and M.~K{\v{r}}{\'{\i}}{\v{z}}ek}, {\em Simplicial
  Partitions with Applications to the Finite Element Method}, Springer-Verlag,
  Cham, 2020, \url{https://doi.org/10.1007/978-3-030-55677-8}.

\bibitem{BKKS09}
{\sc J.~Brandts, S.~Korotov, M.~K{\v{r}}{\'{\i}}{\v{z}}ek, and J.~{\v{S}}olc},
  {\em On nonobtuse simplicial partitions}, SIAM Rev., 51 (2009), pp.~317--335,
  \url{https://doi.org/10.1137/060669073}.

\bibitem{BKK08}
{\sc J.~H. Brandts, S.~Korotov, and M.~K{\v{r}}{\'{\i}}{\v{z}}ek}, {\em The
  discrete maximum principle for linear simplicial finite element
  approximations of a reaction-diffusion problem}, Linear Algebra Appl., 429
  (2008), pp.~2344--2357, \url{https://doi.org/10.1016/j.laa.2008.06.011}.

\bibitem{Bur07}
{\sc E.~Burman}, {\em On nonlinear artificial viscosity, discrete maximum
  principle and hyperbolic conservation laws}, BIT, 47 (2007), pp.~715--733,
  \url{https://doi.org/10.1007/s10543-007-0147-7}.

\bibitem{Bur15}
{\sc E.~Burman}, {\em A monotonicity preserving, nonlinear, finite element
  upwind method for the transport equation}, Applied Mathematics Letters, 49
  (2015), pp.~141--146,
  \url{https://doi.org/https://doi.org/10.1016/j.aml.2015.05.005}.

\bibitem{BE02}
{\sc E.~Burman and A.~Ern}, {\em Nonlinear diffusion and discrete maximum
  principle for stabilized {G}alerkin approximations of the
  convection--diffusion-reaction equation}, Comput. Methods Appl. Mech. Engrg.,
  191 (2002), pp.~3833--3855,
  \url{https://doi.org/10.1016/S0045-7825(02)00318-3}.

\bibitem{BE04}
{\sc E.~Burman and A.~Ern}, {\em Discrete maximum principle for {G}alerkin
  approximations of the {L}aplace operator on arbitrary meshes}, C. R. Math.
  Acad. Sci. Paris, 338 (2004), pp.~641--646,
  \url{https://doi.org/10.1016/j.crma.2004.02.010}.

\bibitem{BE05}
{\sc E.~Burman and A.~Ern}, {\em Stabilized {G}alerkin approximation of
  convection-diffusion-reaction equations: discrete maximum principle and
  convergence}, Math. Comp., 74 (2005), pp.~1637--1652,
  \url{https://doi.org/10.1090/S0025-5718-05-01761-8}.

\bibitem{BH04}
{\sc E.~Burman and P.~Hansbo}, {\em Edge stabilization for {G}alerkin
  approximations of convection-diffusion-reaction problems}, Comput. Methods
  Appl. Mech. Engrg., 193 (2004), pp.~1437--1453,
  \url{https://doi.org/10.1016/j.cma.2003.12.032}.

\bibitem{CH84}
{\sc I.~Christie and C.~Hall}, {\em The maximum principle for bilinear
  elements}, Internat. J. Numer. Methods Engrg., 20 (1984), pp.~549--553,
  \url{https://doi.org/10.1002/nme.1620200312}.

\bibitem{Ciarlet70}
{\sc P.~G. Ciarlet}, {\em Discrete maximum principle for finite-difference
  operators}, Aequationes Math., 4 (1970), pp.~338--352,
  \url{https://doi.org/10.1007/BF01844166}.

\bibitem{CR73}
{\sc P.~G. Ciarlet and P.-A. Raviart}, {\em Maximum principle and uniform
  convergence for the finite element method}, Comput. Methods Appl. Mech.
  Engrg., 2 (1973), pp.~17--31,
  \url{https://doi.org/10.1016/0045-7825(73)90019-4}.

\bibitem{Codina93}
{\sc R.~Codina}, {\em A discontinuity-capturing crosswind-dissipation for the
  finite element solution of the convection-diffusion equation}, Comput.
  Methods Appl. Mech. Engrg., 110 (1993), pp.~325--342,
  \url{https://doi.org/10.1016/0045-7825(93)90213-H}.

\bibitem{Col55}
{\sc L.~Collatz}, {\em Numerische {B}ehandlung von {D}ifferentialgleichungen},
  Die Grundlehren der mathematischen Wissenschaften in Einzeldarstellungen mit
  besonderer Ber\"{u}cksichtigung der Anwendungsgebiete, Bd. LX,
  Springer-Verlag, Berlin, 1955.
\newblock 2te Aufl.

\bibitem{Col60}
{\sc L.~Collatz}, {\em The numerical treatment of differential equations. 3d
  ed}, Die Grundlehren der mathematischen Wissenschaften, Bd. 60,
  Springer-Verlag, Berlin, 1960.
\newblock Translated from a supplemented version of the 2d German edition by P.
  G. Williams.

\bibitem{Dav04}
{\sc T.~A. Davis}, {\em Algorithm 832: {UMFPACK} {V}4.3---an
  unsymmetric-pattern multifrontal method}, ACM Trans. Math. Software, 30
  (2004), pp.~196--199, \url{https://doi.org/10.1145/992200.992206}.

\bibitem{DE12}
{\sc D.~A. Di~Pietro and A.~Ern}, {\em Mathematical aspects of discontinuous
  {G}alerkin methods}, Springer, Heidelberg, 2012,
  \url{https://doi.org/10.1007/978-3-642-22980-0}.

\bibitem{DF15}
{\sc V.~Dolej\v{s}\'{\i} and M.~Feistauer}, {\em Discontinuous {G}alerkin
  method. Analysis and applications to compressible flow}, Springer, Cham,
  2015, \url{https://doi.org/10.1007/978-3-319-19267-3}.

\bibitem{DDS05}
{\sc A.~Dr{\u{a}}g{\u{a}}nescu, T.~F. Dupont, and L.~R. Scott}, {\em Failure of
  the discrete maximum principle for an elliptic finite element problem}, Math.
  Comp., 74 (2005), pp.~1--23,
  \url{https://doi.org/10.1090/S0025-5718-04-01651-5}.

\bibitem{EG13}
{\sc A.~Ern and J.-L. Guermond}, {\em Weighting the edge stabilization}, SIAM
  J. Numer. Anal., 51 (2013), pp.~1655--1677,
  \url{https://doi.org/10.1137/120867482}.

\bibitem{EG21-II}
{\sc A.~Ern and J.-L. Guermond}, {\em Finite elements {II}---{G}alerkin
  approximation, elliptic and mixed {PDE}s}, Springer, Cham, 2021,
  \url{https://doi.org/10.1007/978-3-030-56923-5}.

\bibitem{Eva10}
{\sc L.~C. Evans}, {\em Partial differential equations}, vol.~19 of Graduate
  Studies in Mathematics, American Mathematical Society, Providence, RI,
  second~ed., 2010.

\bibitem{FH06}
{\sc I.~Farag{\'o} and R.~Horv{\'a}th}, {\em Discrete maximum principle and
  adequate discretizations of linear parabolic problems}, SIAM J. Sci. Comput.,
  28 (2006), pp.~2313--2336, \url{https://doi.org/10.1137/050627241}.

\bibitem{FKK12}
{\sc I.~Farag{\'o}, J.~Kar{\'a}tson, and S.~Korotov}, {\em Discrete maximum
  principles for nonlinear parabolic {PDE} systems}, IMA J. Numer. Anal., 32
  (2012), pp.~1541--1573, \url{https://doi.org/10.1093/imanum/drr050}.

\bibitem{floater_2015}
{\sc M.~S. Floater}, {\em Generalized barycentric coordinates and
  applications}, Acta Numerica, 24 (2015), pp.~161--214,
  \url{https://doi.org/10.1017/S0962492914000129}.

\bibitem{FF95}
{\sc L.~P. Franca and C.~Farhat}, {\em Bubble functions prompt unusual
  stabilized finite element methods}, Comput. Methods Appl. Mech. Engrg., 123
  (1995), pp.~299--308, \url{https://doi.org/10.1016/0045-7825(94)00721-X}.

\bibitem{GJM+16}
{\sc S.~Ganesan, V.~John, G.~Matthies, R.~Meesala, S.~Abdus, and U.~Wilbrandt},
  {\em An object oriented parallel finite element scheme for computing {PDE}s:
  {D}esign and implementation}, in IEEE 23rd International Conference on High
  Performance Computing Workshops (HiPCW) Hyderabad, IEEE, 2016, pp.~106--115,
  \url{https://doi.org/10.1109/HiPCW.2016.023}.

\bibitem{Ger30}
{\sc S.~A. {Gershgorin}}, {\em {Fehlerabsch\"atzung f\"ur das
  Differenzenverfahren zur L\"osung partieller Differentialgleichungen}}, {Z.
  Angew. Math. Mech.}, 10 (1930), pp.~373--382,
  \url{https://doi.org/10.1002/zamm.19300100409}.

\bibitem{GT01}
{\sc D.~Gilbarg and N.~S. Trudinger}, {\em Elliptic partial differential
  equations of second order}, Springer-Verlag, Berlin, second~ed., 2001,
  \url{https://doi.org/10.1007/978-3-642-61798-0}.

\bibitem{God59}
{\sc S.~K. Godunov}, {\em A difference method for numerical calculation of
  discontinuous solutions of the equations of hydrodynamics}, Mat. Sb. (N.S.),
  47 (89) (1959), pp.~271--306.

\bibitem{GN14}
{\sc J.-L. Guermond and M.~Nazarov}, {\em A maximum-principle preserving
  {$C^0$} finite element method for scalar conservation equations}, Comput.
  Methods Appl. Mech. Engrg., 272 (2014), pp.~198--213,
  \url{https://doi.org/10.1016/j.cma.2013.12.015}.

\bibitem{GNPY14}
{\sc J.-L. Guermond, M.~Nazarov, B.~Popov, and Y.~Yang}, {\em A second-order
  maximum principle preserving {L}agrange finite element technique for
  nonlinear scalar conservation equations}, SIAM J. Numer. Anal., 52 (2014),
  pp.~2163--2182, \url{https://doi.org/10.1137/130950240}.

\bibitem{GP17}
{\sc J.-L. Guermond and B.~Popov}, {\em Invariant domains and second-order
  continuous finite element approximation for scalar conservation equations},
  SIAM J. Numer. Anal., 55 (2017), pp.~3120--3146,
  \url{https://doi.org/10.1137/16M1106560}.

\bibitem{GPT19}
{\sc J.-L. Guermond, B.~Popov, and I.~Tomas}, {\em Invariant domain preserving
  discretization-independent schemes and convex limiting for hyperbolic
  systems}, Comput. Methods Appl. Mech. Engrg., 347 (2019), pp.~143--175,
  \url{https://doi.org/10.1016/j.cma.2018.11.036}.

\bibitem{GPY17}
{\sc J.-L. Guermond, B.~Popov, and Y.~Yang}, {\em The effect of the consistent
  mass matrix on the maximum-principle for scalar conservation equations}, J.
  Sci. Comput., 70 (2017), pp.~1358--1366,
  \url{https://doi.org/10.1007/s10915-016-0285-7}.

\bibitem{Hajduk21}
{\sc H.~Hajduk}, {\em Monolithic convex limiting in discontinuous {G}alerkin
  discretizations of hyperbolic conservation laws}, Comput. Math. Appl., 87
  (2021), pp.~120--138, \url{https://doi.org/10.1016/j.camwa.2021.02.012}.

\bibitem{HM81}
{\sc W.~H{\"o}hn and H.-D. Mittelmann}, {\em Some remarks on the discrete
  maximum-principle for finite elements of higher order}, Computing, 27 (1981),
  pp.~145--154, \url{https://doi.org/10.1007/BF02243548}.

\bibitem{HM13}
{\sc T.~L. Horv\'{a}th and M.~E. Mincsovics}, {\em Discrete maximum principle
  for interior penalty discontinuous {G}alerkin methods}, Cent. Eur. J. Math.,
  11 (2013), pp.~664--679, \url{https://doi.org/10.2478/s11533-012-0154-z}.

\bibitem{Hua11}
{\sc W.~Huang}, {\em Discrete maximum principle and a {D}elaunay-type mesh
  condition for linear finite element approximations of two-dimensional
  anisotropic diffusion problems}, Numer. Math. Theory Methods Appl., 4 (2011),
  pp.~319--334, \url{https://doi.org/10.4208/nmtma.2011.m1024}.

\bibitem{Ikeda83}
{\sc T.~Ikeda}, {\em Maximum principle in finite element models for
  convection-diffusion phenomena}, North-Holland, Amsterdam, 1983.

\bibitem{Jameson}
{\sc A.~Jameson}, {\em Origins and further development of the
  {J}ameson-{S}chmidt-{T}urkel scheme}, AIAA J., 55 (2017), pp.~1487--1510,
  \url{https://doi.org/10.2514/1.J055493}.

\bibitem{JST81}
{\sc A.~Jameson, W.~Schmidt, and E.~Turkel}, {\em Numerical solution of the
  {E}uler equations by finite volume methods using {R}unge-{K}utta
  time-stepping schemes}, in 14th AIAA Fluid and Plasma Dynamics Conference,
  Palo Alto, CA (USA), 23-25 Jun 1981, AIAA meeting paper 1981-1259, 1981,
  \url{https://doi.org/10.2514/6.1981-1259}.

\bibitem{Jha21}
{\sc A.~Jha}, {\em A residual based a posteriori error estimators for {AFC}
  schemes for convection-diffusion equations}, Comput. Math. Appl., 97 (2021),
  pp.~86--99, \url{https://doi.org/10.1016/j.camwa.2021.05.031}.

\bibitem{JJ19}
{\sc A.~Jha and V.~John}, {\em A study of solvers for nonlinear {AFC}
  discretizations of convection-diffusion equations}, Comput. Math. Appl., 78
  (2019), pp.~3117--3138, \url{https://doi.org/10.1016/j.camwa.2019.04.020}.

\bibitem{JJK22x}
{\sc A.~Jha, V.~John, and P.~Knobloch}, {\em Adaptive grids in the context of
  algebraic stabilizations for convection--diffusion--reaction equations},
  2022, \url{https://arxiv.org/abs/2007.08405}.

\bibitem{JK07}
{\sc V.~John and P.~Knobloch}, {\em On spurious oscillations at layers
  diminishing ({SOLD}) methods for convection--diffusion equations: {P}art {I}
  -- {A} review}, Comput. Methods Appl. Mech. Engrg., 196 (2007),
  pp.~2197--2215, \url{https://doi.org/10.1016/j.cma.2006.11.013}.

\bibitem{JK08}
{\sc V.~John and P.~Knobloch}, {\em On spurious oscillations at layers
  diminishing ({SOLD}) methods for convection--diffusion equations: {P}art {II}
  -- {A}nalysis for {$P_1$} and {$Q_1$} finite elements}, Comput. Methods Appl.
  Mech. Engrg., 197 (2008), pp.~1997--2014,
  \url{https://doi.org/10.1016/j.cma.2007.12.019}.

\bibitem{JK21}
{\sc V.~John and P.~Knobloch}, {\em Existence of solutions of a finite element
  flux-corrected-transport scheme}, Appl. Math. Lett., 115 (2021), p.~Paper No.
  106932, \url{https://doi.org/10.1016/j.aml.2020.106932}.

\bibitem{JK22}
{\sc V.~John and P.~Knobloch}, {\em On algebraically stabilized schemes for
  convection--diffusion--reaction problems}, Numer. Math., 152 (2022),
  pp.~553--585, \url{https://doi.org/10.1007/s00211-022-01325-9}.

\bibitem{JKK21}
{\sc V.~John, P.~Knobloch, and P.~Korsmeier}, {\em On the solvability of the
  nonlinear problems in an algebraically stabilized finite element method for
  evolutionary transport-dominated equations}, Math. Comp., 90 (2021),
  pp.~595--611, \url{https://doi.org/10.1090/mcom/3576}.

\bibitem{JKP22}
{\sc V.~John, P.~Knobloch, and O.~P\'artl}, {\em A numerical assessment of
  finite element discretizations for convection-diffusion-reaction equations
  satisfying discrete maximum principles}, Comput. Methods Appl. Math.,
  (2022), \url{https://doi.org/10.1515/cmam-2022-0125}.

\bibitem{JM01}
{\sc V.~John and G.~Matthies}, {\em {Higher-order finite element
  discretizations in a benchmark problem for incompressible flows.}}, Int. J.
  Numer. Methods Fluids, 37 (2001), pp.~885--903,
  \url{https://doi.org/10.1002/fld.195}.

\bibitem{JMRSTV09}
{\sc V.~John, T.~Mitkova, M.~Roland, K.~Sundmacher, L.~Tobiska, and A.~Voigt},
  {\em Simulations of population balance systems with one internal coordinate
  using finite element methods}, Chemical Engineering Science, 64 (2009),
  pp.~733--741, \url{https://doi.org/10.1016/j.ces.2008.05.004}.

\bibitem{JN12}
{\sc V.~John and J.~Novo}, {\em On (essentially) non-oscillatory
  discretizations of evolutionary convection-diffusion equations}, J. Comput.
  Phys., 231 (2012), pp.~1570--1586,
  \url{https://doi.org/10.1016/j.jcp.2011.10.025}.

\bibitem{Kan78}
{\sc H.~Kanayama}, {\em Discrete maximum principles for salinity distribution
  in a bay: conservation law and maximum principle}, Theoretical Appl. Mech.,
  28 (1978), pp.~559--579.

\bibitem{KKK07}
{\sc J.~Kar{\'a}tson, S.~Korotov, and M.~K{\v{r}}{\'{\i}}{\v{z}}ek}, {\em On
  discrete maximum principles for nonlinear elliptic problems}, Math. Comput.
  Simulation, 76 (2007), pp.~99--108,
  \url{https://doi.org/10.1016/j.matcom.2007.01.011}.

\bibitem{Kik77}
{\sc F.~Kikuchi}, {\em Discrete maximum principle and artificial viscosity in
  finite element approximations to convective diffusion equations}, Institute
  of Space and Aeronautical Science, University of Tokyo, 550 (1977).

\bibitem{Kno06}
{\sc P.~Knobloch}, {\em Improvements of the {M}izukami--{H}ughes method for
  convection--diffusion equations}, Comput. Methods Appl. Mech. Engrg., 196
  (2006), pp.~579--594, \url{https://doi.org/10.1016/j.cma.2006.06.004}.

\bibitem{Knobloch06b}
{\sc P.~Knobloch}, {\em Numerical solution of convection--diffusion equations
  using upwinding techniques satisfying the discrete maximum principle}, in
  Proceedings of {C}zech-{J}apanese {S}eminar in {A}pplied {M}athematics 2005,
  M.~Bene{\v s}, M.~Kimura, and T.~Nakaki, eds., vol.~3 of COE Lect. Note,
  Kyushu Univ., Fukuoka, 2006, pp.~69--76.

\bibitem{Kno07}
{\sc P.~Knobloch}, {\em Application of the {M}izukami--{H}ughes method to
  bilinear finite elements}, in Proceedings of {C}zech-{J}apanese {S}eminar in
  {A}pplied {M}athematics 2006, M.~Bene\v{s}, M.~Kimura, and T.~Nakaki, eds.,
  vol.~6 of COE Lect. Note, Kyushu Univ., Fukuoka, 2007, pp.~137--147.

\bibitem{Knobloch10}
{\sc P.~Knobloch}, {\em Numerical solution of convection-diffusion equations
  using a nonlinear method of upwind type}, J. Sci. Comput., 43 (2010),
  pp.~454--470, \url{https://doi.org/10.1007/s10915-008-9260-2}.

\bibitem{Kno17}
{\sc P.~Knobloch}, {\em On the discrete maximum principle for algebraic flux
  correction schemes with limiters of upwind type}, in Boundary and Interior
  Layers, Computational and Asymptotic Methods BAIL 2016, Z.~Huang, M.~Stynes,
  and Z.~Zhang, eds., vol.~120 of Lect. Notes Comput. Sci. Eng.,
  Springer-Verlag, Cham, 2017, pp.~129--139,
  \url{https://doi.org/10.1007/978-3-319-67202-1_10}.

\bibitem{Kno21}
{\sc P.~Knobloch}, {\em A new algebraically stabilized method for
  convection--diffusion--reaction equations}, in Numerical mathematics and
  advanced applications ENUMATH 2019, F.~Vermolen and C.~Vuik, eds., vol.~139
  of Lect. Notes Comput. Sci. Eng., Springer-Verlag, Cham, 2021, pp.~605--613,
  \url{https://doi.org/10.1007/978-3-030-55874-1_59}.

\bibitem{Kno22x}
{\sc P.~Knobloch}, {\em An algebraically stabilized method for
  convection--diffusion--reaction problems with optimal experimental
  convergence rates on general meshes}, 2022,
  \url{https://arxiv.org/abs/2208.07705}.

\bibitem{KV10}
{\sc S.~{Korotov} and T.~{Vejchodsk\'y}}, {\em A comparison of simplicial and
  block finite elements}, in Numerical mathematics and advanced applications
  2009. Proceedings of ENUMATH 2009, G.~Kreiss, P.~L\"otstedt, A.~M\r{a}lqvist,
  and M.~Neytcheva, eds., Springer-Verlag, Berlin, 2010, pp.~533--541,
  \url{https://doi.org/10.1007/978-3-642-11795-4_57}.

\bibitem{Kuz06}
{\sc D.~Kuzmin}, {\em On the design of general-purpose flux limiters for finite
  element schemes. {I}. {S}calar convection}, J. Comput. Phys., 219 (2006),
  pp.~513--531, \url{https://doi.org/10.1016/j.jcp.2006.03.034}.

\bibitem{Kuz07}
{\sc D.~Kuzmin}, {\em Algebraic flux correction for finite element
  discretizations of coupled systems}, in Proceedings of the Int.~Conf.~on
  Computational Methods for Coupled Problems in Science and Engineering,
  M.~Papadrakakis, E.~O{\~n}ate, and B.~Schrefler, eds., CIMNE, Barcelona,
  2007, pp.~1--5.

\bibitem{Kuz08}
{\sc D.~Kuzmin}, {\em On the design of algebraic flux correction schemes for
  quadratic finite elements}, J. Comput. Appl. Math., 218 (2008), pp.~79--87,
  \url{https://doi.org/10.1016/j.cam.2007.04.045}.

\bibitem{Kuz09}
{\sc D.~Kuzmin}, {\em Explicit and implicit {FEM}-{FCT} algorithms with flux
  linearization}, J. Comput. Phys., 228 (2009), pp.~2517--2534,
  \url{https://doi.org/10.1016/j.jcp.2008.12.011}.

\bibitem{Kuz12a}
{\sc D.~Kuzmin}, {\em Algebraic flux correction {I}. {S}calar conservation
  laws}, in Flux-corrected transport. Principles, algorithms, and applications,
  D.~Kuzmin, R.~L\"ohner, and S.~Turek, eds., Springer, Dordrecht, second~ed.,
  2012, pp.~145--192, \url{https://doi.org/10.1007/978-94-007-4038-9_6}.

\bibitem{Kuz12}
{\sc D.~Kuzmin}, {\em Linearity-preserving flux correction and convergence
  acceleration for constrained {G}alerkin schemes}, J. Comput. Appl. Math., 236
  (2012), pp.~2317--2337, \url{https://doi.org/10.1016/j.cam.2011.11.019}.

\bibitem{Kuz20}
{\sc D.~Kuzmin}, {\em Monolithic convex limiting for continuous finite element
  discretizations of hyperbolic conservation laws}, Comput. Methods Appl. Mech.
  Engrg., 361 (2020), p.~Paper No. 112804,
  \url{https://doi.org/10.1016/j.cma.2019.112804}.

\bibitem{KBS17}
{\sc D.~Kuzmin, S.~Basting, and J.~N. Shadid}, {\em Linearity-preserving
  monotone local projection stabilization schemes for continuous finite
  elements}, Comput. Methods Appl. Mech. Engrg., 322 (2017), pp.~23--41,
  \url{https://doi.org/10.1016/j.cma.2017.04.030}.

\bibitem{KHR20}
{\sc D.~Kuzmin, H.~Hajduk, and A.~Rupp}, {\em Locally bound-preserving enriched
  {G}alerkin methods for the linear advection equation}, Comput. \& Fluids, 205
  (2020), p.~Paper No. 104525,
  \url{https://doi.org/10.1016/j.compfluid.2020.104525}.

\bibitem{KH15}
{\sc D.~Kuzmin and J.~H\"{a}m\"{a}l\"{a}inen}, {\em Finite element methods for
  computational fluid dynamics: {A} practical guide}, Society for Industrial
  and Applied Mathematics (SIAM), Philadelphia, PA, 2015,
  \url{https://doi.org/10.1137/1.9781611973617}.

\bibitem{KS17}
{\sc D.~Kuzmin and J.~N. Shadid}, {\em Gradient-based nodal limiters for
  artificial diffusion operators in finite element schemes for transport
  equations}, Internat. J. Numer. Methods Fluids, 84 (2017), pp.~675--695,
  \url{https://doi.org/10.1002/fld.4365}.

\bibitem{KT02}
{\sc D.~Kuzmin and S.~Turek}, {\em Flux correction tools for finite elements},
  J. Comput. Phys., 175 (2002), pp.~525--558,
  \url{https://doi.org/10.1006/jcph.2001.6955}.

\bibitem{KT04}
{\sc D.~Kuzmin and S.~Turek}, {\em High-resolution {FEM}-{TVD} schemes based on
  a fully multidimensional flux limiter}, J. Comput. Phys., 198 (2004),
  pp.~131--158, \url{https://doi.org/10.1016/j.jcp.2004.01.015}.

\bibitem{LL21}
{\sc D.~Leykekhman and B.~Li}, {\em Weak discrete maximum principle of finite
  element methods in convex polyhedra}, Math. Comp., 90 (2021), pp.~1--18,
  \url{https://doi.org/10.1090/mcom/3560}.

\bibitem{LZ20}
{\sc H.~Li and X.~Zhang}, {\em On the monotonicity and discrete maximum
  principle of the finite difference implementation of {$C^0$}-{$Q^2$} finite
  element method}, Numer. Math., 145 (2020), pp.~437--472,
  \url{https://doi.org/10.1007/s00211-020-01110-6}.

\bibitem{LH10}
{\sc X.~Li and W.~Huang}, {\em An anisotropic mesh adaptation method for the
  finite element solution of heterogeneous anisotropic diffusion problems}, J.
  Comput. Phys., 229 (2010), pp.~8072--8094,
  \url{https://doi.org/10.1016/j.jcp.2010.07.009}.

\bibitem{LH13}
{\sc X.~Li and W.~Huang}, {\em Maximum principle for the finite element
  solution of time-dependent anisotropic diffusion problems}, Numer. Methods
  Partial Differential Equations, 29 (2013), pp.~1963--1985,
  \url{https://doi.org/10.1002/num.21784}.

\bibitem{Loh20}
{\sc C.~Lohmann}, {\em Physics-compatible finite element methods for scalar and
  tensorial advection problems}, Springer Spektrum, Wiesbaden, 2019,
  \url{https://doi.org/10.1007/978-3-658-27737-6}.

\bibitem{LKSM17}
{\sc C.~Lohmann, D.~Kuzmin, J.~N. Shadid, and S.~Mabuza}, {\em Flux-corrected
  transport algorithms for continuous {G}alerkin methods based on high order
  {B}ernstein finite elements}, J. Comput. Phys., 344 (2017), pp.~151--186,
  \url{https://doi.org/10.1016/j.jcp.2017.04.059}.

\bibitem{LMPV87}
{\sc R.~{L\"ohner}, K.~{Morgan}, J.~{Peraire}, and M.~{Vahdati}}, {\em {Finite
  element flux-corrected transport (FEM-FCT) for the Euler and Navier-Stokes
  equations.}}, {Int. J. Numer. Methods Fluids}, 7 (1987), pp.~1093--1109,
  \url{https://doi.org/10.1002/fld.1650071007}.

\bibitem{Lor77}
{\sc J.~Lorenz}, {\em Zur {I}nversmonotonie diskreter {P}robleme}, Numer.
  Math., 27 (1976/77), pp.~227--238, \url{https://doi.org/10.1007/BF01396643}.

\bibitem{LHQ14}
{\sc C.~Lu, W.~Huang, and J.~Qiu}, {\em Maximum principle in linear finite
  element approximations of anisotropic diffusion-convection-reaction
  problems}, Numer. Math., 127 (2014), pp.~515--537,
  \url{https://doi.org/10.1007/s00211-013-0595-8}.

\bibitem{MSK18}
{\sc S.~Mabuza, J.~N. Shadid, and D.~Kuzmin}, {\em Local bounds preserving
  stabilization for continuous {G}alerkin discretization of hyperbolic
  systems}, J. Comput. Phys., 361 (2018), pp.~82--110,
  \url{https://doi.org/10.1016/j.jcp.2018.01.048}.

\bibitem{MH85}
{\sc A.~Mizukami and T.~J.~R. Hughes}, {\em A {P}etrov-{G}alerkin finite
  element method for convection-dominated flows: an accurate upwinding
  technique for satisfying the maximum principle}, Comput. Methods Appl. Mech.
  Engrg., 50 (1985), pp.~181--193,
  \url{https://doi.org/10.1016/0045-7825(85)90089-1}.

\bibitem{OU84}
{\sc K.~Ohmori and T.~Ushijima}, {\em A technique of upstream type applied to a
  linear nonconforming finite element approximation of convective diffusion
  equations}, RAIRO Anal. Num\'er., 18 (1984), pp.~309--332,
  \url{https://doi.org/10.1051/m2an/1984180303091}.

\bibitem{PC86}
{\sc A.~K. Parrott and M.~A. Christie}, {\em {FCT} applied to the 2-{D} finite
  element solution of tracer transport by single phase flow in a porous
  medium}, in Numerical methods for fluid dynamics II, Proc. Conf., Reading/UK
  1985, vol.~7 of Inst. Math. Appl. Conf. Ser., New Ser., 1986, pp.~609--619.

\bibitem{Pazner21}
{\sc W.~Pazner}, {\em Sparse invariant domain preserving discontinuous
  {G}alerkin methods with subcell convex limiting}, Comput. Methods Appl. Mech.
  Engrg., 382 (2021), p.~Paper No. 113876,
  \url{https://doi.org/10.1016/j.cma.2021.113876}.

\bibitem{RT92}
{\sc R.~Rannacher and S.~Turek}, {\em Simple nonconforming quadrilateral
  {S}tokes element}, Numer. Methods Partial Differential Equations, 8 (1992),
  pp.~97--111, \url{https://doi.org/10.1002/num.1690080202}.

\bibitem{RH73}
{\sc W.~Reed and T.~Hill}, {\em Triangular mesh methods for the neutron
  transport equation}, Technical Report LA-UR-73-479, Los Alamos Scientific
  Laboratory, Los Alamos, NM, 1973.

\bibitem{Riv08}
{\sc B.~Rivi\`ere}, {\em Discontinuous {G}alerkin methods for solving elliptic
  and parabolic equations. Theory and implementation}, Society for Industrial
  and Applied Mathematics (SIAM), Philadelphia, PA, 2008,
  \url{https://doi.org/10.1137/1.9780898717440}.

\bibitem{RST08}
{\sc H.-G. Roos, M.~Stynes, and L.~Tobiska}, {\em Robust numerical methods for
  singularly perturbed differential equations. Convection-diffusion-reaction
  and flow problems}, Springer-Verlag, Berlin, second~ed., 2008,
  \url{https://doi.org/10.1007/978-3-540-34467-4}.

\bibitem{Sch80}
{\sc A.~H. Schatz}, {\em A weak discrete maximum principle and stability of the
  finite element method in {$L_{\infty }$} on plane polygonal domains. {I}},
  Math. Comp., 34 (1980), pp.~77--91, \url{https://doi.org/10.2307/2006221}.

\bibitem{STW10}
{\sc A.~H. Schatz, V.~Thom\'{e}e, and L.~B. Wahlbin}, {\em On positivity and
  maximum-norm contractivity in time stepping methods for parabolic equations},
  Comput. Methods Appl. Math., 10 (2010), pp.~421--443,
  \url{https://doi.org/10.2478/cmam-2010-0025}.

\bibitem{Shu16}
{\sc C.-W. Shu}, {\em Discontinuous {G}alerkin methods for time-dependent
  convection dominated problems: basics, recent developments and comparison
  with other methods}, in Building bridges: connections and challenges in
  modern approaches to numerical partial differential equations, G.~R.
  Barrenechea, F.~Brezzi, A.~Cangiani, and E.~H. Georgoulis, eds., vol.~114 of
  Lect. Notes Comput. Sci. Eng., Springer, Cham, 2016, pp.~369--397,
  \url{https://doi.org/10.1007/978-3-319-41640-3_12}.

\bibitem{Stockie11}
{\sc J.~M. Stockie}, {\em The mathematics of atmospheric dispersion modeling},
  SIAM Rev., 53 (2011), pp.~349--372, \url{https://doi.org/10.1137/10080991X}.

\bibitem{SF73}
{\sc G.~Strang and G.~J. Fix}, {\em An analysis of the finite element method},
  Prentice-Hall, Inc., Englewood Cliffs, N. J., 1973.

\bibitem{Tabata77}
{\sc M.~Tabata}, {\em A finite element approximation corresponding to the
  upwind finite differencing}, Mem. Numer. Math., 4 (1977), pp.~47--63.

\bibitem{TW08}
{\sc V.~Thom\'{e}e and L.~B. Wahlbin}, {\em On the existence of maximum
  principles in parabolic finite element equations}, Math. Comp., 77 (2008),
  pp.~11--19, \url{https://doi.org/10.1090/S0025-5718-07-02021-2}.

\bibitem{TsiPet13}
{\sc C.~Tsiotsios and M.~Petrou}, {\em On the choice of the parameters for
  anisotropic diffusion in image processing}, Pattern Recognition, 46 (2013),
  pp.~1369--1381,
  \url{https://doi.org/https://doi.org/10.1016/j.patcog.2012.11.012}.

\bibitem{Tur94}
{\sc S.~Turek}, {\em Tools for simulating nonstationary incompressible flow via
  discretely divergence-free finite element models}, Internat. J. Numer.
  Methods Fluids, 18 (1994), pp.~71--105,
  \url{https://doi.org/10.1002/fld.1650180105}.

\bibitem{Var00}
{\sc R.~S. Varga}, {\em Matrix iterative analysis}, Springer-Verlag, Berlin,
  2000, \url{https://doi.org/10.1007/978-3-642-05156-2}.

\bibitem{Vej10}
{\sc T.~{Vejchodsk\'y}}, {\em {Angle conditions for discrete maximum principles
  in higher-order FEM}}, in Numerical mathematics and advanced applications
  2009. Proceedings of ENUMATH 2009, G.~Kreiss, P.~L\"otstedt, A.~M\r{a}lqvist,
  and M.~Neytcheva, eds., Springer-Verlag, Berlin, 2010, pp.~901--909,
  \url{https://doi.org/10.1007/978-3-642-11795-4_97}.

\bibitem{Vej11}
{\sc T.~{Vejchodsk\'y}}, {\em Discrete Maximum Principles}, habilitation,
  Charles University Prague, Faculty of Mathematics and Physics, 2011.

\bibitem{VS07a}
{\sc T.~Vejchodsk{\'y} and P.~{\v{S}}ol{\'{\i}}n}, {\em Discrete maximum
  principle for a 1{D} problem with piecewise-constant coefficients solved by
  {$hp$}-{FEM}}, J. Numer. Math., 15 (2007), pp.~233--243,
  \url{https://doi.org/10.1515/jnma.2007.011}.

\bibitem{VS07b}
{\sc T.~Vejchodsk{\'y} and P.~{\v{S}}ol{\'{\i}}n}, {\em Discrete maximum
  principle for higher-order finite elements in 1{D}}, Math. Comp., 76 (2007),
  pp.~1833--1846, \url{https://doi.org/10.1090/S0025-5718-07-02022-4}.

\bibitem{VS18}
{\sc F.~J. Vermolen and A.~Segal}, {\em On an integration rule for products of
  barycentric coordinates over simplexes in {$\Bbb R^n$}}, J. Comput. Appl.
  Math., 330 (2018), pp.~289--294,
  \url{https://doi.org/10.1016/j.cam.2017.09.013}.

\bibitem{WSHD07}
{\sc J.~Warren, S.~Schaefer, A.~N. Hirani, and M.~Desbrun}, {\em Barycentric
  coordinates for convex sets}, Adv.Comput. Math., 39 (2007), pp.~319--338,
  \url{https://doi.org/10.1007/s10444-005-9008-6}.

\bibitem{Wes01}
{\sc P.~Wesseling}, {\em Principles of computational fluid dynamics},
  Springer-Verlag, Berlin, 2001,
  \url{https://doi.org/10.1007/978-3-642-05146-3}.

\bibitem{WB_J17}
{\sc U.~Wilbrandt, C.~Bartsch, N.~Ahmed, N.~Alia, F.~Anker, L.~Blank,
  A.~Caiazzo, S.~Ganesan, S.~Giere, G.~Matthies, R.~Meesala, A.~Shamim,
  J.~Venkatesan, and V.~John}, {\em Par{M}oo{N}---{A} modernized program
  package based on mapped finite elements}, Comput. Math. Appl., 74 (2017),
  pp.~74--88, \url{https://doi.org/10.1016/j.camwa.2016.12.020}.

\bibitem{XZ99}
{\sc J.~Xu and L.~Zikatanov}, {\em A monotone finite element scheme for
  convection-diffusion equations}, Math. Comp., 68 (1999), pp.~1429--1446,
  \url{https://doi.org/10.1090/S0025-5718-99-01148-5}.

\bibitem{Zal79}
{\sc S.~T. Zalesak}, {\em Fully multidimensional flux-corrected transport
  algorithms for fluids}, J. Comput. Phys., 31 (1979), pp.~335--362,
  \url{https://doi.org/10.1016/0021-9991(79)90051-2}.

\bibitem{Zal12}
{\sc S.~T. Zalesak}, {\em The design of flux-corrected transport ({FCT})
  algorithms for structured grids}, in Flux-corrected transport. Principles,
  algorithms, and applications, D.~Kuzmin, R.~L\"ohner, and S.~Turek, eds.,
  Springer, Dordrecht, second~ed., 2012, pp.~23--65,
  \url{https://doi.org/10.1007/978-94-007-4038-9_2}.

\bibitem{ZS10}
{\sc X.~Zhang and C.-W. Shu}, {\em On maximum-principle-satisfying high order
  schemes for scalar conservation laws}, J. Comput. Phys., 229 (2010),
  pp.~3091--3120, \url{https://doi.org/10.1016/j.jcp.2009.12.030}.

\bibitem{ZS13}
{\sc X.~Zhang and C.-W. Shu}, {\em Maximum-principle-satisfying and
  positivity-preserving high-order schemes for conservation laws: {S}urvey and
  new developments}, Proc. R. Soc. Lond. Ser. A Math. Phys. Eng. Sci., 467
  (2011), pp.~2752--2776, \url{https://doi.org/10.1098/rspa.2011.0153}.

\bibitem{ZZS13}
{\sc Y.~Zhang, X.~Zhang, and C.-W. Shu}, {\em Maximum-principle-satisfying
  second order discontinuous {G}alerkin schemes for convection-diffusion
  equations on triangular meshes}, J. Comput. Phys., 234 (2013), pp.~295--316,
  \url{https://doi.org/10.1016/j.jcp.2012.09.032}.

\end{thebibliography}

\def\cprime{$'$}

\end{document}